\documentclass[10pt]{amsart}

\usepackage[T1]{fontenc}  
\usepackage[utf8]{inputenc} 
\usepackage{longtable}
\usepackage{amssymb}
\usepackage{amsthm}
\usepackage{amsmath}
\usepackage{wasysym}
\usepackage{dcpic, pictex}
\usepackage{bbm}
\usepackage{enumerate}
\usepackage{amsfonts}
\usepackage{amsthm}
\usepackage{amsmath}
\usepackage{times}
\usepackage{amscd}
\usepackage[mathscr]{eucal}
\usepackage{indentfirst}
\usepackage{pict2e}
\usepackage{epic}
\usepackage{epstopdf} 
\usepackage{verbatim}
\usepackage{cancel}
\usepackage[foot]{amsaddr}
\usepackage{lipsum}
\usepackage{mathrsfs}
\usepackage{bm}

\usepackage{graphicx}

\usepackage{
amsfonts, 
latexsym, 
enumerate,
}
\usepackage[english]{babel}
\allowdisplaybreaks
\makeindex
\newtheorem{theorem}{Theorem}[section]
\newtheorem{lemma}[theorem]{Lemma}
\newtheorem{proposition}[theorem]{Proposition}
\newtheorem{definition}[theorem]{Definition}

\newtheorem{remark}[theorem]{Remark}

\newtheorem{assumption}[theorem]{Assumption}
\usepackage[a4paper,top=3.5cm,bottom=3.5cm,left=3cm,right=3cm]{geometry}
\numberwithin{equation}{section}

\allowdisplaybreaks

\usepackage{ bbold }
\usepackage{xcolor}

\newcommand{\dif}{\mathrm{d}}
\newcommand{\E}{\mathbf{E}}

\newcommand{\R}{\mathbf{R}}
\newcommand{\C}{\mathbf{C}}

\newcommand{\N}{\mathbf{N}}

\newcommand{\X}{\mathcal{X}}

\newcommand{\Cov}{\mathrm{Cov}}
\newcommand{\Var}{\mathrm{Var}}

\newcommand{\Lin}{\mathcal{L}}

\newcommand{\ii}{\mathrm{i}}
\newcommand{\ee}{\mathrm{e}}
\newcommand{\rd}{\mathrm{d}}

\newcommand{\mt}{\mathfrak{t}}

\usepackage{hyperref}

\title[Hyperuniformity]{The eigenvalues of i.i.d. matrices are hyperuniform}

\date{\today}

\begin{document}

\maketitle

\vspace{0.25cm}

\renewcommand{\thefootnote}{\fnsymbol{footnote}}

\noindent
\mbox{}%
\hfill%
\begin{minipage}{0.19\textwidth}
	\centering
	{Giorgio Cipolloni}\footnotemark[1]\\
	\footnotesize{\textit{cipolloni@axp.mat.uniroma2.it}}
\end{minipage}
\hfill%
\begin{minipage}{0.19\textwidth}
	\centering
	{L\'aszl\'o Erd\H{o}s}\footnotemark[2]\\
	\footnotesize{\textit{lerdos@ist.ac.at}}
\end{minipage}
\hfill%
\begin{minipage}{0.19\textwidth}
	\centering
	{Oleksii Kolupaiev}\footnotemark[2]\\
	\footnotesize{\textit{okolupaiev@ist.ac.at}}
\end{minipage}
\hfill%
\mbox{}%
\footnotetext[1]{University of Rome Tor Vergata, Via della Ricerca Scientifica 1, 00133 Rome, Italy.}
\footnotetext[1]{Partially supported by the MUR Excellence Department Project MatMod@TOV awarded to the Department of Mathematics, University of Rome Tor Vergata, CUP E83C18000100006.}
\footnotetext[2]{Institute of Science and Technology Austria, Am Campus 1, 3400 Klosterneuburg, Austria. 
}
\footnotetext[2]{Supported by the ERC Advanced Grant ``RMTBeyond'' No.~101020331.}

\renewcommand*{\thefootnote}{\arabic{footnote}}
\vspace{0.25cm}

\begin{abstract} 
We prove that the point process of the eigenvalues of real or complex non-Hermitian matrices $X$ with independent, identically distributed entries is hyperuniform: the variance of the number of eigenvalues in a subdomain $\Omega$ of the spectrum is much smaller than the volume of $\Omega$. Our main technical novelty is a very precise computation of the covariance between the resolvents of the Hermitization of $X-z_1, X-z_2$, for two distinct complex parameters $z_1,z_2$.
\end{abstract}
\vspace{0.15cm}


\vspace{0.25cm}
\normalsize

\tableofcontents

\section{Introduction}

We consider the \emph{i.i.d. matrix ensemble}, i.e. an ensemble of $N\times N$ non-Hermitian random matrices $X$ with independent, identically distributed (i.i.d.) real or complex centered entries. For convenience, we normalize the entries of $X$ so that $\E |X_{ij}|^2=N^{-1}$. Then eigenvalues $\sigma_i$ of such matrices form a correlated point process on the unit disk of the complex plane; see \cite{borodin2009ginibre, Ginibre65} for the Ginibre ensemble\footnote{We say that $X$ belongs to the real/complex Ginibre ensemble if its entries are (normalized) standard real/complex Gaussian random variables.} and \cite{nonherm_edge_univ, dubova2025bulk, maltsev2023bulk, osman2025bulk} for general i.i.d. matrices (see also the recent \cite{afanasiev2024universality, campbell2025spectral, cipolloni2025non, liu2023critical, liu2025repeated, zhang2024bulk}). These eigenvalues tend to be uniformly distributed over the unit disk $\mathbf{D}$. In particular, for any (sufficiently nice) subdomain of the disk $\Omega\subset\mathbf{D}$, it is well known \cite{bai1997circular, Girko84, gotze2010circular, pan2010circular, tao2008random} that (\emph{circular law})
\begin{equation}
\label{eq:approx}
\mathcal{N}_\Omega:=\#\big\{\sigma_i\in \Omega\big\}=\sum_{i=1}^N \bm1(\sigma_i\in \Omega)=\frac{N|\Omega|}{\pi}+o(N),
\end{equation}
with high probability. In particular, this implies $\E \mathcal{N}_\Omega\sim N$. Note that the lhs. of \eqref{eq:approx} is random, while the leading order term in the rhs. is deterministic. It is then natural to study the fluctuations around this deterministic limit. 
For a system of independent particles, such as a Poisson point process it holds that the \emph{number variance} is proportional to the expectation, i.e. $\mathrm{Var}(\mathcal{N}_\Omega)\sim \E\mathcal{N}_\Omega\sim N$. Thus, the natural question: Is the size of $\mathrm{Var}(\mathcal{N}_\Omega)$ for i.i.d. matrices still proportional to those of the mean? We answer it negatively; 
the number variance is much smaller than the expected size of $\mathcal{N}_\Omega$:

\begin{theorem}[Informal statement]\label{theo:informal}
Let $X$ be an i.i.d. matrix and consider any nice domain $\Omega\subset\mathbf{D}$. Then, there exists a $q>0$ such that
\begin{equation}
\label{eq:hyp}
\mathrm{Var}(\mathcal{N}_\Omega)\le N^{1-q}.
\end{equation} 
\end{theorem}

Our proof explicitly gives $q=1/40$ in the complex case and $q=1/106$ in the real case, see Theorems~\ref{theo:main1} and~\ref{theo:main2} below. Theorem~\ref{theo:informal} establishes a connection between the point process of the eigenvalues of general i.i.d. matrices with \emph{hyperuniformity}, a key concept in condensed matter physics for classifying crystals, quasicrystals, and exotic states of matter. In \cite[Section 1]{Torquato18}, Torquato defines the concept of hyperuniformity as:

\smallskip

\emph{(...) the number variance of particles within a spherical observation window of radius R grows more
slowly than the window volume in the large–R limit.}

\smallskip

For a system of independent particles the number variance is proportional to the volume $\mathrm{Var}(\mathcal{N}_\Omega)\sim N$, hence it is not hyperuniform. Hyperuniformity is typically a signature of strong correlations at large distances in the system that reduce fluctuations.  In random matrix theory, one well-known manifestation of such correlations is the \emph{eigenvalue rigidity}, proven in many Hermitian random matrix models, which asserts that each eigenvalue fluctuates on a  scale only slightly larger (by an $N^\epsilon$ factor) than the typical distance between neighbouring eigenvalues.  In particular, this trivially implies hyperuniformity in the sense of Torquato. In fact, for Hermitian models with eigenvalue rigidity, the variance of $\mathcal{N}_I$, the number of eigenvalues within an interval $I\subset R$, is bounded by $N^\epsilon$, independently of $I$. Note that  here the concept of rigidity uses that the Hermitian spectrum is one dimensional hence the eigenvalues can be ordered. For the same reason, the proof of eigenvalue rigidity easily follows from optimal concentration properties of the resolvent (\emph{optimal local laws}). In the two-dimensional setting of non-Hermitian random matrices, there is no direct concept of eigenvalue rigidity, but hyperuniformity can be interpreted as a higher dimensional version of the same strong correlations that cause rigidity. However, in this case hyperuniformity does not simply follow from optimal local laws (which are well understood in this setting as well); substantial new inputs are needed.

The word hyperuniformity was coined by Torquato and Stillinger \cite{torquato2003local} in the early 2000, analogous concepts appeared in the physics of Coulomb gases much earlier \cite{jancovici1993large, lebowitz1983charge, levesque2000charge, martin1980charge, martin1988sum}. In recent years this phenomenon attracted lots of interest in a variety of models in mathematics and physics (here for concreteness we focus only on two-dimensional objects), including averaged and perturbed lattices \cite{gacs1975problem}, zeroes of random polynomials with i.i.d. Gaussian coefficients \cite{forrester1999exact, shiffman2008number}, invariant point processes in a plane \cite{Weierstr_zeta23}, certain fermionic systems \cite{calabrese2015random, torquato2008point}, Berezin-Toeplitz operators in connection with the Quantum Hall Effect \cite{charles2020entanglement}. A rigorous proof of hyperuniformity in system with two-body interactions is  notoriously difficult. In the prominent model of the two-dimensional Coulomb gas (also known as the two-dimensional one component plasma) hyperuniformity was proved only very recently by Lebl\'e \cite{leble2021two}. This result showed that the number variance  is bounded by $N/(\log N)^a$, for some small $a>0$; so, in particular, it grows more slowly than the volume $N$.

Our result \eqref{eq:hyp} shows that the point process given by the eigenvalues of general non-Hermitian i.i.d. matrices is hyperuniform. It can thus be thought as a random matrix counterpart of \cite{leble2021two}, with a  much stronger control on the number variance: we prove a polynomial  factor $N^q$, while in the Coulomb gas model only a logarithmic factor $(\log N)^a$ was obtained. Previous to our result, in the random matrix setting, hyperuniformity was known only for integrable ensembles such as complex Ginibre matrices \cite{costin1995gaussian, lin2024nonlocal, norm_rep, soshnikov2000gaussian} (see also the related works \cite{abreu2023entanglement, calabrese2015random, forrester2014local, rider2004deviations}, and the \emph{four moment matching} result \cite{TaoVu15}), for real Ginibre matrices away from the real axis \cite{goel2024central} (see also \cite{akemann2023universality} for the symplectic Ginibre ensemble), for elliptic Gaussian matrices and normal matrix models \cite{akemann2024fluctuations} (the latter result was later improved in \cite{fluct_norm_mat26}). In particular, our result is new even for real Ginibre matrices for domains $\Omega$ intersecting the real axis.  

Before explaining our strategy to prove \eqref{eq:hyp}, we point out that in certain systems a stronger version of hyperuniformity is expected. In the so-called “class I hyperuniform” systems, the number variance grows as the perimeter\footnote{We point out that initially this was chosen as the definition of hyperuniformity in \cite{torquato2003local}.}, which is the slowest possible growth (see e.g. \cite{beck1987irregularities}), and not merely slower than the volume. In particular, for point process in two-dimensions, “class I hyperuniform” means that $\mathrm{Var}(\mathcal{N}_\Omega)\sim \sqrt{N}$. This strong version of hyperuniformity has been proven  for Ginibre ensembles using explicit formulas. The two-dimensional Coulomb gas is  also expected to be “class I hyperuniform”, but the result \cite{leble2021two} is still far from catching this. Similarly to these examples, it is expected that the eigenvalues $\sigma_i$ for general i.i.d. matrices form also a “class I hyperuniform” system, i.e. that $\mathrm{Var}(\mathcal{N}_\Omega)\sim \sqrt{N}$. Our result gives the first (and effective) proof of hyperuniformity for the i.i.d. random matrices, however, \eqref{eq:hyp} is also very far from the (expected) optimal $\sqrt{N}$ bound. We leave this for future work. However, we mention that our result also covers certain mesoscopic $N$-dependent domains $\Omega_N\subset \mathbf{D}$ (see Assumption~\ref{ass:Omega} below), showing  that hyperuniformity occurs on mesoscopic scales as well.

To study the variance of $\mathcal{N}_\Omega$, we view it as a special case of (centered) linear statistics\footnote{The linear statistics for $f(\cdot)=\bm1(\cdot\in\Omega)$ are often called \emph{counting statistics}.}
\begin{equation}
L_N(f):=\sum_{i=1}^N f(\sigma_i)-\E \sum_{i=1}^N f(\sigma_i),
\end{equation}
when $f(\cdot)=\bm1(\cdot\in\Omega)$. In particular, $\mathcal{N}_\Omega-\E\mathcal{N}_\Omega=L_N(f)$, for this choice of $f$. We thus need to give a bound on the variance of $L_N(f)$ for $f(\cdot)=\bm1(\cdot\in\Omega)$. 
The behavior of the linear statistics $L_N(f)$ is very well understood for smooth (e.g. $C^2$) test functions \cite{macroCLT_complex, macroCLT_real, Cipolloni_meso, kopel2015linear, rider2006gaussian, rider2007noise} (see also \cite{nguyen2014random}), however extending this study to less regular $f$'s creates substantial difficulties. In particular, the loss of regularity of $f$ does not only present new technical complications, but it substantially changes the answer. In fact, it is expected that the size of the number variance is much bigger than $\mathrm{Var}(L_N(f))$ for smooth $f$. More precisely, as mentioned above, it is expected that $\mathrm{Var}(\mathcal{N}_\Omega)\sim \sqrt{N}$, while it is proven in \cite[Theorem~2.2]{macroCLT_complex} that for smooth test functions the linear statistics have an order one variance
\begin{equation}
\label{eq:varform}
\mathrm{Var}(L_N(f))=V_f+o(1), \qquad\quad V_f:=\frac{1}{4\pi^2}\int_\mathbf{D}|\nabla f(z)|^2\,\dif^2 z+\kappa_4 \left|\frac{1}{\pi}\int_\mathbf{D} f(z)\,\dif^2 z-\frac{1}{2\pi}\int_0^{2\pi}f(e^{\ii\theta})\,\dif \theta\right|^2,
\end{equation}
where $\kappa_4$ is the fourth cumulant of the entries of $X$.
Note that \eqref{eq:varform} is consistent with the fact that if we regularized $\bm1(\cdot\in\Omega)$ on a scale $\sim N^{-1/2}$ (which is the fluctuation scale of the eigenvalues of $X$) and naively plug this regularized $f$ in the formula \eqref{eq:varform} for $V_f$, we would indeed get a size of order $V_f \sim \sqrt{N}$.  
To regularize the indicator function we rely on the Portmanteau principle (see Assumption~\ref{ass:Omega} below), which for the variance states 
\begin{equation}
\label{eq:port}
\mathrm{Var}(\mathcal{N}_\Omega)\lesssim \mathrm{Var}(L_N(f_+))+\mathrm{Var}(L_N(f_-))+\big(\E L_N(f_+)-\E L_N(f_-)\big)^2.
\end{equation}
Here $f_\pm$ are the lower and upper envelope functions with $f_-\le \bm1_\Omega\le f_+$
and they are chosen to be smooth on a scale $N^{-a}$, for a certain $a>1/2$ which we will then optimize.
 While it seems natural to smooth out  $\bm1(\cdot\in\Omega)$ on a scale $\sim N^{-1/2}$ (i.e. at the fluctuation scale of the eigenvalues), in practice we will need to choose a smaller scale $N^{-a}$, with some $a>1/2$, to make sure that the term $(\E L_N(f_+)-\E L_N(f_-))^2$ is negligible compared to the upper bound on the other two terms in the rhs. of \eqref{eq:port}. In particular, this additional term represents the discrepancy in the number of eigenvalues in the support of the functions $\bm1_\Omega$ and $f_\pm$.

Choosing $a$ sufficiently large, 
the bound \eqref{eq:port} enables us to reduce \eqref{eq:hyp} to study the linear statistics of $f_\pm$, which are smooth, even if only on the very small scale $N^{-a}$. To study these eigenvalue statistics we rely on Girko's formula \eqref{eq:girko} below, a backbone in the study of non-Hermitian spectral statistics. The key feature of this formula is that it relates the eigenvalues of any non-Hermitian matrix $X$ with those of a family of Hermitized matrices
\begin{equation}
\label{eq:herm}
H^z:=\left(\begin{matrix}
0 & X-z \\
(X-z)^* & 0
\end{matrix}\right), \qquad\quad z\in\C.
\end{equation}
Then, Girko's formula reads as
\begin{equation}
\label{eq:girko}
\sum_{i=1}^N f(\sigma_i)=\frac{\ii}{4\pi}\int_\C\Delta f(z)\int_0^\infty\mathrm{Tr} G^z(\ii\eta)\, \dif \eta\dif^2 z, \qquad\quad  G^z(\ii\eta):=(H^z-\ii\eta)^{-1},
\end{equation}
for any smooth function $f$. From the appearance of $\Delta f$ in \eqref{eq:girko} it is clear that the smoothness of $f$ is extremely relevant to control the eigenvalue statistics. For our approximating functions $f=f_\pm$ from above we
will use the bound $\lVert \Delta f\rVert_\infty\le N^{2a}$, which grows very quickly as $a>1/2$.

To deduce from \eqref{eq:port} and \eqref{eq:girko} that the eigenvalue process of $X$ is hyperuniform in the sense that  the lhs. of \eqref{eq:hyp} has an upper bound of order $N^{1-\xi}$ for some tiny implicit constant $\xi>0$, it suffices to show with the same implicit level of precision that the spectra of $H^{z_1}$ and $H^{z_2}$ decorrelate on all scales $|z_1-z_2|\gg N^{-1/2}$. Such weak inputs are already available from \cite{macroCLT_complex}. However, along this approach, not only any explicit value of $\xi$ is practically untraceable, but also  the actual proof theoretically would yield an exponent $\xi$ which depends on the model parameters in Assumption~\ref{ass:chi} below, such as the constants in the upper bounds on the moments of the single-entry distribution of $X$. In contrast, in this paper we get the universal quantitative exponent in \eqref{eq:hyp}, which requires substantial new inputs. In particular, we need to quantify all spectral decorrelation rates in terms of powers of $|z_1-z_2|$, which is achieved by performing the analysis that we now explain. For more details see the comment below Proposition~\ref{prop:Var_main}.

To analyze the rhs. of \eqref{eq:girko} it is convenient to split the $\eta$-integral into three regimes: \emph{sub-microscopic regime} ($\eta\ll 1/N$), \emph{microscopic regime} ($\eta\sim 1/N$), and \emph{macro-mesocopic regime} ($\eta\gg 1/N$). They require different methodological approaches. Inspired by \cite{macroCLT_complex}, we can control the size of the sub-microscopic and microscopic regimes using left-tail bounds on the smallest singular value of $X-z$ from \cite{bordenave2012around, tao2010smooth} (see also \cite{cipolloni2020optimal, cipolloni2022condition}) and the Dyson Brownian motion (DBM) from \cite{bourgade2024fluctuations, macroCLT_complex}, respectively. The main novelty of this work lies in the analysis of the regime $\eta\gg 1/N$, where substantially higher precision, compared to previous works, is required.

In order to compute the variance of $L_N(f)$ using \eqref{eq:girko}, we naturally need to compute\footnote{For a matrix $A\in\C^{d\times d}$ we use the short-hand notation $\langle A \rangle=\frac{1}{d}\mathrm{Tr}[A]$.} $\mathrm{Cov}(\langle G_1\rangle, \langle G_2\rangle)$, where $G_i:=G^{z_i}(\ii\eta_i)$. The first natural approach to compute $\mathrm{Cov}(\langle G_1\rangle, \langle G_2\rangle)$ is to rely on similar computations to \cite[Section 6]{macroCLT_complex}. However, the precision of \cite[Eq. (6.31)]{macroCLT_complex} would give (here $\eta_*:=\eta_1\wedge\eta_2$)
\begin{equation}
\label{eq:oldcov}
\mathrm{Cov}(\langle G_1\rangle, \langle G_2\rangle)=\frac{\mathrm{Main \, Term}}{N^2}+\mathcal{O}\left(\frac{1}{\sqrt{N\eta_*}(N\eta_1\eta_2)^2}\right).
\end{equation}
The key drawback of \eqref{eq:oldcov} is that while the "Main Term" in \eqref{eq:oldcov} behaves like $|z_1-z_2|^{-4}$ (see e.g. \cite[Eq. (4.24)]{macroCLT_complex}), this decorrelation
decay in $|z_1-z_2|$ is not reflected in the error term. 
The decorrelation phenomenon is essential to prove a good estimate since the variance of $L_N(f)$
based upon \eqref{eq:girko} involves a double $\dif^2 z_1\dif^2 z_2$ integral of 
$\mathrm{Cov}(\langle G_1\rangle, \langle G_2\rangle)$.

The main technical result of this work is an improved covariance control, which in the Ginibre case gives
\begin{equation}
\label{eq:newcov}
\mathrm{Cov}(\langle G_1\rangle, \langle G_2\rangle)=\frac{\mathrm{Main \, Term}}{N^2}+\mathcal{O}\left(\frac{1}{N[|z_1-z_2|^2+\eta_1+\eta_2](N\eta_1\eta_2)^2}\right).
\end{equation}
Note that the error term in \eqref{eq:newcov} does not only improve \eqref{eq:oldcov} in terms of $|z_1-z_2|$, but also the $N$-power is better (e.g. $N^{-5/2}$ vs. $N^{-3}$). In addition, we prove \eqref{eq:newcov} not only for spectral parameters on the imaginary axis but for all spectral parameters $w_1, w_2$ in the bulk spectrum
of $H^{z_1}, H^{z_2}$ (in which case an additional regularizing term depending on the distance between $w_1$ and $w_2$ appears). In contrast, the previous results in \cite{macroCLT_complex} covered only spectral parameters along the imaginary axis. 

To prove \eqref{eq:newcov}, we use a \emph{chaos expansion} strategy, which consists in performing iterative cumulant expansions in $\mathrm{Cov}(\langle G_1\rangle,\langle G_2\rangle)$. More precisely, we start with a trivial bound of order $(N\eta_1\eta_2)^{-1}$ on the error term in the rhs. of \eqref{eq:newcov}, which simply follows by the single resolvent local law (see e.g.~\eqref{eq:1G_ll_av} below). Then, at each step of the chaos expansion we improve this bound by a (small) factor of $(N\eta_*)^{-1}$. In the mesoscopic regime $\eta_*\ge N^{-1+\epsilon}$, this factor is much smaller than one, however it gains only $N^{-\epsilon}$ for a small $\epsilon>0$. Thus, the factor $(N[|z_1-z_2|^2+\eta_1+\eta_2])^{-1}$ in the rhs. of \eqref{eq:newcov} cannot in general be obtained in one step of this iteration. This forces us to perform iterative expansions, each expansion amounting to a further gain $ N^{-\epsilon}$, until we reach the desired precision in \eqref{eq:newcov}.

This gradual improvement of the bound on the error term to reach the one in the rhs. of \eqref{eq:newcov} involves understanding the size of more complicated covariances, which generalize the lhs. of \eqref{eq:newcov}. For example, the first step generates 
\begin{equation}
\mathrm{Cov}\left(\langle G_1-M_1\rangle\langle (G_1-M_1)M_1^2\rangle,\langle G_2\rangle\right),
\label{eq:Cov_intro}
\end{equation}
where $M_1$ is a certain deterministic approximation to $G_1$ introduced later in \eqref{eq:M}. These covariances are further subjected to the iterative expansions, giving rise to an entire hierarchy of more and more involved covariances. However, they still have a specific structure that plays an important role in the analysis. Though the number of steps in the chaos expansion may be arbitrarily large, we will need to iterate the expansion roughly $1/\epsilon$-times, where $\epsilon$ is the small parameter in the constraint $\eta_1\wedge\eta_2\ge N^{-1+\epsilon}$. However, at the end, the hierarchy of covariances needs to be truncated, meaning that we estimate these quantities arising after the final step of the chaos expansion simply by size, without re-expanding them further. The loss in this size estimate is compensated by the many $1/(N\eta_*)$ factors accumulated along the iteration. This size bound is given by the \emph{multi-resolvent local laws}: concentration bounds for the products of resolvents and deterministic matrices sandwiched in between (see Section~\ref{sec:local_laws} below).

A related strategy was used in \cite{he2020mesoscopic} to compute the covariance $\mathrm{Cov}(\langle G_1\rangle, \langle G_2\rangle)$ for resolvents $G_j:=(W-~\!w_j)^{-1}$, $j=1,2$, of a Wigner matrix $W$. Our analysis of the Hermitized matrix $H^z$ carries three novelties, besides several technical complications.  First, most importantly,  $G_1$ and $G_2$ in \cite{he2020mesoscopic}  have the same spectral resolution, since they are resolvents of the same matrix $W$. This allows for a substantial  simplification of the structure of the hierarchy of covariances using the resolvent identity $G_1G_2=(G_1-G_2)/(w_1-w_2)$, which also automatically gives the decorrelation factor $w_1-w_2$. Since there is no analogous resolvent identity for products of resolvents of $H^{z_1}$ and $H^{z_2}$ with $z_1\neq z_2$, hence the decorrelation factor $(N[|z_1-z_2|^2+\eta_1+\eta_2])^{-1}$ does not appear automatically. Instead, it will be obtained from the very refined bound in a new two-resolvent local law for quantities of the form $\langle G_1B_1G_2B_2\rangle$ with deterministic matrices $B_1,B_2$ (see Proposition~\ref{prop:2G_av} below). Second, the precision in the estimate of the error term in the chaos expansion is governed by the inverse of the so-called \emph{one-body stability operator} $\mathcal{B}_{11}$ introduced later below \eqref{eq:def_B12}. In \cite{he2020mesoscopic} this operator has a bounded inverse, while in our setup it is unstable and the norm of $\mathcal{B}_{11}^{-1}$ blows up as $\eta_1^{-1}$ for small $\eta_1$'s. Thus, showing that this deterioration of the stability bound does not affect the final result \eqref{eq:newcov} requires a separate structural analysis of the hierarchy of covariances. Third, $M_1$ appearing in \eqref{eq:Cov_intro} is a multiple of identity in \cite{he2020mesoscopic}, while now it is an $(2N)\times (2N)$ matrix with a $2\times 2$ block-constant structure. In particular, quantities arising along the chaos expansion contain traces of the form $\langle G_1B_1G_1\cdots G_1B_k\rangle$ for all $k\in\N$ and deterministic matrices $B_j$, $j\in[k]$, with $2\times 2$ block-constant structure, while in \cite{he2020mesoscopic} such traces were automatically simplified to $\langle G_1^k\rangle$.

As a standard consequence of the new averaged two-resolvent local law mentioned in the previous paragraph, we get the following optimal (in $N$) bound on the overlap of singular vectors of $X-z_1$ and $X-z_2$ in the bulk regime, i.e. for $|z_1|, |z_2|\le 1-\delta$,
\begin{equation}
\left\vert \langle \bm{u}_i^{z_1},\bm{u}_j^{z_2}\rangle\right\vert^2 + \left\vert \langle \bm{v}_i^{z_1},\bm{v}_j^{z_2}\rangle\right\vert^2 \lesssim \frac{1}{N}\frac{1}{|z_1-z_2|^2+N^{-1}|i-j|},
\label{eq:overlap_intro}
\end{equation}
with very high probability up to a factor $N^\xi$, for an arbitrary small $\xi>0$. For a precise statement see Proposition~\ref{prop:overlap}. Here $\bm{u}_i^z$ and $\bm{v}_i^z$ are the left and right singular vectors associated to the $i$-th singular value of $X-z$. The main feature of \eqref{eq:overlap_intro} is the identification of the decorrelation in $|i-j|$ for all $i,j$, on top of the $|z_1-z_2|^2$ correlation decay that was known before but only for $|i|,|j|\lesssim N^\xi$, see \cite[Corollary~3.6]{nonHermdecay}. A bound analogous to \eqref{eq:overlap_bound} was obtained in \cite[Theorem~2.6]{eigenv_decorr} for the overlap of eigenvectors of $W+D_1$ and $W+D_2$, where $W$ is a Wigner matrix and $D_j=D_j^*$ are deterministic deformations for $j=1,2$. In~\cite{eigenv_decorr}, $\langle (D_1-D_2)^2\rangle$ played the role of $|z_1-z_2|^2$ in \eqref{eq:overlap_intro}, while further improvements beyond this quadratic bound were encoded in a complicated control parameter coined as the \emph{linear term}, see \cite[Eq.(2,17)]{eigenv_decorr}. In the current paper we encounter a similar control parameter (see the second line of \eqref{eq:def_UV} later). Owing to the more explicit structure of our model compared to \cite{eigenv_decorr}, we manage to analyze precisely this control parameter and to extract from it the $|i-j|$ contribution.
  
We note that the term $|z_1-z_2|^2$ in \eqref{eq:overlap_bound} plays an essential role in our analysis of the critical regime $\eta\sim N^{-1}$ in the Girko's formula \eqref{eq:girko}, and allows us to quantify the error terms coming from the Dyson Brownian motion technique.  Meanwhile, the $|i-j|$ decay in the rhs. of \eqref{eq:overlap_intro} is a result of an independent interest which is not used in the proof of Theorem~\ref{theo:informal}.

\subsection*{Notations and conventions}
We set $[k] := \{1, ... , k\}$ for a positive integer $k \in \N$ and $\langle A \rangle := d^{-1} \mathrm{Tr}(A)$, $d \in \N$, for the normalized trace of a $d \times d$ matrix $A$. Additionally, we denote its transpose by $A^\mt$. The sets of real and complex numbers are denoted by $\R$ and $\C$, respectively, while $\mathbf{D}\subset \C$ stands for the open unit disk. For a finite set $S$ we denote the number of elements in $S$ by $\# S$. For positive quantities $f, g$ we write $f \lesssim g$, $f \gtrsim g$, to denote that $f \le C g$ and $f \ge c g$, respectively, for some $N$-independent constants $c, C > 0$ that depend only on the basic control parameters of the model in Assumptions~\ref{ass:chi} and \ref{ass:Omega} below. In informal explanations we sometimes use the notation $f\ll g$, which indicates that $f$ is "much smaller" than $g$.

We denote vectors by bold-faced lower case Roman letters $\boldsymbol{x}, \boldsymbol{y} \in \C^{d}$, for some $d \in \N$. Moreover, for vectors $\boldsymbol{x}, \boldsymbol{y} \in \C^{d}$ and a matrix $A\in\C^{d\times d}$ we define
 \begin{equation*}
	\langle \boldsymbol{x}, \boldsymbol{y} \rangle := \sum_{i=1}^d \bar{x}_i y_i\,, 
	\qquad\quad A_{\boldsymbol{x} \boldsymbol{y}} := \langle \boldsymbol{x}, A \boldsymbol{y} \rangle\,. 
\end{equation*}
Matrix entries are indexed by lower case Roman letters $a, b, c , ... ,i,j,k,... $ from the beginning or the middle of the alphabet.

The covariance of two complex-valued random variables $\zeta_1$ and $\zeta_2$ is denoted by
\begin{equation*}
\Cov(\zeta_1,\zeta_2):=\E\left[ (\zeta_1-\E\zeta_1)(\overline{\zeta_2}-\E\overline{\zeta_2})\right].
\end{equation*} 
We further denote the variance of a complex-valued random variable $\zeta$ by $\Var[\zeta]:=\Cov(\zeta,\zeta)$. For a pair of real or complex-valued stochastic processes $X=X(t)$ and $Y=Y(t)$, we denote their covariation process by~$[X,Y]_t$.

Finally, we will use the concept  \emph{with very high probability},  meaning that for any fixed $D > 0$, the probability of an $N$-dependent event is bigger than $1 - N^{-D}$ for all $N \ge N_0(D)$. We will use the convention that $\xi > 0$ denotes an arbitrarily small positive exponent, independent of $N$. Moreover, we introduce the common notion of \emph{stochastic domination} (see, e.g., \cite{loc_sc_gen}): For two families
\begin{equation*}
	X = \left(X^{(N)}(u) \mid N \in \N, u \in U^{(N)}\right) \quad \text{and} \quad Y = \left(Y^{(N)}(u) \mid N \in \N, u \in U^{(N)}\right)
\end{equation*}
of non-negative random variables indexed by $N$, and possibly a parameter $u$, we say that $X$ is stochastically dominated by $Y$, if for all $\epsilon, D >0$ we have 
\begin{equation*}
	\sup_{u \in U^{(N)}} \mathbf{P} \left[X^{(N)}(u) > N^\epsilon Y^{(N)}(u)\right] \le N^{-D}
\end{equation*}
for large enough $N \ge N_0(\epsilon, D)$. In this case we write $X \prec Y$. If for some complex family of random variables we have $\vert X \vert \prec Y$, we also write $X = \mathcal{O}_\prec(Y)$. 

\subsection*{Acknowledgment} We thank Leslie Molag for comments on the existing literature and for the references \cite{norm_rep,fluct_norm_mat26, Weierstr_zeta23}.

\section{Main results}

We consider \emph{real} or \emph{complex i.i.d. matrices $X$}, i.e. $N\times N$ matrices whose entries $x_{ab}\stackrel{d}{=}N^{-1/2}\chi$ are independent and identically distributed. We impose the following assumptions on the (possibly $N$-dependent) complex-valued random variable $\chi$.
\begin{assumption}\label{ass:chi} \emph{(i)} The random variable $\chi$ satisfies $\E \chi=0$, $\E |\chi|^2=1$. In addition, in the complex case we also assume $\E\chi^2=0$. We further assume the existence of high moments, i.e. that there exist constants $C_p>0$, for any $p\in \N$, such that  
\begin{equation}
\E |\chi|^p\le C_p.
\label{eq:momass}
\end{equation}

\noindent\emph{(ii)} There exist a (small) $\mathfrak{a}>0$ and a (large) $\mathfrak{b}>0$ such that the probability density $\rho_\chi$ of $\chi$
satisfies
\begin{equation}
\rho_\chi \in L^{1+\mathfrak{a}}(\C),\quad \|\rho_\chi\|_{1+\mathfrak{a}}\le N^\mathfrak{b}.
\label{eq:reg_ass}
\end{equation}
\end{assumption}

\begin{remark} We use the mild regularity assumption \eqref{eq:reg_ass} only to control the $\eta\le N^{-100}$ regime in the Girko's formula \eqref{eq:girko}. This assumption can easily be removed by standard methods as introduced in \cite{TaoVu15} (see also \cite[Remark~2.2]{rightmost_ev}), we omit the details for brevity and to keep the presentation simpler. 
\end{remark}

Denote the eigenvalues of $X$ by $\{\sigma_i\}_{i=1}^N$. For any test function $g:\C\to\C$, the linear eigenvalue statistics of $X$ is given by
\begin{equation*}
\mathcal{L}_N(g):=\sum_{i=1}^N g(\sigma_i).
\end{equation*}
Typically, the smoother is the test function $g$ the simpler is the analysis of $\mathcal{L}_N(g)$ (see e.g. the paragraph below \eqref{eq:girko}). As explained in the introduction, to study the phenomenon of hyperuniformity we need to consider the test function $g$ being the indicator function of a  (possibly $N$-dependent) subset $\Omega=\Omega_N$ of the disk. However, all the existing results on the linear eigenvalue statistics assume at least that $g\in H^1(\C)$, i.e. do not cover our case. We consider domains $\Omega_N$ satisfying the following assumptions:

\begin{assumption}\label{ass:Omega} The domain $\Omega_N\subset\C$ is open and simply connected, with $C^2$ boundary. There exists $\delta>0$ such that $\Omega_N\subset (1-\delta)\mathbf{D}$ for all $N\in\N$. There further exists an exponent $\alpha\in [0,1/2)$, such that $\Omega_N=N^{-\alpha}\widetilde{\Omega}_N$, where 
\begin{equation}
{\rm{diam}}(\widetilde{\Omega}_N):=\sup \{|z_1-z_2|\,:\, z_1,z_2\in\widetilde{\Omega}_N\} \sim 1,
\label{eq:diam}
\end{equation}
and the curvature of the boundary $\partial\widetilde{\Omega}_N$ is bounded from above uniformly in $N$.
\end{assumption}

We are now ready to formulate our main result in the complex case.

\begin{theorem}[Complex case]
\label{theo:main1} Let $X$ be a complex $N\times N$ i.i.d. matrix satisfying Assumption \ref{ass:chi}. Fix $\alpha\in [0,1/2)$ and let $\Omega_N\subset\mathbf{D}$ be a domain satisfying Assumption \ref{ass:Omega} with exponent $\alpha$. Then, there exists $q_0>0$ such that for any fixed $\xi>0$ we have
\begin{equation}
\Var\left[ \#\{\sigma_i: \sigma_i\in \Omega_N\}\right] \le N^{1-2\alpha -q_0(1/2-\alpha)+\xi}
\label{eq:main1}
\end{equation}  
for sufficiently large $N$. Our proof gives\footnote{See Remark~\ref{rem:impr}--(v) for possible improvements on $q_0$.} $q_0=1/20$.
\end{theorem}

To state the analogue of Theorem~\ref{theo:main1} in the real case, we impose the following additional assumption on the domain $\Omega_N$.
\begin{assumption}\label{ass:Omega_real} In the set-up of Assumption~\ref{ass:Omega} we assume
that the angle between the tangent line to $\widetilde{\Omega}_N$ at these points and the real axis is at least $\varphi_0$, for some $N$-independent $\varphi_0>0$. We further assume that there exists a (small) $N$-independent constant $d_0>0$ such that $\partial\widetilde{\Omega}_N\cap \{z:|\Im z|\le d_0\}$ is either empty or consists only of curves intersecting the real axis. 
\end{assumption}

\begin{figure}[h]
\begin{center}
\includegraphics[height=5cm]{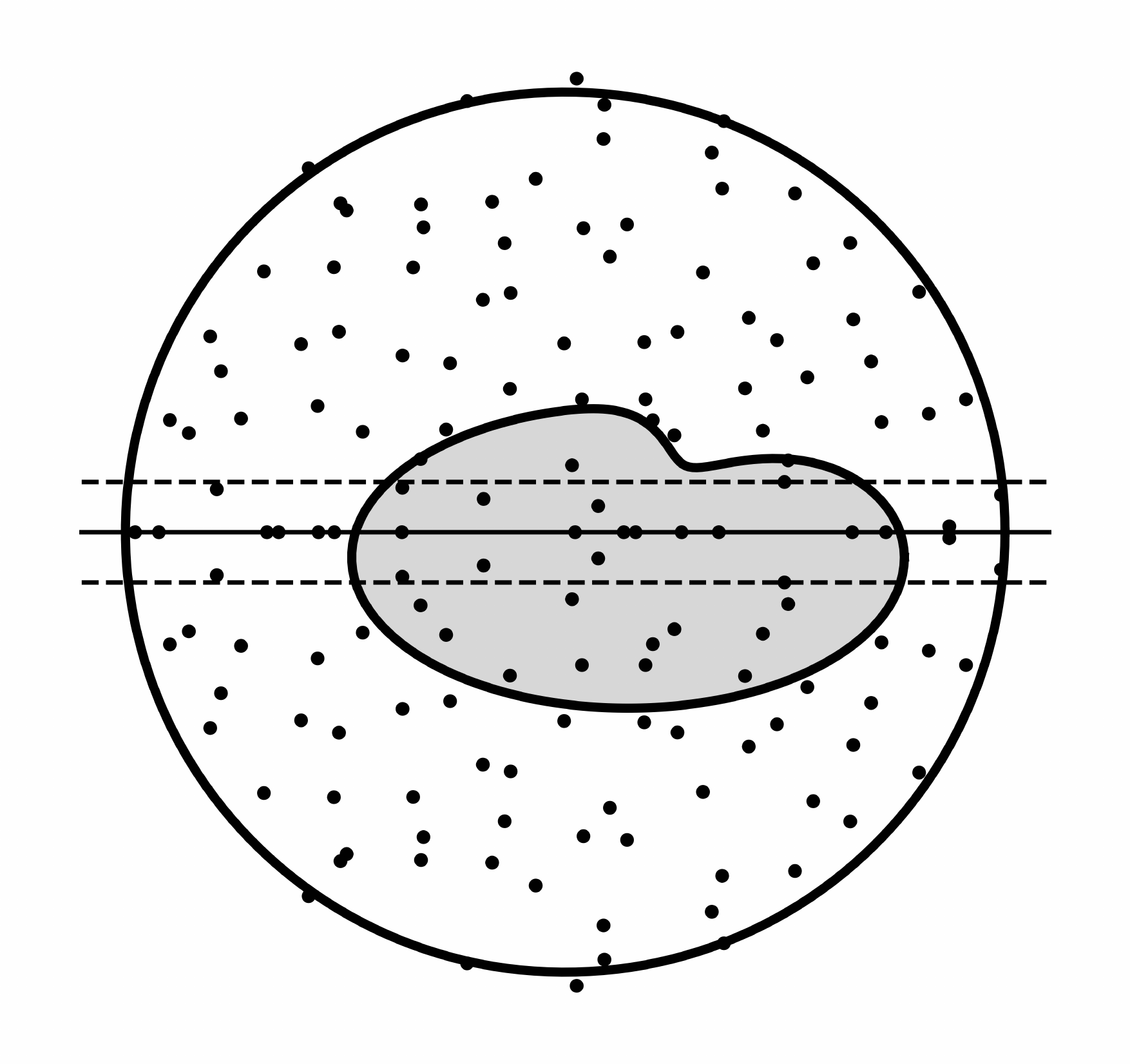}
\end{center}
\caption{In gray, we illustrated a macroscopic ($\alpha=0$) domain $\Omega_N$ in the case when it intersects the real line. The dashed lines correspond to $\Im z=\pm d_0$, where $d_0$ is given by Assumption~\ref{ass:Omega_real}. In this example, the intersection of $\partial\Omega_N$ and the set $|\Im z|\le d_0$ consists of two curves transversely intersecting the real axis. Finally, the scattered dots depict the eigenvalues of a large real i.i.d. matrix $X$.}
\label{fig:real_domain}
\end{figure}

In other words, Assumption~\ref{ass:Omega_real} guarantees that if $\Omega_N$ intersects the real axis, then this happens transversely, cf. Figure \ref{fig:real_domain}. It also prohibits $\partial\Omega_N$ to concentrate in a neighborhood of the real axis of a radius much smaller than the scale of~$\Omega_N$.

\begin{theorem}[Real case]
\label{theo:main2} 
Let $X$ be a real $N\times N$ i.i.d. matrix satisfying Assumption~\ref{ass:chi}. Fix $\alpha\in [0,1/2)$ and let $\Omega_N\subset\mathbf{D}$ be a domain satisfying Assumptions~\ref{ass:Omega} and~\ref{ass:Omega_real} with an exponent $\alpha$. Then, there exist $q_0, q_1>0$ such that for any fixed $\xi>0$ we have
\begin{equation}
\Var\left[ \#\{\sigma_i: \sigma_i\in \Omega_N\}\right] \le N^{1-2\alpha -q_0(1/2-\alpha)+\xi} + N^{1-2\alpha-q_1+\xi}
\label{eq:main2}
\end{equation}  
for sufficiently large $N$. Our proof gives $q_0=1/20$ and $q_1=1/106$.
\end{theorem}

\begin{remark}
\label{rem:impr}

We now discuss several possible extensions of the results in Theorems~\ref{theo:main1} and~\ref{theo:main2}.

\smallskip

\emph{(i) [Relaxation of Assumption \ref{ass:Omega}].} We stated Theorem \ref{theo:main1} under the assumption that the boundary of $\Omega_N$ is a $C^2$ curve. However, a careful examination of our proof shows that this assumption can be relaxed: one can allow $\partial\Omega_N$ to be a piecewise $C^2$ curve. In particular, $\Omega_N$ can be a rectangle whose edges have length of order $N^{-\alpha}$. For more details see Supplementary Section~\ref{app:Omega}. Similarly, Assumption \ref{ass:Omega} can be relaxed in Theorem~\ref{theo:main2}, but only away from the real axis. To avoid technicalities, we only mention that if $\Omega_N$ does not intersect the real axis and lies at a distance of order $N^{-\alpha}$ from it, then the discussion in Supplementary Section~\ref{app:Omega} applies also in the real case.

\smallskip

\emph{(ii)[Relaxation of Assumption \ref{ass:Omega_real}].} In the case when $\Omega_N$ does not intersect the real axis,  Assumption~\ref{ass:Omega_real} requires $\Omega_N$ to be at distance at least of order $N^{-\alpha}$ from the real axis. Our methods allow to treat a more general case when the distance from $\Omega_N$ to $\R$ is at least $N^{-\theta}$ for any fixed $0\le \theta<1/2$, and give a quantitative improvement of the volume law bound $N^{1-2\alpha}$ on the lhs. of \eqref{eq:main2}. However, we refrain from computing the exact exponents in this generalized setting.

\smallskip

\emph{(iii)[Domains overlapping the edge].} In Theorems~\ref{theo:main1} and \ref{theo:main2} we focused on domains $\Omega_N$ well inside the unit disk. However, an analogous proof would also work for domains overlapping the edge of the spectrum (i.e. the unit circle). The main new difficulty is that in this regime the local laws and the covariance expansion in Sections~\ref{sec:local_laws} and~\ref{sec:Cov} would need to be extended to the regime $|z|\approx 1$, when the limiting eigenvalue density of $H^z$ develops a cusp singularity. The methods of Section~\ref{sec:Cov} and Supplementary Section~\ref{sec:multiG_proof} are strong enough to cover this regime as well, however, at the price of more technically involved computations. For these reason we decide to omit the edge regime from the presentation.

\smallskip

\emph{(iv)[Bounds on higher moments].} In Theorems~\ref{theo:main1} and~\ref{theo:main2} we estimate the variance of the number of eigenvalues in $\Omega_N$, since this is the quantity determining whether a point process is hyperuniform or not. Our methods also allow us to show that the $p$-th moment, $p\in\N$, of this random variable is much smaller than $N^{(1/2-\alpha)p}$, and to quantify this smallness.

\smallskip

\emph{(v)[Possible improvements on $q_0, q_1$].} In Theorems~\ref{theo:main1} and~\ref{theo:main2} the exponents $q_0,q_1$ are not optimal. By performing a more careful analysis, one could increase the value of $q_0$ to $1/3$ and remove the term containing $q_1$ in the rhs. of \eqref{eq:main2}. However, this result will be still far from the optimal one corresponding to $q_0=1$. For more details see comments below Propositions~\ref{prop:tail_bound}, \ref{prop:EGG}, and in Section~\ref{sec:real_ideas}.
\end{remark}

\section{Proof strategy}
\label{sec:proofstrat}

To keep the presentation simpler, we now only consider the complex case. We will then explain the minor required differences for the real case in Section~\ref{sec:real_ideas}. In this section we discuss the main ideas of the proof of Theorem~\ref{theo:main1}. To study linear statistics with the test function being a sharp cut-off in Theorem~\ref{theo:main1}, we approximate it by smooth functions. To this end, we consider a smooth bump function~$\omega$ satisfying
\begin{equation}
\omega:\C\to [0,+\infty),\quad \omega\in C_c^\infty(\C),\quad  \int_\C \omega(z)\dif^2 z=1\quad\text{and}\quad \mathrm{supp}(\omega)=\overline{\mathbf{D}},
\label{eq:omega_conditions}
\end{equation}
and rescale it as
\begin{equation*}
\omega_{a,N}(z):=N^{2a}\omega(N^az),\quad \forall \, z\in \C,
\end{equation*} 
for some $a>0$. Note that $\omega_{a,N}$ integrates to 1. We then consider the convolution
\begin{equation}
f_{a,N}:=\phi_N*\omega_{a,N},
\label{eq:def_f}
\end{equation}
for a characteristic function $\phi_N$ of an appropriately chosen subset of the unit disk, which may slightly differ from~$\Omega_N$. Specifically, in the proof of Theorem \ref{theo:main1} we will construct this domain by adding a small neighborhood of $\partial\Omega_N$ to $\Omega_N$, or by removing it from $\Omega_N$ (see~\eqref{eq:Omega_pm} below). However, the exact construction is not important at this moment. By a Portmanteau-type argument (see Lemma~\ref{lem:Portm} below), approximating the characteristic function of $\Omega_N$ both from above and below by functions of the form \eqref{eq:def_f}, we reduce \eqref{eq:main1} to the analysis of expectation and variance of $\Lin_N(f_{a,N})$.

Linear eigenvalue statistics with a test function of the form $\omega(N^a(z-z_0))$, for some $z_0\in\C$, has been extensively studied in the \emph{macroscopic} regime ($a=0$)~\cite{macroCLT_complex,macroCLT_real}, and in the entire \emph{mesoscopic} regime ($0<a<1/2$)~\cite{Cipolloni_meso}. For instance, for $0<a<1/2$ it is known from \cite[Theorem~2.1]{Cipolloni_meso} that for $|z_0|\le 1-\tau$ with any fixed $\tau>0$, the variance of $\Lin_N(\omega(N^a(z-z_0)))$ is of order one. This is a consequence of the fact that the eigenvalues of i.i.d. matrices are “rigid" in the unit disk.

Our setting differs from that of \cite{Cipolloni_meso} in two fundamental aspects. First, to approximate the characteristic function of $\Omega_N$ by $f_{a,N}$ with precision sufficient for \eqref{eq:main1}, we are forced to take $a>1/2$, i.e. we need to consider sub-microscopically localized test functions. In this regime, $\Lin_N(\omega(N^a(z-z_0)))$ is governed by the finitely many eigenvalues of $X$ closest to $z_0$, so the self-averaging mechanism no longer applies. Second, the test function $f_{a,N}$ naturally lives on two scales: a mesoscopic scale $N^{-\alpha}\gg N^{-1/2}$, corresponding to the diameter of $\Omega_N$, and the sub-microscopic smoothing scale $N^{-a}\ll N^{-1/2}$. This two-scale structure requires to understand quite precisely the correlation between eigenvalues of $X$ near $z_1$ and $z_2$, for mesoscopically separated $z_1$, $z_2$.

The difficulties discussed above prevent us from computing the leading order behavior of the variance of $\Lin_N(f_{a,N})$, as it was done for smooth test functions in the mesoscopic scaling in aforementioned \cite[Theorem~2.1]{Cipolloni_meso} and in the macroscopic scaling (i.e. for an $N$-independent test function) in \cite[Theorem~2.2]{macroCLT_complex} and \cite[Theorem~2.2]{macroCLT_real}. Nevertheless, we establish an effective upper bound on the variance of $\Lin_N(f_{a,N})$ in Proposition~\ref{prop:Var_main}. We also compute the expectation of this linear statistics up to the leading order for an explicit range of $a$'s above the critical scale $1/2$ ( i.e. for sub-microscopic scales), see Proposition~\ref{prop:E_omega}. Together Propositions~\ref{prop:E_omega} and~\ref{prop:Var_main} serve as the main ingredients in the proof of Theorem~\ref{theo:main1}. Now we discuss these results in more detail, starting with the computation of the expectation. For simplicity of presentation, we consider the case when $\phi_N$ in \eqref{eq:def_f} is the characteristic function of $\Omega_N$. Since we impose only Assumption~\ref{ass:Omega} on $\Omega_N$, this restriction yields no loss of generality.

To compute $\E\Lin_N(f_{a,N})$, it is more convenient to work not directly with $f_{a,N}$ but with the test function
\begin{equation}
\omega_{a,N}^{(z_0)}(z):=\omega_{a,N}(z-z_0) = N^{2a}\omega(N^a(z-z_0)),\quad z,z_0\in\C,
\label{eq:omega_z0}
\end{equation}
which serves as a more basic object compared to $f_{a,N}$, as the expectation of $f_{a,N}$ can be easily recovered once the expectation of $\omega_{a,N}^{(z_0)}$ is known (see \eqref{eq:E_omega_to_f} below). In particular, \eqref{eq:omega_z0} does not carry any information about $\Omega_N$, which is though encoded in $f_{a,N}$. Once the expectation of the linear statistics with the test function \eqref{eq:omega_z0} has been computed, integrating it over $z_0\in\Omega_N$ yields $\E\Lin_N (f_{a,N})$, i.e.
\begin{equation}
\E \Lin_N (f_{a,N}) = \int_{\Omega_N} \E \Lin_N(\omega_{a,N}^{(z_0)})\dif^2 z_0.
\label{eq:E_omega_to_f}
\end{equation}

Although the assumptions imposed on~$\omega$ in \eqref{eq:omega_conditions} are essential for the proof of Theorem~\ref{theo:main1}, some of them are not needed for the computation of expectation of the linear statistics. In the following proposition we therefore relax~\eqref{eq:omega_conditions} and, in particular, allow~$\omega$ to be complex-valued.

\begin{proposition}\label{prop:E_omega} Let $X$ be a complex $N\times N$ i.i.d. matrix satisfying Assumption \ref{ass:chi}. Fix $\delta>0$, $a\in [1/2,1/2+\nu_0)$ with $\nu_0:=1/14$, and let $\omega\in C^2_c(\C)$ be a complex-valued function with integral 1. Uniformly in $z_0\in\mathbf{D}$ with $|z_0|\le 1-\delta$ it holds that
\begin{equation}
\E\Lin_N\big(\omega_{a,N}^{(z_0)}\big) = \frac{N}{\pi}\left(1+\mathcal{O}(N^{-c})\right)
\label{eq:E_omega}
\end{equation}
for some $c>0$, where $\omega_{a,N}^{(z_0)}$ is defined in \eqref{eq:omega_z0}.

Fix additionally $\alpha\in [0,1/2)$. Let $\Omega_N\subset \mathbf{D}$ be an $N$-dependent domain satisfying Assumption~\ref{ass:Omega} with exponent $\alpha$ and let $f_{a,N}$ be defined as in~\eqref{eq:def_f} with $\phi_N$ equal to the test function of $\Omega_N$. Then we have
\begin{equation}
\E\Lin_N\big(f_{a,N}\big) = \frac{N|\Omega_N|}{\pi}\left(1+\mathcal{O}(N^{-c})\right),
\label{eq:E_f}
\end{equation}
where $|\Omega_N|$ is the volume of $\Omega_N$.
\end{proposition}

We remark that the following analogue of \eqref{eq:E_omega} is known from \cite[Eq.(2.5)]{Cipolloni_meso} in the mesoscopic regime $0< a<1/2$ under slightly more general assumptions on $\omega$:
\begin{equation}
\E\Lin_N\big(\omega_{a,N}^{(z_0)}\big) = \frac{N}{\pi}\int \omega(z)\dif^2 z + \frac{N^{2a}}{8\pi} \int_\C \Delta \omega (z)\dif^2 z + \mathcal{O}(1+N^{2a-c}),
\label{eq:E_omega_meso}
\end{equation}
for some fixed $c>0$. In \eqref{eq:E_omega} there is no analogue of the second term in the rhs. of  \eqref{eq:E_omega_meso}, since the integral of $\Delta \omega$ over the entire complex plane vanishes due to the smoothness and compact support of $\omega$. Meanwhile, the first term in the rhs. of \eqref{eq:E_omega_meso} exactly corresponds to the leading term in the rhs. of \eqref{eq:E_omega}.

As it was discussed around \eqref{eq:E_omega_to_f}, \eqref{eq:E_f} immediately follows from~\eqref{eq:E_omega}, so \eqref{eq:E_f} is stated only for completeness. We also included $a=1/2$ into Proposition~\ref{prop:E_omega}, since the methods developed in this paper for the regime $a>1/2$ also apply to this case. However, Proposition~\ref{prop:E_omega} is not used for $a=1/2$ in the proof of Theorem~\ref{theo:main1}. The same remark applies to the following bound on the variance of $\Lin_N (f_{a,N})$ stated for all~$a\ge 1/2$. 

\begin{proposition}\label{prop:Var_main} Assume the set-up and conditions of Proposition \ref{prop:E_omega}, but instead of $a\in [1/2,1/2+\nu_0)$ consider any fixed $a\ge 1/2$. Then it holds that
\begin{equation}
\Var\left[ \Lin_N(f_{a,N})\right]\lesssim N^{2(a-\alpha) - 2q_0(1/2-\alpha)+\xi} 
\label{eq:Var_main}
\end{equation}
for any fixed $\xi>0$, where $q_0:=1/20$.
\end{proposition}

From the Girko's formula \eqref{eq:girko} we have
\begin{equation}
\Var\left[ \Lin_N(f_{a,N})\right] = \frac{N^2}{(4\pi)^2}\int_\C\int_\C \Delta f(z_1)\Delta f(z_2) \int_0^\infty\int_0^\infty \mathrm{Cov}\left(\langle G^{z_1}(\ii\eta_1)\rangle, \langle G^{z_2}(\ii\eta_2)\rangle\right)\dif\eta_1\dif\eta_2\dif^2 z_1\dif^2 z_2.
\label{eq:Var_explanation}
\end{equation}
The single-resolvent local law (see~\eqref{eq:av_ll_simple} below) gives an upper bound of order $(N^2\eta_1\eta_2)^{-1}$  on the covariance in the rhs. of \eqref{eq:Var_explanation}, up to an $N^\xi$ factor. This leads to the upper bound of order $N^{2(a-\alpha)+\xi}$ on the lhs. of \eqref{eq:Var_main}. Proposition~\ref{prop:Var_main} thus consists in a quantitative improvement of this elementary bound.

The remainder of this section is organized as follows. In Section~\ref{sec:proof_thm1} we prove Theorem~\ref{theo:main1} relying on Propositions~\ref{prop:E_omega} and~\ref{prop:Var_main}. In Section~\ref{sec:results_list} we then summarize  the technical ingredients required for the proofs of these results. Additionally, we explain the origin of the constants $\nu_0$ and $q_0$ appearing in the statements of Propositions~\ref{prop:E_omega} and~\ref{prop:Var_main}, respectively.

\subsection{Proof of Theorem \ref{theo:main1}}\label{sec:proof_thm1} From now on, we take $\omega$ to be a \emph{bump function}, meaning that $\omega$ satisfies \eqref{eq:omega_conditions}.
Fix an exponent $a\in[1/2,1/2+\nu_0)$ with $\nu_0=1/14$, at the end we will optimize the estimates over~$a$. Define
\begin{equation}
\Omega_N^-:=\left\lbrace z\in \Omega_N\,:\, \mathrm{d}(z,\partial{\Omega}_N)> N^{-a}\right\rbrace, \quad\Omega_N^+:=\left\lbrace z\in \C\,:\, \mathrm{d}(z,\Omega_N)< N^{-a}\right\rbrace,
\label{eq:Omega_pm}
\end{equation}
where $\mathrm{d}(\cdot,\cdot)$ is the Euclidean distance in $\C$. It is easy to see that the Assumption \ref{ass:Omega} imposed on $\Omega_N$ implies that $\Omega_N^-$ and $\Omega_N^+$ also satisfy\footnote{Here we use that a boundary of a tubular neighborhood of a $C^2$ curve $\partial\Omega_N$ is a union of two $C^2$ curves.} Assumption \ref{ass:Omega}. Recall the notation $\phi_N$ for the characteristic function of $\Omega_N$ and denote the characteristic functions of $\Omega_N^-$ and $\Omega_N^+$ by $\phi_N^-$ and $\phi_N^+$, respectively. Denote further 
\begin{equation}
f_{a,N}^{\pm}:=\phi_N^\pm *\omega_{a,N}.
\label{eq:def_f_pm}
\end{equation}
Since $\omega$ is supported on the unit disk, the construction \eqref{eq:Omega_pm} implies that
\begin{equation}
f_{a,N}^-(z) \le \phi_N(z)\le f_{a,N}^+(z),\quad \forall z\in\C.
\label{eq:f_phi_order}
\end{equation}

Now we formulate a simple statement allowing us to get an upper bound on the variance of $\mathcal{L}_N(\phi_N)$ in terms of the first two moments of $\mathcal{L}_N(f_{a,N}^\pm)$.

\begin{lemma}[Portmanteau principle]\label{lem:Portm} For any $N\in\N$ it holds that
\begin{equation}
\Var \left[\mathcal{L}_N(\phi_N)\right]\le 2\left(\Var \left[\mathcal{L}_N(f_{a,N}^+)\right] + \Var \left[\mathcal{L}_N(f_{a,N}^-)\right] + \left( \E\mathcal{L}_N(f_{a,N}^+-f_{a,N}^-)\right)^2\right).
\label{eq:Var_pm_bound}
\end{equation}
\end{lemma}

Hence, for an upper bound on $\Var \left[\mathcal{L}_N(\phi_N)\right]$ we need not only an upper bound on $\Var \big[\mathcal{L}_N(f_{a,N}^\pm)\big]$, but also an upper bound on $\big| \E\mathcal{L}_N(f_{a,N}^+-f_{a,N}^-)\big|$. We point out that using Proposition~\ref{prop:E_omega} we actually compute this expectation even up to the leading order, and not only give an upper bound for it. The proof of Lemma~\ref{lem:Portm} is elementary and relies only on the fact that $\Lin_N$ preserves the order~\eqref{eq:f_phi_order}, i.e.
\begin{equation*}
\Lin_N(f_{a,N}^-)\le \Lin_N(\phi_N)\le \Lin_N(f_{a,N}^+),
\end{equation*}
we omit the details for brevity.

We start with computing the expectation in the rhs. of \eqref{eq:Var_pm_bound} by the means of Proposition \ref{prop:E_omega}:
\begin{equation}
\begin{split}
\E \Lin_N(f_{a,N}^+-f_{a,N}^-) &= \E\Lin_N ((\phi_N^+-\phi_N^-)*\omega_{a,N}) = \int_{\C} \left(\phi_N^+(z)-\phi_N^-(z)\right)\E\Lin_N (\omega^{(z)}_{a,N})\dif z\\
&\le |\Omega_N^+\setminus\Omega_N^-| \frac{N}{\pi} \left(1+\mathcal{O}(N^{-c})\right) \lesssim N^{1-a-\alpha}.
\end{split}
\label{eq:E_dif_f}
\end{equation}
To go from the first to the second line we used that $\Omega_N\subset (1-\delta)\mathbf{D}$ by Assumption \ref{ass:Omega}, so that $\Omega_N^+\setminus\Omega_N^-\subset (1-\delta/2)\mathbf{D}$ for sufficiently large $N$, and thus Proposition \ref{prop:E_omega} is applicable to $\omega^{(z)}_{a,N}$. In the second line of \eqref{eq:E_dif_f} we additionally estimated
\begin{equation*}
|\Omega_N^+\setminus\Omega_N^-|\lesssim |\partial\Omega_N|N^{-a}\lesssim N^{-a-\alpha}.
\end{equation*}
Next we apply Proposition \ref{prop:Var_main} to $f_{a,N}^\pm$ and get
\begin{equation}
\Var \left[\mathcal{L}_N(f_{a,N}^+)\right] + \Var \left[\mathcal{L}_N(f_{a,N}^-)\right]\lesssim N^{2(a-\alpha)-2q_0(1/2-\alpha)+\xi}
\label{eq:Var_f_pm}
\end{equation}
for any small $\xi>0$. Combining \eqref{eq:Var_pm_bound}, \eqref{eq:E_dif_f} and \eqref{eq:Var_f_pm} we conclude that
\begin{equation}
\Var \left[\mathcal{L}_N(\phi_N)\right] \lesssim N^{2(a-\alpha)-2q_0(1/2-\alpha)+\xi} + N^{2-2a-2\alpha}.
\end{equation}
Finally, we optimize this bound over $a$ and take
\begin{equation}
a:=\frac{1}{2} + \frac{q_0}{2}\left(\frac{1}{2}-\alpha\right),
\label{eq:complex_a}
\end{equation}
thereby obtaining \eqref{eq:main1}. We can make such a choice, since $a\in [1/2,1/2+\nu_0)$ due to $\alpha\in[0,1/2)$, $\nu_0=1/14$ and $q_0=1/20$. This finishes the proof of Theorem \ref{theo:main1}. \qed

\subsection{Ideas behind the proofs of Propositions \ref{prop:E_omega} and \ref{prop:Var_main}}\label{sec:results_list} For any $z\in\C$ denote
\begin{equation}
H^z:=W-Z,\quad \text{where}\quad 
W=\begin{pmatrix}
0&X\\X^*&0
\end{pmatrix}\in \C^{(2N)\times(2N)},\quad
Z=\begin{pmatrix}
0&z\\\overline{z}&0
\end{pmatrix}\in \C^{(2N)\times(2N)}.
\label{eq:def_hermitization}
\end{equation}
Here $Z$ is a $2\times 2$ block-constant matrix with blocks of size $N\times N$. The matrix $H^z$ is called the \emph{Hermitization} of $X$, while $z$ is known the \emph{Hermitization parameter}. We denote the resolvent of $H^z$ by
\begin{equation}
G^z(w):=(H^z-w)^{-1}=(W-Z-w)^{-1},\quad w\in \C\setminus\R,
\label{eq:def_resolvent}
\end{equation}
and call $w$ a \emph{spectral parameter}.

We denote $f_N:=f_{a,N}$ and express the linear eigenvalue statistics of $X$ in terms of $G^z$ by the means of the Girko's Hermitization formula \cite{Girko84} in the regularized form given in \cite{TaoVu15}:
\begin{equation}
\Lin_N (f_N)=\frac{1}{4\pi} \int_\C \Delta f_N(z) \log\left\vert \det (H^z-\ii T)\right\vert \dif^2 z - \frac{N}{2\pi\ii} \int_\C \Delta f_N(z)\dif^2 z \int_0^T \langle G^z(\ii\eta)\rangle\dif \eta,
\label{eq:Girko}
\end{equation}
where we choose $T:=N^D$ for a large $N$-independent $D>0$. While the first term in the rhs. of \eqref{eq:Girko} will be negligible by a standard argument (see Lemma~\ref{sec:Girko_reduct} below), we focus on the second term and split the $\eta$-integration into four regimes as follows. Consider the scales $\eta_L<\eta_0<\eta_c$, where $\eta_L:=N^{-L}$ for a large $L>0$, $\eta_c:=N^{-1+\delta_c}$ for a small $\delta_c>0$, and $\eta_0<N^{-1}$ is an intermediate scale which we choose at the end to optimize the estimates. For any $0<\eta_l\le \eta_u$, denote
\begin{equation}
\begin{split}
I_{\eta_l}^{\eta_u}=I_{\eta_l}^{\eta_u}(f_N)&:= - \frac{N}{2\pi\ii} \int_\C \Delta f_N(z)\dif^2 z \int_{\eta_l}^{\eta_u} \langle G^z(\ii\eta)\rangle\dif \eta,\\
J_T=J_T(f_N)&:=\frac{1}{4\pi} \int_\C \Delta f_N(z) \log\left\vert \det (H^z-\ii T)\right\vert \dif^2 z.
\end{split}
\label{eq:def_IJ}
\end{equation}
Using this notation, \eqref{eq:Girko} can be written in the form
\begin{equation}
\Lin_N(f) = J_T + I_0^{\eta_L} + I_{\eta_L}^{\eta_0}+I_{\eta_0}^{\eta_c}+I_{\eta_c}^T.
\label{eq:Girko_I}
\end{equation}
Similarly to $J_T$, the term $I_{0}^{\eta_L}$ is negligible, for more details see Section \ref{sec:Girko_reduct}.

Now we explain, why we separated the $\eta$-regimes $I_{\eta_L}^{\eta_0}$, $I_{\eta_0}^{\eta_c}$ and $I_{\eta_c}^T$ in \eqref{eq:Girko_I} and why they exhibit different behavior. Owing to the definition of $H^z$ \eqref{eq:def_hermitization}, the spectrum of $H^z$ is symmetric with respect to zero (\emph{chiral symmetry}). Let $\{\lambda^z_i\}_{i=1}^N$ be the non-negative eigenvalues of $H^z$ labeled in increasing order. The rest of the eigenvalues are given by $\{\lambda^z_{-i}\}_{i=1}^N$ with $\lambda_{-i}^z:=-\lambda_i^z$. Denote the normalized eigenvector of $H^z$ corresponding to $\lambda_{\pm i}^z$ by $\bm{w}_{\pm i}^z\in\C^{2N}$. We have
\begin{equation}
\bm{w}_{\pm i}^z =\begin{pmatrix}
\pm \bm{u}_i^z\\\bm{v}_i^z
\end{pmatrix},
\label{eq:eigenvectors_symmetry}
\end{equation}
where $\bm{u}_i^z, \bm{v}_i^z\in\C^{N}$ are the left and right singular vectors of $X-z$ corresponding to the singular value~$|\lambda_i^z|$:
\begin{equation}
(X-z)\bm{v}_i^z = \lambda_i^z\bm{u}_i^z,\quad (\bm{u}_i^z)^*(X-z) = \lambda_i^z (\bm{v}_j^z)^*,\qquad \|\bm{u}_i^z\|=\|\bm{v}_i^z\|=1/2,\,\, i\in [N].
\label{eq:def_uv}
\end{equation}
Since the eigenvalues of $H^z$ come in pairs $(\lambda_i^z,-\lambda_i^z)$, we have
\begin{equation}
\langle G^z(\ii\eta)\rangle = \ii\langle \Im G^z(\ii\eta)\rangle = \ii\frac{1}{N}\sum_{i=1}^N \frac{\eta}{(\lambda_i^z)^2+\eta^2}.
\label{eq:G_to_ImG}
\end{equation}
We work in the bulk regime $|z|\le 1-\delta$, where the typical size of $\lambda_1^z$ is $N^{-1}$. Therefore, for $\eta\ll N^{-1}$ the leading contribution in the rhs. of \eqref{eq:G_to_ImG} comes from $\lambda_1^z$ and from the atypical regime when $\lambda_1^z$ is much smaller than $N^{-1}$. On the other hand, for $\eta\gg N^{-1}$ the non-trivial contribution comes from the first $N\eta$ eigenvalues, i.e. a collective behavior of an increasing with $N$ number of eigenvalues should be taken into account. This explains why we separate the regimes $[\eta_L,\eta_0]$ and $[\eta_c,T]$.  The remaining regime $\eta\in [\eta_0,\eta_c]$ contains the critical scale $\eta\sim N^{-1}$, where only the first few terms in the rhs. of \eqref{eq:G_to_ImG} have a non-trivial contribution, but the correlations between them cannot be neglected.  

In light of the discussion above, we refer to the $\eta$-regimes $[\eta_L,\eta_0]$, $[\eta_0,\eta_c]$ and $[\eta_c,T]$ as to the \emph{sub-microscopic}, \emph{critical} (or \emph{microscopic}) and \emph{mesoscopic} regimes, respectively. These terms are associated to the spectral resolution of $H^z$, rather than that of $X$, in contrast to the terminology used for the scaling of a test function. 

The idea of decomposing the Girko's formula into several $\eta$-regimes similarly to~\eqref{eq:Girko_I} is well-known in the literature. It has been employed, for example in the analysis of the linear eigenvalue statistics in macroscopic \cite{macroCLT_complex, macroCLT_real} and mesoscopic \cite{Cipolloni_meso} regimes, and in the proof of the non-Hermitian edge universality~\cite{nonherm_edge_univ}. A common feature of these works is that the leading contribution in the Girko's formula arises from the mesoscopic regime $I_{\eta_c}^T$, while the rest of the regimes are negligible. In the current paper we show that this behavior persists for the multiscale test function $f_{a,N}$ with $a\ge 1/2$ once the expectation of $\Lin_N(f_{a,N})$ is considered, as our proof of Proposition~\ref{prop:E_omega} shows. This conclusion, however, does not follow automatically, as establishing the negligibility of $\E I_{\eta_L}^{\eta_0}$ and $\E I_{\eta_0}^{\eta_c}$ requires a very precise analysis of these regimes.  Meanwhile, the analysis of $\Var\left[\Lin_N(f_{a,N})\right]$ is even much more delicate and our methods do not identify a single regime that yields the leading contribution. Instead, we derive an upper bound on the contribution from each of the regimes $I_{\eta_L}^{\eta_0}$, $I_{\eta_0}^{\eta_c}$ and $I_{\eta_c}^T$ to $\Var\left[\Lin_N(f_{a,N})\right]$, and none of these bounds follow from previous works. The final bound in Proposition~\ref{prop:Var_main} is then obtained by optimizing the sum of these contributions in~$\eta_0$. 

The analysis of $I_{\eta_c}^T$, $I_{\eta_L}^{\eta_0}$, and $I_{\eta_0}^{\eta_c}$ is non-trivial and relies on three different sets of tools, which we discuss in detail in Sections~\ref{sec:large_eta}, \ref{sec:small_eta}, and~\ref{sec:interm_eta}, respectively. In the proof of Proposition \ref{prop:Var_main}, the microscopic regime $I_{\eta_0}^{\eta_c}$ restricts the value of $q_0$ in Proposition \ref{prop:Var_main} to $1/20$. Our main contribution is in the mesoscopic regime $I_{\eta_c}^T$ which is discussed first for this reason. To analyze it, we compute the covariance of $\langle G^{z_1}(\ii\eta_1)\rangle$ and $\langle G^{z_2}(\ii\eta_2)\rangle$ for $\eta_1,\eta_2\gg N^{-1}$ with very high precision. This is our main technical result stated in Proposition \ref{prop:Cov}. To avoid unnecessary technical complications, we keep a suboptimal error term in this result, since it is anyway dominated by the error term coming from the critical regime $I_{\eta_0}^{\eta_c}$. We, nevertheless, discuss in detail how to eliminate this suboptimality.

\subsubsection{Mesoscopic regime: $\eta\in [\eta_c,T]$}\label{sec:large_eta} We focus on the ingredients for the analysis of $\Var[I_{\eta_c}^T]$ in the proof of Proposition \ref{prop:Var_main}, while the analysis of $\E I_{\eta_c}^T$ in the proof of Proposition \ref{prop:E_omega} relies on already existing tools. For more details regarding the calculation of the expectation see Section \ref{sec:E}. We write $\Var[I_{\eta_c}^T]$ as
\begin{equation}
\Var[I_{\eta_c}^T]= \left(\frac{N}{2\pi\ii}\right)^2 \int_\C\int_\C \Delta f_N(z_1)\overline{\Delta f_N(z_2)} \int_{\eta_c}^T\int_{\eta_c}^T \Cov\left(\langle G^{z_1}(\ii\eta_1)\rangle, \langle G^{z_2}(\ii\eta_2)\rangle\right)\dif\eta_1\dif\eta_2\dif^2 z_1\dif^2 z_2.
\label{eq:Var_I_ll} 
\end{equation}
The main result of this section is the very precise calculation of the covariance in the rhs. of \eqref{eq:Var_I_ll}, as stated in Proposition \ref{prop:Cov} below, which allows us to compute the lhs. of \eqref{eq:Var_I_ll} with an effective error term. 

To introduce the set-up for Proposition \ref{prop:Cov}, denote
\begin{equation}
E_+:=\begin{pmatrix}
1&0\\0&1
\end{pmatrix},\quad 
E_-:=\begin{pmatrix}
1&0\\0&-1
\end{pmatrix},\quad 
F:=\begin{pmatrix}
0&1\\0&0
\end{pmatrix},\quad E_+,E_-,F\in \C^{(2N)\times (2N)},
\end{equation}
i.e. $E_+,E_-$ and $F$ have a $2\times 2$ block-constant structure with the blocks of size $N\times N$. Recall the definition of $W$ from \eqref{eq:def_hermitization}. The self-energy operator associated to $W$ is defined by its action on the space of $(2N)\times (2N)$ deterministic matrices:
\begin{equation}
\mathcal{S}[R]:=\E[WRW] = \langle RE_+\rangle E_+-\langle RE_-\rangle E_-, \quad \forall R\in \C^{(2N)\times (2N)}.
\label{eq:def_S}
\end{equation}
The \emph{Matrix Dyson Equation} (MDE) is given by
\begin{equation}
-(M^z(w))^{-1} = w+Z+\mathcal{S}[M^z(w)],\quad \Im M^z(w)\Im w>0,\quad z\in\C,\, w\in\C\setminus\R,
\label{eq:MDE}
\end{equation}
where we used the notation $Z$ from \eqref{eq:def_hermitization}. It is known that \eqref{eq:MDE} has a unique solution; see the abstract statement in \cite{Helton07} and its application to \eqref{eq:MDE} in \cite[Lemma 2.2]{loc_circ_law18}. By [Eq.(3.5)--(3.6)]\cite{loc_circ_law18}, this solution can be written in the form 
\begin{equation}
M^z(w)=
\begin{pmatrix}
m^z(w)&-zu^z(w)\\ -\overline{z}u^z(w)&m^z(w)
\end{pmatrix}\in \C^{(2N)\times (2N)},\quad 
u^z(w):=\frac{m^z(w)}{w+m^z(w)},
\label{eq:M}
\end{equation}
where $m^z(w)$ is the unique solution to the scalar equation
\begin{equation}
-\frac{1}{m^z(w)} = w+m^z(w) -\frac{|z|^2}{w+m^z(w)},\quad \Im m^z(w)\Im w>0.
\label{eq:m}
\end{equation}
The \emph{self-consistent density of states} $\rho^z$ is given by 
\begin{equation}
\rho^z(x):=\lim_{\eta\to +0} \rho^z(x+\ii\eta),\quad \text{where}\quad \rho^z(w):=\frac{1}{\pi}\left\vert\Im m^z(w)\right\vert.
\label{eq:def_rho}
\end{equation}
This is an even function on the real line: $\rho^z(x)=\rho^z(-x)$ for any $x\in\R$. Finally, for $\kappa>0$, we define the $\kappa$-bulk of the density $\rho^z$ by
\begin{equation}
\mathbf{B}_\kappa^z:=\{x\in\R\,:\, \rho^z(x)\ge \kappa\}.
\label{eq:def_bulk}
\end{equation}

As $N$ goes to infinity, $G^z(w)$ is well-approximated by $M^z(w)$, such results are known as \emph{single-resolvent local laws}. For example, the \emph{averaged} single-resolvent local law \cite[Theorem~3.1]{macroCLT_real} asserts that
\begin{equation}
\left\vert\langle G^z(w)-M^z(w)\rangle\right\vert\le \frac{N^\xi}{N\eta},\quad \eta:=|\Im w|>0,
\label{eq:av_ll_simple}
\end{equation}
with very high probability for any fixed $z\in\C$. Since the typical size of $|\langle M^z(w)\rangle|$ is of order~1, \eqref{eq:av_ll_simple} manifests the concentration of $\langle G^z(w)\rangle$ around $\langle M^z(w)\rangle$ for $\eta\gg N^{-1}$. For this reason we often call $M^z$ the \emph{deterministic approximation} to $G^z$.

Now we are ready to state the main technical result of this paper: a very precise calculation of the covariance in the rhs. of \eqref{eq:Var_I_ll}. The proof of this proposition is presented in Section \ref{sec:Cov}. 
\begin{proposition}\label{prop:Cov} Let $X$ be a complex $N\times N$ i.i.d. matrix satisfying Assumption~\ref{ass:chi}. Denote $\kappa_4:=\E|\chi|^4-2$. Fix (small) $\delta, \epsilon, \kappa, \xi>0$. Uniformly in $|z_l|\le 1-\delta$ and $w_l\in\C\setminus\R$ with $E_l:=\Re w_l\in \mathbf{B}_\kappa^{z_l}$ and $\eta_l\in [N^{-1+\epsilon}, 1]$, $l=1,2$, it holds that
\begin{equation}
\Cov \left(\langle G^{z_1}(w_1)\rangle,\langle G^{z_2}(w_2)\rangle\right) = \frac{1}{N^2}\cdot \frac{V_{12}+\kappa_4U_1U_2}{2} + \mathcal{O}\left(\left(\frac{1}{N\gamma}+N^{-1/4}\right)\frac{N^\xi}{N^2\eta_1\eta_2}\right),
\label{eq:main_cov}
\end{equation}
where $V_{12}=V_{12}(z_1,z_2,w_1,w_2)$ and $U_l=U_l(z_l,w_l)$ are defined as
\begin{equation}
\begin{split}
\gamma=\gamma(z_1,z_2,w_1,w_2):=&|z_1-z_2|^2 +\mathrm{LT} + ||E_1|-|E_2||^2 +\eta_1+\eta_2,\\
\mathrm{LT}=\mathrm{LT}(z_1,z_2,w_1,w_2):=&|E_1|-|E_2| - \mathrm{sgn}(E_1) \frac{\Im u_1}{\Im m_1} \Re \left[\overline{z}_1(z_1-z_2)\right],\\
V_{12}=V_{12}(z_1,z_2,w_1,w_2):=&-\frac{1}{2}\partial_{w_1}\partial_{w_2}\log \left[ 1+(u_1u_2|z_1||z_2|)^2-m_1^2m_2^2-2u_1u_2\Re [z_1\overline{z_2}]\right],\\
U_l=U_l(z_l,w_l):=&-\frac{1}{\sqrt{2}}\partial_{w_l}m_l^2.
\end{split}
\label{eq:def_UV}
\end{equation}
Here we denoted $m_l=m^{z_l}(w_l)$, $u_l=u^{z_l}(w_l)$ and $\mathrm{sgn}(x)$ is the sign of $x\in \R$ defined with the convention $\mathrm{sgn}(0)=0$. In the Gaussian case the error $N^{-1/4}$ in \eqref{eq:main_cov} can be removed.
\end{proposition}

We now compare \eqref{eq:main_cov} with previous results and discuss the strategy of its proof, as well as briefly discuss the emergence of certain error terms in the rhs. of \eqref{eq:main_cov}.

\smallskip

\textbf{Comparison with the local law estimate.} We start the discussion of Proposition~\ref{prop:Cov} with comparing \eqref{eq:main_cov} with the local law bound on the lhs. of \eqref{eq:main_cov}.  The high probability bound \eqref{eq:av_ll_simple} implies that the lhs. of \eqref{eq:main_cov} is at most of order $(N^2\eta_1\eta_2)^{-1}$, up to a  factor $N^\xi$, and indeed this is the correct size of the leading term in the rhs. of \eqref{eq:main_cov} in the extreme case when $z_1$ and $z_2$ coincide. In contrast, the error term in \eqref{eq:main_cov} improves upon this estimate by a factor of $(N\gamma)^{-1}+N^{-1/4}$, which is effective in the regime $|z_1-z_2|\gg N^{-1/2}$. This factor reflects the decorrelation between $\langle G^{z_1}(\ii\eta_1)\rangle$ and $\langle G^{z_2}(\ii\eta_2)\rangle$ once the separation between the Hermitization parameters $|z_1-z_2|$ exceeds the typical eigenvalue spacing $N^{-1/2}$.

\smallskip

\textbf{Comparison with previous results.} Weaker versions of Proposition~\ref{prop:Cov} have been established earlier in \cite[Proposition~3.3]{macroCLT_complex} and in \cite[Proposition~3.4]{Cipolloni_meso}. In \cite{macroCLT_complex} the leading term in the rhs. of \eqref{eq:main_cov} was computed, but the error term was suboptimal. In fact, it improved only slightly beyond the local law estimate $(N^2\eta_1\eta_2)^{-1}$ by a factor $(N(\eta_1\wedge~\!\eta_2))^{-1/2}$ and only in the case $|z_1-z_2|\sim 1$. Later in \cite{Cipolloni_meso} this bound on the error term was extended to the entire mesoscopic regime $|z_1-z_2|\ge N^{-1/2+\omega}$ for any fixed $\omega>0$. However, the resulting error term was independent of $z_1,z_2$ and did not capture the decay in $|z_1-z_2|$. In contrast, the error term in the rhs. of \eqref{eq:main_cov} improves beyond \cite[Proposition~3.4]{Cipolloni_meso} and explicitly incorporates this dependence. Moreover, in the regime when $\eta_1\wedge\eta_2$ is just slightly above the critical threshold $N^{-1}$ and $|z_1-z_2|\sim 1$ this improvement yields an additional small factor of order $N^{-1/4}$ compared to \cite{Cipolloni_meso}.

\smallskip

\textbf{Discussion of the $N^{-1/4}$ error term.} In fact, the $N^{-1/4}$ term in the rhs. of~\eqref{eq:main_cov} can be eliminated, leaving the error term of order
\begin{equation*}
\frac{1}{N\gamma}\cdot \frac{N^\xi}{N^2\eta_1\eta_2}.
\end{equation*}
This can be achieved within the scope of our methods, however we do not pursue this improvement here, as it is not needed for the proof of Proposition~\ref{prop:Var_main}. The error term $N^{-1/4}$ comes from the third and higher order cumulants of the single-entry distribution of $X$. In particular, this term is not present if $X$ is a Ginibre matrix, as it is stated in the end of Proposition~\ref{prop:Cov}.

\smallskip

\textbf{Proof strategy.} Now we present the strategy of the proof of Proposition~\ref{prop:Cov}. Along the way we also explain how to remove the $N^{-1/4}$ term from the rhs. of~\eqref{eq:main_cov}. The main tool which we use in the proof of Proposition~\ref{prop:Cov} is the \emph{second order (Gaussian) renormalization}, originally introduced in \cite[Eq.(4.2)]{ETH_Wigner} (for earlier works, which did not operate with this formalism but followed a similar approach\footnote{We warn the reader that in \cite{Knowles20} the underline notation is used to denote the centering of a random variable, while the second order renormalization is present in the last term in the rhs. of \cite[Eq.(2.5)]{Knowles20}.} see e.g. \cite[Eq.(2.5)]{Knowles20}). For any differentiable function~$g(W)$ we denote
\begin{equation}
\underline{Wg(W)}:=Wg(W) - \widetilde{\E}\left[ \widetilde{W}\left(\partial_{\widetilde{W}}g\right)(W)\right],
\label{eq:def_under_1}
\end{equation} 
where $\widetilde{W}$ is an independent copy of $W$, the expectation $\widetilde{\E}$ is taken with respect to $\widetilde{W}$ and $\partial_{\widetilde{W}}$ is the derivative in the direction of $\widetilde{W}$. The name second order renormalization arises from the fact that the underline removes the second order terms in the cumulant expansion for $Wg(W)$. In particular, $\E\underline{Wg(W)}=0$ for $W$ with Gaussian entries.

Denote $G_j:=G^{z_j}(\ii\eta_j)$ and $M_j:=M^{z_j}(\ii\eta_j)$, for $j=1,2$. An elementary calculation based on \eqref{eq:def_under_1} shows that (see also \cite[Eq.(5.2)]{macroCLT_complex})
\begin{equation}
G_j=M_j-M_j\underline{WG_j}+M_j\langle G_j-M_j\rangle G_j,
\label{eq:G_exp_intro}
\end{equation}
for $j=1,2$. From \eqref{eq:G_exp_intro} we thus get
\begin{equation}
\begin{split}
\Cov\left(\langle G_1\rangle,\langle G_2\rangle\right) &= \Cov\left(\langle G_1-M_1\rangle \langle (G_1-M_1)A\rangle,\langle G_2\rangle\right) + \frac{1}{4N^2} \sum_{\sigma\in\{\pm\}} \sigma\langle G_1AE_\sigma G_2^2 E_\sigma\rangle\\
& -\E \left[ \langle \underline{WG_1A\rangle (\langle G_2\rangle -\E \langle G_2\rangle}\right],
\end{split}
\label{eq:init_expansion}
\end{equation}
where $A:=(1-\langle M_1^2\rangle)^{-1}M_1$. The derivation and analysis of~\eqref{eq:init_expansion} are presented in Section~\ref{sec:Cov_reduction}. We incorporate the fully underlined term in the second line and the last term in the first line of~\eqref{eq:init_expansion} into the error term. Their treatment will be discussed later in this section, since it is essential for obtaining the bound on the error term in~\eqref{eq:main_cov} and for its refinement mentioned above. For a moment, we ignore these terms and focus on the overall strategy of the proof of~\eqref{eq:main_cov}. By the single-resolvent local law~\eqref{eq:av_ll_simple}, the first term in the rhs. of~\eqref{eq:init_expansion} has an upper bound of order $(N^3\eta_1^2\eta_2)^{-1}$. This improves the trivial $(N^2\eta_1\eta_2)^{-1}$ local law bound on the lhs. of \eqref{eq:init_expansion} by a factor of $(N\eta_1)^{-1}\ll 1$. However, this estimate is insufficient for the proof of Proposition~\ref{prop:Cov} in the regime when $\eta_1,\eta_2$ are only slightly larger than $N^{-1}$; on the other hand this gain would have been already enough for the precision required in \cite[Proposition~3.3]{macroCLT_complex} and in \cite[Proposition~3.4]{Cipolloni_meso}. This necessitates a further expansion of the first term in the rhs. of~\eqref{eq:init_expansion}, in a manner similar to \eqref{eq:init_expansion}. For this reason, we refer to~\eqref{eq:init_expansion} as the~\emph{initial underline expansion}. In fact, we will need to perform such expansions iteratively for the covariances arising along the way to establish \eqref{eq:main_cov}, giving rise to a hierarchy of covariances. 

We refer to the iterative expansion procedure introduced above as the \emph{chaos expansion}, in analogy with the chaos decomposition of a random variable. A similar strategy was employed in~\cite[Theorem 2.1]{Knowles20} to compute the covariance of two resolvent traces in the Hermitian setting. From this point of view, our proof of Proposition~\ref{prop:Cov} can be viewed as a non-Hermitian analogue of the approach developed in~\cite{Knowles20}.

\smallskip

\textbf{Error terms in the chaos expansion.} Two types of error terms arise along the chaos expansion. They are represented by the last term in the first line of~\eqref{eq:init_expansion} (error terms of the \emph{first type}) and the term in the second line of~\eqref{eq:init_expansion} (error terms of the \emph{second type}). For simplicity, we focus here only on these two representatives, however, the proof of Proposition~\ref{prop:Cov} requires the analysis of all error terms generated at each step of the chaos expansion. The error term $N^{-1/4}(N^2\eta_1\eta_2)^{-1}$ in the rhs. of~\eqref{eq:main_cov} originates from the error terms of the first type, while the error terms of the second type contribute only $N^{-1/2}(N^2\eta_1\eta_2)^{-1}$. Consequently, the elimination of the $N^{-1/4}(N^2\eta_1\eta_2)^{-1}$ error term in the  rhs. of~\eqref{eq:main_cov} proceeds in two steps: first we improve it to $N^{-1/2}(N^2\eta_1\eta_2)^{-1}$ by a more precise treatment of the error terms of the first type, and then remove $N^{-1/2}(N^2\eta_1\eta_2)^{-1}$  by using iterative underline expansions for the error terms of the second type.

Although the last term in the first line of~\eqref{eq:init_expansion}, as well as the rest of the error terms of the first type, cannot be controlled using the single-resolvent local law~\eqref{eq:av_ll_simple}, we establish that it still concentrates around its deterministic counterpart, which contributes to the leading term in the rhs. of~\eqref{eq:main_cov}. In the case when $X$ has a Gaussian component of order one, we prove an optimal high probability bound of order
\begin{equation}
\frac{1}{N\gamma}\cdot\frac{N^\xi}{N^2\eta_1\eta_2}
\label{eq:first_type}
\end{equation}
on the fluctuation of $N^{-2}\langle G_1AE_\sigma G_2^2E_\sigma\rangle$. This is a consequence of the new two-resolvent local law, which we discuss in detail in Section~\ref{sec:local_laws}. In fact, the bound of order~\eqref{eq:first_type} holds for general i.i.d. matrices as well, and can be established by a routine Green function comparison (GFT) for the two-resolvent local law. Since this procedure is lengthy, technically involved, and unnecessary for the level of precision of Proposition~\ref{prop:Var_main}, instead we remove the Gaussian component from $X$ directly in $\Cov(\langle G_1\rangle,\langle G_2\rangle)$. This leads to the suboptimal term $N^{-1/4}(N^2\eta_1\eta_2)^{-1}$ in the rhs. of~\eqref{eq:main_cov}. Once~\eqref{eq:first_type} is available for general i.i.d. matrices, the contribution from the error terms of the first type to the rhs. of~\eqref{eq:main_cov} can be bounded by $(N\gamma)^{-1}(N^2\eta_1\eta_2)^{-1}$. Without the treatment of the error terms of the second type, this step alone improves the $N^{-1/4}$ term in the rhs. of \eqref{eq:main_cov} to~$N^{-1/2}$.
 
The treatment of the error terms of the second type is more delicate and requires an adjustment of the chaos expansion procedure. Specifically, one needs to perform the cumulant expansion for these error terms and subject the contribution from the third order cumulants to the iterative underline expansions as well, thereby extending the hierarchy of objects involved into the chaos expansion. A similar refinement in the simpler Hermitian setting was done in~\cite{Knowles20}.

\subsubsection{Sub-microscopic regime: $\eta\in [\eta_L,\eta_0]$}\label{sec:small_eta} As it was mentioned above, in the regime $\eta\in [\eta_L,\eta_0]$ the leading contribution to $\langle G^z(\ii\eta)\rangle$ comes from the least positive eigenvalue in its atypically small position. To control this contribution we prove the following bound on the left tail of the distribution of $\lambda_1^z$.

\begin{proposition}[Left tail of the smallest singular value distribution]\label{prop:tail_bound} Let $X$ be a complex $N\times N$ i.i.d. matrix satisfying Assumption \ref{ass:chi}. Fix (small) $\delta,\xi>0$ and set $\nu_1:=1/10$. Uniformly in $|z|\le 1-\delta$ and $N^{-\nu_1+\xi}\le x\le 1$ it holds that
\begin{equation}
\mathbf{P}\left[ \lambda_1^z \le N^{-1}x\right] \lesssim \log N\cdot x^2.
\label{eq:tail_bound}
\end{equation}
\end{proposition}
The proof of Proposition~\ref{prop:tail_bound} is presented in Section~\ref{sec:tail_bound}. It relies on the comparison with the~\emph{complex Ginibre ensemble} (GinUE), which corresponds to~$\chi$ distributed as a standard complex Gaussian random variable in the set-up introduced above Assumption~\ref{ass:chi}. It is known from~\cite[Corollary~2.4]{least_sing_val_Ginibre} that~\eqref{eq:tail_bound} holds for GinUE. First we transfer this result to i.i.d. random matrices with a small Gaussian component using the relaxation of the \emph{Dyson Brownian Motion} (DBM) from~\cite[Proposition~4.6]{bourgade2024fluctuations}, and then remove the Gaussian component by a short-time \emph{Green function comparison} (GFT) following the strategy introduced in~\cite{EYY12}.

The constant $\nu_1=1/10$ in Proposition~\ref{prop:tail_bound} is not sharp and can be increased by refining the relaxation bound~\cite[Proposition~4.6]{bourgade2024fluctuations} in our specific case when we have rigidity (see \cite[Remark 4.19]{bourgade2024fluctuations}). We also note that the restriction $a<~\! 1/2~\!+~\!\nu_1$ in Proposition~\ref{prop:E_omega} arises from the condition $x\gg N^{-\nu_1}$ imposed in Proposition~\ref{prop:tail_bound}. 

\subsubsection{Critical regime: $\eta\in [\eta_0,\eta_c]$}\label{sec:interm_eta} We treat $\E I_{\eta_0}^{\eta_c}$ and $\Var\left[ I_{\eta_0}^{\eta_c}\right]$ differently in the proofs of Propositions~\ref{prop:E_omega} and~\ref{prop:Var_main}, respectively. First, to show that $\E I_{\eta_0}^{\eta_c}$ is negligible, we compare this quantity with its analogue for the GinUE matrix $\widetilde{X}$, for which this regime can be shown to be negligible using the explicit formula for the one-particle density and the estimates on the rest of the regimes in the Girko's formula. The comparison between~$X$ and~$\widetilde{X}$ is carried out using the following result.

\begin{proposition}\label{prop:EG} Let $X$ be a complex $N\times N$ i.i.d. matrix satisfying Assumption~\ref{ass:chi}, and let $\widetilde{X}$ be an $N\times N$ GinUE matrix. Similarly to \eqref{eq:def_hermitization} denote the Hermitization of $\widetilde{X}$ by $\widetilde{H}^z$, $z\in \C$, and set $\widetilde{G}^z(w):=(\widetilde{H}^z-w)^{-1}$ for any $w\in \C\setminus\R$. Fix (small) $\delta,\epsilon, \omega_*>0$. Then uniformly in $|z|\le 1-\delta$ and $\eta\in [N^{-3/2+\epsilon},1]$ it holds that 
\begin{equation}
\E\langle G^z(\ii\eta)\rangle = \E \langle \widetilde{G}^z(\ii\eta)\rangle + \mathcal{O}(N^\xi\Phi_1(\eta))
\label{eq:EG_comp}
\end{equation}
for any fixed $\xi>0$, where we defined
\begin{equation}
\begin{split}
\Phi_1(\eta)&:=\min_{t\in [N^{-1+\omega_*},N^{-\omega_*}]}\left(N\mathcal{E}_0(t)\left(1+N\eta + \frac{N\mathcal{E}_0(t)}{N\eta}\right) +\sqrt{N}t\left(1+\frac{1}{N\eta}\right)^4\right),\\
\mathcal{E}_0(t)&:= \frac{1}{N}\left(\frac{1}{\sqrt{Nt}}+t\right).
\end{split}
\label{eq:def_Phi}
\end{equation}
\end{proposition}

The proof of Proposition~\ref{prop:EG} is presented in Section~\ref{sec:proof_EG}. Although~\eqref{eq:EG_comp} is stated for completeness for all $\eta\in [N^{-3/2+\epsilon},1]$, we use it only for $\eta\in [\eta_0,\eta_c]$, where $\eta_c$ is just slightly above~$N^{-1}$, as defined above~\eqref{eq:def_IJ}. Moreover, in the key regime $\eta\sim N^{-1}$ the complicated error term in the rhs. of~\eqref{eq:EG_comp} simplifies to
\begin{equation}
\Phi_1(\eta)\sim N^{-1/6},
\label{eq:Phi1_simple}
\end{equation}
as follows by an elementary optimization of~\eqref{eq:def_Phi} in~$t$. This yields an effective improvement of the bound on the error term in~\eqref{eq:EG_comp} compared to the local law bound~\eqref{eq:av_ll_simple}, which only provides an error term of order one in the same regime~$\eta\sim N^{-1}$.

A comparison similar to Proposition~\ref{prop:EG}, but in a substantially weaker form, was used in~\cite[Eq.(4.21)]{macroCLT_complex} for the same purpose. There only the range $\eta\in [N^{-1-\epsilon}, N^{-1+\epsilon}]$ was considered, and the bound on the error term was of order~$N^{-\omega}$ for some small implicit~$\omega>0$. In contrast, we establish a quantitative estimate for~$\omega$. For example,~\eqref{eq:Phi1_simple} implies that~$\omega$ can be taken equal to~$1/6$ for $\eta\sim N^{-1}$. We also extend the range of~$\eta$'s for which the error term in~\eqref{eq:EG_comp} improves upon the local law estimate~$(N\eta)^{-1}$.

Similarly to Proposition~\ref{prop:tail_bound}, the proof of Proposition~\ref{prop:EG} is based on the quantitative relaxation of the DBM from~\cite[Proposition~4.6]{bourgade2024fluctuations}. The resulting bound on the error term in~\eqref{eq:EG_comp} is not optimal and can be improved by refining the bound in~\cite[Proposition~4.6]{bourgade2024fluctuations}.

The analysis of $\Var[I_{\eta_0}^{\eta_c}]$ proceeds self-consistently, i.e. without comparison to the Gaussian case. It is based on the following decorrelation bound for resolvents, which was previously unknown even in the Gaussian case.

\begin{proposition}\label{prop:EGG} Let $X$ be a complex $N\times N$ i.i.d. matrix satisfying Assumption~\ref{ass:chi}. Fix (small) $\delta,\epsilon, \omega_*>0$. Then uniformly in $|z_1|, |z_2|\le 1-\delta$ and $\eta_1,\eta_2\in[ N^{-3/2+\epsilon},1]$ it holds that
\begin{equation}
\left\vert\Cov\left( \langle G^{z_1}(\ii\eta_1)\rangle, \langle G^{z_2}(\ii\eta_2)\rangle\right)\right\vert \lesssim N^\xi \Phi_2(\eta_1,\eta_2, |z_1-z_2|),
\label{eq:EGG}
\end{equation}
for any fixed $\xi>0$, where we defined
\begin{equation}
\begin{split}
\Phi_2(\eta_1,\eta_2,|z_1-z_2|)&:=\min_{t\in [N^{-1+\omega_*}, N^{-\omega_*}]} \min_{R\in [0,N|z_1-z_2|^2]} \\
\Bigg(N(\mathcal{E}_1(t,R)&+\mathcal{E}_0(t))\left(N\eta_1+\frac{1}{N\eta_1}\right)\left(N\eta_2+\frac{1}{N\eta_2}\right) + \sqrt{N}t \left(1+\frac{1}{N(\eta_1\wedge\eta_2)}\right)^3\frac{1}{N^2\eta_1\eta_2}\Bigg),\\ 
\mathcal{E}_1(t,R)&:=\sqrt{\frac{t}{N}}\left(\frac{1}{R^{1/8}}+\sqrt{\frac{R}{N|z_1-z_2|^2}}\right)+\frac{\sqrt{Nt^3}}{R},
\end{split}
\label{eq:def_Phi2}
\end{equation}
and $\mathcal{E}_0(t)$ is defined in \eqref{eq:def_Phi}.
\end{proposition}

The proof of Proposition~\ref{prop:EGG} is given in Section~\ref{sec:proof_EGG}. Similarly to Proposition~\ref{prop:EG}, we later use~\eqref{eq:EGG} only for $\eta_1,\eta_2\in [\eta_0,\eta_c]$, while the regime $\eta_1,\eta_2\gg N^{-1}$ is covered by Proposition~\ref{prop:Cov}, which for sufficiently large~$\eta_1,\eta_2$ gives in fact a much better bound. In the key regime $\eta_1,\eta_2\sim N^{-1}$ and $|z_1-z_2|\gg N^{-1/2}$ the rhs. of \eqref{eq:EGG} simplifies to
\begin{equation}
\Phi_2(\eta_1,\eta_2,|z_1-z_2|)\sim \left(N|z_1-z_2|^2\right)^{-q_0}
\label{eq:Phi2_simple}
\end{equation}
with $q_0=1/20$. In particular,~\eqref{eq:Phi2_simple} is the source of the exponent~$q_0$ appearing in Proposition~\ref{prop:Var_main}. We remark that the local law~\eqref{eq:av_ll_simple} yields an upper bound of order $(N^2\eta_1\eta_2)^{-1}$ on the lhs. of ~\eqref{eq:EGG}, while~\eqref{eq:def_Phi2} improves upon this estimate, at least in the regime $\eta_1,\eta_2\sim N^{-1}$ and $|z_1-z_2|\gg N^{-1/2}$, and explicitly captures the dependence on~$|z_1-z_2|$.

Several earlier works contain results which may be viewed as precursors of Proposition~\ref{prop:EGG}, see \cite[Proposition 3.5]{macroCLT_complex} and \cite[Proposition 3.5]{Cipolloni_meso}. These statements, however, concern only the regime $\eta\in [N^{-1-\epsilon},N^{-1+\epsilon}]$ and yield bounds on the covariance in the lhs. of~\eqref{eq:EGG} of order~$N^{-\omega}$ for $|z_1-z_2|\gg N^{-1/2+\xi}$ with~$\omega$ implicitly depending on~$\xi$. In contrast, our bound~\eqref{eq:EGG} provides a quantitative estimate on~$\omega$.

The proof of Proposition~\ref{prop:EGG} is based on the novel quantitative decorrelation estimate for a pair of Dyson Brownian motions, which we state in Theorem~\ref{eq:maintheoDBM1_main} and establish in Supplementary Section~\ref{sec:DBM}. The value of~$q_0$ in~\eqref{eq:Phi2_simple}, and consequently in Proposition~\ref{prop:Var_main} and Theorem~\ref{theo:main1}, is not optimal and can be increased by refining the estimate in Theorem~\ref{eq:maintheoDBM1_main}, see Supplementary Remark~\ref{rem:DBM}.

\subsection{The real case}\label{sec:real_ideas} The structure of the proof of Theorem~\ref{theo:main2} is the same as the one of Theorem~\ref{theo:main1} discussed above, however, the control on the error terms in Propositions~\ref{prop:E_omega}--\ref{prop:EGG} in the real case becomes weaker. In this section we explain the new technical difficulties causing this deterioration, while the precise statements and their proofs are given in Supplementary Section~\ref{sec:real}. 

First, the difference of the real case from the complex one is reflected in the self-energy operator
\begin{equation}
\mathcal{S}[R]=\langle RE_+\rangle E_+ - \langle RE_-\rangle E_- +\frac{1}{N} \left( E_+R^\mt E_+ - E_-R^\mt E_-\right),\quad \forall\, R\in\C^{(2N)\times (2N)}.
\label{eq:real_S_intro}
\end{equation}
Compared to \eqref{eq:def_S} stated in the complex case, \eqref{eq:real_S_intro} contains two additional \emph{torsion} terms involving $R^\mt$, which we recall that it stands for the transpose of $R$. The emergence of the torsion terms makes the structure of the chaos expansion in the proof of the real version of Proposition~\ref{prop:Cov} slightly more complicated. For instance, the real analogue of \eqref{eq:init_expansion} is as follows:
\begin{equation}
\begin{split}
\mathrm{Cov}\left(\langle G_1\rangle,\langle G_2\rangle\right) =& \mathrm{Cov}\left(\langle G_1-M_1\rangle \langle (G_1-M_1)A\rangle, \langle G_2\rangle\right) + \frac{\sigma}{2N}\mathrm{Cov}\left(\langle G^\mathfrak{t}_1E_\sigma G_1AE_\sigma\rangle,\langle G_2\rangle\right)\\
+&\frac{1}{4N^2}\left\langle G_1AE_\sigma \left(G_2^2 + (G_2^\mathfrak{t})^2\right)E_\sigma\right\rangle - \E\left[\underline{\langle WG_1A\rangle \left(\langle G_2\rangle-\E\langle G_2\rangle\right)}\right].
\end{split}
\label{eq:real_init_expansion_intro}
\end{equation}
We treat all terms in \eqref{eq:real_init_expansion_intro} involving $G_j^\mt$ for some $j\in [2]$ as error terms and do not expand them further. This allows us to use the same hierarchy of covariances as in the complex case at the price of including the error term of order
\begin{equation*}
\left(\frac{1}{N|\Im z_1|^2}+\frac{1}{N|\Im z_2|^2}\right) \frac{N^\xi}{N^2\eta_1\eta_2}.
\end{equation*}
into the rhs. of~\eqref{eq:main_cov}.

To minimize the changes in the proofs of Propositions~\ref{prop:tail_bound}--\ref{prop:EGG} in the real case compared to the complex case, we still use the complex DBM, which requires adding a small complex Ginibre component to a real i.i.d. matrix $X$. This approach was already used in \cite[Appendix~B]{cipolloni2024maximum}. To remove the complex component by a GFT, one needs to also estimate the contribution from the second-order cumulants, which is not present in the purely complex case. This slightly deteriorates the error terms in Propositions~\ref{prop:tail_bound}--\ref{prop:EGG} and it results in the second term in the rhs. of~\eqref{eq:main2}.

\subsection*{Outline of the paper} The remainder of the paper is structured as follows. First, in Section~\ref{sec:local_laws}, we summarize the local laws needed for the proofs of the results stated in Section~\ref{sec:results_list}. In particular, we state the novel two-resolvent averaged local law which is required for the proof of Proposition~\ref{prop:Cov} and explain the ideas behind its proof, postponing the technical details to Supplementary Section~\ref{sec:multiG_proof}. Next, in Section~\ref{sec:Girko_calculations}, we perform calculations in the Girko's formula using Propositions~\ref{prop:Cov}--\ref{prop:EGG}, and conclude the proof of Propositions~\ref{prop:E_omega} and~\ref{prop:Var_main}. In Section \ref{sec:proofprop}, we prove Propositions~\ref{prop:tail_bound}--\ref{prop:EGG} relying on the Dyson Brownian motion, which is justified by a standard argument deferred to Supplementary Section~\ref{sec:DBM}. Finally, in Section~\ref{sec:Cov} we prove Proposition~\ref{prop:Cov} by performing the chaos expansion. The real case, that is the proof of Theorem~\ref{theo:main2}, requires some more technical complications that are addressed in Supplementary Section~\ref{sec:real}.

\section{Local laws: set-up and results}\label{sec:local_laws}

In this section we state our second main technical result: the \emph{optimal averaged two-resolvent local law} for matrices with a Gaussian component, formulated in Proposition~\ref{prop:2G_av}. This result provides a bound for fluctuations of quantities of the form $\langle G^{z_1}(w_1)B_1G^{z_2}(w_2)B_2\rangle$ around their deterministic counterparts. Here $B_1, B_2\in \C^{(2N)\times (2N)}$ are deterministic matrices, and $|\Im w_1|, |\Im w_2|\gg N^{-1}$. Proposition~\ref{prop:2G_av} is one of the key ingredients in the proof of Proposition~\ref{prop:Cov}, allowing for the precise bound on the error term in~\eqref{eq:main_cov}, see the discussion around\footnote{\label{ftn:GG}While \eqref{eq:first_type} concerns the fluctuation of $\langle G_1AE_\sigma G_2^2E_\sigma\rangle$ with $G_j:=G^{z_j}(w_j)$ for $j=1,2$, which is not directly controlled to the two-resolvent local law, the number of resolvents in this quantity can be easily reduced to two using the identity $G_2^2 = \partial_{w_2} G_2$.}~\eqref{eq:first_type}. 

Throughout the paper we refer to concentration bounds of products of resolvents (or, more briefly, \emph{resolvent chains}) as \emph{local laws}.  The proof of Proposition~\ref{prop:Cov} is not the only place where the local laws are used, they appear repeatedly in the proof of Theorem~\ref{theo:main1}, in fact, in each of the technical ingredients listed in Sections~\ref{sec:large_eta}--\ref{sec:interm_eta}. Apart from Proposition~\ref{prop:2G_av} below, two additional new local laws are required. The first one is an averaged two-resolvent local law similar to Proposition~\ref{prop:2G_av}, which provides a weaker bound compared to Proposition~\ref{prop:2G_av}, but for general i.i.d. matrices~$X$, i.e. without imposing the Gauss-divisibility assumption. We state this result in Proposition~\ref{prop:2G_av_suboptimal} and prove it by extending \cite[Theorem~3.4]{nonHermdecay} from the imaginary axis to the bulk regime. The second local law is in fact a full hierarchy of concentration bounds for products of several resolvents and deterministic matrices sandwiched in between. We state these new multi-resolvent local laws in Proposition~\ref{prop:multiG_oneG} for general i.i.d. matrices~$X$. Though multi-resolvent local laws may involve arbitrarily large (but independent of~$N$) number of resolvents, they operate with bounds purely in terms of imaginary parts of spectral parameters and do not take into account the decorrelation decay in $|z_1-z_2|$. This makes the proof of Proposition~\ref{prop:multiG_oneG} substantially easier and fairly standard.

We start this section with recalling the well-established single-resolvent local laws, which are frequently used throughout the proofs of Propositions~\ref{prop:E_omega} and~\ref{prop:Var_main}.  Next, in Section~\ref{sec:2G} we state the new two-resolvent local laws in Propositions~\ref{prop:2G_av} and~\ref{prop:2G_av_suboptimal}. Finally, in Section~\ref{sec:multlonggs} we generalize the concepts introduced in Section~\ref{sec:2G} to longer products of resolvents and state the new multi-resolvent local laws in Proposition~\ref{prop:multiG_oneG}.

\subsection{Single-resolvent local law}

Recall the definition of~$M^z$ from~\eqref{eq:M}. In \eqref{eq:av_ll_simple} we have already indicated that $G^z(w)$ concentrates around $M^z(w)$ in the averaged sense for~$w$ on the imaginary axis. The same concentration result holds also away from the imaginary axis:

\begin{proposition}[Single-resolvent local law, {\cite[Theorem~3.1]{macroCLT_real}}]\label{prop:1G} Fix a (small) $\delta>0$ and a (large) $L>0$. Let $X$ be an $N\times N$ either complex or real i.i.d. matrix satisfying Assumption~\ref{ass:chi}. Uniformly in $|z|\le 1-\delta$ and $w\in \C\setminus\R$ with $|w|\le N^L$ and $\eta:=|\Im w|\ge N^{-L}$, it holds that 
\begin{align}
\left\vert \langle (G^z(w)-M^z(w))A\rangle\right\vert &\prec \frac{\|A\|}{N\eta},\label{eq:1G_ll_av}\\
\left\vert \langle\bm{x}, (G^z(w)-M^z(w))\bm{y}\rangle\right\vert&\prec \left(\frac{1}{\sqrt{N\eta}}+\frac{1}{N\eta}\right)\|\bm{x}\|\|\bm{y}\|,\label{eq:1G_ll_iso} 
\end{align}
for any deterministic matrix $A\in\C^{(2N)\times (2N)}$ and vectors $\bm{x},\bm{y}\in \C^{2N}$.
\end{proposition}

The estimates~\eqref{eq:1G_ll_av} and~\eqref{eq:1G_ll_iso} are known as the \emph{averaged} and \emph{isotropic} single-resolvent local laws, respectively. An important consequence of~\eqref{eq:1G_ll_av} is the bulk \emph{eigenvalue rigidity} from~\cite[Eq.(7.4)]{macroCLT_complex}:
\begin{equation}
|\lambda_i^z-\gamma_i^z|\prec \frac{1}{N},\quad |i|\le (1-\tau)N,
\label{eq:rigidity}
\end{equation}
uniformly in $|z|\le 1-\delta$, for any fixed~$\tau>0$. Here the quantiles $\{\gamma_i^z\}_{|i|\le N}$ of $\rho^z$, from \eqref{eq:def_rho}, are implicitly defined by
\begin{equation}
\frac{i}{N} =\int_0^{\gamma_i^z}\rho^z(E)\dif E,\quad \gamma_{-i}^z:=-\gamma_i^z,\quad 1\le i\le N.
\label{eq:def_quantiles}
\end{equation}
In \cite[Eq.(7.4)]{macroCLT_complex}, \eqref{eq:rigidity} is stated only for $|i|\le cN$ for some small constant $c>0$. However, \eqref{eq:rigidity} holds throughout the bulk regime $|i|\le (1-\tau)N$ as currently formulated. For the detailed derivation of \eqref{eq:rigidity} from Proposition~\ref{prop:1G} see e.g.~\cite[Section~5]{EYY12}. 

Concluding the discussion of eigenvalue rigidity, we record the following estimate on quantiles which we will use occasionally. We have
\begin{equation}
|\gamma_i^z|\sim |i|/N,
\label{eq:gamma_size}
\end{equation}
uniformly in $|i|\le N$ and $|z|\le 1-\delta$ for any fixed~$\delta>0$. This follows directly from \eqref{eq:def_quantiles} together with the fact that $\rho^z(0)>0$ for~$|z|<1$.

\subsection{Novel averaged two-resolvent local law}\label{sec:2G}

Now we introduce the set-up for our new averaged two-resolvent local law. For $z_j\in \C$ and $w_j\in \C\setminus\R$, $j=1,2$, denote $M_j:=M^{z_j}(w_j)$. Consider a resolvent chain $G_1B_1G_2$ for a $(2N)\times (2N)$ deterministic matrix~$B_1$, which we customarily call an \emph{observable}. Unlike one may naively think, $G_1B_1G_2$ does not concentrate around $M_1B_1M_2$ as $N$ goes to infinity.  Instead, the \emph{deterministic approximation} of $G_1B_1G_2$ is given by (see e.g. \cite[Eq.(5.7),(5.8)]{macroCLT_complex})
\begin{equation}
M_{12}^{B_1}=M_{12}^{B_1}(w_1,w_2)=\mathcal{M}[G_1B_1G_2]:=\mathcal{B}_{12}^{-1}[M_1B_1M_2],
\label{eq:def_M12}
\end{equation}
where $\mathcal{B}_{12}$ is the \emph{two-body stability operator} defined by its action on $(2N)\times (2N)$ deterministic matrices:
\begin{equation}
\left(\mathcal{B}_{12}(w_1,w_2)\right)[R]=\mathcal{B}_{12}[R]:=R-M_1\mathcal{S}[R]M_2,\quad \forall R\in\C^{(2N)\times (2N)},
\label{eq:def_B12}
\end{equation}
and $\mathcal{S}$ is defined in~\eqref{eq:def_S}. The dependence of $\mathcal{B}_{12}$ and $M_{12}^{B_1}$ on $z_1, z_2$ is indicated by the subscript. When $z_1=z_2$, we simply denote these objects by $\mathcal{B}_{11}$ and $M_{11}^{B_1}$, respectively. If additionally $w_1$ and $w_2$ coincide, we call $\mathcal{B}_{11}(w_1,w_1)$ the \emph{one-body stability operator}. The invertibility of $\mathcal{B}_{12}$ follows from \cite[Eq.(6.2)]{macroCLT_complex}, so the lhs. of~\eqref{eq:def_M12} is well-defined. We interchangeably use the three notations introduced in~\eqref{eq:def_M12} for the deterministic approximation of $G_1B_1G_2$, leaning towards $\mathcal{M}[G_1B_1G_2]$ in the calculations where longer resolvent chains are involved. 

For simplicity, we restrict our attention to observables of size $(2N)\times (2N)$ which are linear combinations of $E_+,E_-,F$ and~$F^*$, as this setting is sufficient for the proof Proposition~\ref{prop:Cov}. The general case can be managed without any additional technical difficulties by projecting observables on $\mathrm{span}\{E_\pm,F^{(*)}\}$. Our restriction is equivalent to saying that observables posses $2\times 2$ block-constant structure. Since the matrices~$M_j$, $j=1,2$, share this structure, and the operator~$\mathcal{B}_{12}$ preserves it, we can identify observables, deterministic approximations to resolvents and to the 2-resolvent chains~\eqref{eq:def_M12} with $2\times 2$ matrices. As it will become apparent later, this identification extends to deterministic approximations to longer resolvent chains as well. From now on, we will use this identification interchangeably, viewing the objects mentioned in this paragraph as $(2N)\times (2N)$ matrices with $2\times 2$ block-constant structure and as $2\times 2$ matrices.

We use the following control parameter to quantify the concentration of $G_1B_1G_2$ around $M_{12}^{B_1}$: 
\begin{equation}
\widehat{\beta}_{12}=\widehat{\beta}_{12}(w_1,w_2):=\min\left\lbrace \left\lVert \left(\mathcal{B}_{12}\big(w_1^{(*)},w_2^{(*)}\big)\right)[R]\right\rVert \,:\,R\in \mathrm{span}\{E_\pm,F^{(*)}\}, \|R\|=1\right\rbrace,
\label{eq:def_beta_hat}
\end{equation}
where $w_j^{(*)}$ indicates both choices~$w_j$ and~$\overline{w}_j$. The dependence of~$\widehat{\beta}_{12}$ on $z_1,z_2$ is recorded in the subscript similarly to $\mathcal{B}_{12}$. An important property of~$\widehat{\beta}_{12}$ which immediately follows from its definition and~\eqref{eq:def_M12} is that
\begin{equation}
\|M_{12}^{B_1}\|\le \big(\widehat{\beta}_{12}\big)^{-1},\quad\forall B_1\in \mathrm{span}\{E_\pm,F^{(*)}\}.
\label{eq:M_12_beta_bound}
\end{equation}
Moreover, this bound is optimal among the ones which are uniform in~$B_1$. We further have that $\widehat{\beta}_{12}(w_1,w_2)$ is insensitive to the complex conjugations of $w_1, w_2$ and that
\begin{equation}
0<\widehat{\beta}_{12}(w_1,w_2)\le 1,\quad \widehat{\beta}_{12}(w_1,w_2)=\widehat{\beta}_{21}(w_2,w_1),\qquad \forall \, z_1,z_2\in\C,\, w_1,w_2\in\C\setminus\R.
\label{eq:beta_hat_trivial}
\end{equation}
Here the positivity of $\widehat{\beta}_{12}$ follows from invertibility of $\mathcal{B}_{12}$, and to verify the bound $\widehat{\beta}_{12}\le 1$ it suffices to take $R:=F$ in the rhs. of \eqref{eq:def_beta_hat} and use that $\mathcal{B}_{12}[F]=F$. The second part of~\eqref{eq:beta_hat_trivial} directly follows from the identity
\begin{equation}
\left(\left(\mathcal{B}_{12}(w_1,w_2)\right)[R]\right)^* = \left(\mathcal{B}_{21}(\overline{w}_2,\overline{w}_1)\right)[R^*],\quad \forall R\in \mathrm{span}\{E_\pm,F^{(*)}\}.
\end{equation}

Now we are ready to state our second main technical result, the optimal two-resolvent local law for matrices with a Gaussian component, postponing the proof to Section~\ref{sec:multiG_proof}.

\begin{proposition}[Optimal averaged two-resolvent local law in the bulk, Gauss-divisible case]\label{prop:2G_av} Let $X_0$ be a complex $N\times N$ i.i.d. matrix satisfying Assumption \ref{ass:chi} and let $\widetilde{X}$ be a complex Ginibre matrix independent of $X_0$. For a possibly $N$-dependent $\mathfrak{s}=\mathfrak{s}_N\in (0,1)$ set $X:=\sqrt{1-\mathfrak{s}^2}X_0+\mathfrak{s}\widetilde{X}$. Fix $b\in [0,1]$ and a (large) $L>0$, and assume that $\mathfrak{s}\ge L^{-1} N^{-b}$.

For spectral parameters $w_1,w_2\in\C\setminus\R$, denote $\eta_l:=|\Im w_l|$ and $\eta_*:=\eta_1\wedge \eta_2\wedge 1$.  Then for any fixed $\delta,\epsilon,\kappa>0$ we have
\begin{equation}
\left\vert \left\langle \left(G^{z_1}(w_1)B_1G^{z_2}(w_2)-M_{12}^{B_1}(w_1,w_2)\right)B_2\right\rangle\right\vert \prec \frac{1}{N\eta_*\left(\widehat{\beta}_{12}(w_1,w_2)\wedge N^{-b} + \eta_*\right)},
\label{eq:2G_av}
\end{equation}  
uniformly in $B_1,B_2\in\mathrm{span}\{E_\pm,F^{(*)}\}$, $|z_l|\le 1-\delta$, $\Re w_l\in\mathbf{B}_\kappa^{z_l}$ and $N^{-1+\epsilon}\le \eta_l\le N^{100}$ for $l=1,2$.
\end{proposition}

We now discuss the result in \eqref{eq:2G_av}, relate it to previous results, and give a sketch of its proof.

\smallskip

\textbf{Gauss divisibility assumption.} We start the discussion of Proposition~\ref{prop:2G_av} by commenting on the Gauss divisibility assumption. First, we note that the word \emph{optimal} in the title of Proposition~\ref{prop:2G_av} refers only to the case $b=0$, i.e. when $X$ contains a Gaussian component of order one. In this case~\eqref{eq:2G_av} simplifies to
\begin{equation}
\left\vert \left\langle \left(G^{z_1}(w_1)B_1G^{z_2}(w_2)-M_{12}^{B_1}(w_1,w_2)\right)B_2\right\rangle\right\vert \prec \frac{1}{N\eta_*\widehat{\beta}_{12}(w_1,w_2)}.
\label{eq:2G_av_b=0}
\end{equation}  
Here we omitted $\eta_*$ in the second factor in the rhs. of~\eqref{eq:2G_av}; this does not change the bound, since $\widehat{\beta}_{12}\gtrsim \eta_*$ by \cite[Lemma 6.1]{macroCLT_real}. The bound~\eqref{eq:2G_av_b=0} is optimal in the sense that, for certain observables $B_1,B_2$, the variance of the lhs. of~\eqref{eq:2G_av_b=0} is of the same order as the square of the rhs. of \eqref{eq:2G_av_b=0}. At the opposite extreme $b=1$,~\eqref{eq:2G_av} yields an upper bound of order $(N\eta_*^2)^{-1}$, which is already known from \cite[Theorem~3.4]{nonHermdecay} for general i.i.d. matrices $X$ in the regime when $w_1,w_2$ are on the imaginary axis. For intermediate values $b\in (0,1)$,~\eqref{eq:2G_av} interpolates between \eqref{eq:2G_av_b=0} and this $(N\eta_*^2)^{-1}$ bound. In fact, the magnitude of fluctuation of the lhs. of~\eqref{eq:2G_av} is not sensitive to the size of Gaussian component, and the apparent deterioration of~\eqref{eq:2G_av} as $b$ increases from~0 to~1 is purely technical. In principle, one could remove the Gaussian component from $X$ without weakening~\eqref{eq:2G_av_b=0}, and prove that~\eqref{eq:2G_av_b=0} holds for general i.i.d. matrices~$X$. This could be achieved by a routine, though technically involved, GFT argument. Since~\eqref{eq:2G_av} is already sufficient for the proof of Proposition~\ref{prop:Cov}, we do not pursue these additional refinements.

\smallskip

\textbf{Relation to the upper bound on $\|M_{12}^B\|$.} We further focus on the case $b=0$ in Proposition~\ref{prop:2G_av}, which is our main contribution. An important feature of~\eqref{eq:2G_av_b=0} is that this bound on the fluctuation of $\langle G_1B_1G_2B_2\rangle$ is  stronger than the bound on its deterministic approximation $\langle M_{12}^{B_1}B_2\rangle$  \eqref{eq:M_12_beta_bound} by a factor of $(N\eta_*)^{-1}\ll 1$. Thus,~\eqref{eq:2G_av} should be viewed as a concentration bound. The emergence of the $(N\eta_*)^{-1}$ improvement in multi-resolvent averaged local laws compared to the size estimate is well-expected and was first observed in the Hermitian setting for Wigner matrices in~\cite{Multi_res_llaws}.

\smallskip

\textbf{Relation to the previous results.} Several weaker versions of~\eqref{eq:2G_av_b=0} appeared earlier in the literature, see in chronological order~\cite[Theorem~5.2]{macroCLT_complex}, \cite[Theorem~3.3]{Cipolloni_meso} and \cite[Theorem~3.4]{nonHermdecay}. We remark that none of these works assumes Gauss divisibility of~$X$. In~\cite{macroCLT_complex} the bound in the averaged two-resolvent local law is worse than~\eqref{eq:2G_av_b=0} by a typically large factor~$\eta_*^{-11/12}$, while in~\cite{Cipolloni_meso}~$\widehat{\beta}_{12}$ is replaced by a weaker control parameter and is not gained in the optimal power. Finally, for $w_1,w_2$ on the imaginary axis, \cite[Theorem~3.4]{nonHermdecay} asserts that
\begin{equation}
\left\vert \left\langle \left(G^{z_1}(w_1)B_1G^{z_2}(w_2)-M_{12}^{B_1}(w_1,w_2)\right)B_2\right\rangle\right\vert \prec \frac{1}{\sqrt{N\eta_*}\left(|z_1-z_2|^2 + (\eta_1+\eta_2)\wedge 1\right)}\wedge \frac{1}{N\eta_*^2}.
\label{eq:2G_av_borrowed}
\end{equation}
Besides the fact that the local law is needed for the proof of Proposition~\ref{prop:Cov} for $w_1,w_2$ not only on the imaginary axis, but also in a neighborhood of it (as already mentoned above Proposition~\ref{prop:2G_av}), there is an additional loss in~\eqref{eq:2G_av_borrowed} compared to~\eqref{eq:2G_av_b=0}. Although the $\widehat{\beta}_{12}$ factor appears in the rhs. of~\eqref{eq:2G_av_borrowed}, as by \cite[Lemma~3.3]{nonHermdecay} it is of the same order as $|z_1-z_2|^2 + (\eta_1+\eta_2)\wedge 1$ for spectral parameters on the imaginary axis, there is still a substantial loss of a large factor $\sqrt{N\eta_*}$ in the denominator. In particular, even if~\eqref{eq:2G_av_borrowed} was extended to the bulk regime, the replacement of Proposition~\ref{prop:2G_av} by \eqref{eq:2G_av_borrowed} in the proof of Proposition~\ref{prop:Cov} would considerably weaken the bound on the error term in~\eqref{eq:main_cov}.

\smallskip

\textbf{Proof strategy.} The proof of Proposition~\ref{prop:2G_av} relies on the \emph{method of characteristics}, which consists in studying the evolution of the resolvent along a stochastic flow which enables us to transfer information from spectral parameters with large imaginary part to those with much smaller imaginary part. In \cite{pastur1972spectrum}, it was observed that the resolvent of a Wigner matrix evolving along the  Ornstein-Uhlenbeck flow (see \eqref{eq:OU_ll_discussion} below) solves complex Burgers-like equations, and it can thus be analyzed along the characteristics (which have the convenient property described above). More recently this idea was used to prove universality at the edge of certain deformed Hermitian matrices \cite{lee2015edge}, and, closer to our setting, to prove local laws for single resolvents \cite{adhikari2020dyson, HL_rig, von2019random}. Only more recently it was shown that this \emph{method of characteristics} is powerful to also consider products of resolvents \cite{bourgade2022liouville, Cipolloni_meso}. We refer the interested reader to \cite[Section 1.4]{univ_extr} for a detailed description of the history and related results.

Specifically, consider the Ornstein-Uhlenbeck flow
\begin{equation}
\dif X_t = -\frac{1}{2}X_t\dif t + \frac{\dif B_t}{\sqrt{N}},
\label{eq:OU_ll_discussion}
\end{equation}
where $B_t$ is an $N\times N$ matrix composed of $N^2$ independent standard complex-valued Brownian motions, and choose the initial condition in such a way that $X\stackrel{d}{=}X_T$ for some $T\sim N^{-b}$. We complement \eqref{eq:OU_ll_discussion} by the \emph{characteristic flow}
\begin{equation}
\frac{\dif}{\dif t} z_{j,t}=-\frac{1}{2} z_{j,t},\qquad\frac{\dif}{\dif t} w_{j,t} = -\frac{1}{2}w_{j,t} - \langle M^{z_{j,t}}(w_{j,t})\rangle,\quad t\in [0,T],\, j=1,2,
\label{eq:char_flow_discussion}
\end{equation} 
with the final condition $z_{j,T}:=z_j$, $w_{j,T}:=w_j$ (instead of the customary initial condition). This flow is well-defined by \cite[Lemma~5.2]{Cipolloni_meso}. Define $G_{j,t}:=G^{z_{j,t}}(w_{j,t})$ and set
\begin{equation}
\mathcal{Y}_{\sigma_1,\sigma_2,t}:=\left\langle \left(G_{1,t}E_{\sigma_1}G_{2,t}-M_{12,t}^{E_{\sigma_1}}\right)E_{\sigma_2}\right\rangle,\quad \mathcal{Y}_t:=\left(\mathcal{Y}_{+,+,t},\mathcal{Y}_{+,-,t},\mathcal{Y}_{-,+,t},\mathcal{Y}_{-,-,t}\right)^\mt\in\C^4, 
\end{equation}
where $\sigma_1,\sigma_2\in\{\pm\}$ and $M_{12,t}^B$ is the time-dependent deterministic approximation. An elementary application of It\^{o} calculus gives
\begin{equation}
\dif \mathcal{Y}_{[2],t} = \left(I+\mathcal{A}_{[2],t}\right)\mathcal{Y}_{[2],t}\dif t +\mathcal{F}_{[2],t}\dif t+\dif\mathcal{E}_{[2],t},
\label{eq:Y_syst_intro}
\end{equation}
for a similar calculation see e.g. \cite[Eq.(5.42)--(5.47)]{univ_extr} and \cite[Eq.(4.67)--(4.69)]{nonHermdecay}. In \eqref{eq:Y_syst_intro} the \emph{generator} $\mathcal{A}_{[2],t}\in\C^{4\times 4}$ is an explicit $4\times 4$ matrix constructed from the deterministic approximations to two-resolvent chains, $\mathcal{F}_{[2],t}$ is a \emph{forcing term} and $\mathcal{E}_{[2],t}$ is a \emph{martingale term}. The exact form of $\mathcal{A}_{[2],t}$ is given in Appendix~\ref{app:propagator} later (see also Supplementary Eq.~\eqref{eq:def_general_A}). 

We will solve \eqref{eq:Y_syst_intro} by the Duhamel formula, so need to understand the propagator of this ODE. To this end, we prove the following bound on the propagator in Appendix~\ref{app:propagator}. 

\begin{lemma}\label{lem:propag_bound2} Denote $f_{[2],r}:=\left[ \max \mathrm{spec}(\Re \mathcal{A}_{[2],r})\right]_+$ for $r\in [0,T]$, where $[\cdot]_+$ stands for the positive part of a real number. Then for any $t\in [0,T]$ we have
\begin{equation}
\exp\left\lbrace \int_s^t f_{[2],r}\dif r\right\rbrace \lesssim \left(\frac{\widehat{\beta}_{12,s}}{\widehat{\beta}_{12,t}}\right)^2,
\label{eq:propag_bound2}
\end{equation}
with $\widehat{\beta}_{12,t}$ equal to $\widehat{\beta}_{12}$ (see \eqref{eq:def_beta_hat} for its definition) evaluated at $z_{j,t}$ and $w_{j,t}$, $j=1,2$.
\end{lemma}

This lemma is the analogue of Lemma~5.7 and Eq.~(5.58) from \cite{univ_extr}, where similar estimates on the propagator were derived in the regime when $z_1,z_2$ are close to the non-Hermitian spectral edge and $w_1,w_2$ are purely imaginary. A major simplification in \cite{univ_extr} is that the off-diagonal entries of $\mathcal{A}_{[2],t}$ are purely imaginary when $w_1,w_2$ are on the imaginary axis, so the real part of the generator is diagonal. In our setting this simplification is no longer available, so we separately treat the off-diagonal entries of $\mathcal{A}_{[2],t}$.

The remaining details in the proof of Proposition~\ref{prop:2G_av} are standard, and thus we present them in Supplementary Section~\ref{sec:multiG_proof}.

\smallskip

\textbf{Discussion of the control parameter $\widehat{\beta}_{12}$.} Although \eqref{eq:def_beta_hat} explicitly defines $\widehat{\beta}_{12}$, it does not make transparent the dependence of $\widehat{\beta}_{12}$ on $z_1,z_2$ and $w_1,w_2$. This more detailed information is in fact not needed for the proof of Proposition \ref{prop:2G_av}, where we mostly rely on the bound \eqref{eq:M_12_beta_bound}. However, effective estimates on $\widehat{\beta}_{12}$ in terms of $z_1,z_2$ and $w_1,w_2$ become essential when applying the two-resolvent local law \eqref{eq:2G_av} in the proof of Proposition \ref{prop:Cov}, as well as in the potential future applications.  The following statement fully covers this question by providing precise asymptotics for $\widehat{\beta}_{12}$ expressed directly in these simpler terms.

\begin{proposition}\label{prop:stab_bound} Fix (small) $\delta,\kappa>0$. Uniformly in Hermitization parameters $z_1,z_2\in (1-\delta)\mathbf{D}$ and spectral parameters $w_1,w_2\in \C\setminus\R$ such that $E_j:=\Re w_j\in\mathbf{B}_\kappa^{z_j}$ and $\eta_j:=|\Im w_j|\in (0,1]$ for $j=1,2$, it holds that
\begin{equation}
\widehat{\beta}_{12}(w_1,w_2)\sim \gamma(z_1,z_2,w_1,w_2),
\label{eq:stab_LT}
\end{equation}  
where 
$\gamma$ is defined in \eqref{eq:def_UV}. Moreover, uniformly in $|z_j|\lesssim 1$ and $w_j\in\C\setminus\R$, $j=1,2$, we have
\begin{equation}
\widehat{\beta}_{12}\gtrsim \eta_*,\quad \text{where}\quad\eta_*:=\eta_1\wedge\eta_2\wedge 1.
\label{eq:eta_stab_bound}
\end{equation}
\end{proposition}

Our main contribution in Proposition \ref{prop:stab_bound} is the identification of the term linear in $z_1-z_2$ and $w_1-w_2$ in $\gamma$, which is captured by $\mathrm{LT}$ defined in~\eqref{eq:def_UV}. The remaining terms in $\gamma$ have already appeared as a lower bound on a quantity closely related to $\widehat{\beta}_{12}$ in \cite[Lemma~6.1]{macroCLT_real}, together with the estimate~\eqref{eq:eta_stab_bound}. We remark that the linear term is analogous to the quantity introduced under the same name in \cite[Eq.(2.17)]{eigenv_decorr}, where it played a crucial role in \cite[Proposition~3.1]{eigenv_decorr} in improving the lower bound on the two-body stability operator from quadratic to linear dependence on the distance between spectral parameters in the Hermitian setting. For further details and the proof of Proposition \ref{prop:stab_bound} see Appendix~\ref{app:stab_bound}.

Since Proposition~\ref{prop:Cov} concerns general i.i.d. matrices, while Proposition~\ref{prop:2G_av} is stated for i.i.d. matrices with a Gaussian component, a GFT argument is still required in the proof of Proposition~\ref{prop:Cov}. However, a Gaussian component is removed from $X$ directly at the level of covariance (\emph{direct} GFT), rather than in the local law estimate, as it is technically simpler. There is one more instance in the proof of Proposition~\ref{prop:Var_main}, where an averaged two-resolvent local law is needed (though with less precision than in~\eqref{eq:2G_av_b=0}), and where we prefer to perform GFT in the local law instead of working out a direct GFT, since the former one is more standard. Further details will be given later in Proposition~\ref{prop:overlap}. In such a way, we establish a (suboptimal) bound in the two-resolvent averaged local law for general i.i.d. random matrices. We accomplish this in the following proposition by extending~\eqref{eq:2G_av_borrowed} to the entire bulk regime. 

\begin{proposition}\label{prop:2G_av_suboptimal} Let $X$ be a complex $N\times N$ i.i.d. matrix satisfying Assumption \ref{ass:chi}. For any fixed $\delta,\epsilon,\kappa>0$ it holds that
\begin{equation}
\left\vert \left\langle \left(G^{z_1}(w_1)B_1G^{z_2}(w_2)-M_{12}^{B_1}(w_1,w_2)\right)B_2\right\rangle\right\vert \prec \frac{1}{\sqrt{N\eta_*}\widehat{\beta}_{12}}
\label{eq:2G_av_suboptimal}
\end{equation}  
uniformly in $B_1,B_2\in\mathrm{span}\{E_\pm,F^{(*)}\}$, $|z_l|\le 1-\delta$, $\Re w_l\in\mathbf{B}_\kappa^{z_l}$ and $N^{-1+\epsilon}\le \eta_l\le N^{100}$ for $l=1,2$.
\end{proposition}

The proof of Proposition~\ref{prop:2G_av_suboptimal} is presented in Supplementary Section~\ref{sec:2G_av_subopt}. It is obtained by a minor modification of the proof of~\cite[Theorem~3.4]{nonHermdecay} using the new propagator bound from Lemma~\ref{lem:propag_bound2} developed for the proof of Proposition~\ref{prop:2G_av}.

As a standard corollary of Proposition~\ref{prop:2G_av_suboptimal}, we get the following overlap bound for the left and right singular vectors of $X-z_1$ and $X-z_2$ associated to the bulk eigenvalues\footnote{Here we recall the relation between singular vectors of $X-z$ and eigenvectors of $H^z$ discussed in~\eqref{eq:eigenvectors_symmetry}--\eqref{eq:def_uv}} of $H^{z_1}$, $H^{z_2}$. We cover the entire bulk regime, which is possible since we prove Proposition \ref{prop:2G_av_suboptimal} throughout the bulk. Meanwhile~\eqref{eq:2G_av_suboptimal} restricted to the spectral parameters on the imaginary axis would imply an effective bound only on the singular vector overlaps associated to the eigenvalues of $H^{z_1}, H^{z_2}$ in the interval $[-N^{-1+\omega}, N^{-1+\omega}]$ for a small~$\omega>0$.

\begin{proposition}[Singular vector overlap]\label{prop:overlap} Let $X$ be a complex $N\times N$ i.i.d. matrix satisfying Assumption~\ref{ass:chi}. Fix (small) $\delta, \tau>0$. Recall the notations introduced in \eqref{eq:M}, \eqref{eq:def_uv} and \eqref{eq:def_quantiles}. It holds that
\begin{equation}
\left\vert \langle \bm{u}_i^{z_1},\bm{u}_j^{z_2}\rangle\right\vert^2 + \left\vert \langle \bm{v}_i^{z_1},\bm{v}_j^{z_2}\rangle\right\vert^2 \le \frac{N^\xi}{N}\frac{1}{|z_1-z_2|^2+N^{-1}|i-j|}\wedge 1,
\label{eq:overlap_bound}
\end{equation}
with very high probability for any fixed $\xi>0$, uniformly in $i,j\in [(1-\tau)N]$ and $|z_1|,|z_2|\le 1-\delta$.
\end{proposition}

The proof of Proposition~\ref{prop:overlap} is presented in Appendix~\ref{app:overlap}. It relies on Proposition~\ref{prop:2G_av_suboptimal} along with the spectral decomposition, and on a precise analysis of the linear term in the second line of \eqref{eq:def_UV}.

\subsection{Multi-resolvent local laws for longer chains}
\label{sec:multlonggs}

In this section we state local laws (both in the averaged and isotropic sense) for resolvent chains of the form
\begin{equation}
G_1B_1G_2\cdots B_{k-1}G_k,\quad\text{with}\quad G_j:=G^{z_j}(w_j),\, z_j\in\C,\, w_j\in\C\setminus\R,\, j\in [k],
\label{eq:kG_chain}
\end{equation}
for any $k\in \N$ and observables $B_1,\ldots,B_{k-1}$ with the $2\times 2$ block-constant structure. There are two instances where these resolvent chains appear in the proof of Theorem~\ref{theo:main1}. First, they arise as building blocks of covariances from the hierarchy obtained from the iterative underline expansions similar to \eqref{eq:init_expansion}. Those ones, however, contain only resolvents of the same type\footnote{Since Proposition~\ref{prop:Cov} concerns two resolvents $G_1$ and $G_2$, one may think that both of them may appear in the same product, which has to be further expanded similarly to \eqref{eq:init_expansion}. However, this does not happen and only the error terms contain the products of $G_1$ and $G_2$ under the same trace. These terms are reduced to two-resolvent quantities (see Footnote \ref{ftn:GG}) and then estimated by Proposition \ref{prop:2G_av} without involving local laws for products of more than two resolvents.} , i.e. $G_j=G_1$ for all $j\in [k]$. Second, local laws for \eqref{eq:kG_chain} with $k\le 4$ and at most two types of resolvents are used as an input for the proof of Proposition \ref{prop:2G_av} for $b>0$. In particular, we could restrict attention to the case when there are at most two distinct Hermitization parameters in~\eqref{eq:kG_chain}, but this would not yield any simplifications. However, importantly, in both applications of the concentration bounds for \eqref{eq:kG_chain} it is unnecessary to include the refined control parameter~$\widehat{\beta}_{12}$ in the local law estimates. Instead, we carry out all estimates purely in terms of
\begin{equation}
\eta_*:=\eta_1\wedge\cdots \wedge\eta_k\wedge 1,\quad\text{where}\quad \eta_j:=|\Im w_j|,\, j\in [k].
\label{eq:def_eta*_multiG}
\end{equation}
This is a major simplification compared to the set-up of Proposition~\ref{prop:2G_av}, even though the number of resolvents increased.

To state local laws for resolvent chains of the form \eqref{eq:kG_chain}, we first introduce the notion of the deterministic approximation of \eqref{eq:kG_chain}. As in~\eqref{eq:def_M12}, we use interchangeably the notations
\begin{equation*}
M_{[k]}^{B_1,\ldots,B_k},\quad M_{[k]}^{B_1,\ldots,B_k}(w_1,\ldots, w_k)\quad \text{and}\quad\mathcal{M}\left[G_1B_1\cdots B_{k-1}G_k\right]
\end{equation*}
for the deterministic approximation of $G_1B_1\cdots B_{k-1}G_k$, and by induction on $k$ define
\begin{equation}
M_{[k]}^{B_1,\ldots,B_{k-1}}:=\mathcal{B}_{1k}^{-1}\bigg[M_1B_1 M_{[2,k]}^{B_2,\ldots,B_{k-1}} + \sum_{j=2}^{k-1} M_1\mathcal{S}\left[ M_{[1,j]}^{B_1,\ldots,B_{j-1}}\right] M_{[j,k]}^{B_j,\ldots,B_{k-1}}\bigg],
\label{eq:def_multiM}
\end{equation}
using \eqref{eq:def_M12} as the starting point. The size of the lhs. of \eqref{eq:def_multiM} is governed by the following bound.

\begin{proposition}\label{prop:multiM_bound} Fix $\delta,\kappa>0$ and $k\in\N$. Recall the definition of $\eta_*$ from \eqref{eq:def_eta*_multiG}. It holds that
\begin{equation}
\left\lVert M_{[k]}^{B_1,\ldots,B_{k-1}}\right\rVert \lesssim \frac{1}{\eta_*^{k-1}},
\label{eq:multiM_bound}
\end{equation}
uniformly in Hermitization parameters $z_l\in (1-\delta)\mathbf{D}$, spectral parameters $w_j\in\C\setminus\R$ with $\Re w_l\in \mathbf{B}_\kappa^z$, $l\in[k]$, and observables $B_l\in \mathrm{span}\{E_\pm,F^{(*)}\}$, $l\in[k-1]$.
\end{proposition}

Now we present the multi-resolvent local laws for general i.i.d. matrices $X$ in the bulk regime.

\begin{proposition}[Multi-resolvent local laws in the bulk regime]\label{prop:multiG_oneG}
Let $X$ be a complex $N\times N$ i.i.d. matrix satisfying Assumption~\ref{ass:chi}. Fix $\delta,\epsilon,\kappa>0$. Then we have the \emph{averaged} local law
\begin{equation}
\left\vert \left\langle \left(G^{z_1}(w_1)B_1\cdots G^{z_k}(w_k) - M_{[k]}^{B_1,\ldots,B_{k-1}}(w_1,\ldots,w_k)\right)B_k\right\rangle\right\vert \prec \frac{1}{N\eta_*^k},
\label{eq:multiG_av}
\end{equation}
and the \emph{isotropic} local law
\begin{equation}
\left\vert \left\langle \bm{x},\left(G^{z_1}(w_1)B_1\cdots G^{z_k}(w_k) - M_{[k]}^{B_1,\ldots,B_{k-1}}(w_1,\ldots,w_k)\right)\bm{y}\right\rangle\right\vert \prec \frac{1}{\sqrt{N\eta_*}\eta_*^{k-1}},
\label{eq:multiG_iso}
\end{equation}
uniformly in the Hermitization parameters $z_l\in( 1-\delta)\mathbf{D}$, spectral parameters $w_l\in\C\setminus\R$ with $\eta_l:=|\Im w_l|\in [N^{-1+\epsilon}, N^{100}]$ and $\Re w_l\in\mathbf{B}_\kappa^z$, $l\in [k]$, observables $B_l\in\mathrm{span}\{E_\pm,F^{(*)}\}$ for $l\in[k]$, and deterministic vectors $\bm{x},\bm{y}\in\C^{2N}$ with $\|\bm{x}\|=\|\bm{y}\|=1$.
\end{proposition}

The proofs of Propositions~\ref{prop:multiM_bound} and \ref{prop:multiG_oneG} are presented in Section~\ref{sec:multiG_oneG}. We remark that similarly to~\eqref{eq:2G_av}, the bound~\eqref{eq:multiG_av} improves upon the deterministic size bound \eqref{eq:multiM_bound} by the small factor $(N\eta_*)^{-1}$. Meanwhile, the bound in the isotropic local law \eqref{eq:multiG_iso} yields an improvement upon \eqref{eq:multiM_bound} by a factor of $(N\eta_*)^{-1/2}$. These gains are optimal and were first observed in the context of Wigner matrices in~\cite[Theorem 2.5]{Multi_res_llaws}. 

The special case of Proposition~\ref{prop:multiG_oneG} when all Hermitization parameters coincide and the spectral parameters lie in a narrow cone with vertex at origin is proven in \cite[Theorem 1.3]{sparse_universality_Osman}, however our proof does not rely on this work.


\section{Calculations with the Girko's formula}\label{sec:Girko_calculations}

In this section we prove Propositions~\ref{prop:E_omega} and~\ref{prop:Var_main}. First, in Section~\ref{sec:Girko_reduct} we establish the bounds on the regularization terms $J_T, I_0^{\eta_L}$ and on the term $I_{\eta_L}^{\eta_0}$ in~\eqref{eq:Girko_I}. This is done for a general smooth test function to fit the set-up both of Propositions~\ref{prop:E_omega} and~\ref{prop:Var_main}, which are subsequently proven in Sections~\ref{sec:E} and~\ref{sec:Var}, respectively.

Throughout the entire Section~\ref{sec:Girko_calculations} we use~$\xi>0$ to denote an exponent which can be taken arbitrarily small, the exact value of~$\xi$ may change from line to line.

\subsection{Preliminary reductions in the Girko's formula}\label{sec:Girko_reduct}

We present a standard bound on $J_T$ and $I_0^{\eta_L}$, implying that these terms are negligible in the rhs. of~\eqref{eq:Girko_I}. Additionally we prove an upper bound on the moments of $I_{\eta_L}^{\eta_c}$ relying on the tail bound estimate for $\lambda_1^z$ from Proposition \ref{prop:tail_bound}. Though for the proof of Propositions \ref{prop:E_omega} and \ref{prop:Var_main} one needs only the bound on the first two moments of $J_T$, $I_0^{\eta_L}$ and $I_{\eta_L}^{\eta_0}$, we consider the general $p$-th moment of these random variables, since this does not cause any additional difficulties.
\begin{lemma}\label{lem:Girko_reduction} Let $X$ be a complex $N\times N$ i.i.d. matrix satisfying Assumption \ref{ass:chi}. Fix $\delta>0$ and let $g_N\in C^2(\C)$ be a complex-valued function with ${\rm{supp}}(g_N)\subset (1-\delta)\mathbf{D}$. Suppose that there exists $D>0$ such that $\|g_N\|_\infty\!+\|\Delta g_N\|_1\lesssim N^{D}$ for all $N\in\N$. Then for any fixed $p\in\N$ and $K>0$ there exist $L,D'>0$ such that
\begin{equation}
\E |J_T(g_N)|^p + \E |I_0^{\eta_L}(g_N)|^p \lesssim N^{-K},
\label{eq:J}
\end{equation}
where $\eta_L=N^{-L}$ and $T=N^{D'}$. We further have
\begin{equation}
\E\left\vert I_{\eta_L}^{\eta_0}(g_N)\right\vert^p \lesssim \left((N\eta_0)^2 + N^{-2\nu_1}\right) \|\Delta g_N\|_1^pN^\xi \label{eq:I_small_eta}
\end{equation}
for any $\eta_0\in (\eta_L,N^{-1})$, where $\nu_1:=1/10$ is given by Proposition~\ref{prop:tail_bound}. The implicit constants in~\eqref{eq:J} and~\eqref{eq:I_small_eta} depend only on $\delta$, $p$, $K$ and the model parameters in Assumption~\ref{ass:chi}.
\end{lemma}

\begin{proof}[Proof of Lemma \ref{lem:Girko_reduction}] We start with the proof of \eqref{eq:J} and closely follow the proof of similar estimates in \cite[Lemma~3]{nonherm_edge_univ}. From \cite[Eq.(5.27)--(5.28)]{loc_circ_law18} we have that
\begin{equation}
|J_T(g_N)|=\left\vert \frac{1}{4\pi} \int_\C \Delta g_N(z)\log \vert \det (H^z-\ii T)\vert\dif^2 z\right\vert \prec \frac{N\lVert \Delta g_N\rVert_{L^1}}{T^2}
\label{eq:J_prec}
\end{equation}
with very high probability. By choosing $T:=N^{D'}$ for a sufficiently large fixed $D'>0$, we make the rhs. of \eqref{eq:J_prec} smaller than $N^{-K}$. On the small probability event where \eqref{eq:J_prec} does not hold, we use the trivial bound
\begin{equation}
|J_T(g_N)| \lesssim N\lVert \Delta g_N\rVert_{L^1}\log T,
\label{eq:J_a_priori}
\end{equation}
Together with \eqref{eq:J_prec} this gives $\E |J_T(g_N)|^p \lesssim N^{-K}$.

To prove the upper bound on $\E\left\vert I_0^{\eta_L}(g_N)\right\vert^p$ from \eqref{eq:J}, we first estimate $\langle \Im G^z(\ii\eta)\rangle\le 1/\eta$, use \eqref{eq:J_a_priori} and rearrange \eqref{eq:Girko} to express $I_0^{\eta_L}(g_N)$. We get the following bound on $I_0^{\eta_L}(g_N)$ which holds on the entire probability space:
\begin{equation*}
|I_0^{\eta_L}(g_N)|\! \lesssim\! \left\vert \mathcal{L}_N(g_N)\right\vert + N\|\Delta g_N\|_1\log N + N\|\Delta g_N\|_1\!\!\int_{N^{-L}}^T \frac{\dif\eta}{\eta}\!\lesssim\! N\|g_N\|_\infty + N\|\Delta g_N\|_1\log N\lesssim N^{D+2}.
\end{equation*}
Therefore, it is sufficient to prove that $\E \left\vert I_0^{\eta_L}(g_N)\right\vert \lesssim N^{-K'}$ for some fixed  $K'>K+(p-1)(D+2)$. We have
\begin{equation}
\E \left\vert I_0^{\eta_L}(g_N)\right\vert = \E\int_\C |\Delta g_N(z)|\int_0^{\eta_L}\sum_{i=1}^N \frac{\eta}{(\lambda_i^z)^2+\eta^2}\dif\eta\dif^2 z\lesssim N\|g_N\|_1\sup_{|z|\le 1-\delta} \E\log \left(1+(N^{L}\lambda^z_1)^{-2}\right).
\label{eq:tiny_eta}
\end{equation}
By \cite[Lemma~4.12]{bordenave2012around} (see also \cite[Proposition~5.7]{loc_circ_law18}), under the regularity assumption \eqref{eq:reg_ass} it holds that 
\begin{equation}
\mathbf{P}[\lambda_1^z<N^{-1}u] \le C_{\mathfrak{a}} u^{2\mathfrak{a}/(1+\mathfrak{a})}N^{\mathfrak{b}+1},
\label{eq:bordenave_tail}
\end{equation}
uniformly in $|z|\le 1$ and $u>0$, for some constant $C_{\mathfrak{a}}>0$. We proceed by a standard dyadic argument and decompose the probability space into the events $2^kN^{-L}<\lambda_1^z\le 2^{k+1}N^{-L}$, $k\in\mathbf{Z}$, whose probabilities are then estimated by \eqref{eq:bordenave_tail}.  We thus conclude that one can choose $L>0$ sufficiently large in the definition of $\eta_L=N^{-L}$, so that the rhs. of~\eqref{eq:tiny_eta} is smaller than $N^{-K'}$. This finishes the proof of~\eqref{eq:J}. 

Now we prove \eqref{eq:I_small_eta}. Denote $I:=I_{\eta_L}^{\eta_0}(g_N)$ and fix a (small) $\xi_0>0$. We use \eqref{eq:G_to_ImG} to go from $\langle G^z(\ii\eta)\rangle$ to the sum over eigenvalues in $I$ and split this summation into two regimes:
\begin{equation}
I= I_1+I_2 := - \frac{N}{2\pi}\int_\C \Delta g_N(z) \dif^2 z \int_{\eta_L}^{\eta_0} \frac{1}{N} \bigg(\sum_{i\le N^{\xi_0}}+\sum_{i> N^{\xi_0}}\bigg) \frac{\eta}{(\lambda_i^z)^2+\eta^2} \dif \eta.
\label{eq:small_eta_splitting}
\end{equation}
In \eqref{eq:small_eta_splitting}, as well as everywhere further in this proof, the summation index $i$ takes values only from the set $[1,N]$ even if this is not mentioned explicitly. To estimate $\E |I_1|^p$, we observe that for any random variable $h^z$ and any $p\in\N$ it holds that
\begin{equation}
\begin{split}
\E \left\vert \int_\C \Delta g_N(z)h^z\dif^2 z\right\vert^p &\le \int_{\C^p} \prod_{j=1}^p |\Delta g_N(z_j)|\E\Bigg[\prod_{j=1}^p |h^{z_j}|\Bigg]\prod_{j=1}^p \dif^2 z_j\\
&\le \sum_{l=1}^p \int_{\C^p} \prod_{j=1}^p |\Delta g_N(z_j)|\E\left[ |h^{z_l}|^p\right]\prod_{j=1}^p \dif^2 z_j \le p \|\Delta g_N\|_1^p \sup_{|z|\le 1-\delta} \E |h^z|^p, 
\end{split}
\label{eq:pth_moment}
\end{equation}
where in the last bound we used that $\mathrm{supp}(\Delta g_N) \subset (1-\delta)\mathbf{D}$. We use \eqref{eq:pth_moment} for
\begin{equation}
h^z_1:= \int_{\eta_L}^{\eta_0} \sum_{i\le N^{\xi_0}}\frac{\eta}{(\lambda_i^z)^2+\eta^2} \dif \eta.
\label{eq:def_h1}
\end{equation}
Performing integration in $\eta$ in the rhs. of~\eqref{eq:def_h1} we get
\begin{equation}
h^z_1 \le \sum_{i\le N^{\xi_0}} \log \left(1+\frac{(N\eta_0)^2}{(N\lambda_i^z)^2+(N\eta_L)^2}\right)\le N^{\xi_0} \log \left(1+\frac{(N\eta_0)^2}{(N\lambda_1^z)^2+(N\eta_L)^2}\right).
\label{eq:h1_bound}
\end{equation}
To estimate the $p$-th moment of the logarithm in the rhs. of \eqref{eq:h1_bound} we decompose the probability space into the events $P_0:=\{\lambda_1^z\le N\eta_L\}$ and $P_k:=\{2^{k-1} N\eta_L<\lambda_1^z\le 2^kN\eta_L\}$, $k\in \N$. Since $\eta_L=N^{-L}$, only $k\lesssim \log N$ are involved in this decomposition. We have from Proposition~\ref{prop:tail_bound} that  
\begin{equation}
\mathbf{P}\left[P_k\right] \lesssim N^\xi \left((2^kN\eta_L)^2 + N^{-2\nu_1}\right),\quad \forall k\ge 0,
\label{eq:tail_bound_modified}
\end{equation}
uniformly in $|z|\le 1-\delta$. An elementary calculation based on \eqref{eq:tail_bound_modified} shows that
\begin{equation}
\E \left\vert \log \left(1+\frac{(N\eta_0)^2}{(N\lambda_1^z)^2+(N\eta_L)^2}\right)\right\vert^p\!\lesssim\!\! \sum_{0\le k\lesssim \log N}\!\! \mathbf{P}[P_k] \left\vert \log \left(1+\frac{(N\eta_0)^2}{(2^kN\eta_L)^2}\right)\right\vert^p\!\!\lesssim N^\xi\left((N\eta_0)^2 + N^{-2\nu_1}\right)\!.
\label{eq:X_int_bound}
\end{equation}
Together with \eqref{eq:pth_moment} applied to $h^z_1$ and \eqref{eq:h1_bound} this implies
\begin{equation}
\E|I_1|^p \lesssim N^{p\xi_0+\xi}((N\eta_0)^2+N^{-2\nu_1})\lVert \Delta g_N\rVert_1^p.
\label{eq:I1_bound}
\end{equation}

Next we estimate $\E|I_2|^p$. Similarly to~\eqref{eq:def_h1} we take
\begin{equation}
h^z_2:= \int_{\eta_L}^{\eta_0} \sum_{i> N^{\xi_0}}\frac{\eta}{(\lambda_i^z)^2+\eta^2} \dif \eta.
\label{eq:def_h2}
\end{equation}
However, instead of immediately performing integration in the rhs. of~\eqref{eq:def_h2}, we first bound the sum over~$i$ for any~$\eta>0$:
\begin{equation}
\sum_{|i|\ge N^{\xi_0}}\frac{\eta}{(\lambda_i^z)^2+\eta^2} \lesssim \sum_{|i|\ge N^{\xi_0}}\frac{\eta N^\xi}{(\lambda_i^z)^2 + (N^{-1+\xi_0})^2} \le \eta N^{2+\xi-\xi_0}\langle \Im G^z(\ii N^{-1+\xi_0})\rangle\lesssim \eta N^{2+\xi-\xi_0},
\label{eq:h2_bound}
\end{equation}
with very high probability uniformly in $|z|\le 1-\delta$. Here we used the eigenvalue rigidity \cite[Eq.(7.4)]{macroCLT_complex} in the first step and the averaged single-resolvent local law \cite[Eq.(5.4)]{loc_circ_law18} in the last step to bound $|\langle \Im G^z\rangle|\lesssim 1$. Together with the trivial bound $h^z_2\lesssim N\log N$ which holds with probability~1, \eqref{eq:h2_bound} implies
\begin{equation}
\E |h_2^z|^p \lesssim N^{p(\xi-\xi_0)}\left\vert N^2\int_{\eta_L}^{\eta_0}\eta\dif \eta\right\vert^p\lesssim N^{p(\xi-\xi_0)} (N\eta_0)^{2p}.
\label{eq:Eh2}
\end{equation}
Using~\eqref{eq:Eh2} along with~\eqref{eq:pth_moment} applied to~$h^z_2$ we get
\begin{equation}
\E|I_2|^p \lesssim N^{p(\xi-\xi_0)} (N\eta_0)^{2p} \|\Delta g_N\|_1^p.
\label{eq:I2_bound}
\end{equation}
Finally, we combine~\eqref{eq:I1_bound} with~\eqref{eq:I2_bound} and choose $\xi_0>0$ to be sufficiently small. This finishes the proof of Lemma~\ref{lem:Girko_reduction}.
\end{proof}

\subsection{Calculation of the expectation: proof of Proposition \ref{prop:E_omega}}\label{sec:E} In this section we denote $\omega_N:=\omega_{a,N}^{(z_0)}$ and omit the arguments from the notations introduced in~\eqref{eq:def_IJ}, which are meant to be equal to $\omega_N$ throughout the proof. Since $|z_0|<1-\delta$, it holds that ${\rm supp}(\omega_N)\subset (1-\delta/2)\mathbf{D}$ for sufficiently large~$N$. From the definition of~$\omega_N$ given in~\eqref{eq:omega_z0} we have that 
\begin{equation}
\|\omega_N\|_\infty\lesssim N^{2a}\quad\text{and}\quad \|\Delta\omega_N\|_1\lesssim N^{2a}.
\label{eq:omega_bounds}
\end{equation}
Let $\delta_0\in (0,1/2+\nu_0-a)$ and $\delta_c>0$ be positive exponents which will be taken sufficiently small at the end. We take 
\begin{equation}
\eta_L:=N^{-L},\quad \eta_0:=N^{-1/2-a-\delta_0},\quad \eta_c:=N^{-1+\delta_c}\quad\text{and}\quad T:=N^{D'}
\label{eq:E_eta_scales}
\end{equation}
for sufficiently large $L,D'>0$. Actually, the choice of $\eta_L$, $\eta_c$ and~$T$ was already made in Section~\ref{sec:results_list} and is simply recalled here. However, the choice of~$\eta_0$ is specific to the proof of Proposition~\ref{prop:E_omega}. We have from Lemma~\ref{lem:Girko_reduction} for~$p=1$ that 
\begin{equation}
|\E J_T| + |\E I_{0}^{\eta_L}|+|\E I_{\eta_L}^{\eta_0}|\lesssim ((N\eta_0)^2+N^{-2\nu_1})N^{2a+\xi}\lesssim N^{1-\delta_0+\xi}.
\label{eq:E_cut_tails}
\end{equation}
Therefore, the contribution from the terms in the lhs. of~\eqref{eq:E_cut_tails} to $\E\Lin(\omega_N)$ can be incorporated into the error term in~\eqref{eq:E_omega}. We further analyze separately $I_{\eta_c}^T$ and $I_{\eta_0}^{\eta_c}$. First we compute $I_{\eta_c}^T$ up to the leading order. As we will see, this is the main contributing regime. Afterwards we show that $I_{\eta_0}^{\eta_c}$ is negligible by comparison with GinUE.

\medskip

\noindent\underline{Analysis of $I_{\eta_c}^T$.} Performing integration by parts in $I_{\eta_c}^T$, first with respect to $z$ and then to $\overline{z}$, we get 
\begin{equation}
\begin{split}
I_{\eta_c}^T &= -\frac{N}{2\pi \ii}\int_\C \Delta \omega_N(z) \int_{\eta_c}^T \langle G^z(\ii\eta)\rangle\dif\eta = \frac{2N}{\pi\ii} \int_\C\partial_{\bar{z}} \omega_N(z)\dif^2 z\int_{\eta_c}^T \partial_z \langle G^z(\ii\eta)\rangle\dif\eta\\
& = \frac{2N}{\pi}\int_\C \partial_{\bar{z}} \omega_N(z) \langle G^z(\ii\eta_c)F\rangle\dif^2 z + \mathcal{O}\left(\frac{N^{1+a}}{T}\right)\\
&=-\frac{2N}{\pi}\int_\C \omega_N(z) \langle G^z(\ii\eta_c)F^*G^z(\ii\eta_c)F\rangle \dif^2 z +\mathcal{O}\left(\frac{N}{T}\right).
\end{split}
\label{eq:E_by_parts}
\end{equation}
Here to go from the first to the second line we used that $\partial_z \langle G^z(\ii\eta)\rangle = -\ii\partial_\eta \langle G^z(\ii\eta)F\rangle$, $\lVert \partial_{\bar{z}} \omega_N\rVert_{1}\lesssim N^a$ and $\| G^z(\ii T)\|\le T^{-1}$. Since $F$ and $F^*$ are off-diagonal matrices, they are regular in the sense of \cite[Definition 3.1]{non-herm_overlaps}, so the two-resolvent local law with regular observables \cite[Theorem~4.4]{non-herm_overlaps} implies that
\begin{equation}
\left\langle \big(G^z(\ii\eta_c)F^*G^z(\ii\eta_c)- M_{11}^{F^*}\big)F\right\rangle\prec (N\eta_c)^{-1/2},
\label{eq:2G_reg}
\end{equation}
where $M_{11}^{F^*}$ is the deterministic approximation to $G^z(\ii\eta_c)F^*G^z(\ii\eta_c)$ (for the definition see~\eqref{eq:def_M12} with $M_1=M_2=M^z(\ii\eta_c)$). By a straightforward explicit calculation using \eqref{eq:M}, \eqref{eq:m}, \eqref{eq:def_M12} and~\eqref{eq:def_B12} we have
\begin{equation}
\langle M_{11}^{F^*}F\rangle = \frac{m^2}{2}\cdot \frac{1-m^2+|z|^2u^2}{1-m^2-|z|^2u^2} =-\frac{1}{2}+\mathcal{O}(\eta_c),
\label{eq:MFMF_expl}
\end{equation}
with $m=m^z(\ii\eta_c)$ and $u=u^z(\ii\eta_c)$. Therefore, combining \eqref{eq:E_by_parts}, \eqref{eq:2G_reg}, and \eqref{eq:MFMF_expl} we conclude that
\begin{equation}
I_{\eta_c}^T = \frac{N}{\pi}\int_\C \omega_N(z)\dif^2 z \left(1+\mathcal{O}(N^{-\delta_c/2})\right)=\frac{N}{\pi} \left(1+\mathcal{O}(N^{-\delta_c/2})\right),
\label{eq:E_large_eta}
\end{equation}
where we also used that the integral of $\omega_N(z)$ over $\C$ equals to~1 due to the scaling~\eqref{eq:omega_z0}.

\medskip

\noindent\underline{Analysis of $I_{\eta_0}^{\eta_c}$.} Let $\widetilde{X}$ be an $N\times N$ GinUE matrix. Denote the $\widetilde{X}$-counterparts of $\Lin_N$, $I_{\eta_l}^{\eta_u}$, and $G^z(w)$ by $\widetilde{\Lin}_N$, $\widetilde{I}_{\eta_l}^{\eta_u}$, and $\widetilde{G}^z(w)$, respectively, where $0\le \eta_l\le\eta_u$. Since $\mathrm{supp}(\omega_N)\subset (1-\delta/2)\mathbf{D}$ and the one-particle density of the eigenvalue process of GinUE converges to $(2\pi)^{-1}$ with exponential speed in the bulk \cite[Eq.(1.44)]{Ginibre65}, i.e. for  $|z|\le (1-\delta/2)$, \eqref{eq:E_omega} holds for $\E \widetilde{\Lin}_N(\omega_N)$ with an error term bounded by~$\ee^{-cN}$ for some~$c>0$ independent of~$N$. Together with~\eqref{eq:E_cut_tails} and~\eqref{eq:E_large_eta} applied to~$\widetilde{X}$, this gives 
\begin{equation}
\left\vert\E \widetilde{I}_{\eta_0}^{\eta_c}\right\vert\le \left\vert \E\widetilde{\Lin}_N(\omega_N)-\E \widetilde{I}_{\eta_c}^T\right\vert + \left\vert\E \widetilde{J}_T\right\vert + \left\vert\E \widetilde{I}_{0}^{\eta_L}\right\vert+\left\vert\E \widetilde{I}_{\eta_L}^{\eta_0}\right\vert \lesssim N^{1-\delta_0+\xi} + N^{1-\delta_c/2}.
\label{eq:I_Gin_bound}
\end{equation}
Finally, we compare $\E I_{\eta_0}^{\eta_c}$ with $\E \widetilde{I}_{\eta_0}^{\eta_c}$ by the means of Proposition \ref{prop:EG}, which is applicable provided that $\eta_0\ge N^{-3/2+\epsilon}$ for some fixed $\epsilon>0$:
\begin{equation}
\left\vert\E \left[ I_{\eta_0}^{\eta_c}-\widetilde{I}_{\eta_0}^{\eta_c}\right]\right\vert \!\lesssim\! N\!\!\int_\C |\Delta \omega_N(z)|\dif^2 z\! \int_{\eta_0}^{\eta_c} \left\vert \E \left[\langle G^z(\ii\eta)\rangle-\langle \widetilde{G}^z(\ii\eta)\rangle\right]\right\vert\dif \eta\lesssim N^{1+\xi}\|\Delta \omega_N\|_1\int_{\eta_0}^{\eta_c}\Phi_1(\eta)\dif\eta.
\label{eq:I_Gin_comparison}
\end{equation}
In the regime $\eta\in [N^{-1}, \eta_c]$ we simply estimate $\Phi_1(\eta)$ from above by taking $t:=N^{-2/3}$ in \eqref{eq:def_Phi}:
\begin{equation}
\Phi_1(\eta)\lesssim N^{-1/6}N\eta \le N^{-1/6+\delta_c}.
\label{eq:Phi1_bound1}
\end{equation}
For $\eta\in [N^{-1-\nu_0},N^{-1})$ an elementary optimization of the rhs. of the first line of \eqref{eq:def_Phi} in $t$ shows that this quantity attains its minimum for $t\sim N^{-2/3}(N\eta)^{8/3}$, which gives
\begin{equation}
\Phi_1(\eta)\sim N^{-1/6}(N\eta)^{-4/3}.
\label{eq:Phi1_bound2}
\end{equation}
In fact, $\Phi_1(\eta)$ remains much smaller than the local law estimate of order $(N\eta)^{-1}$ on the error term in~\eqref{eq:EG_comp} for a certain range of $\eta<N^{-1-\nu_0}$. However, we do not explore this regime, since the rhs. of~\eqref{eq:I_Gin_comparison} is of order~$N$ when $\eta_0\sim N^{-1-\nu_0}$, as we will now see. Thus, for $\eta_0\ll N^{-1-\nu_0}$ the rhs. of \eqref{eq:I_Gin_comparison} is even much larger than~$N$, which would contribute an error term to the rhs. of~\eqref{eq:E_omega} exceeding the leading order. Together with the choice of~$\eta_0$ in~\eqref{eq:E_eta_scales} this explains the constraint $a<1+\nu_0$ imposed in Proposition~\ref{prop:E_omega}.

Using \eqref{eq:Phi1_bound1} for $\eta\in [N^{-1},\eta_c]$, \eqref{eq:Phi1_bound2} for $\eta\in [N^{-1-\nu_0},N^{-1}]$, and recalling the definition of~$\eta_0$ from~\eqref{eq:E_eta_scales}, we get
\begin{equation}
\int_{\eta_0}^{\eta_c} \Phi_1(\eta)\dif \eta\lesssim N^{-7/6+2\delta_c} + N^{-1/6-4/3}\eta_0^{-1/3}\le N^{-1+2\delta_c+\delta_0/3}\left(N^{-1/6}+N^{-1/3(1-a)}\right)
\label{eq:Phi1_integrated}
\end{equation}
Finally, we combine \eqref{eq:I_Gin_comparison}, \eqref{eq:Phi1_integrated}, \eqref{eq:omega_bounds} ,and use that $a\ge 1/2$:
\begin{equation}
\left\vert\E \left[ I_{\eta_0}^{\eta_c}-\widetilde{I}_{\eta_0}^{\eta_c}\right]\right\vert \lesssim N^{1+(a-1/2-\nu_0)7/3 +2\delta_c+\delta_0/3+\xi}.
\label{eq:I_comparison_fin}
\end{equation}
Importantly, $a<1/2+\nu_0$, so the rhs. of \eqref{eq:I_comparison_fin} is much smaller than $N$ for sufficiently small $\delta_c,\delta_0,\xi>0$.

Using \eqref{eq:Girko_I}, collecting the error terms from \eqref{eq:E_cut_tails}, \eqref{eq:E_large_eta}, \eqref{eq:I_Gin_bound}, \eqref{eq:I_comparison_fin}, and taking $\delta_0$ and $\delta_c$ sufficiently small, we finish the proof of Proposition \ref{prop:E_omega}.

\subsection{Upper bound on the variance: proof of Proposition \ref{prop:Var_main}}\label{sec:Var}

Throughout this section we abbreviate $f_N:=f_{a,N}$. We start with introducing some additional notation. For a set $S\subset\C$ and~$h>0$ denote the $h$-neighborhood of~$S$ by
\begin{equation}
{\mathfrak{N}}_h(S):= \{ z\in \C\, :\, \mathrm{d}(z,S)<h\}
\label{eq:def_Nh}
\end{equation}
where $\mathrm{d}$ stands for the Euclidean distance. In particular, for $z\in\C$ we denote an open ball with the center at~$z$ and radius~$r>0$ by $B_r(z):={\mathfrak{N}}_r(\{z\})$. Since~$\omega$ is compactly supported, there exists~$\mathcal{R}>0$ such that ${\rm{supp}}(\omega)\subset B_\mathcal{R}(0)$. Denote the \emph{tubular neighborhood} of $\partial\Omega_N$ of radius $\mathcal{R}N^{-a}$ by\footnote{We point out that here $\Delta$ denotes the tubular neighborhood of a curve, and it should not be confused with the Laplacian.}
\begin{equation}
\Delta\Omega_N:= {\mathfrak{N}}_{\mathcal{R}N^{-a}}(\partial\Omega_N).
\end{equation}
By construction of $f_N$ in \eqref{eq:def_f}, $f_N(z)=1$ for $z\in\Omega_N\setminus\Delta \Omega_N$ and $f_N(z)=0$ for  $z\in\C\setminus\left(\Omega_N\cup\Delta \Omega_N\right)$. Therefore, only $z\in \Delta\Omega_N$ contribute to the rhs. of the Girko's formula \eqref{eq:Girko}. For the further reference we note that
\begin{equation}
|\Delta\Omega_N|\lesssim |\partial \Omega_N|\mathcal{R}N^{-a}\sim N^{-\alpha-a},
\label{eq:Delta_Omega_area}
\end{equation}
since $\Omega_N$ satisfies Assumption \ref{ass:Omega}. We additionally observe that $\lVert f_N\rVert_\infty\lesssim 1$ and
\begin{equation}
\lVert \Delta f_N\rVert_{1} \le |\Delta \Omega_N|\|\Delta f_N\|_\infty \lesssim N^{-a-\alpha}N^{2a}=N^{a-\alpha}.
\label{eq:Delta_f_L1}
\end{equation}

As the last preparation for the proof of Proposition~\ref{prop:Var_main}, we present a bound on the integrals of regularized singularities of the form $|z_1-z_2|^{-q}$, $q> 0$, over $(\Delta\Omega_N)^2$. Later these estimates will be used for $q~\!=~\!2q_0,2,4$, where $q_0=1/20$ is defined in Proposition~\ref{prop:Var_main}. 
\begin{lemma}\label{lem:int} Suppose that $\Omega_N$ satisfies Assumption~\ref{ass:Omega}. Then for any fixed exponent $q>0$ and a sequence $\{r_N\}\subset (0,+\infty)$ it holds that 
\begin{equation}
\int_{\Delta\Omega_N}\int_{\Delta\Omega_N} \frac{1}{|z_1-z_2|^q + r_N}\dif^2 z_1\dif^2 z_2 \lesssim N^{-2a}\begin{cases}
N^{-\alpha}r_N^{-\frac{q-1}{q}} \left(|\log r_N|+\log N\right),&\text{if}\,\, q\ge 1,\\
N^{-\alpha(2-q)},&\text{if}\,\, 0<q< 1,
\end{cases}
\label{eq:int}
\end{equation}
where the implicit constant depends only on the model parameters from Assumption \ref{ass:Omega}.
\end{lemma}

The proof of Lemma \ref{lem:int} is an elementary calculus exercise using that $\Delta \Omega_N$ is a tubular neighborhood of a $C^2$ curve $\partial\Omega_N$, we omit further details.

Now we are ready to prove Proposition~\ref{prop:Var_main}. The structure of the proof is similar to the one presented  in Section~\ref{sec:E}. First we note that ${\rm supp}(f_N)\subset (1-\delta/2)\mathbf{D}$ for sufficiently large~$N$, so applying Lemma~\ref{lem:Girko_reduction} for $p=2$ along with the bound~\eqref{eq:Delta_f_L1} we get
\begin{equation}
\E |J_T|^2 + \E |I_0^{\eta_L}|^2 + \E |I_{\eta_L}^{\eta_0}|^2 \lesssim \left((N\eta_0)^2+N^{-2\nu_1}\right) \|\Delta f_N\|_1^2 N^\xi\lesssim \left((N\eta_0)^2+N^{-2\nu_1}\right) N^{2(a-\alpha)+\xi}
\label{eq:Var_tails}
\end{equation}
for sufficiently large $L,D'>0$ and $\eta_L=:N^{-L}, T:=N^{D'}$. We further take $\eta_c:=N^{-1+\delta_c}$ for a small $\delta_c>0$, and $\eta_0\in [N^{-1-\nu_1},N^{-1})$. The exact choice of~$\eta_0$ will be specified at the end. The Girko's formula~\eqref{eq:Girko_I} implies that
\begin{equation}
\Var\left[\Lin_N(f_N)\right] \lesssim \E|J_T|^2 + \E|I_0^{\eta_L}|^2 + \E|I_{\eta_L}^{\eta_0}|^2+\Var\left[I_{\eta_0}^{\eta_c}\right]+\Var\left[I_{\eta_c}^T\right].
\label{eq:Var_I_bound}
\end{equation}
The bound \eqref{eq:Var_tails} takes care of the first three terms in \eqref{eq:Var_I_bound}, now we proceed to the analysis of the last two terms in the rhs. of \eqref{eq:Var_I_bound}.

\medskip

\noindent\underline{Analysis of $\Var\left[I_{\eta_c}^{T}\right]$.} We write the variance of $I_{\eta_c}^T$ as in \eqref{eq:Var_I_ll} and use the formula for the covariance of resolvent traces from \eqref{eq:main_cov}. To bound the contribution from the error term in \eqref{eq:main_cov} to $\Var\left[I_{\eta_c}^{T}\right]$, we estimate $\gamma\ge |z_1-z_2|^2 +\eta_c$ and perform the integration over $\eta_1,\eta_2$:
\begin{equation}
\begin{split}
&\left(\frac{N}{2\pi}\right)^2 \int_\C\int_\C |\Delta f_N(z_1)||\Delta f_N(z_2)| \left(\int_{\eta_c}^T\int_{\eta_c}^T \left(\frac{1}{N\gamma}+N^{-1/4}\right)\frac{N^\xi}{N^2\eta_1\eta_2}\dif\eta_1\dif\eta_2\right)\dif^2 z_1\dif^2 z_2\\
&\quad\lesssim \|\Delta f_N\|_\infty^2 (\log N)^2 N^\xi \int_{\Delta \Omega_N}\int_{\Delta \Omega_N} \left(\frac{1}{N}\frac{1}{|z_1-z_2|^2 + \eta_c}+N^{-1/4}\right)\dif^2 z_1\dif^2 z_2\\
&\quad \lesssim N^{4a+\xi}\left(N^{-\alpha -2a-1/2-\delta_c/2}+N^{-2\alpha -2a-1/4}\right). 
\end{split}
\label{eq:Cov_err_int}
\end{equation}
To go from the second to the third line we used that $\|\Delta f_N\|_\infty\lesssim N^{2a}$, estimated the area of $\Delta\Omega_N$ by \eqref{eq:Delta_Omega_area}, and applied Lemma \ref{lem:int} to $q=2$ and $r_N=\eta_c$. Next, compute the contribution from the leading term in the rhs. of \eqref{eq:main_cov} using the following explicit formula from \cite[Lemma 4.8, Lemma 4.10]{macroCLT_complex}:
\begin{equation}
\begin{split}
&\left(\frac{N}{2\pi\ii}\right)^2 \int_\C\int_\C \Delta f_N(z_1)\Delta f_N(z_2) \left(\int_0^\infty\int_0^\infty \frac{1}{N^2}\cdot \frac{V_{12}+\kappa_4U_1U_2}{2} \dif\eta_1\dif\eta_2\right)\dif^2 z_1\dif^2 z_2\\
&\quad = \frac{1}{4\pi}\int_C |\nabla f_N(z)|^2\dif^2 z +\frac{\kappa_4}{\pi^2} \left\vert \int_\C f_N(z)\dif^2 z\right\vert^2\lesssim N^{a-\alpha} + N^{-4\alpha}.
\end{split}
\label{eq:expl_UV_calc}
\end{equation}
In the first line of~\eqref{eq:expl_UV_calc} the trivial regimes $\eta_l\in [0,\eta_c]\cup[T,+\infty)$ are included, since the identity stated in~\eqref{eq:expl_UV_calc} is available only for the $\eta$-integrals over the entire regime $[0,+\infty)$. However, these trivial regimes are not present in the rhs. of~\eqref{eq:Var_I_ll}, so we need to eliminate them from~\eqref{eq:expl_UV_calc}. An elementary calculation using \cite[Eq. (4.24)]{macroCLT_complex} shows that 
\begin{equation}
\left\vert\int_\C\int_\C \Delta f_N(z_1)\Delta f_N(z_2) \left(\left(\int_0^{\eta_c}+\int_T^\infty\right)^2 \frac{V_{12}+\kappa_4U_1U_2}{2}\dif\eta_1\dif\eta_2\right)\dif^2 z_1\dif^2 z_2\right\vert \lesssim N^{2a-\alpha -1/2-\delta_c/2+\xi},
\label{eq:tails_UV}
\end{equation}
for the proof see Supplementary Section~\ref{app:cut_tails}. Combining~\eqref{eq:Var_I_ll}, \eqref{eq:Cov_err_int}, \eqref{eq:expl_UV_calc}, \eqref{eq:tails_UV}, we conclude that
\begin{equation}
\Var\left[I_{\eta_c}^T\right] \lesssim N^{2a-\alpha-1/2+\xi} + N^{2a-2\alpha-1/4+\xi},
\label{eq:Var_ll_bound}
\end{equation}
where we also used that $a>1/2$ to show that the error terms coming from the rhs. of~\eqref{eq:expl_UV_calc} are suppressed by the ones in the rhs. of~\eqref{eq:Var_ll_bound}.

\medskip

\noindent\underline{Analysis of $\Var\left[I_{\eta_0}^{\eta_c}\right]$.} First we write the variance of $I_{\eta_0}^{\eta_c}$ similarly to \eqref{eq:Var_I_ll}. In the regime $\eta_1,\eta_2\in [\eta_0,\eta_c]$ the covariance of $\langle G^{z_1}(\ii\eta_1)\rangle$ and $\langle G^{z_2}(\ii\eta_2)\rangle$ is controlled by the means of Proposition \ref{prop:EGG}. However, for $|z_1-z_2|\le N^{-1/2}$ this control is not effective. Indeed, in this case $R$ varies from $0$ to $N|z_1-z_2|^2\le 1$ in the definition of $\Phi_2$ in \eqref{eq:def_Phi2}, so $N\mathcal{E}_1(t,R)\ge \sqrt{Nt}$, which gives that $\Phi_2\gg (N^2\eta_1\eta_2)^{-1}$. Moreover, $\Phi_2$ exceeds $(N^2\eta_1\eta_2)^{-1}$ by a non-negligible power of $N$ as $|z_1-z_2|$ decreases further. So, instead of applying Proposition \ref{prop:EGG} in the regime $|z_1-z_2|\le N^{-1/2}$, we simply use the averaged single-resolvent local law \eqref{eq:av_ll_simple} and get
\begin{equation}
\int_{\eta_0}^{\eta_c}\int_{\eta_0}^{\eta_c}\left\vert\Cov\left(\langle G^{z_1}(\ii\eta_1)\rangle,\langle G^{z_2}(\ii\eta_2)\rangle\right)\right\vert\dif\eta_1\dif\eta_2\lesssim N^\xi \int_{\eta_0}^{\eta_c}\int_{\eta_0}^{\eta_c} \frac{1}{N^2\eta_1\eta_2}\dif\eta_1\dif\eta_2\lesssim N^{-2+\xi},
\label{eq:Cov_ll_eta_int}
\end{equation}
uniformly in $z_1,z_2\in\Delta\Omega_N$. We note that the contribution of the regime $|z_1-z_2|\le N^{-1/2}$ to the analogue of~\eqref{eq:Var_I_ll} for $I_{\eta_0}^{\eta_c}$ gains smallness from the volume of the integration domain, rather than from~\eqref{eq:Cov_ll_eta_int}. In particular,~\eqref{eq:Cov_ll_eta_int} is not affordable for general~$z_1,z_2$.

In the complementary regime $|z_1-z_2|>N^{-1/2}$, Proposition \ref{prop:EGG} already improves beyond the trivial local law bound, and gives that
\begin{equation}
\left\vert\Cov\left(\langle G^{z_1}(\ii\eta_1)\rangle,\langle G^{z_2}(\ii\eta_2)\rangle\right)\right\vert\lesssim N^\xi \Phi_2(\eta_1,\eta_2,|z_1-z_2|)\lesssim  \left(N(\mathcal{E}_1(t,R) + \mathcal{E}_0(t))+\frac{\sqrt{N}t}{(N\eta_*)^3}\right)\frac{N^{6\delta_c+\xi}}{N^2\eta_1\eta_2}, 
\label{eq:Phi2_bound}
\end{equation}
with $\eta_*:=\eta_1\wedge\eta_2$, uniformly in $z_1,z_2\in\Delta\Omega_N$, $\eta_1,\eta_2\in [\eta_0,\eta_c]$ and $t\in [N^{-1+\omega_*},N^{-\omega_*}]$, $R\in [0,N|z_1-z_2|^2]$, for any fixed $\omega_*>0$. We integrate \eqref{eq:Phi2_bound} over $\eta_1,\eta_2\in [\eta_0,\eta_c]$ and obtain
\begin{equation}
\int_{\eta_0}^{\eta_c}\int_{\eta_0}^{\eta_c}\left\vert\Cov\left(\langle G^{z_1}(\ii\eta_1)\rangle,\langle G^{z_2}(\ii\eta_2)\rangle\right)\right\vert\dif\eta_1\dif\eta_2\lesssim N^{-2+8\delta_c+\xi} \left(N(\mathcal{E}_1(t,R)+\mathcal{E}_0(t)) + \frac{\sqrt{N}t}{(N\eta_0)^3}\right).
\label{eq:Phi2_bound_eta_int}
\end{equation}
Next, we optimize the rhs. of \eqref{eq:Phi2_bound_eta_int} over $R$ and $t$, which leads to the choice
\begin{equation}
R:=\left(N|z_1-z_2|^2\right)^{4/5},\quad t:=N^{-1}\left(N|z_1-z_2|^2\right)^{1/10}.
\end{equation}
This gives
\begin{equation}
\int_{\eta_0}^{\eta_c}\int_{\eta_0}^{\eta_c}\left\vert\Cov\left(\langle G^{z_1}(\ii\eta_1)\rangle,\langle G^{z_2}(\ii\eta_2)\rangle\right)\right\vert\dif\eta_1\dif\eta_2\lesssim N^{-2+8\delta_c+\xi} \left(N|z_1-z_2|^2\right)^{-1/20},
\label{eq:Phi2_bound_eta_int_optimized}
\end{equation}
where we additionally used that $N\eta_0\gtrsim (N|z_1-z_2|^2)^{1/20}N^{-1/6}$, since $|z_1-z_2|\lesssim 1$ and $\eta_0\ge N^{-1-\nu_1}$, with $\nu_1=1/10$.

Using the analogue of \eqref{eq:Var_I_ll} for $I_{\eta_0}^{\eta_c}$, and employing \eqref{eq:Cov_ll_eta_int} for $|z_1-z_2|\le N^{-1/2}$ and \eqref{eq:Phi2_bound_eta_int_optimized} for $|z_1-z_2|>N^{-1/2}$, we get 
\begin{equation}
\Var\left[I_{\eta_0}^{\eta_c}\right] \lesssim N^{8\delta_c+\xi}\int_{\Delta \Omega_N}\int_{\Delta \Omega_N} |\Delta f_N(z_1)||\Delta f_N(z_2)| (N|z_1-z_2|^2 +1)^{-1/20}\dif^2 z_1\dif^2 z_2.
\label{eq:Var_I_optimized}
\end{equation}
This integral is further estimated by Lemma \ref{lem:int} with $q:=2q_0=1/10$ and $r_N:=N^{-1/2}$:
\begin{equation}
\Var\left[I_{\eta_0}^{\eta_c}\right] \lesssim N^{8\delta_c+\xi} \|\Delta f_N\|_\infty ^2 N^{-2a-\alpha(2-2q_0)-q_0}\lesssim N^{2a-2\alpha -2q_0(1/2-\alpha) +8\delta_c+\xi},
\label{eq:Var_interm_fin}
\end{equation}
where we also used that $\|\Delta f_N\|_\infty\lesssim N^{2a}$.

Note that the bound \eqref{eq:Var_I_optimized} is insensitive to the exact choice of $\eta_0\in [N^{-1-\nu_1},N^{-1})$, so we choose $\eta_0:=N^{-1-\nu_1}$ to minimize the rhs. of \eqref{eq:Var_tails}. Taking $\delta_c$ sufficiently small, using \eqref{eq:Var_I_bound} and collecting the error terms from \eqref{eq:Var_tails}, \eqref{eq:Var_ll_bound}, \eqref{eq:Var_I_optimized}, we get
\begin{equation}
\Var\left[\mathcal{L}_N(f_N)\right] \lesssim N^{2a-2\alpha+\xi}\left(N^{-2\nu_1} + N^{-1/4} + N^{-2q_0(1/2-\alpha)}\right) + N^{2a-\alpha-1/2+\xi}.
\label{eq:Var_bound_fin}
\end{equation}
Recalling that $\nu_1=1/10$, $q_0=1/20$ and $\alpha\in [0,1/2)$, we obtain that $N^{2a-2\alpha-2q_0(1/2-\alpha)+\xi}$ dominates the rest of the terms in the rhs. of \eqref{eq:Var_bound_fin} and complete the proof of Proposition \ref{prop:Var_main}.

\section{Proofs of the ingredients required for the analysis of the regime $\eta\lesssim N^{-1}$}\label{sec:proofprop}

In this section we prove various technical results needed to estimate the contribution of the regime~$\eta\lesssim N^{-1}$ in Girko's formula. Our arguments rely on the Dyson Brownian motion techniques: in the proofs of Propositions~\ref{prop:tail_bound} and~\ref{prop:EG} we use the relaxation of the DBM from \cite[Proposition~4.6]{bourgade2024fluctuations}, while in the proof of Proposition~\ref{prop:EGG} we employ the quantitative decorrelation of two DBMs from Theorem~\ref{eq:maintheoDBM1_main} stated below. DBM plays a central role in this section, since it regularizes the eigenvalue distribution of $H^z$ on small scales at the cost of adding a small Gaussian component to~$X$. 
This component is subsequently removed via a separate Green function comparison argument (GFT), whose technical implementation in the proof of Proposition~\ref{prop:tail_bound} differs from the one in the proofs of Propositions~\ref{prop:EG} and~\ref{prop:EGG}.

First, we introduce the set-up for the quantitative decorrelation of two DBMs. Consider the flow 
\begin{equation}
\label{eq:matDBM1_main}
\dif X_t=\frac{\dif B_t}{\sqrt{N}}, \qquad\quad X_0=X,
\end{equation}
where the entries $(B_t)_{ab}$ are independent standard complex Brownian motions. Let $H_t^z$ be the Hermitization of $X_t-z$ defined as in \eqref{eq:def_hermitization} and denote its eigenvalues by $\lambda_i^z(t)$, $|i|\le N$.  Here and also further in this section we implicitly exclude zero from the set of indices $|i|\le N$. It is known that the eigenvalues of $H_t^z$ satisfy the following SDE called the Dyson Brownian motion
\begin{equation}
\label{eq:DBM1_main}
\dif \lambda_i^z(t)=\frac{\dif b_i^z(t)}{\sqrt{2N}}+\frac{1}{2N}\sum_{j\ne i}\frac{1}{\lambda_i^z(t)-\lambda_j^z(t)}\dif t,
\end{equation}
though the exact form of this SDE is not important in this section. In \eqref{eq:DBM1_main}, $\{b_i^z(t)\}_{i=1}^N$ are independent standard real Brownian motions, and $b_{-i}^z(t)=-b_i^z(t)$. We consider the evolution of eigenvalues of $H_t^{z_1}$ and $H_t^{z_2}$ for some $z_1,z_2\in\C$ and prove the following result.

\begin{theorem}\label{eq:maintheoDBM1_main} Let $X$ be a complex i.i.d. matrix satisfying Assumption~\ref{ass:chi}(i), and let $X_t$ be the solution of \eqref{eq:matDBM1_main} with initial condition $X_0=X$. Fix any small $\omega_*, \delta>0$, $z_1,z_2\in\C$ with $|z_1|,|z_2|<1-\delta$, $T>N^{-1+\omega_*}$, and a possibly $N$-dependent $0<R<N|z_1-z_2|^2$. Then there exist two diffusion processes $\mu^{(l)}_{i}(t)$, $|i|\le N$, $t\in [0,T]$ for $l=1,2$, adapted to the filtration induced by the Brownian motion in \eqref{eq:matDBM1_main}, which satisfy the following properties.
\begin{enumerate}
\item[(i)] For any $t\in [0,T]$ and $l=1,2$, the $\{\mu^{(l)}_i(t)\}_{|i|\le N}$ are distributed as the eigenvalues of the Hermitization of $\widetilde{X}^{(l)}_t$, where $\widetilde{X}^{(l)}_t$ is defined by the flow \eqref{eq:matDBM1_main} with the initial condition given by a complex Ginibre matrix $\widetilde{X}^{(l)}_0:=\widetilde{X}^{(l)}$.
\item[(ii)] The processes $\{(\mu^{(1)}_i(t))_{t\in [0,T]}\}_{|i|\le N}$ and $\{(\mu^{(2)}_i(t))_{t\in [0,T]}\}_{|i|\le N}$ are independent.
\item[(iii)] For $l=1,2$ and any $\xi>0$ it holds that
\begin{equation}
\left\vert \rho_t^{z_l}(0)\lambda_i^{z_l}-\rho_{sc,t}(0)\mu_i^{(l)}\right\vert \le N^\xi \left(\mathcal{E}_1(t,R) + |i|\mathcal{E}_0(t) + \frac{|i|^2}{N^2}\right),
\label{eq:DBM_bound2}
\end{equation}
with very high probability simultaneously for all $t\in [N^{-1+\omega_*},T]$ and indices $|i|\le R$. Here $\mathcal{E}_0$ and $\mathcal{E}_1$ are defined in \eqref{eq:def_Phi} and \eqref{eq:def_Phi2}, respectively, and $\rho_t^z$, $\rho_{\mathrm{sc},t}$ are the evolutions of $\rho^z$ and $\rho_{\mathrm{sc}}$ along the semicircular flow\footnote{Here we denoted $\mu_{\mathrm{sc},t}(\dif x)=\rho_{\mathrm{sc},t}(x)\,\dif x:= t^{-1/2}\rho_{\mathrm{sc}}(t^{-1/2}x)\dif x$,  where $\rho_{\mathrm{sc}}(x):=(2\pi)^{-1}\sqrt{[4-x^2]_+}$. The density $\rho_t^{z_l}$ is defined as follows. Let $\mu^{z_l}(\dif x):=\rho_0^{z_l}(x)\dif x$, then
\[
\rho_t^{z_l}(x):=\lim_{\Im z\to 0^+}\frac{1}{\pi}\Im\int_\R \frac{\mu_t^{z_l}(\dif x)}{x-z},
\]
where $\mu_t^{z_l}:=\mu^{z_l}\boxplus \mu_{\mathrm{sc},t}$ is given by the free additive convolution between $\mu^{z_l}$ and the time-rescaled semicircular measure $\mu_{\mathrm{sc},t}(\dif x)$. We refer to \cite[Section~4.4]{bourgade2024fluctuations} for a detailed description of the semicircular flow evolution.}, respectively.
\end{enumerate}
\end{theorem}

Theorem~\ref{eq:maintheoDBM1_main} establishes decorrelation of the eigenvalues of $H^{z_1}_t$ and $H^{z_2}_t$ in the regime $|z_1-z_2|\gg N^{-1/2}$ by coupling them to two independent processes $\{\mu^{(1)}_i(t)\}_{|i|\le N}$ and $\{\mu^{(2)}_i(t)\}_{|i|\le N}$. This decorrelation has been obtained in various forms \cite{bourgade2024fluctuations, macroCLT_complex, cipolloni2024maximum} (see also \cite{cipolloni2023quenched}) relying on the analysis of weakly correlated Dyson Brownian motions (DBM). Our proof of Theorem \ref{eq:maintheoDBM1_main} relies on the recent proof of \cite[Theorem~4.1]{bourgade2024fluctuations} to give a quantitative estimate of the decorrelation exponent, and is presented in Supplementary Section~\ref{sec:DBM}. We remark that one of the main inputs required for the proof of Theorem~\ref{eq:maintheoDBM1_main} is the overlap bound from Proposition~\ref{prop:overlap}.

For ease of reference, we now state \cite[Proposition~4.6]{bourgade2024fluctuations}, whose set-up is similar but much simpler than the one of Theorem~\ref{eq:maintheoDBM1_main}. Namely, it concerns a single Hermitization parameter $z\in\C$ instead of two $z_1, z_2\in\C$. Fix (small) $\omega_*,\delta,\xi>0$. By \cite[Proposition~4.6]{bourgade2024fluctuations} there exists a diffusion process $\{\mu_i(t)\}_{|i|\le N}$ satisfying Theorem~\ref{eq:maintheoDBM1_main}(i), such that
\begin{equation}
\left\vert\rho_t^{z}(0)\lambda^z_i(t)-\rho_{sc,t}(0)\mu_i(t)\right\vert \le N^\xi \left(|i|\mathcal{E}_0(t)+\frac{|i|^2}{N^2}\right),\quad\text{with}\quad \mathcal{E}_0(t):=\frac{1}{N}\left(\frac{1}{\sqrt{Nt}}+t\right),
\label{eq:DBM_bound1}
\end{equation}
with very high probability, uniformly in $|z|\le 1-\delta$, $t\in [N^{-1+\omega_*},1]$ and $|i|\le N$. Here $\rho_t^z$ and $\rho_{sc,t}$ are defined as in Theorem~\ref{eq:maintheoDBM1_main}, and $\mathcal{E}_0$ is defined in \eqref{eq:def_Phi}.

We further assume in this section that $\rho_t^z(0)=\rho_{sc,t}(0)$ to simplify the application of both~\eqref{eq:DBM_bound1} and Theorem~\ref{eq:maintheoDBM1_main}. This assumption can be easily removed by a simple time-rescaling. 

The proofs of Propositions~\ref{prop:tail_bound}, \ref{prop:EG}, and~\ref{prop:EGG} share the same strategy consisting of two steps. First, we fix a small~$\omega_*>0$ and for any $N^{-1+\omega_*}\le t\le 1$ establish the desired bounds for the i.i.d. matrix distributed as 
\begin{equation}
\frac{1}{\sqrt{1+t}}X_t\stackrel{d}{=} \frac{1}{\sqrt{1+t}}X + \frac{\sqrt{t}}{\sqrt{1+t}}\widetilde{X}, 
\label{eq:distr_t}
\end{equation}
where $X_t$ follows \eqref{eq:matDBM1_main} and~$\widetilde{X}$ is a complex Ginibre matrix independent from~$X$. Here we normalized the lhs. of \eqref{eq:distr_t} by $(1+t)^{-1/2}$ to ensure that the second-order moment structure of $X$ coincides with the one of the lhs. of \eqref{eq:distr_t}. In the second step we embed~$X$ into the matrix-valued Ornstein-Uhlenbeck flow
\begin{equation}
\dif X^t = -\frac{1}{2}X^t +\frac{\dif B_t}{\sqrt{N}},\quad X^0=X, 
\label{eq:OU_GFT}
\end{equation}
with $B_t$ defined in the same way as in~\eqref{eq:matDBM1_main}. Since
\begin{equation}
X^s\stackrel{d}{=} \ee^{-s/2}X+\sqrt{1-\ee^{-s}}\widetilde{X}
\label{eq:distr_GFT}
\end{equation}
for any~$s>0$, it holds that
\begin{equation}
X^{s_0}\stackrel{d}{=}\frac{1}{\sqrt{1+t}}X_t\quad\text{for}\quad s_0=\log(1+t).
\label{eq:strat_concl}
\end{equation}
This enables us to perform a GFT using the flow~\eqref{eq:OU_GFT} and remove the Gaussian component from~$X_t$. These two steps (adding and removing a Gaussian component) are afterwards complemented by the optimization in~$t$ of the sum of two error terms obtained one at each step.

We denote the analogue of~$H^z$ defined in~\eqref{eq:def_hermitization} for $X_t, X^t$ and~$\widetilde{X}$ by~$H^z_t$, $H^{t,z}$ and~$\widetilde{H}^z$, respectively. The same notational convention is used for the analogues of~$W$ and~$G^z(w)$ defined in~\eqref{eq:def_hermitization} and~\eqref{eq:def_resolvent}, respectively. We also use~$\xi$ to denote a positive $N$-independent exponent which can be taken arbitrarily small and whose exact value may change from line to line. Now we separately prove each of the Propositions~\ref{prop:tail_bound}, \ref{prop:EG}, and~\ref{prop:EGG} using the strategy explained above.

\subsection{Proof of Proposition \ref{prop:tail_bound}}\label{sec:tail_bound} Fix a small $\omega_*>0$. First we prove that
\begin{equation}
\mathbf{P}\left[\lambda_1^z(t) <N^{-1}x\right] \lesssim \log N\cdot x^2 + N^\xi (N\mathcal{E}_0(t))^2
\label{eq:tail_t}
\end{equation}
for any fixed $\xi>0$ uniformly in $t\in [N^{-1+\omega_*},1]$ and $x\in (0,1]$. Let $\{\mu_i(t)\}_{|i|\le N}$ be given by \cite[Proposition~4.6]{bourgade2024fluctuations}. Applying \eqref{eq:DBM_bound1} for $i=1$, we get that
\begin{equation}
\mathbf{P}\left[\lambda_1^z(t) \le N^{-1}x\right] \le \mathbf{P}\left[\mu_1(t)\le N^{-1}x + N^\xi \mathcal{E}_0(t) \right] + \mathcal{O}(N^{-D})
\label{eq:lambda_to_mu_tail}
\end{equation}
for any (large) fixed $D>0$. Since $\mu_1(t)$ is distributed as the least positive eigenvalue of the Hermitization of $\widetilde{X}_t$, where $(1+t)^{-1/2}\widetilde{X}_{t}$ is a complex Ginibre matrix, \cite[Corollary 2.4]{least_sing_val_Ginibre} for $z=0$ implies that
\begin{equation}
\mathbf{P}\left[\mu_1(t) \le N^{-1}y\right]\lesssim (|\log y|+1)y^2 
\label{eq:tail_Gin}
\end{equation}
uniformly in $y>0$. Here we additionally used that $1+t\sim 1$. Combining~\eqref{eq:tail_Gin} for $y\sim x + N^{1+\xi}\mathcal{E}_0$ with~\eqref{eq:lambda_to_mu_tail}, we obtain~\eqref{eq:tail_t}.

Now we compare the probability in the lhs. of \eqref{eq:tail_t} at times $t$ and $0$ by closely following the GFT strategy introduced in \cite{EYY12} and later employed in \cite{small_dev, LFL_Wigner}. Recall the definition of $X^s$ from \eqref{eq:OU_GFT} and denote the least positive eigenvalue of $H^{s,z}$ by $\lambda_1^{s,z}$ for $s\ge 0$. By \eqref{eq:strat_concl} we have that
\begin{equation}
\lambda_1^{s_0,z}\stackrel{d}{=}(1+t)^{-1/2}\lambda_1^z(t)\quad\text{for}\quad s_0=\log (1+t).
\label{eq:F_GFT_s_to_t}
\end{equation}
For $E>0$, let $\bm{1}_E$ be the characteristic function of $[-E,E]$. Trivially, it holds that
\begin{equation}
\mathbf{P}[\lambda_1^{s,z}\le E] = \mathbf{P}\left[\mathrm{Tr}\bm{1}_E(H^{s,z})\ge 1\right],\quad s\ge 0. 
\label{eq:P_to_chi}
\end{equation}
In order to represent the rhs. of \eqref{eq:P_to_chi} in terms of $G^{s,z}$, we consider a non-decreasing function $F\in C^\infty(\R)$ such that
\begin{equation*}
F(y)=0\,\, \text{for}\,\, 0\le y\le 1/9;\qquad F(y)=1\,\, \text{for}\,\, y\ge 2/9.
\end{equation*}
Since $\mathrm{Tr}\bm{1}_E(H^{s,z})$ is integer-valued, \eqref{eq:P_to_chi} implies that
\begin{equation}
\mathbf{P}[\lambda_1^{s,z}\le E] = \E\left[F(\mathrm{Tr}\bm{1}_E(H^{s,z}))\right],\quad s\ge 0. 
\label{eq:P_to_F}
\end{equation}
Similarly to \cite[Lemma 2.3]{small_dev}, for any $\eta, l>0$ such that $\eta\le l\le E$, we have
\begin{equation}
\mathrm{Tr}\bm{1}_{E-l}*\theta_\eta(H^{s,z}) - N^\xi \frac{\eta}{l}\le \mathrm{Tr} \bm{1}_E(H^{s,z}) \le \mathrm{Tr}\bm{1}_{E+l}*\theta_\eta(H^{s,z}) + N^\xi \frac{\eta}{l}\quad\text{with}\quad \theta_\eta(y):=\frac{\eta}{\pi (y^2+\eta^2)}
\label{eq:lemF1}
\end{equation}
for any fixed $\xi>0$, with probability at least $1-N^{-D}$.

Given $x\in(0, 1]$, we fix a (small) $\epsilon>0$ and choose 
\begin{equation}
E:=N^{-1}x,\quad l:=N^{-\epsilon}E,\quad \eta:=N^{-2\epsilon}E.
\label{eq:F_GFT_param}
\end{equation}
In particular $\eta\ll l$, which in combination with \eqref{eq:lemF1} implies that
\begin{equation}
\E\left[F\left(\mathrm{Tr}\bm{1}_{E-l}*\theta_\eta(H^{s,z})\right)\right] - N^{-D}\le \E\left[F\left(\mathrm{Tr} \bm{1}_E(H^{s,z})\right)\right] \le \E\left[F\left(\mathrm{Tr}\bm{1}_{E+l}*\theta_\eta(H^{s,z})\right)\right] + N^{-D}
\label{eq:lemF2}
\end{equation}
for any fixed $D>0$ uniformly in $0\le s\lesssim 1$. Observe that
\begin{equation}
\mathrm{Tr}\bm{1}_{E+l}*\theta_\eta(H^{s,z}) = \frac{2N}{\pi}\int_{-(E+l)}^{E+l} \langle \Im G^{s,z}(y+\ii\eta)\rangle\dif y.
\label{eq:Tr_to_int}
\end{equation}
From now on we restrict to the case $\eta\ge N^{-3/2+\epsilon}$, which is sufficient for the proof of Proposition \ref{prop:tail_bound}. We claim that
\begin{equation}
\frac{\dif}{\dif s} \E \left[ F\left(\mathrm{Tr}\bm{1}_{E+l}*\theta_\eta(H^{s,z})\right)\right]\lesssim N^{1/2+\xi} (N\eta)^{-3}  \E \left[ F\left(\mathrm{Tr}\bm{1}_{E+3l}*\theta_\eta(H^{s,z})\right)\right] + N^{-D},
\label{eq:dF}
\end{equation}
uniformly in $s\in [0,1]$. This is a direct analogue of \cite[Eq.(2.62)]{small_dev}. The proof of \eqref{eq:dF} follows the same strategy as that of \cite[Eq.(2.62)]{small_dev} and thus is omitted. However, one aspect of the proof simplifies. While \cite{small_dev} exploits the concentration property of resolvent entries, the so-called isotropic single-resolvent local law, resolvents appearing in the proof of \cite[Eq.(2.62)]{small_dev} are decomposed into their deterministic counterparts and the fluctuations around them, and then these two terms are analyzed separately. In our situation $\eta$ is below the scale of the spectral resolution and the fluctuation given by the isotropic local law exceeds the size of the deterministic approximation. So, we simply have from \cite[Theorem~3.4]{singleG_circ} (see also \cite[Theorem~3.1]{nonHermdecay}) that
\begin{equation}
\label{eq:llawbelsc}
\left\vert(G^{s,z}(y+\ii\eta))_{ab}\right\vert\prec \frac{1}{N\eta}
\end{equation}
uniformly in $s\in [0,1]$, $|z|\le 1-\delta$, $|y|\le \tau'$, $\eta\in (0,N^{-1})$ and $a,b\in [2N]$, for some (small) constant $\tau'>0$. In particular, in the proof of \eqref{eq:dF} we never decompose $G^{s,z}_{ab}$ into $M^{s,z}_{ab}$ and $(G^{s,z}-M^{s,z})_{ab}$, where $M^{s,z}$ is the deterministic approximation to $G^{s,z}$, but just use the bound \eqref{eq:llawbelsc}.

To complete the GFT, we fix $x$ and $t$ such that
\begin{equation}
N^{-1/2+2\epsilon+\xi}\le x\le 1\quad \text{and}\quad N^{-1+\omega_*}\le t\le N^{-1/2-2\xi}(N\eta)^3=N^{-1/2-6\epsilon-2\xi}x^3.
\label{eq:xt_constraints}
\end{equation}
Here the lower bound on $x$ is a consequence of the constraint $\eta\ge N^{-3/2+\xi}$ with $\eta$ defined in \eqref{eq:F_GFT_param}. The lower bound on $t$ in \eqref{eq:xt_constraints} comes from the first step of this proof (see below~\eqref{eq:tail_t}), while the upper bound ensures that
\begin{equation}
N^{1/2+\xi}(N\eta)^{-3}t\le N^{-\xi},
\label{eq:F_GFT_propag}
\end{equation}
which is crucial for the application of a Gronwall-type argument to \eqref{eq:dF}. To make the set of $t$'s defined by \eqref{eq:xt_constraints} non-empty, we further restrict the range of $x$ in \eqref{eq:xt_constraints} to
\begin{equation}
N^{-1/6+\omega_*+2\epsilon +\xi}\le x\le 1.
\end{equation}
Finally, we take $s_0:=\log (1+t)$ and observe that $s_0\sim t$. Combining \eqref{eq:lemF2}, \eqref{eq:P_to_F}, \eqref{eq:F_GFT_s_to_t}, and \eqref{eq:tail_t} we get that
\begin{equation}
\begin{split}
&\E \left[ F\left(\mathrm{Tr}\bm{1}_{E+kl}*\theta_\eta(H^{s_0,z})\right)\right]\le \mathbf{P}\left[\lambda_1^{s_0,z}\le E+(k+1)l\right]+N^{-D}\\
&\quad = \mathbf{P}\left[\lambda_1^{z}(t)\le (1+t)^{1/2}\left(E+(k+1)l\right)\right]+N^{-D}\lesssim \log N\cdot x^2 + N^\xi (N\mathcal{E}_0(t))^2 
\end{split}
\label{eq:tail_s0}
\end{equation}
for any fixed $k\in\N$.

We derive from \eqref{eq:dF} that 
\begin{equation}
\begin{split}
&\E \left[ F\left(\mathrm{Tr}\bm{1}_{E+l}*\theta_\eta(H^{0,z})\right)\right]\lesssim \E \left[ F\left(\mathrm{Tr}\bm{1}_{E}*\theta_\eta(H^{s_0,z})\right)\right]\\
&\quad + \sum_{m=1}^{k-1} \left(N^{1/2+\xi}(N\eta)^{-3}t\right)^m \E \left[ F\left(\mathrm{Tr}\bm{1}_{E+(2m+1)l}*\theta_\eta(H^{s_0,z})\right)\right] + \left(N^{1/2+\xi}(N\eta)^{-3}t\right)^k
\end{split}
\label{eq:F_iteration}
\end{equation}
for any fixed $k\in \N$. The proof of~\eqref{eq:F_iteration} is identical to \cite[Eq.(3.19)--(3.20)]{LFL_Wigner} and is thus omitted. Taking sufficiently large $N$-independent $k\in\N$ in~\eqref{eq:F_iteration} and using \eqref{eq:F_GFT_propag}, \eqref{eq:tail_s0} as an input, we get
\begin{equation}
\E \left[ F\left(\mathrm{Tr}\bm{1}_{E+l}*\theta_\eta(H^{0,z})\right)\right] \lesssim \log N\cdot x^2 + N^\xi (N\mathcal{E}_0(t))^2.
\end{equation}
Together with \eqref{eq:lemF2}, \eqref{eq:P_to_F}, and \eqref{eq:F_GFT_param} this implies
\begin{equation}
\mathbf{P}\left[\lambda_1^z \le N^{-1}x\right] \lesssim \log N\cdot x^2 + N^\xi (N\mathcal{E}_0(t))^2,
\label{eq:tail_bound_parametrized}
\end{equation}
where $\lambda_1^z=\lambda_1^z(0)$.

It remains to optimize~\eqref{eq:tail_bound_parametrized} over the range of~$t$ given in~\eqref{eq:xt_constraints} using the explicit form of~$\mathcal{E}_0(t)$ from~\eqref{eq:DBM_bound1}. An elementary calculation shows that~$\mathcal{E}_0(t)$ achieves its minimum on $[N^{-1+\omega_*}, N^{-1/2-6\epsilon-2\xi}x^3]$ at the rightmost point of this interval. Choosing $\omega_*,\epsilon,\xi>0$ sufficiently small we thus get
\begin{equation}
\mathbf{P}\left[\lambda_1^z \le N^{-1}x\right] \lesssim \log N\cdot x^2 + N^{-1/2+C\xi}x^{-3},\quad x\in [N^{-1/6+\xi},1],
\label{eq:tail_bound_optimized}
\end{equation}
for any fixed $\xi>0$ and for some constant~$C>0$. Since for $x\ge N^{-\nu_1+C\xi}$ with $\nu_1:=1/10$ the first term in the rhs. of~\eqref{eq:tail_bound_optimized} dominates the second one,~\eqref{eq:tail_bound_optimized} completes the proof of Proposition~\ref{prop:tail_bound}.

\subsection{Proof of Proposition \ref{prop:EG}}\label{sec:proof_EG} Let $\widetilde{X}_t$ be the solution to \eqref{eq:matDBM1_main} with the initial condition given by a complex Ginibre matrix. Denote the resolvent of the Hermitization of $\widetilde{X}_t-z$ by $\widetilde{G}^z_t$. First we show that
\begin{equation}
\E \langle G_t^z(\ii\eta)\rangle = \E \langle \widetilde{G}^z_t(\ii\eta)\rangle + \mathcal{O}\left(N^{1+\xi}\mathcal{E}_0(t)\left(1+N\eta + \frac{N\mathcal{E}_0(t)}{N\eta} \right)\right)
\label{eq:EG_t}
\end{equation}
for any fixed $\xi>0$, uniformly in $\eta\in (0,1]$, $|z|\le 1-\delta$ and $t\in [N^{-1+\omega_*},1]$.

Throughout the proof of \eqref{eq:EG_t} the parameters $\eta,z$, and $t$ remain fixed, so we often omit them from notations. In particular, we simply denote $\lambda_i:=\lambda_i^z(t)$ and $\mu_i:=\mu_i(t)$, where $\{\mu_i(t)\}_{|i|\le N}$ is given by \cite[Proposition~4.6]{bourgade2024fluctuations}. Introduce the notation
\begin{equation}
\mathscr{L}_i=\mathscr{L}_i(t,\eta):= \frac{\eta}{( \lambda_i(t))^2 +\eta^2},\quad \mathscr{M}_i=\mathscr{M}_i(t,\eta):= \frac{\eta}{(\mu_i(t))^2 +\eta^2},\quad \text{for}\quad i\in [N].
\label{eq:def_T}
\end{equation}
We claim that 
\begin{equation}
\E\frac{1}{N}\sum_{i=1}^N \mathscr{L}_i(t,\eta) = \E\frac{1}{N}\sum_{i=1}^N \mathscr{M}_i(t,\eta) + \mathcal{O}\left(N^{1+\xi}\mathcal{E}_0(t)\left(1+N\eta + \frac{N\mathcal{E}_0(t)}{N\eta}\right)\right).
\label{eq:Lambda_to_M}
\end{equation}
Once \eqref{eq:Lambda_to_M} is obtained, the rest of the proof of \eqref{eq:EG_t} goes as follows. Denote the eigenvalues of $\widetilde{H}^z_t$ by $\{\widetilde{\lambda}_i^z(t)\}_{|i|\le N}$ and let $\{\widetilde{\mu}_i(t)\}_{|i|\le N}$ be the comparison process given by \cite[Proposition~4.6]{bourgade2024fluctuations}. Define $\widetilde{\mathscr{L}}_i$ and $\widetilde{\mathscr{M}}_i$ as the $\widetilde{X}$-counterparts of the quantities defined in \eqref{eq:def_T}. We apply \eqref{eq:Lambda_to_M} to $\widetilde{\mathscr{L}}_i, \widetilde{\mathscr{M}}_i$ and note that $\widetilde{\mathscr{M}}_i$ equals in distribution to $\mathscr{M}_i$ by Theorem~\ref{eq:maintheoDBM1_main}(i) (we mentioned above \eqref{eq:DBM_bound1} that this part of Theorem~\ref{eq:maintheoDBM1_main} is applicable in the set-up of \eqref{eq:DBM_bound1}). Together with \eqref{eq:Lambda_to_M} this gives
\begin{equation}
\E\frac{1}{N}\sum_{i=1}^N \mathscr{L}_i(t,\eta) = \E\frac{1}{N}\sum_{i=1}^N \widetilde{\mathscr{L}}_i(t,\eta)+\mathcal{O}\left(N^{1+\xi}\mathcal{E}_0(t)\left(1+N\eta + \frac{N\mathcal{E}_0(t)}{N\eta}\right)\right).
\label{eq:Lambda_to_Lambda_tilde}
\end{equation}
By \eqref{eq:G_to_ImG}, the lhs. of \eqref{eq:Lambda_to_Lambda_tilde} equals to $-\ii \E\langle G^z_t(\ii\eta)\rangle$, while the sum in the rhs. of \eqref{eq:Lambda_to_Lambda_tilde} equals to $-\ii \E\langle \widetilde{G}^z_t(\ii\eta)\rangle$. This finishes the derivation of \eqref{eq:EG_t} from \eqref{eq:Lambda_to_M}.

Now we prove~\eqref{eq:Lambda_to_M}. In the calculations below the index $i$ stands for a positive integer smaller than~$N$. Recall the definition of the densities $\rho_t^z$ and $\rho_{sc,t}$ from Theorem \ref{eq:maintheoDBM1_main}. Similarly to~\eqref{eq:def_quantiles}, we denote by $\{\gamma_i\}$ and $\{\widetilde{\gamma}_i\}$ the quantiles of~$\rho^z_t$ and $\rho_{sc,t}$, respectively. Although~$\widetilde{\gamma}_i$ depends on $t$ and $\gamma_i$ depends on both~$z$ and~$t$, we suppress this dependence in the notation for brevity. Fix a (small) exponent $\xi_0>0$. We distinguish between the regimes $i> N^{\xi_0}$ and $i\le N^{\xi_0}$ in the analysis of~$|\mathscr{L}_i-\mathscr{M}_i|$.

Consider at first the case $i>N^{\xi_0}$, where we have that $\lambda_i\sim \gamma_i$ with very high probability. In the regime $i\le (1-\tau)N$, for any small fixed $\tau>0$, this follows from the eigenvalue rigidity \eqref{eq:rigidity} and \eqref{eq:gamma_size}, while for $i>(1-\tau)N$ we simply use that $\lambda_i\sim 1$ and $\gamma_i\sim 1$. Similarly we have $\mu_i\sim \widetilde{\gamma}_i$, since $\{\mu_i\}_{|i|\le N}$ are distributed as the eigenvalues of $\widetilde{H}^0_t$. Using additionally that $\gamma_i\sim\widetilde{\gamma}_i$ by \eqref{eq:gamma_size}, we obtain
\begin{equation}
\left\vert \mathscr{L}_i-\mathscr{M}_i\right\vert =  \frac{\eta\left\vert \lambda_i-\mu_i\right\vert(\lambda_i+\mu_i)}{\big(\lambda_i^2 +\eta^2\big)\big(\mu_i^2 +\eta^2\big)} \sim \frac{\eta\left\vert \lambda_i-\mu_i\right\vert\gamma_i }{\big(\gamma_i^2 +\eta^2\big)^2}\lesssim \frac{\eta \left\vert \lambda_i-\mu_i\right\vert}{\gamma_i^3} 
\label{eq:T-T_large_i}
\end{equation}
with very high probability. Combining \eqref{eq:T-T_large_i} with \eqref{eq:DBM_bound1} and \eqref{eq:gamma_size}, we get
\begin{equation}
\frac{1}{N}\sum_{i> N^{\xi_0}} \left\vert \mathscr{L}_i-\mathscr{M}_i\right\vert \lesssim N^{2+\xi}\eta\sum_{i> N^{\xi_0}}\frac{1}{|i|^3}\left(|i|\mathcal{E}_0(t) + \frac{|i|^2}{N^2}\right)\lesssim (N\eta)(N\mathcal{E}_0(t))
\label{eq:T-T_large_i_fin}
\end{equation}
with very high probability, where in the last step we used that $\mathcal{E}_0(t)\ge N^{-2+\xi'}$ for some $\xi'>0$.

For $i\le N^{\xi_0}$ we perform a more subtle analysis of the lhs. of \eqref{eq:T-T_large_i}. On the event 
\begin{equation*}
\mathfrak{E}_i := \bm{1}\left\lbrace \lambda_i \le N^{2\xi_0}\mathcal{E}_0(t)\right\rbrace
\end{equation*}
we estimate $\lambda_i, \mu_i$ trivially from below by zero and bound the probability of $\mathfrak{E}_i$ from above first by the probability of $\mathfrak{E}_1$ and then by \eqref{eq:tail_t} applied to $x:=N^{1+2\xi_0}\mathcal{E}_0(t)$:
\begin{equation}
\E\big[|\mathscr{L}_i-\mathscr{M}_i|\bm{1}\{\mathfrak{E}_i\}\big]\le \E\big[\left(|\mathscr{L}_i|+|\mathscr{M}_i|\right)\bm{1}\{\mathfrak{E}_1\}\big]\lesssim \eta^{-1}\mathbf{P}\left[\mathfrak{E}_1\right]\lesssim \eta^{-1}N^{4\xi_0+\xi}(N\mathcal{E}_0(t))^2.
\label{eq:small_lambda_event}
\end{equation}
On the complementary event, \eqref{eq:DBM_bound1} implies that $\lambda_i\sim\mu_i$ and we perform the same calculation as in the first step in \eqref{eq:T-T_large_i}:
\begin{equation}
\E \big[\left\vert \mathscr{L}_i-\mathscr{M}_i\right\vert\left(1-\bm{1}\{\mathfrak{E}_i\}\right)\big]\!\lesssim\! \E \left[\!\frac{\eta |\lambda_i-\mu_i|\lambda_i}{\left(\lambda_i^2 +\eta^2\right)^2}(1-\bm{1}\{\mathfrak{E}_i\})\!\right]\!\lesssim\! N^{\xi_0+\xi}\eta\mathcal{E}_0(t)\E \left[\!\frac{\lambda_i}{\lambda_i^4+\eta^4}(1-\bm{1}\{\mathfrak{E}_i\})\!\right]\!,
\label{eq:T-T_small_i}
\end{equation}
where in the second bound we also used that $|\lambda_i-\mu_i|\lesssim N^{\xi_0+\xi}\mathcal{E}_0(t)$ for $i\le N^{\xi_0}$ by~\eqref{eq:DBM_bound1}. To estimate the expectation in the rhs. of \eqref{eq:T-T_small_i}, consider the following dyadic decomposition
\begin{equation}
1-\bm{1}\{\mathfrak{E}_i\} \le \sum_k \bm{1}\left\lbrace 2^{k-1}\eta< \lambda_i\le 2^k\eta\right\rbrace,
\label{eq:dyad_decomp}
\end{equation}
where the summation runs from $k\ge \log_2(\eta^{-1}\mathcal{E}_0(t))$ to $k\lesssim \log N$. We use~\eqref{eq:tail_t} together with the trivial bound $\lambda_1\le \lambda_i$ to estimate the probabilities of events in the rhs. of~\eqref{eq:dyad_decomp}. This yields
\begin{equation}
\E \left[\frac{\lambda_i}{\lambda_i^4+\eta^4}(1-\bm{1}\{\mathfrak{E}_i\})\right]\lesssim\frac{N^{2+\xi}}{\eta}.
\label{eq:lambda_int_bound}
\end{equation}
Finally, combining \eqref{eq:small_lambda_event}--\eqref{eq:lambda_int_bound} we obtain
\begin{equation}
\frac{1}{N} \sum_{i\le N^{\xi_0}} \E |\mathscr{L}_i-\mathscr{M}_i|\lesssim \left(1 + \frac{N\mathcal{E}_0(t)}{N\eta}\right)N^{1+4\xi_0+\xi}\mathcal{E}_0(t).
\label{eq:T-T_small_i_final}
\end{equation}
We take $\xi_0>0$ sufficiently small and collect error terms from~\eqref{eq:T-T_large_i_fin} and~\eqref{eq:T-T_small_i_final}. This finishes the proof of~\eqref{eq:Lambda_to_M}.

Now we proceed to the second step of the proof of Proposition \ref{prop:EG} and remove the Gaussian component from $X_t$ by the Green function comparison argument. Specifically, we show that
\begin{equation}
\left\vert\E \left\langle G^{s_0,z}(\ii\eta)\right\rangle-\E\left\langle G^{0,z}(\ii\eta)\right\rangle\right\vert \lesssim N^{1/2+\xi}s_0 \left(1+\frac{1}{N\eta}\right)^4
\label{eq:EG_GFT}
\end{equation}
uniformly in $s_0\in [0,1]$, $\eta\in [N^{-3/2+\epsilon},1]$ and $|z|\le 1-\delta$, for any fixed $\epsilon,\delta>0$. Denote $G^s:=G^{s,z}(\ii\eta)$ for $s\ge 0$. The proof of \eqref{eq:EG_GFT} is based on the explicit formula
\begin{equation}
\frac{\dif}{\dif s} \E \left\langle G^{s}\right\rangle = \frac{1}{2} \E\left\langle W_s\left(G^{s}\right)^2\right\rangle +\E\left\langle \mathcal{S}[G^{s}]G^{s}\right\rangle,\quad \forall \, s\ge 0,
\label{eq:dG}
\end{equation}
where $\mathcal{S}$ is defined in \eqref{eq:def_S}. We estimate the rhs. of \eqref{eq:dG} in absolute value from above. To avoid carrying the dependence on $s$ in notations, we demonstrate this bound only for $s=0$, while for general $s\ge 0$ the estimates are exactly the same. 

For a set of index pairs
\begin{equation}
\bm{\alpha}=\left(\alpha_1,\ldots,\alpha_k)\in([1,N]\times [N+1,2N]\cup [N+1,2N]\times [1,N]\right)^k
\label{eq:def_index_set}
\end{equation}
denote the normalized cumulant of the corresponding elements of $W$ by
\begin{equation}
\kappa(\bm{\alpha})=\kappa(\alpha_1,\ldots,\alpha_k):=\kappa(\sqrt{N}w_{\alpha_1},\ldots,\sqrt{N}w_{\alpha_k}).
\label{eq:def_cum}
\end{equation}
Since $X$ has i.i.d. entries, $\kappa(\bm{\alpha})$ does not vanish only when $\alpha_1,\ldots,\alpha_k\in\{ab,ba\}$ for some $a\in[1,N]$, $b\in [N+1,2N]$. We further denote $\partial_{\bm{\alpha}}:=\partial_{\alpha_1}\cdots\partial_{\alpha_k}$, where $\partial_{\alpha}$ denotes the directional derivative in the direction $w_{\alpha}$, for an index pair~$\alpha$. Performing the cumulant expansion \cite[Lemma~3.2]{Knowles20} (see also \cite[Lemma~3.1]{HeKnowles} and \cite[Section~II]{Khorunzhy96}) in the first term in the rhs. of \eqref{eq:dG} we obtain
\begin{equation}
\E \langle WG^2\rangle =\frac{1}{N}\sum_{a,b} \E\left[w_{ab}(G^2)_{ba}\right]= \frac{1}{N}\sum_{a,b}\sum_{\ell=2}^L\sum_{\bm{\alpha}\in\{ab,ba\}^{\ell-1}}\frac{\kappa(ab,\bm{\alpha})}{N^{\ell/2}(\ell-1)!}\E\partial_{\bm{\alpha}}\left[ (G^2)_{ba}\right] + \mathcal{O}(N^{-D}),
\label{eq:cum_exp_first}
\end{equation}
for any fixed $D>0$. In \eqref{eq:cum_exp_first} the upper index of summation $L\in \N$ depends only on $D$, the model parameters from Assumption~\ref{ass:chi}, and $\epsilon>0$ from the constraint $\eta\ge N^{-3/2+\epsilon}$. The $(a,b)$-summation in~\eqref{eq:cum_exp_first} runs over all $a,b$ such that either $a\in [N]$, $b\in [N+1,2N]$ or $a\in [N+1,2N]$, $b\in [N]$. 

First, consider the second order ($\ell=2$) terms in \eqref{eq:cum_exp_first}. Performing in these terms differentiation in $\partial_{\bm{\alpha}}(G^2)_{ba}$ and summing up the resulting quantities over $a,b$, we conclude by a straightforward calculation that this sum exactly cancels out with the second term in the rhs. of \eqref{eq:dG}. Thus, we only need to estimate the contribution from the third and higher order terms ($\ell\ge 3$) in the rhs. of \eqref{eq:cum_exp_first}, which we do by the means of the entry-wise bound 
\begin{equation}
\left\vert(G^{s})_{ab}\right\vert\prec 1+\frac{1}{N\eta},
\label{eq:1G_ll_iso_below_scale}
\end{equation}
that holds uniformly in $a,b\in [2N]$, $|z|\le 1-\delta$ and $\eta\in (N^{-L},1]$ for any fixed $L>0$. This bound immediately follows from the single-resolvent isotropic local law \eqref{eq:1G_ll_iso}. Using \eqref{eq:1G_ll_iso_below_scale}, we obtain that the contribution from the third order cumulant terms to the lhs. of \eqref{eq:EG_GFT} is bounded by the rhs. of \eqref{eq:EG_GFT}, while the contribution from the rest of the terms is even smaller. We omit the rest of the details, since the proof of \eqref{eq:EG_GFT} is fairly standard.

To complete the proof of Proposition~\ref{prop:EG}, we fix $t\in [N^{-1+\omega_*}, N^{-\omega_*}]$ and apply \eqref{eq:EG_t} with $(1+t)^{1/2}z$ and $(1+t)^{1/2}\eta$ instead of $z$ and $\eta$, respectively. Here $z$ and $\eta$ satisfy the conditions of Proposition~\ref{prop:EG}. Due to the upper bound $t\le N^{-\omega_*}$, we have that $(1+t)^{1/2}|z|\le 1-\delta/2$ for sufficiently large~$N$. Thus, after rescaling,~\eqref{eq:EG_t} compares the resolvents of $(1+t)^{-1/2} X_t$  and $(1+t)^{-1/2} \widetilde{X}_t$. Using additionally that $(1+t)^{-1/2}\widetilde{X}_t$ and~$\widetilde{X}$ have the same distribution, we get
\begin{equation}
\E \left\langle\left((1+t)^{-1/2}X_t - Z -\ii\eta\right)^{-1}\right\rangle = \E \left\langle\left(\widetilde{X} - Z -\ii\eta\right)^{-1}\right\rangle + \mathcal{O}\left(N^{1+\xi}\mathcal{E}_0(t)\left(1+N\eta + \frac{N\mathcal{E}_0(t)}{N\eta} \right)\right)
\label{eq:EG_rescaled}
\end{equation}
Combining \eqref{eq:EG_rescaled} with \eqref{eq:EG_GFT} applied to $s_0:=\log(1+t)$ and recalling \eqref{eq:strat_concl}, we obtain
\begin{equation}
\E\langle G^z(\ii\eta)\rangle = \E \langle \widetilde{G}^z(\ii\eta)\rangle + N^\xi\mathcal{O}\left(N\mathcal{E}_0(t)\left(1+N\eta + \frac{N\mathcal{E}_0(t)}{N\eta}\right) +\sqrt{N}t\left(1+\frac{1}{N\eta}\right)^4  \right)
\label{eq:EG_comp_parametrized}
\end{equation}
uniformly in $t\in [N^{-1+\omega_*}, N^{-\omega_*}]$, $|z|\le 1-\delta$ and $\eta\in [N^{-3/2+\epsilon},1]$. Optimizing the error term in the rhs. of \eqref{eq:EG_comp_parametrized} over $t$, we finish the proof of Proposition \ref{prop:EG}.

\subsection{Proof of Proposition \ref{prop:EGG}}\label{sec:proof_EGG} In this section we work in the set-up of Theorem~\ref{eq:maintheoDBM1_main}. Fix a small $\omega_*>0$, take $T:=1$ and a possibly $N$-dependent $0<R<N|z_1-z_2|^2$. First, we show that
\begin{equation}
\left\vert \Cov \left(\langle G^{z_1}_t(\ii\eta_1)\rangle, \langle G^{z_2}_t(\ii\eta_2)\rangle\right)\right\vert \lesssim N^{1+\xi}\left(\mathcal{E}_1(t,R) + \mathcal{E}_0(t)\right)\left(N\eta_1+\frac{1}{N\eta_1}\right)\left(N\eta_2+\frac{1}{N\eta_2}\right)
\label{eq:Cov_t}
\end{equation}
uniformly in $\eta_1,\eta_2\in (0,1]$, $|z_1|,|z_2|\le 1-\delta$, $t\in [N^{-1+\omega_*},1]$ and $0<R<N|z_1-z_2|^2$. To prove \eqref{eq:Cov_t}, we introduce the following notation similar to \eqref{eq:def_T}:
\begin{equation}
\mathscr{L}_i^{z_l}=\mathscr{L}_i^{z_l}(t,\eta):= \frac{\eta}{( \lambda_i^{z_l}(t))^2 +\eta^2},\quad \mathscr{M}_i^{(l)}=\mathscr{M}_i^{(l)}(t,\eta):= \frac{\eta}{(\mu_i^{(l)}(t))^2 +\eta^2},\quad \text{for}\quad i\in [N],\, l=1,2,
\label{eq:def_T_l}
\end{equation}
Since $\mu_i^{(1)}(t)$ is independent of $\mu_j^{(2)}(t)$ for any $|i|,|j|\le N$ by Theorem~\ref{eq:maintheoDBM1_main}(ii), we have
\begin{equation}
\begin{split}
\left\vert \Cov\left(\langle G^{z_1}(\ii\eta_1)\rangle,\langle G^{z_2}(\ii\eta_2)\rangle\right)\right\vert \le& \bigg\vert\Cov \bigg(\frac{1}{N}\sum_{i=1}^N \left(\mathscr{L}^{z_1}_i-\mathscr{M}^{(1)}_i\right), \langle G^{z_2}(\ii\eta_2)\rangle\bigg)\bigg\vert\\
 + &\bigg\vert\Cov \bigg(\frac{1}{N}\sum_{i=1}^N \mathscr{M}^{(1)}_i, \frac{1}{N}\sum_{i=1}^N \left(\mathscr{L}^{z_2}_i-\mathscr{M}^{(2)}_i\right)\bigg)\bigg\vert,
\end{split}
\label{eq:Cov_indep_bound}
\end{equation}
where we additionally used \eqref{eq:G_to_ImG} to express $\langle G^{z_l}\rangle$ in terms of $\lambda_i^{(l)}$'s, for $l=1,2$. We further focus on the second term in the rhs. of \eqref{eq:Cov_indep_bound}, as the analysis of the first term is identical. From Theorem~\ref{eq:maintheoDBM1_main}(i) and the averaged single-resolvent local law \eqref{eq:av_ll_simple} we have
\begin{equation}
\frac{1}{N}\sum_{i=1}^N \mathscr{M}^{(1)}_i\stackrel{d}{=}\langle \Im \widetilde{G}^0_t(\ii\eta_1)\rangle \lesssim N^{\xi}\left(1+\frac{1}{N\eta_1}\right),
\end{equation}
for any $\xi>0$ with very high probability. Therefore, the second term in the rhs. of \eqref{eq:Cov_indep_bound} has an upper bound of order
\begin{equation}
\E \left[\frac{1}{N}\sum_{i=1}^N \left\vert\mathscr{L}^{z_2}_i-\mathscr{M}^{(2)}_i\right\vert\right] \left(1+\frac{1}{N\eta_1}\right)N^\xi.
\label{eq:Cov_bound_2term}
\end{equation}
Fix a (small) $\xi_0>0$. Arguing similarly to \eqref{eq:T-T_large_i}--\eqref{eq:T-T_large_i_fin} and \eqref{eq:T-T_small_i}--\eqref{eq:T-T_small_i_final}, we get
\begin{equation}
\begin{split}
\frac{1}{N}\sum_{i> N^{\xi_0}} \left\vert \mathscr{L}_i^{z_2}-\mathscr{M}_i^{(2)}\right\vert &\lesssim  (N\eta_2)\left(N(\mathcal{E}_1(t,R) + \mathcal{E}_0(t)) + \frac{1}{R}\right),\\
\frac{1}{N} \sum_{i\le N^{\xi_0}} \E \left\vert\mathscr{L}_i^{z_2}-\mathscr{M}_i^{(2)}\right\vert&\lesssim \left(1 + \frac{N\left(\mathcal{E}_1(t,R)+\mathcal{E}_0(t)\right)}{N\eta_2}\right)N^{1+4\xi_0+\xi}\left(\mathcal{E}_1(t,R)+\mathcal{E}_0(t)\right),
\end{split}
\label{eq:Cov_T-T}
\end{equation}
where the first line of \eqref{eq:Cov_T-T} holds with very high probability. The first and the second line of \eqref{eq:Cov_T-T} coincide with \eqref{eq:T-T_large_i_fin} and \eqref{eq:T-T_small_i_final} modulo the replacement of $\mathcal{E}_0$ with $\mathcal{E}_1+\mathcal{E}_0$, and the first line of \eqref{eq:Cov_T-T} features a new term $R^{-1}$ compared to \eqref{eq:T-T_large_i_fin}, which comes from $i>R$. For these indices \eqref{eq:DBM_bound2} does not hold, so we simply bound
\begin{equation}
|\lambda_i^{z_2}(t)-\mu_i^{(2)}(t)| \le \lambda_i^{z_2}(t) + \mu_i^{(2)}(t)\lesssim \frac{|i|}{N},
\end{equation}
where in the second bound we argued similarly to the discussion above~\eqref{eq:T-T_large_i}. However, this additional~$R^{-1}$ term does not play any role since it is smaller than $N\mathcal{E}_1(t,R)$. Combining~\eqref{eq:Cov_indep_bound},~\eqref{eq:Cov_bound_2term}, and~\eqref{eq:Cov_T-T}, we finish the proof of~\eqref{eq:Cov_t}.

Next, similarly to \eqref{eq:EG_GFT} we show that
\begin{equation}
\begin{split}
&\left\vert\Cov\left(\langle G^{s_0,z_1}(\ii\eta_1)\rangle,\langle G^{s_0,z_2}(\ii\eta_2)\rangle\right) - \Cov\left(\langle G^{0,z_1}(\ii\eta_1)\rangle,\langle G^{0,z_2}(\ii\eta_2)\rangle\right)  \right\vert\\
&\quad \lesssim N^{1/2+\xi}s_0 \left(1+\frac{1}{N(\eta_1\wedge\eta_2)}\right)^3 \frac{1}{N^2\eta_1\eta_2},
\end{split}
\label{eq:EGG_GFT}
\end{equation}
for any fixed $\xi>0$, uniformly in $s_0\in [0,1]$, $\eta_1,\eta_2\in [N^{-3/2+\epsilon},1]$ and $z_1,z_2\in (1-\delta)\mathbf{D}$. The proof of \eqref{eq:EGG_GFT} is standard and thus is omitted. Finally, we combine \eqref{eq:Cov_t} and \eqref{eq:EGG_GFT} in the way similar to \eqref{eq:EG_rescaled}--\eqref{eq:EG_comp_parametrized} and complete the proof of Proposition \ref{prop:EGG}.

\section{Chaos expansion for covariance: Proof of Proposition \ref{prop:Cov}}\label{sec:Cov}

Unlike in the proof of Proposition~\ref{prop:2G_av}, no dynamics is used in the proof of Proposition~\ref{prop:Cov}. By this we mean that throughout the argument the random matrix~$X$ and the parameters $z_1,z_2,w_1,w_2$ satisfying assumptions of Proposition~\ref{prop:Cov} remain unchanged. We denote $G_l:=G^{z_l}(w_l)$ for $l=1,2$, and set
\begin{equation}
\eta_*:=\eta_1\wedge\eta_2\wedge 1.
\label{eq:def_eta*}
\end{equation}
Since in Proposition~\ref{prop:Cov} we assume that $\eta_1,\eta_2\le 1$, here~$\eta_*$ simply equals to $\eta_1\wedge\eta_2$. Nevertheless, we present the definition of $\eta_*$ in the form \eqref{eq:def_eta*} to be consistent with~\eqref{eq:beta_hat_trivial}. In the proof of Proposition \ref{prop:Cov} we frequently use that $\widehat{\beta}_{12}\sim \gamma$ by Proposition~\ref{prop:stab_bound}, where $\widehat{\beta}_{12}$ and $\gamma$ are defined in~\eqref{eq:def_beta_hat} and~\eqref{eq:def_UV}, respectively. For brevity, the reference to Proposition~\ref{prop:stab_bound} is often omitted. Moreover, whenever $\widehat{\beta}_{12}$ and $\gamma$ appear in the proof with suppressed arguments, they are meant to be evaluated at $z_1,z_2,w_1,w_2$.

First, we formulate the following version of Proposition \ref{prop:Cov} for matrices with a Gaussian component.

\begin{proposition}\label{prop:Cov_b} Assume that $X$ satisfies conditions of Proposition \ref{prop:2G_av} for some fixed $b\in [0,1]$. Fix (small) $\delta, \epsilon,\kappa,\xi>0$ and recall the notations introduced in Proposition \ref{prop:Cov}. Then, we have
\begin{equation}
\Cov \left(\langle G^{z_1}(w_1)\rangle,\langle G^{z_2}(w_2)\rangle\right) = \frac{1}{N^2}\cdot \frac{V_{12}+\kappa_4U_1U_2}{2} + \mathcal{O}\left(\left(\frac{1}{N(\gamma\wedge N^{-b} + \eta_*)} + N^{-1/2} \right)\frac{N^\xi}{N^2\eta_1\eta_2}\right),
\label{eq:main_cov_b}
\end{equation}
uniformly in $|z_l|\le 1-\delta$, $E_l\in\mathbf{B}_\kappa^{z_l}$ and $\eta_l\in [N^{-1+\epsilon}, 1]$, for $l=1,2$. In the Gaussian case the error $N^{-1/2}$ in \eqref{eq:main_cov_b} can be removed.
\end{proposition}

Now we prove Proposition~\ref{prop:Cov} relying on Proposition~\ref{prop:Cov_b}.

\begin{proof}[Proof of Proposition~\ref{prop:Cov}] Let~$X$ be a general i.i.d. matrix satisfying assumptions of Proposition~\ref{prop:Cov_b}. We embed~$X$ into the matrix-valued Ornstein-Uhlenbeck process~\eqref{eq:OU_GFT} with $X^0:=X$, and recall the notation $G^{t,z}(w)$ for the resolvent of Hermitization of $X^t-z$, introduced below~\eqref{eq:strat_concl}. Fix $b\in [0,1]$ and take $t_b:=N^{-b}$. As it follows from~\eqref{eq:distr_GFT}, $X^{t_b}$ contains a Gaussian component of order~$N^{-b}$, so it satisfies assumptions of Proposition~\ref{prop:Cov_b}. We apply~\eqref{eq:main_cov_b} to $X^{t_b}$ and note that the first term in the rhs. of~\eqref{eq:main_cov_b} does not depend on~$t_b$, so it coincides with the first term in the rhs. of~\eqref{eq:main_cov}. Therefore,~\eqref{eq:main_cov} holds for~$X$ with the error term of order
\begin{equation}
\begin{split}
&\left(\frac{1}{N(\gamma\wedge N^{-b} + \eta_*)} + N^{-1/2} \right)\frac{N^\xi}{N^2\eta_1\eta_2}\\
&+ \left\vert \Cov \left(\langle G^{t,z_1}(w_1)\rangle,\langle G^{t,z_2}(w_2)\rangle\right) - \Cov \left(\langle G^{0,z_1}(w_1)\rangle,\langle G^{0,z_2}(w_2)\rangle\right)\right\vert. 
\end{split}
\label{eq:Cov_b_removal}
\end{equation} 
Arguing similarly to~\eqref{eq:EGG_GFT} (see also the proof of~\eqref{eq:EG_GFT}), we get that the term in the second line of~\eqref{eq:Cov_b_removal} has an upper bound of order
\begin{equation}
N^{1/2-b}\frac{N^\xi}{N^2\eta_1\eta_2}.
\label{eq:EGG_GFT_b}
\end{equation}
Finally, we optimize over $b\in [0,1]$ the sum of the terms in the first line of~\eqref{eq:Cov_b_removal} and in ~\eqref{eq:EGG_GFT_b}, which leads to the choice $b:=3/4$. This finishes the proof of Proposition~\ref{prop:Cov}.  
\end{proof}

The remainder of this section is devoted to the proof of Proposition~\ref{prop:Cov_b}. To simplify presentation, we focus on the case $b=0$, which simplifies the error term in the rhs. of~\eqref{eq:main_cov_b} to 
\begin{equation*}
\mathcal{O}\left(\left(\frac{1}{N\gamma} + N^{-1/2} \right)\frac{N^\xi}{N^2\eta_1\eta_2}\right).
\end{equation*}
However, this does not lead to any loss of generality, since the only adjustment of the proof below needed for the general case $b\in [0,1]$ is that the local law from Proposition \ref{prop:2G_av} is applied not for $b=1$, but for the given~$b$.

The proof of Proposition~\ref{prop:Cov_b} is structured as follows. In Section~\ref{sec:Cov_reduction} we prove the initial expansion~\eqref{eq:init_expansion} and derive Proposition~\ref{prop:Cov_b} from a suboptimal bound on the first term in the rhs. of~\eqref{eq:init_expansion} formulated in Proposition~\ref{prop:psi_22}. Next, in Section~\ref{sec:H_Cov} we state a hierarchy of so called \emph{master inequalities} for covariances arising along the chaos expansion, which consists in iterative expansions of the first term in the rhs. of~\eqref{eq:init_expansion}. We solve this hierarchy and deduce Proposition~\ref{prop:psi_22} from it. Finally, in Section~\ref{sec:H_Cov} we prove master inequalities.

\subsection{Initial expansion}\label{sec:Cov_reduction}  We derive \eqref{eq:main_cov_b} from the following bound on the first term in the rhs. of~\eqref{eq:init_expansion} and prove~\eqref{eq:init_expansion} along the way.

\begin{proposition}\label{prop:psi_22} Assume the set-up and conditions of Proposition~\ref{prop:Cov_b} with $b=0$. For any fixed $\xi>0$ it holds that
\begin{equation}
\left\vert\Cov(\langle G_1-M_1\rangle\langle (G_1-M_1)M_1\rangle,\langle G_2\rangle)\right\vert \lesssim \left(\frac{1}{N\gamma}+\frac{1}{\sqrt{N}}\right)\frac{N^\xi}{N^2\eta_*\eta_2},
\label{eq:psi_22}
\end{equation}
where $M_1=M^{z_1}(w_1)$, $G_l=G^{z_l}(w_l)$ for $l=1,2$, and~$\eta_*$ is defined in~\eqref{eq:def_eta*}.
\end{proposition}

The proof of Proposition~\ref{prop:psi_22} is postponed to Section~\ref{sec:H_Cov}. As we have mentioned below~\eqref{eq:init_expansion}, the averaged single-resolvent local law~\eqref{eq:1G_ll_av} implies that the lhs. of~\eqref{eq:psi_22} has an upper bound of order $(N^3\eta_1^2\eta_2)^{-1}$, which is sharper than~\eqref{eq:psi_22} already for $\eta_1,\eta_2>N^{-1/2}$. However, our main regime of interest is when $\eta_1,\eta_2\sim N^{-1+\epsilon}$ for a small $\epsilon>0$, where \eqref{eq:psi_22} improves upon the local law bound by a factor $\eta_*/\gamma +\sqrt{N}\eta_*$, which is much smaller than 1 already for $|z_1-z_2|\gg N^{-1/2}$. The estimate~\eqref{eq:psi_22} is sufficient for the proof of Proposition~\ref{prop:Cov_b}, as we now demonstrate, so we do not pursue optimality in Proposition~\ref{prop:psi_22}.  

\begin{proof}[Proof of Proposition \ref{prop:Cov_b}] We start with the proof of \eqref{eq:init_expansion}. It holds that 
\begin{equation}
\langle G_1-M_1\rangle = \langle G_1-M_1\rangle \langle (G_1-M_1)A\rangle - \langle \underline{WG_1}A\rangle, \quad A=\left(\left(\mathcal{B}^*_{11}\right)^{-1}[E_+]\right)^*M_1,
\label{eq:<G-M>}
\end{equation}
where $\mathcal{B}^*_{11}$ is defined with respect to the scalar product $\langle R,S\rangle:=\langle R^*S\rangle$ for $R,S\in\C^{(2N)\times (2N)}$. This identity can be easily derived from \eqref{eq:def_under_1} and \eqref{eq:MDE}, for a detailed derivation see \cite[Eq.(6.9)]{macroCLT_complex}. An explicit calculation based on~\eqref{eq:def_B12} shows that
\begin{equation}
\left(\mathcal{B}^*_{11}\right)^{-1}[R] = R + \frac{\langle M^*_1RM^*_1\rangle}{1-\langle (M^*_1)^2\rangle} - \frac{\langle M^*_1RM^*_1E_-\rangle}{1+\langle (M^*_1E_-)^2\rangle}E_-,\quad \forall R\in\C^{(2N)\times (2N)}.
\label{eq:B*_inverse}
\end{equation}
Therefore, $A=(1-\langle M_1^2\rangle)^{-1}M_1$ and, in particular, $\|A\|\lesssim 1$. Using \eqref{eq:<G-M>} we express the lhs. of \eqref{eq:main_cov} as follows:
\begin{equation}
\Cov(\langle G_1\rangle, \langle G_2\rangle) = \Cov\left( \langle G_1-M_1\rangle \langle (G_1-M_1)A\rangle, \langle G_2\rangle\right) - \E \left[ \langle \underline{WG_1}A\rangle \left(\langle G_2\rangle - \E\langle G_2\rangle\right)\right].
\label{eq:Cov_05_exp}
\end{equation}
In the second term in the rhs. of \eqref{eq:Cov_05_exp} we extend the underline on the second factor:
\begin{equation}
\begin{split}
\E \left[ \langle \underline{WG_1}A\rangle \left(\langle G_2\rangle - \E\langle G_2\rangle\right)\right] &= \E \left[\underline{ \langle WG_1A\rangle \left(\langle G_2\rangle - \E\langle G_2\rangle\right)}\right] - \E \widetilde{\E} \left[ \langle \widetilde{W}G_1A\rangle \langle \widetilde{W}G_2^2\rangle\right]\\
&=\E \left[\underline{ \langle WG_1A\rangle \left(\langle G_2\rangle - \E\langle G_2\rangle\right)}\right] -\frac{1}{4N^2}\sum_{\sigma\in\{\pm\}} \sigma\langle G_1AE_\sigma G_2^2E_\sigma\rangle,
\end{split}
\label{eq:under_ext}
\end{equation}
where $\widetilde{W}$ is an independent copy of $W$. To go from the first to the second line in \eqref{eq:under_ext} we additionally used the following simple identity:
\begin{equation}
\E \langle WR\rangle \langle WS\rangle = \frac{1}{4N^2}\left(\langle RE_+SE_+\rangle - \langle RE_-SE_-\rangle\right),\quad \forall R,S\in\C^{(2N)\times (2N)}.
\label{eq:EWW_id}
\end{equation}
Combining \eqref{eq:Cov_05_exp} with \eqref{eq:under_ext} we complete the proof of \eqref{eq:init_expansion}.

We decompose the last term in the first line of \eqref{eq:init_expansion} into the deterministic approximation and the fluctuation around it, and represent $G_2^2$ as a contour integral:
\begin{equation}
\begin{split}
&\left\langle G_1AE_\sigma G_2^2 E_\sigma\right\rangle - \left\langle M_{122}^{AE_\sigma,E_+}(w_1,w_2,w_2)E_\sigma\right\rangle\\
&\quad = \frac{1}{2\pi\ii}\int_\mathscr{C} \frac{1}{(\zeta-w_2)^2} \left(\left\langle G_1AE_\sigma G_2(\zeta) E_\sigma\right\rangle - \left\langle M_{12}^{AE_\sigma}(w_1,\zeta)E_\sigma\right\rangle\right)\dif\zeta,
\end{split}
\label{eq:G^2_rep}
\end{equation}
where $\mathscr{C}$ is a circle with the center at~$w_2$ and radius $\eta_2/2$, and $G_2(\zeta)=G^{z_2}(\zeta)$. To go from the first to the second line of~\eqref{eq:G^2_rep}, we used the standard meta-argument to show that the deterministic approximation also admits the integral representation (for more details on the meta-argument see e.g. \cite[Section~2.6]{Najim16}, \cite[Proof of Lemma~D.1]{non-herm_overlaps} and also the proof of Supplementary Lemma~\ref{lem:M_bounds_weak}). For $\eta_2\sim 1$, the contour $\mathscr{C}$ may exit the bulk regime, however for $\zeta \in\mathscr{C}$ with $\Re\zeta\notin \mathbf{B}_{\kappa/2}^{z_2}$ we have $|\Im\zeta|\sim 1$, in which case the fluctuation of the two-resolvent chain in the rhs. of~\eqref{eq:G^2_rep} has an upper bound of order $N^{-1+\xi}$ by the global law (for more details see Supplementary Eq.~\eqref{eq:global_2G}). Thus, applying~\eqref{eq:2G_av_b=0} for $\Re\zeta\in\mathbf{B}_{\kappa/2}^{z_2}$ and the global law in the complementary regime, we get from~\eqref{eq:G^2_rep} that
\begin{equation}
\left\vert \left\langle G_1AE_\sigma G_2^2 E_\sigma\right\rangle - \left\langle M_{122}^{AE_\sigma,E_+}(w_1,w_2,w_2)E_\sigma\right\rangle\right\vert \lesssim \frac{N^\xi}{N\eta_*\eta_2\gamma}, 
\label{eq:V_error}
\end{equation}
with very high probability for any fixed $\xi>0$. Here we additionally used \eqref{eq:beta_hat_trivial} and \eqref{eq:stab_LT} to show that
\begin{equation*}
\widehat{\beta}_{12}(w_1,\zeta) \sim \gamma(z_1,z_2,w_1,\zeta)\sim\gamma(z_1,z_2,w_1,w_2).
\end{equation*}
Combining~\eqref{eq:init_expansion}, \eqref{eq:psi_22}, and~\eqref{eq:V_error} we conclude that
\begin{equation}
\begin{split}
\Cov\left(\langle G_1\rangle,\langle G_2\rangle\right) &=\frac{1}{4N^2}\sum_{\sigma\in\{\pm\}} \sigma\left\langle M_{122}^{AE_\sigma,E_+}E_\sigma\right\rangle -\E \left[\underline{ \langle WG_1A\rangle \left(\langle G_2\rangle - \E\langle G_2\rangle\right)}\right]\\
&+\mathcal{O}\left(\left(\frac{1}{N\gamma}+\frac{1}{\sqrt{N}}\right)\frac{N^\xi}{N^2\eta_*\eta_2}\right),
\label{eq:Cov_1_exp}
\end{split}
\end{equation}
where we suppressed the arguments of $M_{122}$ for brevity.

We now estimate the contribution from the fully underlined term in \eqref{eq:Cov_1_exp}. Recall the notations introduced around \eqref{eq:def_index_set}--\eqref{eq:def_cum}. Using the cumulant expansion (see e.g. \cite[Lemma 3.2]{Knowles20}) similarly to \eqref{eq:cum_exp_first}, we get
\begin{equation}
\begin{split}
&\E\underline{\langle WG_1A\rangle\left(\langle G_2\rangle - \E\langle G_2\rangle\right)}\\
&\quad = \frac{1}{N}\sum_{a,b}\sum_{\ell=3}^L\sum_{\bm{\alpha}\in\{ab,ba\}^{\ell-1}}\frac{\kappa(ab,\bm{\alpha})}{N^{\ell/2}(\ell-1)!}\E\partial_{\bm{\alpha}}\left[ \left(G_1A\right)_{ba}\left(\langle G_2\rangle - \E\langle G_2\rangle\right)\right] + \mathcal{O}(N^{-D})
\end{split}
\label{eq:Cov_cum_exp}
\end{equation}
for any fixed $D>0$ and for some $L\in\N$ which depends on $D$, but does not depend on $N$. Here we used that due to the definition of the underline \eqref{eq:def_under_1}, the second order cumulants are absent in the expansion. In \eqref{eq:Cov_cum_exp} and throughout the proof, if not specified, summation runs over indices $a\in [1,N], b\in [N+1,2N]$ and $a\in [N+1,2N], b\in [1,N]$.

We call the parameter $\ell$ in \eqref{eq:Cov_cum_exp} the order of the corresponding group of terms, and treat separately the terms of order $\ell\ge 5$, $\ell = 4$, and $\ell=3$. We start with the high order terms ($\ell\ge 5$). Observe that for any $r\in\N$ and $\widehat{\bm{\alpha}}\in \{ab,ba\}^r$ it holds that
\begin{equation}
\left\vert \partial_{\widehat{\bm{\alpha}}} \left(G_1A\right)_{ba}\right\vert \prec 1,\quad \left\vert \partial_{\widehat{\bm{\alpha}}}\left(\langle G_2\rangle - \E\langle G_2\rangle\right) \right\vert \prec \frac{1}{N\eta_2}.
\label{eq:cum_dif_simple}
\end{equation}
Let us prove the first part of \eqref{eq:cum_dif_simple}, while the second part follows analogously. Performing the differentiation in $\partial_{\widehat{\bm{\alpha}}} \left(G_1A\right)_{ba}$ we get a sum of at most $r!$ terms of the form
\begin{equation}
(G_1)_{\beta_1}\cdots (G_1)_{\beta_r}(GA)_{\beta_{r+1}}\quad\text{with}\quad \beta_1,\ldots,\beta_{r+1}\in\{ab,ba,aa,bb\}. 
\label{eq:GA_iso_dif_term}
\end{equation}
We decompose each of the factors into the deterministic approximation and the fluctuation around it, estimate the fluctuations from above by $(N\eta_1)^{-1/2}$ by the single-resolvent isoropic local law~\eqref{eq:1G_ll_iso}, and bound $\|M_1\|\lesssim 1$. Multiplying the resulting bounds we prove the first part of \eqref{eq:cum_dif_simple}. From \eqref{eq:cum_dif_simple} we immediately conclude that for any fixed $\ell\ge 5$ the sum of all $\ell$-th order terms in \eqref{eq:Cov_cum_exp} has an upper bound of order
\begin{equation*}
\frac{1}{N}N^2 \frac{1}{N^{\ell/2}} \cdot  \frac{1}{N\eta_2}N^\xi\le \frac{1}{\sqrt{N}}\cdot\frac{N^\xi}{N^2\eta_*\eta_2}.
\end{equation*}

Now we consider the terms of order $\ell=4$ in \eqref{eq:Cov_cum_exp}. Each of these terms contains 3 derivatives, which we distribute according to the Leibniz rule over the factors $\left(G_1A\right)_{ba}$ and $\langle G_2\rangle -\E\langle G_2\rangle$. Consider at first the case when all of them hit $\left(G_1A\right)_{ba}$. We proceed by a slight refinement of the argument presented around~\eqref{eq:GA_iso_dif_term}. Performing the differentiation, we arrive to the sum consisting of terms of the form~\eqref{eq:GA_iso_dif_term}, each multiplied by $\langle G_2\rangle -\E\langle G_2\rangle$. Decomposing each isotropic factor in the way discussed below~\eqref{eq:GA_iso_dif_term}, we see that the product of deterministic approximations does not contribute to the expectation of the entire term because of the centered $\langle G_2\rangle$ factor. Therefore, we gain from at least one fluctuation and get
\begin{equation}
\left\vert \E\left[(G_1)_{\beta_1}(G_1)_{\beta_2}(G_1)_{\beta_3}(GA)_{\beta_{4}}(\langle G_2\rangle -\E\langle G_2\rangle)\right]\right\vert \prec \frac{1}{\sqrt{N\eta_1}}\cdot\frac{1}{N\eta_2}\le \frac{\sqrt{N}}{N^2\eta_*\eta_2}, 
\label{eq:kappa4_init_exp1}
\end{equation}
where we additionally used the single-resolvent averaged local law~\eqref{eq:1G_ll_av} to bound the last factor in the lhs. of~\eqref{eq:kappa4_init_exp1}. Thus, the fourth order terms in~\eqref{eq:Cov_cum_exp} where all three derivatives hit $(G_1A)_{ba}$ contribute at most $N^{-1/2}(N^2\eta_*\eta_2)^{-1}$.

In the case when at least one out of three derivatives in the $\ell=4$ term hits $\langle G_2\rangle-\E\langle G_2\rangle$, all terms produced by the differentiation are of the form
\begin{equation}
\frac{1}{N}(G_2^2)_{\beta_1}(G_1A)_{\beta_2}G_{\beta_3}G_{\beta_4}\quad\text{with}\quad \beta_1,\ldots,\beta_{4}\in\{ab,ba,aa,bb\},
\label{eq:init_kappa_4_2}
\end{equation}
where each of $G$'s without subscript equals either to $G_1$ or to $G_2$. Again decomposing each of the factors into the deterministic approximation and the fluctuation and arguing as in \eqref{eq:kappa4_init_exp1}, we conclude that if at least one factor in \eqref{eq:init_kappa_4_2} is replaced by its fluctuation, then the entire term is bounded by $N^{1/2}(N^2\eta_*\eta_2)^{-1}$, e.g.   
\begin{equation}
\frac{1}{N}\left(G_2^2-M_{22}^{E_+}\right)_{\beta_1}(M_1A)_{\beta_2}M_{\beta_3}M_{\beta_4}\prec \frac{1}{N}\cdot \frac{1}{\sqrt{N\eta_2}\eta_2}\lesssim \frac{\sqrt{N}}{N^2\eta_*\eta_2},
\label{eq:kappa4_init_exp2}
\end{equation}
where we used the convention that $M$ without a subscript coincides either with $M_1$ or with $M_2$, and employed \eqref{eq:multiG_iso} for $k=2$ in the first bound. 

In order to complete the analysis of the 4th order terms in~\eqref{eq:Cov_cum_exp} it remains to compute the contribution  from~\eqref{eq:init_kappa_4_2}, where each of the factors is replaced by its deterministic counterpart, to the rhs. of~\eqref{eq:Cov_1_exp}. In contrast to~\eqref{eq:kappa4_init_exp2}, the contribution from the sum of these terms cannot be fully incorporated into the error term in the rhs. of~\eqref{eq:Cov_1_exp}, as it contributes to the leading order in the rhs. of~\eqref{eq:main_cov}. We start with the identity
\begin{equation}
N^{-1}(M_{22}^{E_+})_{\beta_1}(M_1A)_{\beta_2}M_{\beta_3}M_{\beta_4}=N^{-1}\partial_{w_2}(M_2)_{\beta_1}\partial_{w_1}(M_1)_{\beta_2}M_{\beta_3}M_{\beta_4},
\label{eq:kappa4_determ}
\end{equation}
where we computed
\begin{equation}
M_{22}^{E_+}=\frac{M_2^2}{1-\langle M_2^2\rangle} = \partial_{w_2} M_2\quad\text{and}\quad M_1A = \frac{M_1^2}{1-\langle M_1^2\rangle} =\partial_{w_1}M_1.
\label{eq:kappa4_id}
\end{equation}
The first identity in~\eqref{eq:kappa4_id} follows from~\eqref{eq:def_M12}, \eqref{eq:def_B12}, and \eqref{eq:MDE}, while the second one is an immediate consequence of~\eqref{eq:MDE} and the explicit formula for~$A$ given below \eqref{eq:B*_inverse}. If at least one of $\beta_j$, $j\in\!~ [4]$, in~\eqref{eq:kappa4_determ} is off-diagonal, then due to the block-constant structure of $M^z(w)$, the rhs. of~\eqref{eq:kappa4_determ} may be non-zero only for $a\in [N], b=N+a$ and $b\in[N],a=N+b$. Therefore, terms of the form~\eqref{eq:kappa4_determ} with at least one off-diagonal factor contribute at most~$N^{-3}$ to the sum in~\eqref{eq:Cov_cum_exp} and thus can be incorporated into the error term in the rhs. of~\eqref{eq:Cov_1_exp}. Here we additionally used that $\|\partial_w M^z(w)\|\lesssim 1$ for $\Re w$ in the bulk of~$\rho^z$.

Consider now the case when all factors in the rhs. of~\eqref{eq:kappa4_determ} are diagonal, i.e. $\beta_j\in\{aa,bb\}$ for $j~\!\in~\![4]$. Such terms arise in~\eqref{eq:Cov_cum_exp} only when $(G_1A)_{ba}$ is differentiated once in~$w_{ab}$ and $\langle G_2\rangle-\E\langle G_2\rangle$ is differentiated both in~$w_{ab}$ and~$w_{ba}$. In this case $\kappa(ab,\bm{\alpha})=\kappa_4$ and exactly one~$M$ without a subscript in the rhs. of~\eqref{eq:kappa4_determ} equals to~$M_1$ and one to~$M_2$. Together with the fact that all diagonal entries of $M^z(w)$ are equal to $m^z(w)$ by~\eqref{eq:M}, this gives that the rhs. of~\eqref{eq:kappa4_determ} equals to
\begin{equation}
N^{-1} \left(\partial_{w_2}\langle M_2\rangle\right) \left(\partial_{w_1}\langle M_1\rangle\right) \langle M_1\rangle\langle M_2\rangle = (4N)^{-1} \partial_{w_1}\langle M_1\rangle^2\partial_{w_2}\langle M_2\rangle^2 = (2N)^{-1}U_1U_2,
\label{eq:kappa4_contrib} 
\end{equation}  
where $U_1,U_2$ are defined in~\eqref{eq:def_UV}. As it is easy to see, this term arises $6N^2$ times in the rhs. of~\eqref{eq:Cov_cum_exp}, and each time it carries the factor $-\kappa_4/(3!N^3)$. Therefore, the deterministic approximations to the 4th order terms in the rhs. of~\eqref{eq:Cov_cum_exp} with all diagonal factors contribute $\kappa_4U_1U_2/(2N^2)$ to the rhs. of~\eqref{eq:Cov_1_exp}.

Finally, we analyze the 3rd order terms ($\ell=3$) in \eqref{eq:Cov_cum_exp}. We will show that their contribution to $\Cov(\langle G_1\rangle, \langle G_2\rangle)$ is negligible and can be incorporated into the error term. Each of the third order terms contains two derivatives. Let us start with the case when both of them hit $(G_1A)_{ba}$ factor. Performing the differentiation, we arrive to a product of 3 isotropic terms. Counting the number of $a$'s in the indices we see that there is at least one off-diagonal factor in this product. Consider the exemplary case when there is exactly one off-diagonal factor: 
\begin{equation}
\frac{1}{N^{5/2}} \E\left(\sum_{a,b} \frac{\kappa(ab,ab,ba)}{2}  (G_1)_{bb}(G_1A)_{ba}(G_1)_{aa}\right)\left(\langle G_2\rangle - \E\langle G_2\rangle\right).
\label{eq:init_kappa_3}
\end{equation}
Let us restrict the summation to $a\in [1,N]$, $b\in [N+1,2N]$, since the estimates in the regime $a>b$ are similar. Decompose resolvents in the diagonal factors as $G_1=M_1+(G_1-M_1)$. In the term with both deterministic terms we perform two isotropic resummations
\begin{equation*}
\sum_{a=1}^N\sum_{b=N+1}^{2N}(M_1)_{bb}(G_1A)_{ba}(M_1)_{aa} = (G_1A)_{\bm{x}\bm{y}}, \quad \bm{x},\bm{y}\in\C^{2N},
\end{equation*}
where $\bm{x}_b=(M_1)_{bb}$ for $a\in[N+1,2N]$ and $\bm{x}_b=0$ for $a\le N$ ($\bm{y}$ is defined similarly with last $N$ being equal to zero). Using the bound $(G_1A)_{\bm{x}\bm{y}}\prec N/\eta_1$ we get
\begin{equation*}
\frac{1}{N^{5/2}} \E\left(\sum_{a=1}^N\sum_{b=N+1}^{2N} \frac{\kappa(ab,ab,ba)}{2}  (M_1)_{bb}(G_1A)_{ba}(M_1)_{aa}\right)\left(\langle G_2\rangle - \E\langle G_2\rangle\right)\lesssim \frac{1}{\sqrt{N}}\cdot \frac{N^\xi}{N^2\eta_1\eta_2}.
\end{equation*}
In the case when only one $G_1$ is replaced by the deterministic counterpart in~\eqref{eq:init_kappa_3}, we perform one isotropic resummation and estimate the fluctuation factor by the isotropic local law~\eqref{eq:1G_ll_iso}
\begin{equation*}
\begin{split}
&\sum_{a=1}^N\sum_{b=N+1}^{2N}(M_1)_{bb}(G_1A)_{ba}(G_1-M_1)_{aa} = \sum_{a=1}^N(G_1A)_{\bm{x}a}(G_1-M_1)_{aa}\\
&\quad\prec \sum_{a=1}^N \left( \left\vert(M_1A)_{\bm{x}a}\right\vert + \frac{\sqrt{N}}{\sqrt{N\eta_1}}\right)\left(\frac{1}{\sqrt{N\eta_1}}\wedge \frac{1}{\sqrt{N}\eta_1}\right)\lesssim \frac{\sqrt{N}}{\eta_1}. 
\end{split}
\end{equation*}
In the case when both $G_1$ are replaced by the fluctuation (i.e. by $G_1-M_1$), we use isotropic local law~\eqref{eq:1G_ll_iso} for these factors and do not perform isotropic resummation. In such a way we conclude that the sum in \eqref{eq:init_kappa_3} has an upper bound of order $N^{-5/2}(\eta_1\eta_2)^{-1}$. The rest of the terms coming from the third order cumulant are estimated in the same way. For more details see \cite[Eq. (6.35)-(6.36)]{macroCLT_complex}.

Combining~\eqref{eq:Cov_1_exp} with the analysis of the cumulant expansion~\eqref{eq:Cov_cum_exp} performed above, we get that
\begin{equation}
\Cov\left(\langle G_1\rangle,\langle G_2\rangle\right) = \frac{1}{2N^2}\left(\frac{1}{2}\sum_{\sigma\in\{\pm\}} \sigma\left\langle M_{122}^{AE_\sigma,E_+}E_\sigma\right\rangle + \kappa_4U_1U_2\right) +\mathcal{O}\left(\left(\frac{1}{N\gamma}+\frac{1}{\sqrt{N}}\right)\frac{N^\xi}{N^2\eta_*\eta_2}\right).
\label{eq:Cov_F1}
\end{equation}
Let us explicitly compute the sum over $\sigma\in\{\pm\}$ in the rhs. of~\eqref{eq:Cov_F1}. First, for $\sigma\in\{\pm\}$ we have  
\begin{equation}
\left\langle M_{122}^{AE_\sigma,E_+}E_\sigma\right\rangle = \partial_{w_2}\left\langle M_{12}^{AE_\sigma}E_\sigma\right\rangle = \partial_{w_2}\left\langle M_{21}^{E_\sigma}AE_\sigma\right\rangle = \partial_{w_2}\left\langle\mathcal{B}_{21}^{-1}[M_2E_\sigma M_1] AE_\sigma\right\rangle.
\label{eq:V_comp1}
\end{equation}
To get the first two identities in~\eqref{eq:V_comp1} we used the meta-argument twice, while the last identity follows from~\eqref{eq:def_M12}. Next, inverting $\mathcal{B}_{21}$ on $M_2E_\pm M_1$ using~\eqref{eq:def_B12}, we get
\begin{equation}
\begin{split}
\mathcal{B}_{21}^{-1}[M_2E_+M_1] &= \frac{(1+\langle M_1E_-M_2E_-\rangle)M_2E_+M_1 -\langle M_1E_-M_2\rangle M_2E_-M_1}{(1+\langle M_1E_-M_2E_-\rangle)(1-\langle M_1M_2\rangle) + \langle M_1E_-M_2\rangle\langle M_2E_-M_1\rangle},\\
\mathcal{B}_{21}^{-1}[M_2E_-M_1] &= \frac{\langle M_2E_-M_1\rangle M_2E_+M_1 +(1-\langle M_1M_2\rangle) M_2E_-M_1}{(1+\langle M_1E_-M_2E_-\rangle)(1-\langle M_1M_2\rangle) + \langle M_1E_-M_2\rangle\langle M_2E_-M_1\rangle}.
\end{split}
\label{eq:V_comp2}
\end{equation}
Combining \eqref{eq:V_comp1}, \eqref{eq:V_comp2}, the explicit formula for~$A$ given below~\eqref{eq:B*_inverse}, and the last identity in~\eqref{eq:kappa4_id}, we obtain
\begin{equation}
\sum_{\sigma\in\{\pm\}}\!\! \sigma\left\langle M_{122}^{AE_\sigma,E_+}E_\sigma\right\rangle\! = \!- \partial_{w_1}\partial_{w_2}\log \Big((1+\langle M_1E_-M_2E_-\rangle)(1-\langle M_1M_2\rangle) \!+ \!\langle M_1E_-M_2\rangle\langle M_2E_-M_1\rangle\Big).
\label{eq:V_comp3}
\end{equation}
We use \eqref{eq:M} for $M_1,M_2$, and conclude that the rhs. of~\eqref{eq:V_comp3} equals to $2V_{12}$, where~$V_{12}$ is defined in~\eqref{eq:def_UV}. Together with~\eqref{eq:Cov_F1} this verifies~\eqref{eq:main_cov_b} in a slightly weaker version, specifically with $\eta_1$ being replaced by $\eta_*$. Interchanging $G_1$ and $G_2$ and performing the same argument as above, we see that~\eqref{eq:Cov_F1} holds also with $\eta_2$ being replaced by $\eta_*$. Since $\eta_*\eta_2+\eta_1\eta_*\sim \eta_1\eta_2$, this finishes the proof of Proposition~\ref{prop:Cov_b}. 
\end{proof}

\subsection{Hierarchy of covariances: Proof of Proposition \ref{prop:psi_22}}\label{sec:H_Cov}

To prove Proposition~\ref{prop:psi_22}, we iteratively expand the lhs. of~\eqref{eq:psi_22} in a manner analogous to~\eqref{eq:init_expansion}. Since our goal in~\eqref{eq:psi_22} is to obtain a size bound, rather than a leading order term, we do not need to compute the leading order terms in the rhs. of these expansions, as it was done e.g. in \eqref{eq:kappa4_determ}--\eqref{eq:kappa4_contrib} and \eqref{eq:V_comp1}--\eqref{eq:V_comp3} above. This allows us to work with a hierarchy of inequalities deduced from the iterative underline expansions instead of more complicated  identities obtained by these expansions. 

Before formulating this hierarchy of inequalities in Proposition~\ref{prop:master_av}, we introduce some notation. For observables $B_i\in \C^{(2N)\times (2N)}$, $i\in[k]$, we define
\begin{equation}
\Delta(\bm{B}):=\langle \left(G_1 B_1 G_1\cdots B_{k-1}G_1 - \mathcal{M}[G_1 B_1 G_1\cdots B_{k-1}G_1]\right) B_k\rangle,\quad \bm{B}=(B_1,\ldots,B_k),
\label{eq:def_delta}
\end{equation}
For brevity, the dependence of $\Delta(\bm{B})$ on $z_1,w_1$ and $N$ is suppressed. Note that $\Delta(\bm{B})$ is a fluctuation of a resolvent chain constructed solely from $G_1$, while $G_2$ does not appear in this definition. For a collection of matrix tuples $\bm{B}^{(j)}=\big(B_1^{(j)},\ldots, B_{k_j}^{(j)}\big)$, $j\in [n]$, we further set
\begin{equation*}
\Delta\left(\overline{\bm{B}}\right)=\Delta\big(\bm{B}^{(1)};\ldots;\bm{B}^{(n)}\big):=\prod_{j=1}^n \Delta\big(\bm{B}^{(j)}\big),\quad \overline{\bm{B}} = \big(\bm{B}^{(1)},\ldots, \bm{B}^{(n)}\big).
\end{equation*}
Here the bar in $\overline{\bm{B}}$ is used to distinguish between a vector of matrix tuples and a matrix tuple, which we denote simply by $\bm{B}$. Next, denote
\begin{equation}
\mathcal{C}\left(\overline{\bm{B}}\right)=\mathcal{C}\big(\bm{B}^{(1)};\ldots;\bm{B}^{(n)}\big):=\Cov \left(\Delta\left(\overline{\bm{B}}\right),\langle G_2\rangle\right).
\end{equation}
For example, this notation allows us to abbreviate
\begin{equation*}
\Cov(\langle G_1\rangle,\langle G_2\rangle)=\mathcal{C}(E_+)\quad\text{and}\quad  \Cov\left(\langle G_1-M_1\rangle\langle (G_1-M_1)M_1\rangle, \langle G_2\rangle\right)=\mathcal{C}\left(E_+; M_1\right).
\end{equation*}
In the last identity we denoted $\mathcal{C}\left(E_+; M_1\right):=\mathcal{C}\left((E_+); (M_1)\right)$, that is, we omit parentheses in the notation for a matrix tuple whenever the tuple consists of a single matrix. We will further adopt this convention when referring to $\Delta(\overline{\bm{B}})$ and $\mathcal{C}(\overline{\bm{B}})$.

To introduce our main control parameters, which are the maxima of $|\mathcal{C}(\overline{\bm{B}})|$ over appropriate sets of $\overline{\bm{B}}$, we restrict attention to observables $B\in\{E_+,E_-,F,F^*\}$, and for
\begin{equation}
\overline{\bm{B}}= \big(\bm{B}^{(1)},\ldots, \bm{B}^{(n)}\big)\,\,\text{with}\,\, \bm{B}^{(i)}=\big(B^{(i)}_1,\ldots, B^{(i)}_{k_i}\big),i\in[n],\quad\text{denote}\quad \mathfrak{n}(\overline{\bm{B}}):=k_1+\ldots+k_n,
\end{equation}
i.e. $\mathfrak{n}(\overline{\bm{B}})$ counts the number of matrices in $\overline{\bm{B}}$. For $n,S\in\N$, let $\mathfrak{S}_{n,S}$ be the set of all vectors $\overline{\bm{B}} = \big(\bm{B}^{(1)},\ldots, \bm{B}^{(n)}\big)$ consisting of $n$ matrix tuples with $\mathfrak{n}(\overline{\bm{B}})\le S$. Our main control parameters are given by
\begin{equation}
\Psi_{n,S}:=\max \left\lbrace N^n\eta_*^{\mathfrak{n}(\overline{\bm{B}})}N\eta_2\left\vert \mathcal{C}\left(\overline{\bm{B}}\right)\right\vert\,:\,\overline{\bm{B}}\in \mathfrak{S}_{n,S}\right\rbrace.
\label{eq:def_psi}
\end{equation}
In the case when $S<n$ or $n\le 0$, the set $\mathfrak{S}_{n,S}$ is empty, and $\Psi_{n,S}$ is set to be equal to zero. Note that $S$ is an upper bound on the number of $G_1$'s in $\mathcal{C}(\overline{\bm{B}})$, and that there is one more resolvent $G_2$ in the definition of $\mathcal{C}(\overline{\bm{B}})$. The normalization in \eqref{eq:def_psi} is chosen in such a way that for any $\xi>0$ we have
\begin{equation}
\Psi_{n,S}\lesssim N^\xi
\label{eq:psi_a_priori}
\end{equation} 
by the multi-resolvent local laws from Proposition~\ref{prop:multiG_oneG}.

Now we state the hierarchy of inequalities, which relates the control parameters $\Psi_{n,S}$ for different values of~$n$ and~$S$.

\begin{proposition}[Master inequalities]\label{prop:master_av} Assume the set-up and conditions of Proposition~\ref{prop:Cov_b} with $b=0$. Fix a (small) $\xi>0$ and denote
\begin{equation}
\Upsilon=\Upsilon(z_1,z_2;\eta_1,\eta_2):=\left(\sqrt{N}+\frac{1}{\gamma}\right)\eta_*N^\xi,
\label{eq:def_Upsilon}
\end{equation}
where $\gamma$ is defined in \eqref{eq:def_UV}. Then for any $n,S\in\N$ it holds that
\begin{equation}
\Psi_{n,S}\lesssim \Psi_{n-2,S-2} + \Psi_{n,S-1}+\frac{1}{N\eta_*}\left(\Psi_{n-1,S+1}+\Psi_{n+1,S+1}\right) + \left(1+\bm{1}_{S>n}(N\eta_*)^{1/2}\right)\Upsilon.\label{eq:master_S}
\end{equation}
\end{proposition}

We refer to the inequalities in Proposition \ref{prop:master_av} as \emph{master inequalities}, adopting the terminology from \cite{Multi_res_llaws}, where a similar hierarchy \cite[Eq. (3.20a), (3.20b)]{Multi_res_llaws} was employed to prove the analogues of the multi-resolvent local laws \eqref{eq:multiG_av}, \eqref{eq:multiG_iso} for Wigner matrices. The key difference is that in \cite{Multi_res_llaws} the hierarchy was used to control the magnitude of fluctuations of a product of resolvents, while now we control the expectation.


The proof of Proposition~\ref{prop:master_av} is postponed to Section~\ref{subsec:master_proof}. Having~\eqref{eq:master_S} in hand, we complete the proof of Proposition~\ref{prop:psi_22} by iterating~\eqref{eq:master_S}.

\begin{proof}[Proof of Proposition \ref{prop:psi_22}] First note that \eqref{eq:psi_22} is equivalent to the bound
\begin{equation}
\Psi_{2,2} \lesssim \Upsilon.
\label{eq:psi_22_ups}
\end{equation}
To prove \eqref{eq:psi_22_ups} we show that for any $n,S\in\N$ it holds that
\begin{equation}
\Psi_{n,S}\lesssim (N\eta_*)^{1/2}\Upsilon.
\label{eq:av_impr}
\end{equation}
Once \eqref{eq:av_impr} is obtained, we apply \eqref{eq:master_S} to $n=S=2$ and then use \eqref{eq:av_impr} for $\Psi_{1,3}$ and $\Psi_{3,3}$, which gives
\begin{equation}
\Psi_{2,2}\lesssim \frac{1}{N\eta_*} \left(\Psi_{1,3}+\Psi_{3,3}\right) + \Upsilon \lesssim \frac{1}{N\eta_*} (N\eta_*)^{1/2}\Upsilon +\Upsilon \lesssim \Upsilon,
\end{equation} 
i.e. \eqref{eq:psi_22_ups} holds.

Now we prove \eqref{eq:av_impr}. For $n>S$ this bound is a triviality since in this case $\Psi_{n,S}=0$. Hence, whenever convenient, we will assume that $n\le S$ without mentioning this explicitly. We claim that
\begin{equation}
\Psi_{n,S}\lesssim \frac{1}{N\eta_*}\sum_{r=1}^{n+1}\Psi_{r,S+1}+(N\eta_*)^{1/2}\Upsilon,\quad \forall n,S\in\N.
\label{eq:1_iter}
\end{equation}
We derive \eqref{eq:1_iter} from \eqref{eq:master_S}, without exploiting the improvement in the last term in the rhs. of \eqref{eq:master_S} occurring for $n\le S$, and simply estimating this term from above by $(N\eta_*)^{1/2}\Upsilon$.
The proof of \eqref{eq:1_iter} proceeds by induction on $n$ and $S$ which we now describe. As the base case we consider $n=S=1$, where \eqref{eq:1_iter} directly follows from \eqref{eq:master_S} applied to $\Psi_{1,1}$. Next we fix $n:=1$ and proceed by induction on $S$. To perform the induction step from $S\le S_0$ to $S=S_0+1$ for some $S_0\in\N$, we use \eqref{eq:master_S} for $n=1$ and $S=S_0+1$:
\begin{equation}
\begin{split}
\Psi_{1,S_0+1}&\lesssim \Psi_{1,S_0} + \frac{1}{N\eta_*} \Psi_{2,S_0+2}+(N\eta_*)^{1/2} \Upsilon\lesssim \frac{1}{N\eta_*}\left(\Psi_{1,S_0+1}+\Psi_{2,S_0+1}+\Psi_{2,S_0+2}\right) +(N\eta_*)^{1/2} \Upsilon\\
&\lesssim \frac{1}{N\eta_*}\left(\Psi_{1,S_0+2}+\Psi_{2,S_0+2}\right) +(N\eta_*)^{1/2} \Upsilon.
\end{split}
\label{eq:n=1_ind_step}
\end{equation}
Here in the second step we used \eqref{eq:1_iter} for $(n,S_0)$ and to go from the first to the second line we used that $\Psi_{n,S_1}\le \Psi_{n,S_2}$, for any $n,S_1,S_2\in \N$ with $S_1\le S_2$, which follows from \eqref{eq:def_psi}. Once \eqref{eq:1_iter} has been established for $n=1$ and all $S\in\N$, we fix $n:=2$ and repeat the same argument, continuing this process inductively for all $n\in\N$. The induction step from $S\le S_0$ to $S=S_0+1$ for general $n,S\in\N$ follows similarly to \eqref{eq:n=1_ind_step} from \eqref{eq:master_S} for $(n,S_0+1)$ together with \eqref{eq:1_iter} for $\Psi_{n-2,S_0-1}$ and $\Psi_{n,S_0}$, and the monotonicity of $\Psi_{n,S}$ in $S$. This completes the verification of~\eqref{eq:1_iter}.

We use \eqref{eq:1_iter} for $\Psi_{r,S+1}$ for all $r\in [n+1]$ and combine these bounds with \eqref{eq:1_iter} for $\Psi_{r,S}$, arriving to 
\begin{equation}
\Psi_{n,S}\lesssim \left(\frac{1}{N\eta_*}\right)^k\sum_{r=1}^{n+k}\Psi_{r,S+k}+(N\eta_*)^{1/2}\Upsilon.
\label{eq:2_iter}
\end{equation}
for $k=2$. Iterating this procedure several times, we obtain \eqref{eq:2_iter} for any fixed $k\in\N$. Finally, we choose $k>1/\epsilon$ in \eqref{eq:2_iter} and use the a priori bound \eqref{eq:psi_a_priori}  for all $\Psi$ terms in the rhs. of \eqref{eq:2_iter}. This finishes the proof of \eqref{eq:av_impr} and thereby of Proposition~\ref{prop:psi_22}.
\end{proof}

\subsection{Proof of the master inequalities in Proposition \ref{prop:master_av}}\label{subsec:master_proof}

In this section we prove Proposition~\ref{prop:master_av}, which is the last technical ingredient in the proof of Proposition~\ref{prop:Cov}. To do so, we first derive the system of identities similar to~\eqref{eq:init_expansion}, relating the quantities of the form $\mathcal{C}(\overline{\bm{B}})$ for different values of $\overline{\bm{B}}$. This is the statement of Lemma~\ref{lem:av_C_dif}. Then we prove Lemma~\ref{lem:av_C_dif}, and finally turn it into the system of master inequalities, thereby proving Proposition~\ref{prop:master_av}.

To streamline the presentation, we introduce some additional short-hand notations. For $n\in\N$ and for a set of distinct indices $\{i_1,\ldots,i_s\}\subset [n]$ denote
\begin{equation*}
\overline{\bm{B}}^{[i_1,\ldots,i_s]}:=\overline{\bm{B}}\setminus\left\lbrace\bm{B}^{(i_1)},\ldots,\bm{B}^{(i_s)}\right\rbrace,\quad\text{where}\quad \overline{\bm{B}}=(\bm{B}^{(1)},\ldots,\bm{B}^{(n)}).
\end{equation*}
We also set
\begin{equation*}
m(\bm{B})=m(B_1,\ldots,B_k):=\langle \mathcal{M}[G_1B_1G_1\cdots B_{k-1}G_1] B_k\rangle.
\end{equation*}
Throughout the calculations below, whenever the index $\sigma$ appears, it is meant that the corresponding expression is summed over $\sigma\in\{\pm\}$. Additionally we adopt the convention that a summation with a lower limit exceeding the upper one is considered to be empty.

The following statement is the analogue of the initial expansion \eqref{eq:init_expansion} for the general covariance $\mathcal{C}(\overline{\bm{B}})$. In its formulation, sequences of matrices of the form $B_i, B_{i+1}, \ldots, B_j$ frequently appear as arguments of $\mathcal{C}$ for some $i,j\in\N\cup\{0\}$. By this notation we mean the sequence of matrices indexed from $i$ to $j$ when $i<j$, the single matrix $B_i$ when $i=j$, and the empty sequence when $i>j$. 

\begin{lemma}\label{lem:av_C_dif} Assume the set-up and assumptions of Proposition~\ref{prop:Cov_b} with $b=0$. Let $\bm{B}^{(j)}=\big(B_1^{(j)},\ldots, B_{k_j}^{(j)}\big)$ for $j\in[n]$ be a vector consisting of $(2N)\times (2N)$ deterministic matrices. Denote 
\begin{equation}
A=A_{k_1}^{(1)}:=\left(\left(\mathcal{B}_{11}^*\right)^{-1}\left[ \left(B_{k_1}^{(1)}\right)^*\right]\right)^*M_1.
\label{eq:def_A}
\end{equation}
Let ${\rm cycl}(k)$ be the set of all cyclic permutations of length $k$. Then we have
\begin{align}
\mathcal{C}(\overline{\bm{B}}) \!&= \! \frac{\sigma}{4N^2}\!\E\left[\langle G_1B_1^{(1)}G_1\cdots G_1B_{k_1-1}^{(1)}G_1AE_\sigma G_2^2 E_\sigma\rangle\Delta\left(\overline{\bm{B}}^{[1]}\right)\right]\!+\!\mathcal{C}\!\left(\!\left(AB_1^{(1)},B_2^{(1)},\ldots,B_{k_1-1}^{(1)}\right)\!; \overline{\bm{B}}^{[1]}\right) \label{eq:av_C_dif}\\ 
&+\sum_{j=1}^{k_1-1} \sigma m\left(E_\sigma, B_1^{(1)},\ldots, B_j^{(1)}\right)\mathcal{C}\left(\left(AE_\sigma, B_{j+1}^{(1)},\ldots, B_{k_1-1}^{(1)}\right); \overline{\bm{B}}^{[1]}\right) \nonumber\\ 
&+\sum_{j=0}^{k_1-2} \sigma m\left(AE_\sigma, B_{j+1}^{(1)},\ldots, B_{k_1-1}^{(1)}\right)\mathcal{C}\left(\left(E_\sigma,B_1^{(1)},\ldots, B_j^{(1)}\right); \overline{\bm{B}}^{[1]}\right)\nonumber\\
&+ \sum_{j=0}^{k_1-1} \sigma \mathcal{C}\left(\left(E_\sigma,B_1^{(1)},\ldots,B_j^{(1)}\right);\left(AE_\sigma,B_{j+1}^{(1)},\ldots,B_{k_1-1}^{(1)}\right);\overline{\bm{B}}^{[1]}\right)\nonumber\\
&+\frac{1}{4N^2}\sum_{i=2}^n \sum_{c\in {\rm{cycl}}(k_i)}\sigma m\left(AE_\sigma,c(\bm{B}^{(i)}),E_\sigma,B_1^{(1)},\ldots, B_{k_1-1}^{(1)}\right)\mathcal{C}\left(\overline{\bm{B}}^{[1,i]}\right)\nonumber\\
&+\frac{1}{4N^2}\sum_{i=2}^n \sum_{c\in {\rm{cycl}}(k_i)} \sigma\mathcal{C}\left(\left(AE_\sigma,c(\bm{B}^{(i)}),E_\sigma, B_1^{(1)},\ldots, B_{k_1-1}^{(1)}\right);\overline{\bm{B}}^{[1,i]}\right)\nonumber\\
&-\E\left[\underline{\langle AWG_1B_1^{(1)}G_1\cdots G_1B_{k_1-1}^{(1)}G_1\rangle\Delta\left(\overline{\bm{B}}^{[1]}\right)\left(\langle G_2\rangle - \E\langle G_2\rangle\right)}\right].\nonumber
\end{align}
\end{lemma}
We choose to invert $\mathcal{B}_{11}$ on $\big(B_{k_1}^{(1)}\big)^*$ in Lemma~\ref{lem:av_C_dif} only for definiteness, while a similar expansion holds with $B_{k_1}^{(1)}$ replaced in~\eqref{eq:def_A} by $B_i^{(l)}$ for any $l\in [n]$ and $i\in [k_l]$.

Before delving into the proof of Lemma~\ref{lem:av_C_dif}, let us demonstrate~\eqref{eq:av_C_dif} with some simple examples. For $n=1$ and $k_1=1$ only a few terms in the rhs. of~\eqref{eq:av_C_dif} are left: the first term in the rhs. of the first line, $j=0$ term in the fourth line, and the fully underlined term in the last line. Thus,~\eqref{eq:av_C_dif} implies that
\begin{equation*}
\mathcal{C}(B) =  \frac{\sigma}{4N^2} \E\langle G_1AE_\sigma G_2^2E_\sigma\rangle + \mathcal{C}(E_+;A) - \E\underline{\langle WG_1A\rangle \left(\langle G_2\rangle -\E \langle G_2\rangle\right)},
\end{equation*}
which generalizes \eqref{eq:init_expansion} formulated in the spacial case $B=E_+$. Here we additionally used that $\langle G_1E_-\rangle=0$, so the $\sigma=-$ term in the fourth line of \eqref{eq:av_C_dif} vanishes. Meanwhile, for $n=2$ and $\bm{B}^{(1)}$, $\bm{B}^{(2)}$ consisting of single matrices $B^{(1)}$ and $B^{(2)}$, \eqref{eq:av_C_dif} is of the form
\begin{equation}
\begin{split}
\mathcal{C}(B^{(1)};B^{(2)})& = \frac{\sigma}{4N^2}\E\langle G_1AE_\sigma G_2^2 E_\sigma\rangle\langle (G_1-M_1)B^{(2)}\rangle + \mathcal{C}(E_+;A;B^{(2)})\\
&+\frac{\sigma}{4N^2}\mathcal{C}(AE_\sigma, B^{(2)}, E_\sigma) - \E\underline{\langle WG_1A\rangle \langle (G_1-M_1)B^{(2)}\rangle\left(\langle G_2\rangle -\E \langle G_2\rangle\right)},
\end{split}
\label{eq:C_exp_example}
\end{equation}
where the first term in the rhs. of~\eqref{eq:C_exp_example} comes from the first term in the rhs. of~\eqref{eq:av_C_dif}, while the second and third terms come from the $j=0$ term in the fourth line of~\eqref{eq:av_C_dif} and $i=2$ term in the sixth line of~\eqref{eq:av_C_dif}, respectively. 

\begin{proof}[Proof of Lemma \ref{lem:av_C_dif}] At first we prove \eqref{eq:av_C_dif} for $n=1$. Denote for brevity $\bm{B}^{(1)}=\bm{B}=(B_1,\ldots,B_k)$, $G:=G_1$ and $M:=M_1$.  The case $k=1$ with $B_1=E_+$ is already covered by \eqref{eq:<G-M>}--\eqref{eq:under_ext}, and for general $B_1$ the argument is identical, so we further focus on the case $k\ge 2$. From \cite[Eq.(5.2)]{macroCLT_complex} and the fact that $\mathcal{S}[G-M]=\langle G-M\rangle$ (see also Lemma~\ref{lem:E_-} later) we have
\begin{equation}
G=M-M\underline{WG}+M\langle G-M\rangle G.
\label{eq:G_exp}
\end{equation}
We multiply this identity on the right by $B_1G\cdots B_{k-1}G$ and in the term containing $\underline{WG}$ extend the underline onto the entire product using \eqref{eq:def_under_1}:
\begin{equation}
\begin{split}
GB_1G\cdots B_{k-1}G &= MB_1G\cdots B_{k-1}G +M\langle G-M\rangle G B_1G\cdots B_{k-1}G\\
& + M\sum_{j=1}^{k-1}\mathcal{S}[GB_1G\cdots B_jG]GB_{j+1}G\cdots B_{k-1}G-M\underline{WGB_1G\cdots GB_{k-1}G}.
\end{split}
\label{eq:GBG_exp}
\end{equation}
Next, in the $j=k-1$ term in \eqref{eq:GBG_exp} we decompose the last $G$ into $G-M$ and $M$, and move the term containing $M$ to the lhs. of this identity, obtaining $\mathcal{B}_{11}[GB_1G\cdots B_{k-1}G]$. Multiplying the resulting identity from the left by $AM^{-1}$ with $A$ defined in \eqref{eq:def_A}, and taking the normalized trace we get
\begin{equation}
\begin{split}
&\langle GB_1G\cdots B_{k-1}GB_k\rangle = \langle AB_1G\cdots B_{k-1}G\rangle  + \langle G-M\rangle \langle GAGB_1\cdots GB_{k-1}\rangle\\
&\quad + \sum_{j=1}^{k-2}\sigma\langle GE_\sigma GB_1\cdots GB_j\rangle\langle G AE_\sigma G B_{j+1} \cdots GB_{k-1}\rangle\\
&\quad+\sigma\langle GE_\sigma GB_1\cdots GB_{k-1}\rangle \langle (G-M)AE_\sigma\rangle - \langle A\underline{WGB_1G\cdots GB_{k-1}G}\rangle.
\end{split}
\label{eq:Delta_exp_aux}
\end{equation}
We subtract the deterministic approximations from both sides of \eqref{eq:Delta_exp_aux} using \eqref{eq:def_multiM} and arrive to the following underline expansion for $\Delta(\bm{B})$: 
\begin{align}
\Delta(\bm{B})&=\Delta(AB_1,B_2,\ldots,B_{k-1})-\langle A\underline{WGB_1G\cdots GB_{k-1}G}\rangle\label{eq:Delta_exp}\\
&+\sum_{j=1}^{k-1} \sigma m(E_\sigma, B_1,\ldots, B_j)\Delta (AE_\sigma, B_{j+1},\ldots, B_{k-1})\nonumber\\
&+\sum_{j=0}^{k-2} \sigma m(AE_\sigma, B_{j+1},\ldots, B_{k-1})\Delta ( E_\sigma, B_1,\ldots, B_j)\nonumber\\
&+ \sum_{j=0}^{k-1} \sigma \Delta(E_\sigma, B_1,\ldots,B_j)\Delta(AE_\sigma,B_{j+1},\ldots,B_{k-1}).\nonumber
\end{align}
To derive \eqref{eq:av_C_dif} from \eqref{eq:Delta_exp}, we multiply \eqref{eq:Delta_exp} by $\langle G_2\rangle -\E\langle G_2\rangle$ and then take the expectation. In this way, the first term in the rhs. of \eqref{eq:Delta_exp} becomes the second term in the rhs. of \eqref{eq:av_C_dif}, while the sums in the second, third and fourth lines of \eqref{eq:Delta_exp} become the sums in the corresponding lines of \eqref{eq:av_C_dif}. The terms in the fifth and sixth lines of \eqref{eq:av_C_dif} are absent for $n=1$, and the remaining term in the first line of \eqref{eq:av_C_dif} arises from extending the second order renormalization in the second term in the rhs. of \eqref{eq:Delta_exp} to the additional factor $\langle G_2\rangle -\E\langle G_2\rangle$. This extension is performed similarly to \eqref{eq:under_ext} by the means of \eqref{eq:EWW_id}. 

In the case $n\ge 2$ we multiply at first \eqref{eq:Delta_exp} by $\Delta\big(\overline{\bm{B}}^{[1]}\big)$ and then take covariance with $\langle G_2\rangle$. The rest follows by a straightforward calculation.
\end{proof}

\begin{proof}[Proof of Proposition \ref{prop:master_av}] Throughout the proof we fix $n, S\in\N$ with $n\le S$, in the complementary regime $\Psi_{n,S}=0$ and \eqref{eq:master_S} is trivial. Recall the definition of $\mathfrak{S}_{n,S}$ introduced above \eqref{eq:def_psi}. To prove \eqref{eq:master_S} one needs to show that properly rescaled $|\mathcal{C}^{\rm{av}}(\overline{\bm{B}})|$ is bounded by the rhs. of \eqref{eq:master_S} for any $\overline{\bm{B}}\in\mathfrak{S}_{n,S}$. We further focus on the critical case when $\overline{\bm{B}}$ contains exactly $S$ matrices, for $\overline{\bm{B}}$ with less than $S$ matrices estimates are simpler and thus are omitted.

We begin with the following simple bound:
\begin{equation}
\|(\mathcal{B}_{11}^*)^{-1}[R]\|\lesssim 1\,\,\text{for}\,\,\ R\in \{E_+,F,F^*\},
\label{eq:B*E-}
\end{equation}
which immediately follows from \eqref{eq:B*_inverse} and \cite[Eq.(3.1),(3.4)]{eigenv_decorr} (see also the justification of applicability of \cite{eigenv_decorr} around Supplementary Eq.~\eqref{eq:MDE_deformed}).
In contrast to \eqref{eq:B*E-}, $\|(\mathcal{B}_{11}^*)^{-1}[E_-]\|$ may be large, e.g.~\eqref{eq:B*_inverse} together with~\eqref{eq:M} imply that this norm scales as $\eta_1^{-1}$ for $w_1=\ii\eta_1$ with $\eta_1\ll 1$. While the exact size of $\|(\mathcal{B}_{11}^*)^{-1}[E_-]\|$ does not matter for our proof, its possible blow-up for small~$\eta_1$ makes it necessary to distinguish between the cases when all matrices in $\overline{\bm{B}}$ coincide with~$E_-$ and when at least one of them differs from $E_-$, since the inversion of $\mathcal{B}_{11}^*$ acting on one of the matrices from $\overline{\bm{B}}$ naturally appears in~\eqref{eq:def_A}.

Consider at first the case when each matrix in $\overline{\bm{B}}$ equals to $E_-$. We do not use Lemma \ref{lem:av_C_dif} in this case, instead, our proof relies on the following identities.
\begin{lemma}\label{lem:E_-}  For any $n\in\N$ and for $G=G^z(w)$ with $z\in\C$ and $w\in\C\setminus\R$, it holds that
\begin{align}
\langle (GE_-)^{2n-1}\rangle &= 0,\label{eq:E_-_identities1}\\
\langle (GE_-)^{2n}\rangle &= 2\sum_{s=0}^{n-1} (-1)^{s+1} \binom{2n-2-s}{n-1} \frac{\langle G^{s+1}\rangle}{(2w)^{2n-s-1}}. \label{eq:E_-_identities2}
\end{align}
Moreover, \eqref{eq:E_-_identities1}--\eqref{eq:E_-_identities2} remain valid after taking the deterministic approximations to all quantities appearing in these identities.
\end{lemma}

\begin{proof}[Proof of Lemma \ref{lem:E_-}] We start with the proof of~\eqref{eq:E_-_identities1}. Inverting the $2\times 2$-block matrix $W-Z-w$ defined in~\eqref{eq:def_hermitization} we get
\begin{equation}
G\!=\!\begin{pmatrix}
w \widetilde{H}_z(w)& X_z H_z(w)\\ H_z(w)X_z^*& w H_z(w)
\end{pmatrix}\!\!,\,\, \text{where}\,\, X_z\!:=X-z,\,
H_z(w)\!:=\! (X_z^*X_z-w^2)^{-1}\!\!,\, \widetilde{H}_z(w)\!:=\! (X_zX_z^*-w^2)^{-1}\!\!.
\label{eq:G_expl}
\end{equation}
Here the matrices $X_z^*X_z-w^2$ and $X_zX_z^*-w^2$ are invertible, since either $\Re w=0$, in which case $-w^2=(\Im w)^2>0$, or $\Re w\neq 0$, which implies $\Im (w^2)\neq 0$. Denote additionally $\widetilde{G}:=G^z(-w)$. From~\eqref{eq:G_expl} and the symmetry relations $\widetilde{H}_z(w)=\widetilde{H}_z(-w)$ and $H_z(w)=H_z(-w)$ we have that
\begin{equation}
GE_-=-E_-\widetilde{G}.
\label{eq:GE}
\end{equation}
Thus, the lhs. of~\eqref{eq:E_-_identities1} can be written as
\begin{equation}
\begin{split}
\langle (GE_-)^{2n-1}\rangle &= \langle (GE_-GE_-)^{n-1} GE_-\rangle = (-1)^{n-1} \langle (G\widetilde{G})^{n-1} GE_-\rangle\\
& = (-1)^{n-1} (2w)^{-n+1}\langle (G-\widetilde{G})^{n-1} GE_-\rangle=w\langle (\widetilde{H}_z(w))^{n}-(H_z(w))^{n}\rangle.
\end{split}
\label{eq:odd_k}
\end{equation}
To go from the first to the second line we employed the resolvent identity and in the last identity we used~\eqref{eq:G_expl}. Observing additionally that $\widetilde{H}_z(\eta)$ and $H_z(\eta)$ have identical spectra, we conclude that the rhs. of~\eqref{eq:odd_k} equals to zero. This finishes the proof of~\eqref{eq:E_-_identities1}.

To prove~\eqref{eq:E_-_identities2}, we compute similarly to~\eqref{eq:odd_k} that
\begin{equation}
\langle (GE_-)^{2n}\rangle = (-1)^n \langle (G\widetilde{G})^n\rangle =(-1)^n \langle G^n\widetilde{G}^n\rangle = -\frac{\partial_{w_1}^{n-1}\partial_{w_2}^{n-1} \langle G^z(w_1)G^z(-w_2)\rangle\big|_{w_1=w_2=w}}{((n-1)!)^2}.
\label{eq:even_k}
\end{equation}
Combining the resolvent identity with \eqref{eq:GE} we get
\begin{equation}
\langle G^z(w_1)G^z(-w_2)\rangle = \frac{\langle G^z(w_1)\rangle - \langle G^z(-w_2)\rangle}{w_1+w_2} = \frac{\langle G^z(w_1)\rangle}{w_1+w_2} + \frac{\langle G^z(w_2)\rangle}{w_1+w_2}.
\label{eq:even_k_res_id}
\end{equation}
Finally, we differentiate the first and the second term in the rhs. of~\eqref{eq:even_k_res_id} $n-1$ times in $w_2$ and $w_1$, respectively, and combine the result with \eqref{eq:even_k}:
\begin{equation*}
\langle (GE_-)^{2n}\rangle = \frac{(-1)^n2}{(n-1)!}\partial_y^{n-1} \frac{\langle G^z(y)\rangle}{(y+w)^n}\big|_{y=w}= \frac{2}{(n-1)!} \sum_{s=0}^{n-1}(-1)^{s+1} \frac{(2n-2-s)!}{(n-1-s)!} \cdot \frac{\langle G^{s+1}\rangle}{(2w)^{2n-1-s}}.
\end{equation*}
This finishes the proof of \eqref{eq:E_-_identities2}.

The proof of the analogues of~\eqref{eq:E_-_identities1} and~\eqref{eq:E_-_identities2} for the deterministic approximations follows by a standard meta-argument and thus is omitted.
\end{proof}

If $\overline{\bm{B}}$ has at least one block of odd size, then $\mathcal{C}^{\rm{av}}(\overline{\bm{B}})$ equals to zero by~\eqref{eq:E_-_identities1}. Assume now that each block has even size. Let $2k$ be the size of the first block. Applying~\eqref{eq:E_-_identities2} to this block we get
\begin{equation*}
\left\vert\mathcal{C}\left(\overline{\bm{B}}\right)\right\vert \lesssim \sum_{r=0}^{k-1} \eta_1^{-2k+r+1}\left\vert\mathcal{C}\left((E_+\times (r+1));\overline{\bm{B}}^{[1]}\right)\right\vert,
\end{equation*}
where $(E_+\times (r+1))$ is a block consisting of $r+1$ matrices $E_+$. Combined with \eqref{eq:def_psi}, this bound leads~to
\begin{equation*}
N^{n+1} \eta^S_*\eta_2\left\vert\mathcal{C}\left(\overline{\bm{B}}\right)\right\vert \!\lesssim\! \sum_{r=0}^{k-1} N^{n+1}\eta_*^{S-2k+r+1}\eta_2\left\vert\mathcal{C}\left((E_+\times (r+1));\overline{\bm{B}}^{[1]}\right)\right\vert\!\le\! \sum_{r=0}^{k-1} \Psi_{n, S-2k+r+1}\!\lesssim\! \Psi_{n,S-1},
\end{equation*} 
which is the desired bound on the rescaled covariance $\mathcal{C}(\overline{\bm{B}})$.

Consider now the case when at least one matrix in $\overline{\bm{B}}$ is different from $E_-$. Without loss of generality one may assume that $B^{(1)}_{k_1}\neq E_-$. We apply \eqref{eq:av_C_dif} to $\mathcal{C}\left(\overline{\bm{B}}\right)$ and note that $\|A\|\lesssim 1$ by \eqref{eq:B*E-}. For the first term in the rhs. of \eqref{eq:av_C_dif} we have the following bound
\begin{equation}
\begin{split}
&N^n\eta^S_* N\eta_2 \frac{1}{4N^2}\left\vert\langle G_1B_1^{(1)}G_1\cdots G_1B_{k_1-1}^{(1)}G_1AE_\sigma G_2^2 E_\sigma\rangle\Delta\left(\overline{\bm{B}}^{[1]}\right)\right\vert\prec \left(1+\bm{1}_{k_1>1}(N\eta_*)^{1/2}\right)\frac{\eta_*}{\gamma}.
\end{split}
\label{eq:multiG_size}
\end{equation}
The proof of \eqref{eq:multiG_size} relies on the local laws from Propositions~\ref{prop:2G_av} and \ref{prop:multiG_oneG}, and is presented in Supplementary Section~\ref{app:multiG_size}. Clearly, the rhs. of \eqref{eq:multiG_size} is bounded by the last term in the rhs. of \eqref{eq:master_S} for $n<S$. Moreover, this bound also holds in the case $n=S$, since then each block of $\overline{\bm{B}}$ consists of a single matrix, so we have $k_1=1$.

The second term in the first line of \eqref{eq:av_C_dif} is present only for $k_1\ge 2$, in which case it contains $n$ blocks and $S-1$ observables. Therefore,
\begin{equation*}
N^n\eta_*^SN\eta_2\left\vert\mathcal{C}\left(\left(AB_1^{(1)},B_2^{(1)},\ldots,B_{k_1-1}^{(1)}\right); \overline{\bm{B}}^{[1]}\right)\right\vert\le \eta_*\Psi_{n,S-1}
\end{equation*}
by the definition of $\Psi_{n,S-1}$ \eqref{eq:def_psi}. In the second line of \eqref{eq:av_C_dif} we estimate by \eqref{eq:multiM_bound} that
\begin{equation*}
\left\vert m\left(B_1^{(1)},\ldots, B_j^{(1)}, E_\sigma\right)\right\vert \lesssim \eta^{-j}_*,
\end{equation*}
and counting powers of $N$ and $\eta_*$ conclude that the second line contributes at most $\Psi_{n,S-1}$. Similarly, the contribution from the third line is bounded by $\Psi_{n,S-1}$. In the fourth line the numbers of blocks and observables equal to $n+1$ and $S+1$, respectively, so the corresponding sum is bounded by $(N\eta_*)^{-1}\Psi_{n+1,S+1}$.  The sixth line is only present for $n\ge 3$, in which case it contributes at most $\Psi_{n-2,S-2}$. Finally, we see that the last but one line in \eqref{eq:av_C_dif} is bounded by $(N\eta_*)^{-1}\Psi_{n-1,S+1}$.

In order to conclude the proof of Proposition \ref{prop:master_av} it is left to estimate the fully underlined term in the last line of \eqref{eq:av_C_dif}. Performing the cumulant expansion as in \eqref{eq:Cov_cum_exp}, we get that this term equals to
\begin{equation}
\frac{1}{N}\sum_{a,b} \sum_{\ell=3}^L\sum_{\bm{\alpha}\in \{ab,ba\}^{\ell-1}}\frac{\kappa(ab,\bm{\alpha})}{(\ell-1)!N^{\ell/2}}\E \partial_{\bm{\alpha}} \left[\left(G_1B_1^{(1)}G_1\cdots G_1B_{k_1-1}^{(1)}G_1A\right)_{ba}\!\!\Delta\left(\overline{\bm{B}}^{[1]}\right)\left(\langle G_2\rangle - \E\langle G_2\rangle\right)\right]
\label{eq:av_underline_exp}
\end{equation}
up to an error term of order $N^{-D}$, for any (large) fixed $D>0$. In contrast to the analysis of \eqref{eq:Cov_cum_exp} in the proof of Proposition \ref{prop:Cov}, we do not distinguish between the terms with $\ell\ge 5$, $\ell=4$ and $\ell=3$, but treat all these cases in the same robust way. Similarly to \eqref{eq:cum_dif_simple}, observe that for any $r\in\N$ and $\widehat{\bm{\alpha}}\in \{ab,ba\}^r$ it holds that
\begin{equation}
\left\vert \partial_{\widehat{\bm{\alpha}}} \left(G_1B_1^{(1)}G_1\cdots G_1B_{k_1-1}^{(1)}G_1A\right)_{ba}\right\vert \prec \frac{1}{\eta_1^{k_1-1}},
\label{eq:cum_iso_dif}
\end{equation}
\begin{equation}
\left\vert \partial_{\widehat{\bm{\alpha}}} \Delta \left(\bm{B}^{(j)}\right)\right\vert \prec \frac{1}{N\eta_1^{k_j}},\qquad \left\vert \partial_{\widehat{\bm{\alpha}}}\left(\langle G_2\rangle -\E\langle G_2\rangle\right)\right\vert\prec \frac{1}{N\eta_2}.
\label{eq:cum_av_dif}
\end{equation}
The proof of \eqref{eq:cum_iso_dif}--\eqref{eq:cum_av_dif} is analogous to the argument around \eqref{eq:GA_iso_dif_term} and relies of the local laws from Proposition \ref{prop:multiG_oneG} and on the deterministic estimate from Proposition \ref{prop:multiM_bound}.  
Applying \eqref{eq:cum_iso_dif}--\eqref{eq:cum_av_dif} to \eqref{eq:av_underline_exp} we obtain that for any $\ell\ge 3$ the sum of all terms of order $\ell$ in \eqref{eq:av_underline_exp} has an upper bound of order
\begin{equation*}
\frac{1}{N}N^2 \frac{1}{N^{\ell/2}}\cdot \frac{1}{\eta_1^{k_1-1}} \left(\prod_{j=2}^n \frac{1}{N\eta_1^{k_j}}\right) \frac{1}{N\eta_2}N^\xi = \frac{1}{N^n \eta_*^S N\eta_2}\cdot\frac{\eta_*}{N^{\ell/2-2}}N^\xi\le \frac{1}{N^n \eta_*^S N\eta_2}\Upsilon.
\end{equation*}
Therefore, the term in the last line of \eqref{eq:av_C_dif} multiplied by $N^n \eta_*^S N\eta_2$ is bounded by the last term in the rhs. of \eqref{eq:master_S}. This finishes the proof of Proposition \ref{prop:master_av}.
\end{proof}

\appendix

\section{Additional analysis of the matrix Dyson equation} 

\subsection{2-body stability analysis: Proof of Proposition \ref{prop:stab_bound}}\label{app:stab_bound}

We begin by introducing the notation used throughout this section. Unless stated otherwise, we work under the assumptions $z_1,z_2\in (1-\delta)\mathbf{D}$ and $w_1,w_2\in \C\setminus\R$. We denote $E_j:=\Re w_j$, $\eta_j:=|\Im w_j|$ and $\eta_*:=\eta_1\wedge\eta_2\wedge 1$, $j=1,2$. Whenever a statement involves index $j$, it is understood to hold for $j=1,2$. We further recall the notations introduced in \eqref{eq:M}, \eqref{eq:m} and set $M_j:=M^{z_j}(w_j)$, $m_j:=m^{z_j}(w_j)$ and $u_j:=u^{z_j}(w_j)$. Finally, we denote
\begin{equation}
\begin{split}
\beta_{12,\pm} =\beta_{12,\pm}(w_1,w_2) :&=1-\Re \left[z_1\overline{z}_2\right] u_1u_2\pm \sqrt{m_1^2m_2^2 - \left(\Im \left[z_1\overline{z}_2\right]\right)^2u_1^2u_2^2},\\
\beta_{12,*}=\beta_{12,*}(w_1,w_2):&=\min\{|\beta_{12,+}|,|\beta_{12,-}|\},
\end{split}
\label{eq:def_beta_pm_main}
\end{equation}
It is known from \cite[Appendix B]{Cipolloni_meso} that $\beta_{12,\pm}$ are the eigenvalues of the two-body stability operator $\mathcal{B}_{12}$ defined in \eqref{eq:def_B12}, while the remaining two eigenvalues equal to 1, see e.g. \cite[Appendix B]{Cipolloni_meso}.

The crucial step in the proof of Proposition \ref{prop:stab_bound} is to reduce the analysis of $\widehat{\beta}_{12}$ to that of $\beta_{12,*}$ defined in \eqref{eq:def_beta_pm_main}. The latter quantity is much more accessible due its definition by the explicit formula. An elementary calculation based on the explicit inversion of $\mathcal{B}_{12}$ gives that
\begin{equation}
\widehat{\beta}_{12}(w_1,w_2)\sim\min\left\lbrace\beta_{12,*}\big(w_1^{(*)}, w_2^{(*)}\big)\wedge 1\right\rbrace,
\label{eq:beta_hat_to_beta_app}
\end{equation}
uniformly in $z_j\in (1-\delta)\mathbf{D}$ and $w_j\in \C\setminus\R$, for $j=1,2$. In \eqref{eq:beta_hat_to_beta_app} the minimum is taken over all four choices of stars. For more details on the proof of \eqref{eq:beta_hat_to_beta_app} see Supplementary Section~\ref{app:beta_hat}.

First we derive \eqref{eq:eta_stab_bound}, which in the view of \eqref{eq:beta_hat_to_beta_app} is equivalent to
\begin{equation}
\beta_{12,*}\gtrsim \eta_*,\quad \forall z_j\in (1-\delta)\mathbf{D},\, w_j\in \C\setminus\R.
\label{eq:eta_bound_app}
\end{equation}
We have from \cite[Lemma 6.1]{macroCLT_real} that
\begin{equation}
\beta_{12,*}(w_1,w_2) \gtrsim \left\vert|E_1|-|E_2|\right\vert^2 +|z_1-z_2|^2 +\eta_1+\eta_2
\label{eq:stab_quad}
\end{equation} 
uniformly in $|z_j|, |w_j|\lesssim 1$. This immediately implies \eqref{eq:eta_bound_app} in the regime $|w_j|\lesssim 1$. To cover the complementary regime where $\max\{|w_1|,|w_2|\}$ is large, we observe that 
\begin{equation}
\beta_{*,12}(w_1,w_2)\ge 1-\lVert M^{z_1}(w_1)\rVert\lVert M^{z_2}(w_2)\rVert,
\label{eq:beta_star_trivial}
\end{equation}
for any $z_j\in \C$, $w_j\in \C\setminus\R$, as it trivially follows from the eigenvalue equation for $\mathcal{B}_{12}$ and the bound $\lVert\mathcal{S}[R]\rVert\le \lVert R\rVert$, for all $R\in\mathrm{span}\{E_\pm, F^{(*)}\}$.  The MDE \eqref{eq:MDE} implies that $\lVert M^{z}(w)\rVert\lesssim |w|^{-1}$ uniformly in $|z|\le 1$ and $|w|\ge C$ for some (large) implicit constant $C>0$. Combined with the trivial bound $\|M^z(w)\|\lesssim 1$ from \cite[Eq.(3.5)]{macroCLT_real} and \eqref{eq:beta_star_trivial} this yields the lower bound of order one on the lhs. of \eqref{eq:eta_bound_app} in the regime $\max\{|w_1|,|w_2|\}\ge C$ and thus finishes the proof of \eqref{eq:eta_stab_bound}.

Now we proceed to the proof of \eqref{eq:stab_LT}, and start with a small simplification. Due to the symmetry relation 
\begin{equation}
M^z(w)E_-=-E_-M^z(-w)
\label{eq:M_symmetry}
\end{equation}
it holds that
\begin{equation*}
\mathcal{B}_{12}(w_1,-w_2)[RE_-] = \mathcal{B}_{12}(w_1,w_2)[R]E_-,
\end{equation*}
so $\widehat{\beta}_{12}(w_1,w_2) = \widehat{\beta}_{12}(w_1,-w_2)$ for any $z_j\in \C$ and $w_j\in \C\setminus\R$. Similarly we have $\widehat{\beta}_{12}(w_1,w_2) = \widehat{\beta}_{12}(-w_1,w_2)$. Therefore, it is sufficient to prove \eqref{eq:stab_LT} in the case when $E_1, E_2\ge 0$.

We deduce \eqref{eq:stab_LT} from the following estimates on $\beta_{12,*}$.

\begin{lemma}\label{lem:beta_estimates} Fix (small) $\delta, \kappa>0$. Uniformly in $z_j\in (1-\delta)\mathbf{D}$ and $w_j\in \C\setminus\R$ with $E_j:=\Re w_j\in \mathbf{B}_\kappa^{z_j}$, $E_j\ge 0$, $\eta_j:=|\Im w_j|\in (0,1]$ for $j=1,2$, it holds that
\begin{align}
\beta_{12,*}(w_1,w_2)\sim \left\vert E_1-E_2 + \frac{\langle (Z_1-Z_2)\Im M_1\rangle}{\langle \Im M_1\rangle}\right\vert + |E_1-E_2|^2 + |z_1-z_2|^2 +\eta_1+\eta_2,
\label{eq:beta_opposite}
\end{align}
if $w_1,w_2$ are in the opposite complex half-planes, and
\begin{equation}
\beta_{12,*}(w_1,w_2)\sim E_1+E_2 + |z_1-z_2|^2 +\eta_1+\eta_2,
\label{eq:beta_same}
\end{equation}
if $w_1,w_2$ are in the same half-plane.
\end{lemma}

We postpone the proof of Lemma \ref{lem:beta_estimates} to the end of this section, and proceed to the proof of \eqref{eq:stab_LT} relying on \eqref{eq:beta_opposite}--\eqref{eq:beta_same}.

\begin{proof}[Proof of \eqref{eq:stab_LT}] Using the explicit form of $M_1$ from \eqref{eq:M} we get
\begin{equation}
E_1-E_2 + \frac{\langle (Z_1-Z_2)\Im M_1\rangle}{\langle \Im M_1\rangle} = E_1-E_2 - \frac{\Im u_1\Re[\overline{z}_1(z_1-z_2)]}{\Im m_1} = \mathcal{O}(E_1+E_2),
\label{eq:LT_same}
\end{equation}
where we additionally used that $|\Im m_1|\sim 1$ due to the conditions $E_1\in\mathbf{B}_\kappa^{z_1}$ and $\eta_1\lesssim 1$, and computed
\begin{equation}
\Im u_1=\Im u^{z_1}(w_1) = \Im [u^{z_1}(w_1)-u^{z_1}(\ii\Im w_1)] = \mathcal{O}(E_1), 
\end{equation} 
since $u^z(w)$ is real for $w$ on the imaginary axis. It follows from \eqref{eq:beta_opposite}, \eqref{eq:beta_same} and \eqref{eq:LT_same} that
\begin{equation}
\beta_{12,*}(w_1,w_2)\gtrsim \beta_{12,*}(w_1,\overline{w}_2) + \beta_{12,*}(\overline{w}_1,w_2),
\end{equation}
for all $w_1, w_2$ in the same half-plane satisfying assumptions of Lemma \ref{lem:beta_estimates}. Combined with \eqref{eq:beta_hat_to_beta_app} and \eqref{eq:beta_opposite}, this immediately finishes the proof of \eqref{eq:stab_LT}. 
\end{proof}

\begin{proof}[Proof of Lemma \ref{lem:beta_estimates}] Instead of working with $\beta_{12,*}$ directly, we use that
\begin{equation}
\beta_{12,*}\sim \left\vert (1-\langle M_1M_2\rangle)(1+\langle M_1E_-M_2E_-\rangle) + \langle M_1E_-M_2\rangle \langle M_2E_-M_1\rangle \right\vert,
\label{eq:beta_to_det}
\end{equation}
and further analyze the rhs. of \eqref{eq:beta_to_det}. For the proof of \eqref{eq:beta_to_det} see Supplementary equations \eqref{eq:B_inverse_aux} and~\eqref{eq:detP}.

\medskip

\noindent\underline{1st case ($\Im w_1\Im w_2<0$): proof of \eqref{eq:beta_opposite}.} We expand $M_2$ around $M_1^*$ up to the first order:
\begin{equation}
M_2 = M_1^* + M_1^*DM_1^* + \frac{\left\langle D(M_1^*)^2\right\rangle}{1-\langle (M_1^*)^2\rangle} (M_1^*)^2 + \mathcal{O}(|\overline{w}_1-w_2|^2 + |z_1-z_2|^2),\quad D:=w_2-\overline{w}_1 + Z_2-Z_1,
\label{eq:M_expansion1}
\end{equation}
where the error term is a $2\times 2$ matrix bounded in operator norm by $|\overline{w}_1-w_2|^2 + |z_1-z_2|^2$, up to a multiplicative constant. A simple calculation based on \eqref{eq:M} shows that $\langle M_1E_-M_1^*\rangle=\langle M_1^*E_-M_1\rangle=0$, which in combination with \eqref{eq:M_expansion1} gives
\begin{equation}
\langle M_1E_-M_2\rangle \langle M_2E_-M_1\rangle = \mathcal{O}\left(|\overline{w}_1-w_2|^2 + |z_1-z_2|^2\right).
\label{eq:offdiag_terms}
\end{equation}
Denote $m_1:=m^{z_1}(w_1)$ and $u_1:=u^{z_1}(w_1)$. We further compute from \eqref{eq:M} and \eqref{eq:M_expansion1} that
\begin{equation}
\begin{split}
1+ \langle M_1E_-M_2E_-\rangle &= 1+ |m_1|^2 - |z_1|^2 |u_1|^2 + \mathcal{O}\left(|\overline{w}_1-w_2| + |z_1-z_2|\right)\\
1-\langle M_1M_2\rangle &= \frac{\eta_1}{\eta_1+|\Im m_1|} + \langle DX\rangle + \mathcal{O}\left(|\overline{w}_1-w_2|^2 + |z_1-z_2|^2\right),
\end{split}
\label{eq:same_side_factors_exp}
\end{equation}
where $X\in \C^{2\times 2}$ is given by
\begin{equation}
X:=M_1^*M_1M_1^* + \frac{\langle M_1(M_1^*)^2\rangle}{1-\langle (M_1^*)^2\rangle}(M_1^*)^2.
\end{equation}
In the second line of \eqref{eq:same_side_factors_exp} we additionally used that
\begin{equation}
M_1M_1^* = \frac{\Im M_1}{\Im w_1 + \Im m_1}
\label{eq:ImM}
\end{equation}
to compute $1-\langle M_1M_1^*\rangle$. We further use \eqref{eq:ImM} to compute $X$:
\begin{equation}
\begin{split}
X&=\frac{1}{2\ii (\Im w_1+\Im m_1)} \left( M_1^*(M_1-M_1^*) + \frac{\langle M_1^*(M_1-M_1^*)\rangle}{1-\langle (M_1^*)^2\rangle}(M_1^*)^2 \right)\\
&=\frac{1}{2\ii (\Im w_1+\Im m_1)} \left( M_1M_1^* - \frac{1-\langle M_1M_1^*\rangle}{1-\langle (M^*)^2\rangle}(M_1^*)^2\right)= \frac{\Im M_1}{2\ii (\Im m_1)^2} + \mathcal{O}(\eta_1).
\end{split}
\label{eq:X_calculation}
\end{equation}
Combining \eqref{eq:beta_to_det}, \eqref{eq:offdiag_terms}, \eqref{eq:same_side_factors_exp}, \eqref{eq:X_calculation}, and additionally using that
\begin{equation*}
1+|m_1|^2 -|z_1|^2 |u_1|^2 \sim 1,
\end{equation*}
as it follows from the bound $|u_1|<1$ from \cite[Eq. (3.5)]{macroCLT_real}, we arrive to 
\begin{equation}
\begin{split}
\beta_{12,*}\sim & \left\vert \frac{\eta_1}{\eta_1+|\Im m_1|} + \frac{\langle D\Im M_1\rangle}{2\ii (\Im m_1)^2} + \mathcal{O}\left(|\overline{w}_1-w_2|^2 + |z_1-z_2|^2 + \eta_1^2\right)\right\vert\\
=& \left\vert \left(E_1-E_2 + \frac{\langle (Z_1-Z_2)\Im M_1\rangle}{\langle \Im M_1\rangle}\right) + \mathcal{O}\left(|E_1-E_2|^2 + |z_1-z_2|^2 + \eta_1+\eta_2\right)\right\vert.
\end{split}
\label{eq:opposite_side_beta_exp}
\end{equation}
To go from the first to the second line in \eqref{eq:opposite_side_beta_exp} we used that $|\Im m_1|\sim 1$ since $E_1\in \mathbf{B}^{z_1}_\kappa$ and $\eta_1\lesssim 1$, included $\eta_1/(\eta_1+|\Im m_1|)$ into the error term and recalled the definition of $D$ from \eqref{eq:M_expansion1}. Finally, to conclude \eqref{eq:beta_opposite} from \eqref{eq:opposite_side_beta_exp}, we make the following trivial observation. Let $a,b,c\in\R$ satisfy conditions $|c|\sim |a+b|$ and $|c|\gtrsim |b|$. Then it holds that $|c|\sim |a|+|b|$. The proof of this fact proceeds by distinguishing between the cases $|a|\ge L|b|$ and $|a|<L|b|$, for sufficiently large fixed $L>0$. Applying this observation to $c:=\beta_{12,*}$ and $a,b$ given by the first and second terms in the second line of \eqref{eq:opposite_side_beta_exp}, respectively, and recalling the lower bound \eqref{eq:stab_quad}, we finish the proof of \eqref{eq:beta_opposite}.

\medskip

\noindent\underline{2nd case ($\Im w_1\Im w_2>0$): proof of \eqref{eq:beta_same}.} We argue similarly to the previous case, and analyze the rhs. of \eqref{eq:beta_to_det}. Expanding $M_2$ around $M_1$ we get analogously to \eqref{eq:offdiag_terms} that
\begin{equation}
\langle M_1E_-M_2\rangle \langle M_2E_-M_1\rangle = \mathcal{O}\left(|w_1-w_2|^2 + |z_1-z_2|^2\right).
\label{eq:offdiag_terms_same}
\end{equation}
We further have from \cite[Eq.(3.4)]{eigenv_decorr} that
\begin{equation}
|1-\langle M_1M_2\rangle|\sim 1,
\label{eq:same_stab}
\end{equation}
for the applicability of \cite{eigenv_decorr} see Supplementary equation \eqref{eq:MDE_deformed} and the discussion around. Next we compute $1+\langle M_1E_-M_2E_-\rangle$ using \eqref{eq:M} both for $M_1$ and $M_2$:
\begin{equation}
\begin{split}
&1+\langle M_1E_-M_2E_-\rangle = 1 +m_1m_2-\Re [z_1\overline{z}_2]u_1u_2\\
&\quad = 1+ \frac{m_1^2 + m_2^2 - (m_1-m_2)^2}{2} - \frac{|z_1|^2 +|z_2|^2 - |z_1-z_2|^2}{2}\frac{u_1^2 +u_2^2-(u_1-u_2)^2}{2}\\
&\quad=1+\frac{m_1^2+m_2^2}{2} - \frac{|z_1|^2+|z_2|^2}{2}\frac{u_1^2+u_2^2}{2} +\mathcal{O}(|z_1-z_2|^2 +|w_1-w_2|^2)\\
&\quad=\frac{1+m_1^2 -|z_1|^2u_1^2}{2}+\frac{1+m_2^2 -|z_2|^2u_2^2}{2} + \frac{(|z_1|^2-|z_2|^2)(u_1^2-u_2^2)}{4} + \mathcal{O}(|z_1-z_2|^2 +|w_1-w_2|^2)\\
&\quad= \frac{w_1}{2(w_1+m_1)}+\frac{w_2}{2(w_2+m_2)} + \mathcal{O}(|z_1-z_2|^2 +|w_1-w_2|^2).
\end{split}
\label{eq:MEM_new}
\end{equation}
Here to go from the second to the third and from the fourth to the fifth line we used that
\begin{equation*}
|u_1-u_2|=\mathcal{O}(|z_1-z_2| +|w_1-w_2|).
\end{equation*}
Additionally, in the last step we used \eqref{eq:m} twice: for $m_1$ and $m_2$. 

Combining \eqref{eq:beta_to_det}, \eqref{eq:offdiag_terms_same}, \eqref{eq:same_stab} and \eqref{eq:MEM_new} we get
\begin{equation}
\beta_{12,*}\sim \left\vert \frac{E_1}{2(w_1+m_1)}+\frac{E_2}{2(w_2+m_2)} + \mathcal{O}(|z_1-z_2|^2 +|E_1-E_2|^2+\eta_1+\eta_2)\right\vert.
\label{eq:same_side_beta_exp}
\end{equation}
We repeat the argument presented below \eqref{eq:opposite_side_beta_exp} using \eqref{eq:stab_quad} and \eqref{eq:same_side_beta_exp} as an input, and obtain
\begin{equation}
\beta_{12,*}\sim \left\vert \frac{E_1}{2(w_1+m_1)}+\frac{E_2}{2(w_2+m_2)}\right\vert + |z_1-z_2|^2 +|E_1-E_2|^2+\eta_1+\eta_2.
\label{eq:same_side_beta_fin}
\end{equation}
Finally, we notice that the first term in the rhs. of \eqref{eq:same_side_beta_fin} is of order $E_1+E_2$, since the imaginary parts of $w_1+m_1$, $w_2+m_2$ share the same sign and have an absolute value of order 1. This finishes the proof of~\eqref{eq:beta_same}.
\end{proof}

\subsection{Bound on the propagator: proof of Lemma \ref{lem:propag_bound2}}\label{app:propagator}

In this section we incorporate the notations and conventions introduced in the beginning of Appendix~\ref{app:stab_bound}. We also frequently use \eqref{eq:beta_hat_to_beta_app} without mentioning this further. The time-dependent version of $M_j$ is denoted by $M_{j,r}:=M^{z_{j,r}}(w_{j,r})$, for $r\in [0,T]$. The time-dependent versions of the related quantities such as $m_j$ and $u_j$ are denoted similarly by replacing $z_j,w_j$ with $z_{j,r}, w_{j,r}$. Denote for short $\beta_{\sigma,r}:=\beta_{12,\sigma,r}$ and $\beta_{*,r}:=\beta_{12,*,r}$ for $\sigma\in\{\pm\}$ and $r\in [0,T]$. Finally, we set
\begin{equation}
a_{12}^\sigma :=\sigma \langle M_{12}^{E_\sigma}E_\sigma\rangle,\,\sigma\in\{\pm\},\quad d_{12}:=\langle M_{12}^{E_+}E_-\rangle = -\langle M_{12}^{E_-}E_+\rangle.
\end{equation}
The last identity is established for $w_i, w_j$ on the imaginary axis in \cite[Eq.(5.26)]{univ_extr}. Since the functions $\langle M_{ij}^{E_+}E_-\rangle$ and $-\langle M_{ij}^{E_-}E_+\rangle$ are analytic in $w_i,w_j$ in the upper and lower complex half-planes, this identity holds for all $w_i,w_j\in\C\setminus\R$.

Recall that $N^{-b}\lesssim T\lesssim 1$, $|z_{j,T}|\le 1-\delta$, $\Re w_{j,T}\in \mathbf{B}_\kappa^{z_j}$ and $|\Im w_{j,T}|\ge N^{-1+\epsilon}$ by the set-up of Lemma~\ref{lem:propag_bound2}, for some fixed $\delta,\kappa,\epsilon>0$. We assume that $z_{j,r}, w_{j,r}$ satisfy the same conditions for all $t\in [0,T]$ with $\delta$ and $\kappa$ decreased to $\delta/2$ and $\kappa/2$, respectively. This can be achieved by choosing $T$ sufficiently small, for more details see the set-up of Supplementary Lemma~\ref{lem:propag_bound}. 

An explicit calculation discussed around \eqref{eq:Y_syst_intro} gives that
\begin{equation}
\mathcal{A}_{[2],t}=
\begin{pmatrix}
2a_{12}^+&d_{12}&-d_{12}&0\\
d_{12}&a_{12}^++a_{12}^-&0&-d_{12}\\
-d_{12}&0&a_{12}^-+a_{12}^+&d_{12}\\
0&-d_{12}&d_{12}&2a_{12}^-
\end{pmatrix}.
\end{equation}
Estimating the off-diagonal part of $\Re\mathcal{A}_{[2],t}$ simply by operator norm and recalling the definition of $f_{[2],r}$ given in Lemma~\ref{lem:propag_bound2}, we get 
\begin{equation}
f_{[2],r}\le 2 \max\{a_{12,r}^+,a_{12,r}^-,0\} + d_{12,r},\quad r\in [0,T].
\label{eq:f2_bound}
\end{equation}
Thus, to verify \eqref{eq:propag_bound2}, it suffices to show that 
\begin{equation}
\int_s^t \max\{\Re a_{12,r}^+,\Re a_{12,r}^-,0\}\dif r=\int_s^t \max\{\Re a_{12,r}^+,\Re a_{12,r}^-\}\dif r + \mathcal{O}(1) = \log  \frac{\beta_{12,*,s}}{\beta_{12,*,t}} + \mathcal{O}(1),\label{eq:a_bound_main}
\end{equation}
\begin{equation}
\int_s^t |\Re d_{12,r}|\dif r=\mathcal{O}(1), \label{eq:d_bound_main}
\end{equation}
for any $0\le s\le t\le T$. We note that both identities in \eqref{eq:a_bound_main} are non-trivial and require verification. Further in the proof of \eqref{eq:a_bound_main} and \eqref{eq:d_bound_main} the time parameters $s,r,t\in [0,T]$ are ordered as $s\le r\le t$, unless stated otherwise.

We start with the preliminary analysis of $\beta_{\pm,r}$. First we observe that there exists $\sigma \in \{\pm\}$ such that 
\begin{equation}
|\beta_{\sigma, r}|\sim 1,\quad \forall r\in [0,T].
\label{eq:one_beta}
\end{equation}
Indeed, from \eqref{eq:char_flow_discussion} we have that
\begin{equation}
z_{j,r}=\ee^{-r/2}z_{j,0},\quad m_{j,r}=\ee^{r/2} m_{j,0},\quad u_{j,r}=\ee^r u_{j,0},\qquad \forall \,r\in [0,T].
\label{eq:time_scaling_main}
\end{equation}
Therefore, the evolution in time of the square root in the rhs. of \eqref{eq:def_beta_pm_main} is given by
\begin{equation}
\Re\sqrt{m_{1,r}^2m_{2,r}^2 - \Im \left[z_{1,r}\overline{z}_{2,r}\right]u_{1,r}^2u_{2,r}^2} = \ee^r \Re\sqrt{m_{1,0}^2m_{2,0}^2 - \Im \left[z_{1,0}\overline{z}_{2,0}\right]u_{1,0}^2u_{2,0}^2}. 
\label{eq:sqr_rescale_main}
\end{equation}
In particular, this quantity preserves the sign for all $r\in[0,T]$. By \cite[Eq. (3.5)]{macroCLT_real} it holds that $|u_{j,r}|<1$, so
\begin{equation}
\Re\left[ 1-\Re [z_{1,r}\overline{z}_{2,r}]u_{1,r}u_{2,r}\right] \ge 1-|z_{1,r}z_{2,r}|\ge \delta,
\label{eq:Re_beta_main}
\end{equation}
where in the last step we estimated $|z_{j,r}|\le 1-\delta/2$. Combining \eqref{eq:def_beta_pm_main} with \eqref{eq:sqr_rescale_main} and \eqref{eq:Re_beta_main}, we finish the verification of~\eqref{eq:one_beta}.

Without loss of generality we may assume that
\begin{equation}
|\beta_{+,r}|\sim 1,\quad |\beta_{-,r}|\sim \beta_{*,r},\qquad \forall r\in [0,T].
\label{eq:beta_relations}
\end{equation}
We further observe that \eqref{eq:time_scaling_main} implies
\begin{equation}
\beta_{12,\sigma,r} = 1-\ee^{r}(1-\beta_{12,\sigma,0}),\quad \sigma\in\{\pm\},\,\, r\in [0,T],
\label{eq:beta_time_dep_main}
\end{equation}
so the trajectory of $\beta_{\sigma,r}$ is a part of a straight line. In particular, $\beta_{\sigma,r}$ does not wind around the origin, so we have
\begin{equation}
\log\frac{\beta_{\sigma,s}}{\beta_{\sigma,t}} = \log\left\vert \frac{\beta_{\sigma,s}}{\beta_{\sigma,t}}\right\vert +\mathcal{O}(1),\quad \forall\, s,t\in [0,T], 
\label{eq:Im_log}
\end{equation} 
which we will use in the remaining part of the proof.

Next, we proceed to the proof of \eqref{eq:a_bound_main}. We borrow the following identity for $a_{12}^\sigma$, $\sigma\in\{\pm\}$, from \cite[Eq.(A.13)-(A.14)]{univ_extr}:
\begin{equation}
a_{12}^\sigma = \frac{\sigma m_1m_2 -|z_1z_2|^2u_1^2u_2^2 +m_1^2m_2^2 +\Re [z_1\overline{z}_2]u_1u_2}{1+|z_1z_2|^2u_1^2u_2^2 -m_1^2m_2^2 -2\Re [z_1\overline{z}_2]u_1u_2},
\label{eq:a_explicit}
\end{equation}
Using the definition of $\beta_\pm$ given in \eqref{eq:def_beta_pm_main} we compute
\begin{equation}
\beta_+\beta_- = 1+|z_1z_2|^2u_1^2u_2^2 -m_1^2m_2^2 -2\Re [z_1\overline{z}_2]u_1u_2.
\label{eq:beta_prod}
\end{equation}
Combining \eqref{eq:time_scaling_main}, \eqref{eq:a_explicit} and \eqref{eq:beta_prod} we obtain that
\begin{equation}
a_{12,r}^\sigma= -\frac{1}{2}\partial_r \log \beta_{+,r}\beta_{-,r} + \sigma \frac{m_{1,r}m_{2,r}}{\beta_{+,r}\beta_{-,r}}.
\label{eq:a_deriv_beta}
\end{equation}
Therefore, \eqref{eq:beta_relations} implies that a half of the rhs. of \eqref{eq:a_bound_main} comes from the first term in the rhs. of \eqref{eq:a_deriv_beta}, while the second half should come from the second term, i.e. the second identity in \eqref{eq:a_bound_main} is equivalent to
\begin{equation}
\int_s^t \left\vert \Re \frac{m_{1,r}m_{2,r}}{\beta_{+,r}\beta_{-,r}}\right\vert\dif r = \frac{1}{2}\log \frac{\beta_{*,s}}{\beta_{*,t}} + \mathcal{O}(1).
\label{eq:m_over_beta_bound1}
\end{equation}

Let us show that the real part of $m_{1,r}m_{2,r}/(\beta_{+,r}\beta_{-,r})$ changes the sign at most 3 times\footnote{The exact number of sign changes does not matter, but it is important that this happens only finitely many times independently of the trajectory of the characteristic flow} on $[0,T]$. To see this, we use \eqref{eq:beta_prod}, \eqref{eq:time_scaling_main} and obtain
\begin{equation}
\Re \frac{m_{1,r}m_{2,r}}{\beta_{+,r}\beta_{-,r}} =  \frac{\Re \left[\ee^r \overline{m}_{1,0}\overline{m}_{2,0}\left( 1+\ee^{2r}|z_{1,0}z_{2,0}|^2 u_{1,0}^2u_{2,0}^2 - \ee^{2r} m_{1,0}^2m_{2,0}^2 -\ee^r 2\Re [z_{1,0}\overline{z}_{2,0}]u_{1,0}u_{2,0}\right)\right]}{|\beta_{+,r}\beta_{-,r}|^2}. 
\label{eq:Re_m_over_beta}
\end{equation}
The numerator in the rhs. of \eqref{eq:Re_m_over_beta} is a polynomial with real coefficients in variable $\ee^r$ of degree at most 3. In particular, it either vanishes identically or has at most 3 zeros. Therefore, the lhs. of \eqref{eq:Re_m_over_beta} changes the sign at most 3 times for $r\in [0,T]$. This allows us to split the integration in the lhs. of \eqref{eq:m_over_beta_bound1} into at most 4 intervals such that on each of them the sign of the lhs. of \eqref{eq:Re_m_over_beta} is preserved. Since the first term in the rhs. of \eqref{eq:m_over_beta_bound1} is an additive function of the time interval, i.e.
\begin{equation*}
\log \frac{\beta_{*,s}}{\beta_{*,t}}=\log\frac{\beta_{*,s}}{\beta_{*,r}} +\log \frac{\beta_{*,r}}{\beta_{*,t}}
\end{equation*}
for any $s,r,t\in [0,T]$, $s\le r\le t$, to verify the second identity in \eqref{eq:a_bound_main} it is sufficient to show that 
\begin{equation}
\left\vert\Re\int_{s_1}^{t_1} \frac{m_{1,r}m_{2,r}}{\beta_{+,r}\beta_{-,r}}\dif r\right\vert = \frac{1}{2}\log \frac{\beta_{*,s_1}}{\beta_{*,t_1}} + \mathcal{O}(1)
\label{eq:m_over_beta_bound2}
\end{equation}
holds for any $s_1,t_1\in [0,T]$ such that $s_1\le t_1$ and the integrand in the lhs. of \eqref{eq:m_over_beta_bound2} preserves the sign on~$[s_1,t_1]$. Actually, the last condition on $s_1, t_1$ is not needed for the proof of \eqref{eq:m_over_beta_bound2}, so we will only assume that $0\le s_1\le t_1\le T$.

Now we prove \eqref{eq:m_over_beta_bound2} and start with a small simplification. Let $c_0>0$ be a constant which will be taken sufficiently small later. If $z_{1,0}$, $z_{2,0}$ satisfy the inequality $|z_{1,0}-z_{2,0}|\ge c_0$, then
\begin{equation*}
\beta_{*,r}\gtrsim |z_{1,r}-z_{2,r}|^2=\ee^{-r}|z_{1,0}-z_{2,0}|^2\ge \ee^{-T}c_0^2,
\end{equation*}
where we used \eqref{eq:stab_LT} in the first step and \eqref{eq:time_scaling_main} in the second. Thus, $|\beta_{+,r}|$ and $|\beta_{-,r}|$ are of order 1 for all $r\in [s_1,t_1]$ and the lhs. of \eqref{eq:m_over_beta_bound2} has an upper bound of order 1, i.e. \eqref{eq:m_over_beta_bound2} holds. Therefore, we may further assume that $|z_{1,0}-z_{2,0}|<c_0$. Using \eqref{eq:beta_time_dep_main} along with \eqref{eq:time_scaling_main}, we explicitly compute the integral in the lhs. of \eqref{eq:m_over_beta_bound2}:
\begin{equation}
\int_{s_1}^{t_1} \frac{m_{1,r}m_{2,r}}{\beta_{-,r}\beta_{+,r}}\dif r = \frac{m_{1,0}m_{2,0}}{\beta_{+,0}-\beta_{-,0}}\left(\log \frac{\beta_{-,s_1}}{\beta_{-,t_1}}-\log \frac{\beta_{+,s_1}}{\beta_{+,t_1}}\right).
\label{eq:int_m_over_beta_expl}
\end{equation}
We further compute the first factor in the rhs. of \eqref{eq:int_m_over_beta_expl} by expressing $\beta_{+,0}-\beta_{-,0}$ from \eqref{eq:def_beta_pm_main}:
\begin{equation}
\frac{m_{1,0}m_{2,0}}{\beta_{+,0}-\beta_{-,0}} = \pm \frac{1}{2} \left( 1- \frac{\Im [z_{1,0}\overline{z}_{2,0}]u_{1,0}u_{2,0}}{m_{1,0}m_{2,0}}\right)^{-1/2}.
\label{eq:m_over_beta_dif}
\end{equation}
The sign in the rhs. of \eqref{eq:m_over_beta_dif} appears due to the fact that $\sqrt{m_{1,0}^2m_{2,0}^2}=\pm m_{1,0}m_{2,0}$ depending on the phase of $m_{1,0}m_{2,0}$. We further estimate
\begin{equation}
\begin{split}
&\left\vert\frac{\Im [z_{1,0}\overline{z}_{2,0}]u_{1,0}u_{2,0}}{m_{1,0}m_{2,0}}\right\vert = \left\vert\frac{\Im [(z_{1,T}-z_{2,T})\overline{z}_{2,T}]}{(w_{1,T}+m_{1,T})(w_{1,T}+m_{1,T})}\right\vert\\
&\quad\le \frac{|z_{1,T}-z_{2,T}|}{\left(|\Im w_{1,T}| + |\kappa+\mathcal{O}(|\Im w_{1,T}|)|\right)\left(|\Im w_{2,T}| + |\kappa+\mathcal{O}(|\Im w_{2,T}|)|\right)}\lesssim \frac{|z_{1,T}-z_{2,T}|}{\kappa^2}.
\end{split}
\label{eq:dif_z_term}
\end{equation}
In the first line of \eqref{eq:dif_z_term} we used \eqref{eq:time_scaling_main}, while to go from the first to the second line we estimated each factor in the denominator from below by the corresponding imaginary part and used that $\Re w_{j,T}\in\mathbf{B}_\kappa^{z_{j,T}}$. In particular, taking $c_0$ sufficiently small and recalling that $|z_{1,T}-z_{2,T}|\le |z_{1,0}-z_{2,0}|\le c_0$, we make the lhs. of \eqref{eq:dif_z_term} smaller than $1/2$. Thus the rhs. of \eqref{eq:m_over_beta_dif} can be linearized, and we get
\begin{equation}
\frac{m_{1,0}m_{2,0}}{\beta_{+,0}-\beta_{-,0}} = \pm \frac{1}{2} + \mathcal{O}(|z_{1,T}-z_{2,T}|),
\end{equation}
where we used \eqref{eq:m_over_beta_dif} and \eqref{eq:dif_z_term}. Together with \eqref{eq:int_m_over_beta_expl}, \eqref{eq:Im_log} and \eqref{eq:beta_relations} this implies
\begin{equation}
\left\vert \Re\int_{s_1}^{t_1} \frac{m_{1,r}m_{2,r}}{\beta_{-,r}\beta_{+,r}}\dif r\right\vert = \left(\frac{1}{2} + \mathcal{O}(|z_{1,T}-z_{2,T}|)\right) \left\vert \log \frac{\beta_{*,s_1}}{\beta_{*,t_1}} + \mathcal{O}(1) \right\vert.
\label{eq:a_bound_prefin}
\end{equation}
Finally, notice that from \eqref{eq:stab_LT} we have
\begin{equation*}
|z_{1,T}-z_{2,T}|\left\vert \log \frac{\beta_{*,s_1}}{\beta_{*,t_1}}\right\vert \lesssim \left\vert(z_{1,T}-z_{2,T})\log |z_{1,T}-z_{2,T}|\right\vert \lesssim 1.
\end{equation*}
Thus, \eqref{eq:m_over_beta_bound2} holds with the logarithm replaced by its absolute value in the right hand side. To deal with the case when this logarithm is negative, we observe that
\begin{equation}
\beta_{*,s}\sim \beta_{*,t}+|t-s|,\quad \forall \, 0\le s\le t\le T,
\label{eq:gamma_monot_main}
\end{equation}
for the proof see Supplementary Lemma~\ref{lem:beta_hat}. This implies an upper bound of order one on the absolute value of the logarithm, so it can be absorbed into the error term. This finishes the proof of \eqref{eq:m_over_beta_bound2} and thereby of the second identity in~\eqref{eq:a_bound_main}.

To complete the proof of \eqref{eq:a_bound_main}, it remains to show that
\begin{equation}
\int_s^t \max\{\Re a_{12,r}^+,\Re a_{12,r}^-,0\}\dif r=\log  \frac{\beta_{*,s}}{\beta_{*,t}} + \mathcal{O}(1).\label{eq:a_plus_bound}
\end{equation}
Since we have already established the second identity in \eqref{eq:a_bound_main}, the lhs. of \eqref{eq:a_plus_bound} is automatically greater or equal to the rhs., so it suffices to prove the reverse inequality. Proceeding as in \eqref{eq:Re_m_over_beta} and using \eqref{eq:a_explicit}, we get that $\max\{\Re a_{12,r}^+,\Re a_{12,r}^-\}$ changes the sign only finitely many times for $r\in [s,t]$. Let $\{[s_i,t_i]\}_{i=1}^k$ be the intervals of positivity of this quantity, ordered so that
\begin{equation*}
s_1\le t_1\le s_2\le\cdots\le s_k\le t_k.
\end{equation*}
We may assume that $s_1=s$ and $t_k=t$, since $s_1$ and $t_1$ ($s_k$ and $t_k$, respectively) may coincide. Applying the second identity in \eqref{eq:a_bound_main} to each of the intervals $[s_i,t_i]$, we get
\begin{equation}
\int_s^t \max\{\Re a_{12,r}^+,\Re a_{12,r}^-,0\}\dif r = \sum_{i=1}^k \log \frac{\beta_{*,s_i}}{\beta_{*,t_i}} + \mathcal{O}(1) = \log \frac{\beta_{*,s}}{\beta_{*,t}} - \sum_{i=1}^{k-1} \log \frac{\beta_{*,t_i}}{\beta_{*,s_{i+1}}}+\mathcal{O}(1).
\label{eq:a_plus_aux}
\end{equation}
By \eqref{eq:gamma_monot_main}, each of the subtracted logarithms in the rhs. of \eqref{eq:a_plus_aux} is bounded from below by $-C$ for some constant $C>0$. Therefore, the rhs. of \eqref{eq:a_plus_aux} is bounded from above by $\log  (\beta_{*,s}/\beta_{*,t})+ \mathcal{O}(1)$. This finishes the proof of \eqref{eq:a_plus_bound}.

The proof of \eqref{eq:d_bound_main} is analogous to the proof of \eqref{eq:a_bound_main} up to several simplifications, which we now discuss. First, by a simple calculation analogous to \cite[Eq. (A.14)]{univ_extr} we have
\begin{equation}
d_{12,r}=\langle M_{12,r}^{E_+}E_-\rangle = \frac{\ii \Im \left[z_{1,r}\overline{z}_{2,r}\right] u_{1,r}u_{2,r}}{\beta_{-,r}\beta_{+,r}}.
\label{eq:d_expl}
\end{equation} 
Similarly to \eqref{eq:Re_m_over_beta} one can show that $\Re d_{12,r}$ changes the sign at most 3 times for $r\in [0,T]$, so it is sufficient to show that
\begin{equation}
\left\vert\Re\int_{s_1}^{t_1}d_{12,r}\dif r\right\vert\lesssim 1
\label{eq:d_bound2}
\end{equation}
for any $s_1,t_1\in [0,T]$, $s_1\le t_1$. Similarly to the argument above \eqref{eq:int_m_over_beta_expl}, we may assume that there exists $r_0\in [s_1,t_1]$ such that $|\beta_{*,r_0}|<\delta/2$. Computing explicitly the integral in the lhs. of \eqref{eq:d_bound2} by the means of \eqref{eq:time_scaling_main}, \eqref{eq:beta_time_dep_main} and \eqref{eq:d_expl}, we get
\begin{equation}
\int_{s_1}^{t_1}d_{12,r}\dif r = \frac{\ii \Im \left[z_{1,0}\overline{z}_{2,0}\right] u_{1,0}u_{2,0}}{\beta_{+,0}-\beta_{-,0}}\left(\log \frac{\beta_{-,s_1}}{\beta_{-,t_1}}-\log \frac{\beta_{+,s_1}}{\beta_{+,t_1}}\right).
\label{eq:int_d_expl}
\end{equation}
Note that by \eqref{eq:time_scaling_main}, the time in the first factor in the rhs. of \eqref{eq:int_d_expl} can be changed from 0 to $r_0$:
\begin{equation*}
\frac{\ii \Im \left[z_{1,0}\overline{z}_{2,0}\right] u_{1,0}u_{2,0}}{\beta_{+,0}-\beta_{-,0}} = \frac{\ii \Im \left[z_{1,r_0}\overline{z}_{2,r_0}\right] u_{1,r_0}u_{2,r_0}}{\beta_{+,r_0}-\beta_{-,r_0}} = \mathcal{O}(|z_{1,r_0}-z_{2,r_0}|),
\end{equation*}
where in the last step we used that $|\beta_{+,r_0}-\beta_{-,r_0}|\ge \delta/2$ and $|u_{1,r_0}u_{2,r_0}|< 1$ by \cite[Eq. (3.5)]{macroCLT_real}. The rest of the proof follows similarly to the argument below \eqref{eq:a_bound_prefin}.

\subsection{Proof of Proposition \ref{prop:overlap}}\label{app:overlap}

Recall the definition of the quantiles $\{\gamma_i^z\}_{i=-N}^N$ of density $\rho^z$, $z\in\C$, from \eqref{eq:def_quantiles}. Since $i,j\in [(1-\tau)N]$ and $|z_1|,|z_2|\le 1-\delta$ in Proposition~\ref{prop:overlap}, there exists a (small) $\kappa>0$ which depends only on $\tau$ and $\delta$, such that $\gamma_i^{z_1}\in\mathbf{B}^{z_1}_\kappa$ and $\gamma_j^{z_2}\in\mathbf{B}^{z_2}_\kappa$. We frequently use this fact throughout the proof and do not refer to it for brevity.

Fix a (small) $\epsilon>0$ and set $\eta:=N^{-1+\epsilon}$. By spectral decomposition and eigenvalue rigidity~\eqref{eq:rigidity} we get
\begin{equation}
N\!\left(\left\vert \langle \bm{u}_i^{z_1},\bm{u}_j^{z_2}\rangle\right\vert^2 + \left\vert \langle \bm{v}_i^{z_1},\bm{v}_j^{z_2}\rangle\right\vert^2\right)\!\!\prec\! (N\eta)^2 \left\langle \Im G^{z_1}(\gamma_i^{z_1}+\ii\eta)\Im G^{z_2}(\gamma_j^{z_2}+\ii\eta) \right\rangle\!\!\prec\! \frac{(N\eta)^2}{\widehat{\beta}_{12}(\gamma_i^{z_1}+\ii\eta,\gamma_j^{z_2}+\ii\eta)}.
\label{eq:overlap_proof} 
\end{equation} 
Here in the first step we additionally used the relation~\eqref{eq:eigenvectors_symmetry} between eigenvectors of~$H^z$ and singular vectors of $X-z$, and in the second step employed~\eqref{eq:2G_av_suboptimal} together with the bound~\eqref{eq:M_12_beta_bound}. Using \eqref{eq:stab_LT} and the fact that $|\partial_w m^z(w)|\lesssim 1$ uniformly in $\Re w\in\mathbf{B}_\kappa^z$ and $0<\Im w\lesssim 1$, we obtain 
\begin{equation}
\begin{split}
\widehat{\beta}_{12}(\gamma_i^{z_1}+\ii\eta,\gamma_j^{z_2}+\ii\eta)&\gtrsim |z_1-z_2|^2 + \mathrm{LT}+\eta,\,\,\,\text{where}\,\,\,\mathrm{LT}=\left\vert \gamma_i^{z_1}-\gamma_j^{z_2} - \frac{\Im u_1}{\Im m_1}\Re \left[\overline{z}_1(z_1-z_2)\right]\right\vert,
\end{split}
\label{eq:beta_quant}
\end{equation}
and $m_1=m^{z_1}(\gamma^{z_1}_i+\ii 0)$, $u_1:=m^{z_1}(\gamma^{z_1}_i+\ii 0)$.

It remains to show that
\begin{equation}
\mathrm{LT}\sim N^{-1}|i-j|+\mathcal{O}(|z_1-z_2|^2).
\label{eq:LT_asymp}
\end{equation}
Indeed, \eqref{eq:overlap_bound} immediately follows from \eqref{eq:overlap_proof}, \eqref{eq:beta_quant} and \eqref{eq:LT_asymp} once $\epsilon>0$ is taken sufficiently small. Denote for short $\Delta z:=z_2-z_1$, $\zeta:=\Delta z/|\Delta z|$ and let $\partial_\zeta$ be derivative in the direction of $\zeta$ in the complex variable $z$, i.e. for a function $f(z)$ we set $\partial_\zeta f(z):=\lim_{h\to 0}(f(z+h\zeta)-f(z))/h$. In the sequel, in case of a function of several variables, $\partial_\zeta$ will always act on the $z$ or $z_1$ variables and not on $w$ and $E$. Since $|\partial_\zeta^2 m^z(w)|\lesssim 1$ uniformly in $\Re w\in\mathbf{B}_\kappa^z$ and $0<\Im w\lesssim 1$, \eqref{eq:def_quantiles} together with \eqref{eq:def_rho} imply that 
\begin{equation}
\gamma_i^{z_1}-\gamma_i^{z_2} = -|\Delta z|\partial_\zeta \gamma_i^{z_1}  + \mathcal{O}(|\Delta z|^2).
\label{eq:quant_dif}
\end{equation}
Differentiating \eqref{eq:def_quantiles} in the direction of $\zeta$ we compute
\begin{equation}
\partial_\zeta \gamma_i^{z_1} = - \frac{1}{\Im m_1} \Im \int_0^{\gamma_i^{z_1}}\partial_\zeta m^{z_1}(E+\ii 0)\dif E.
\label{eq:dif_gamma_init}
\end{equation}
We observe that $\Im u^{z_1}(+\ii 0)=0$ and write $\Im u_1$ as an integral of $\partial_E \Im u^{z_1}(E+\ii 0)$ over $E\in [0,\gamma_i^{z_1}]$. Together with the second part of \eqref{eq:beta_quant}, \eqref{eq:quant_dif} and \eqref{eq:dif_gamma_init} this gives
\begin{equation}
\mathrm{LT} = \left\vert \gamma_i^{z_2}-\gamma_j^{z_2} +\frac{1}{\Im m_1}\Im\int_0^{\gamma_i^{z_1}}\Big(|\Delta z|\partial_\zeta m^{z_1}(E+\ii 0) + \Re \left[\overline{z}_1\Delta z\right] \partial_E u^{z_1}(E+\ii 0)\Big)\dif E\right\vert + \mathcal{O}(|\Delta z|^2).
\label{eq:LT_int_rep}
\end{equation}
Now we show that the function integrated in \eqref{eq:LT_int_rep} vanishes for every $E\in [0,\gamma^{z_1}_i]$. Since along this argument $E$ and $z_1$ remain fixed, we simply denote $m:=m^{z_1}(E+\ii 0)$, $u:=u^{z_1}(E+\ii 0)$ and $M:=M^{z_1}(E+\ii 0)$ (these notations should not be confused with $m_1$ and $u_1$, which correspond to $E:=\gamma^{z_1}_i$).  We differentiate \eqref{eq:MDE} in the direction of $\zeta$ and use \eqref{eq:M}, arriving to
\begin{equation}
\partial_\zeta m = \partial_\zeta \langle M\rangle = -2\Re [\overline{z}_1\zeta]mu(1-\langle M^2\rangle)^{-1}.
\label{eq:dif_zeta_M}
\end{equation}
To compute $\partial_E u$, we first differentiate \eqref{eq:MDE} in $E$ and get
\begin{equation}
\partial_E m = \partial_E\langle M\rangle = \langle M^2\rangle(1-\langle M^2\rangle)^{-1}.
\label{eq:dif_w_M}
\end{equation}
We further recall the definition of $u^z$ from \eqref{eq:M} and apply \eqref{eq:dif_w_M}:
\begin{equation}
\partial_E u =\partial_E \frac{m}{E+m} = \frac{mu}{1-\langle M^2\rangle}\left(\frac{\langle M^2\rangle}{m^2}-\frac{1}{m(E+m)}\right) = \frac{2mu}{1-\langle M^2\rangle},
\label{eq:dif_w_u}
\end{equation}
where in the last step we used \eqref{eq:m}. Combining \eqref{eq:dif_zeta_M} with \eqref{eq:dif_w_u}, we obtain
\begin{equation}
|\Delta z|\partial_\zeta m + \Re \left[\overline{z}_1\Delta z\right] \partial_E u=0,
\end{equation}
which together with \eqref{eq:LT_int_rep} implies
\begin{equation}
\mathrm{LT}=|\gamma^{z_2}_i-\gamma^{z_2}_j|+\mathcal{O}(|\Delta z|^2).
\end{equation}
It remains to notice that by \eqref{eq:def_quantiles}, $|\gamma^{z_2}_i-\gamma^{z_2}_j|\sim N^{-1}|i-j|$. This finishes the proof of \eqref{eq:LT_asymp}, and thereby of Proposition~\ref{prop:overlap}.\hfill $\qed$

\bibliographystyle{plain} 
\bibliography{refs}

\clearpage

\section*{Supplementary Material}

\setcounter{section}{0}
\renewcommand{\thesection}{S\arabic{section}}
\renewcommand{\sectionname}{Section}

\section{Decorrelation for the Dyson Brownian motion}\label{sec:DBM}

In this section we prove Theorem~\ref{eq:maintheoDBM1_main} by explicitly constructing the comparison processes $\{\mu^{(l)}_i(t)\}_{|i|\le N}$ for $l=1,2$, and relying on the recent proof of \cite[Theorem~4.1]{bourgade2024fluctuations}. For completeness, we remind the reader the set-up of Theorem~\ref{eq:maintheoDBM1_main}. The complex i.i.d. matrix $X$ is embedded into the flow
\begin{equation}
\label{eq:matDBM1}
\dif X_t=\frac{\dif B_t}{\sqrt{N}}, \qquad\quad X_0=X,
\end{equation}
where the entries $(B_t)_{ab}$ are independent standard complex Brownian motions.
Let $H_t^z$ be the Hermitization of $X_t-z$ defined as in \eqref{eq:def_hermitization}, denote its eigenvalues by $\lambda_i^z(t)$, and by ${\bf w}_i^z(t)\in\C^{2N}$ the corresponding eigenvectors. We recall that due to the chiral symmetry of $H_t^z$ its spectrum is symmetric around zero, i.e. $\lambda_{-i}(t)=-\lambda_i(t)$ for any $i\in [N]$ and $t\ge 0$. As a consequence the eigenvectors corresponding to $\lambda_{\pm i}(t)$ are given by\footnote{Here $\mathfrak{t}$ denotes the transpose.} ${\bf w}_{\pm i}^z(t)=({\bf u}_i^z(t),\pm {\bf v}_i^z(t))^\mathfrak{t}\in\C^{2N}$, for $i\in [N]$, where ${\bf u}_i^z(t), {\bf v}_i^z(t)$ are the left/ right singular vectors of $X_t-z$, respectively. It is well known \cite[Eq.(5.8)]{loc_relax_flow} that $\lambda_i^z(t)$ are the unique strong solution of the Dyson Brownian motion: 
 \begin{equation}
\label{eq:DBM1}
\dif \lambda_i^z(t)=\frac{\dif b_i^z(t)}{\sqrt{2N}}+\frac{1}{2N}\sum_{j\ne i}\frac{1}{\lambda_i^z(t)-\lambda_j^z(t)}\dif t,
\end{equation}
where, for $i\in [N]$, the $b_i^z(t)$ are independent standard real Brownian motions, and $b_{-i}^z(t)=-b_i^z(t)$. We point out that here, and throughout this section, the summation $\sum_{j\ne i}$ is performed over indices $|j|\le N$. However, for different $z$'s the driving Brownian motions have a non-trivial correlation depending on the overlap of the singular vectors
\begin{equation}
\label{eq:corrbms11}
\dif \big[b_i^{z_1},  b_j^{z_2}\big]_t=4\Re\big[\langle {\bf u}_i^{z_1}(t),{\bf u}_j^{z_2}(t)\rangle\langle {\bf v}_j^{z_2}(t), {\bf v}_i^{z_1}(t)\rangle \big]\dif t.
\end{equation}

We construct the comparison processes $\{\mu_i^{(l)}(t)\}_{|i|\le N}$ for $l=1,2$ as follows. Consider two independent sets $\{\mu_i^{(l)}\}_{|i|\le N}$, for $l=1,2$, of eigenvalues of the Hermitization of two independent Ginibre matrices, and let $\mu_i^{(l)}(t)$ be the solution of \eqref{eq:DBM1} with initial condition $\mu_i^{(l)}$ and with the ${\bm \lambda}$'s being replaced by the ${\bm \mu}$'s.  More precisely, we define the $\mu_i^{(l)}(t)$ as the unique strong solutions of
\begin{equation}
\label{eq:coupleDBM1}
\dif \mu_i^{(l)}(t)=\frac{\dif \beta_i^{(l)}(t)}{\sqrt{2N}}+\frac{1}{2N}\sum_{j\ne i}\frac{1}{\mu_i^{(l)}(t)-\mu_j^{(l)}(t)}\dif t.
\end{equation}
Here $\{(\beta_i^{(l)}(t))_{i=1}^N, l=1,2\}$ is a family of i.i.d. standard real Brownian motions, and it is extended to negative indices by symmetry similarly to the $b_i^z(t)$'s. In particular, we stress that the two processes $\{\mu_i^{(l)}(t),\, |i|\le N\}$, for $l=1,2$, are fully independent. In order to show the independence of the solutions of \eqref{eq:DBM1} for $z_1,z_2$ with $|z_1-z_2|\gg N^{-1/2}$, we will (almost completely) couple the evolution of the $\lambda_i^{z_l}(t)$ with the two fully independent processes $\mu_i^{(l)}(t)$.

Fix $T, R>0$. We will now couple the driving martingales in \eqref{eq:DBM1}--\eqref{eq:coupleDBM1} for times $t\in [0,T]$ and indices $|i|\le R$. To simplify the notation, here and throughout this section we assume that $R$ is an integer. We point out that this coupling idea first originated in \cite{macroCLT_complex} and it was then used in several further works (see e.g. \cite{minor, bao2025numerical, bourgade2024fluctuations, cipolloni2023quenched, Cipolloni_meso, cipolloni2024maximum}).  By the martingale representation theorem \cite[Theorem~18.12]{Kallenberg} we can write
\begin{equation}
\label{eq:coupl1}
\left(\begin{matrix}
\rd \mathbf{b}_t^{z_1} \\
\rd \mathbf{b}_t^{z_2}
\end{matrix}\right)=\sqrt{C(t)}\left(\begin{matrix}
\rd \mathbf{\bm \beta}_t^{(1)} \\
\rd \mathbf{\bm \beta}_t^{(2)}
\end{matrix}\right),
\end{equation}
where $\mathbf{b}_t^{z_l}:=\big(b_1^{z_l}(t), b_2^{z_l}(t), \dots, b_R^{z_l}(t)\big)$, for $l=1,2$, and ${\bm \beta}_t^{(l)}$ is defined similarly.  In the following we may use the short-hand notations $\mathbf{b}_t:=(\mathbf{b}_t^{z_1}, \mathbf{b}_t^{z_2})^t$ and ${\bm \beta}_t:=({\bm \beta}_t^{(1)},{\bm \beta}_t^{(2)})^t$. Here $C(t)$ is the  $(2R) \times (2R)$ covariance matrix of the vector $\mathbf{b}_t$ (see \eqref{eq:corrbms11}). Additionally, we stress that the $\beta_i^{(l)}(t)$ for larger indices are independent from $\{(b_i^{z_l}(t))_{i=1}^N, l=1,2\}$, and that $\{(b_i^{z_l}(t))_{i=R+1}^N, l=1,2\}$ are independent of the $\{(\beta_i^{(l)}(t))_{i=1}^N, l=1,2\}$.  

\begin{theorem}
\label{eq:maintheoDBM1}
Let $X$ be a complex i.i.d. matrix, and let $X_t$ be the solution of \eqref{eq:matDBM1} with initial condition $X_0=X$. Fix any small $\omega_*, \delta>0$, $z_1,z_2\in\C$ with $|z_1|,|z_2|<1-\delta$, $T>N^{-1+\omega_*}$, and a possibly $N$-dependent $0<R<N|z_1-z_2|^2$. For $l=1,2$, let $\lambda_i^{z_l}(t)$, $\mu_i^{(l)}(t)$ be the solutions of \eqref{eq:DBM1} and \eqref{eq:coupleDBM1}, respectively. Assume that for any small $\xi>0$ and any large $D>0$ we have
\begin{equation}
\label{eq:evectors1}
\mathbf{P}\left(\big|\langle {\bf u}_i^{z_1}(t),{\bf u}_j^{z_2}(t)\rangle\big|^2+\big|\langle {\bf v}_j^{z_2}(t), {\bf v}_i^{z_1}(t)\rangle\big|^2\le \frac{N^\xi}{N|z_1-z_2|^2} \,\, \forall\, t\in [0,T] \, \mathrm{and} \, |i|,|j|\le R\right)\ge 1-N^{-D}.
\end{equation}
Then \eqref{eq:DBM_bound2} holds for $l=1,2$ with very high probability simultaneously for all $t\in [N^{-1+\omega_*},T]$ and indices $|i|\le R$, i.e. we have
\begin{equation}
\label{eq:desbnew1}
\big|\rho_t^{z_l}(0)\lambda_i^{z_l}(t)-\rho_{\mathrm{sc},t}(0)\mu_i^{(l)}(t)\big|\le N^\xi\!\!\left[\sqrt{\frac{t}{N}}\left(\frac{1}{R^{1/8}}\!+\!\sqrt{\frac{R}{N|z_1-z_2|^2}}\right)\!+\!\frac{\sqrt{Nt^3}}{R}\!+\!\frac{|i|}{N}\left(\frac{1}{\sqrt{Nt}}\!+\!\frac{|i|}{N}+t\right)\!\right]\!.
\end{equation}
\end{theorem}

We point out that \eqref{eq:desbnew1} improves on the previous results in \cite{bourgade2024fluctuations, macroCLT_complex}, giving an explicit error in terms of~$t$, $|z_1-z_2|$, and the number of coupled Brownian motions $R$. However, at the technical level no new inputs are needed compared to \cite{bourgade2024fluctuations}. We just put together several estimates presented in \cite{bourgade2024fluctuations}. For the convenience of the reader we present the main steps of the proof. At the end of this section we comment on a possible improvement of the bound in the rhs. of \eqref{eq:desbnew1}, which, however, we do not pursue for brevity.

Relying on Theorem~\ref{eq:maintheoDBM1}, we now prove Theorem~\ref{eq:maintheoDBM1_main}.
\begin{proof}[Proof of Theorem~\ref{eq:maintheoDBM1_main}] It suffices to verify \eqref{eq:evectors1}, which is equivalent to proving that the overlap bound \eqref{eq:overlap_bound} holds simultaneously for all $t\in [0,T]$ with very high probability. Since the eigenvectors of $H^{z_1}, H^{z_2}$ are highly unstable, we instead show that
\begin{equation}
\left\vert\langle \Im G^{z_1}_t(w_1)\Im G^{z_2}_t(w_2)\rangle\right\vert \lesssim N^\xi \left(\widehat{\beta}_{12}(w_1,w_2)\right)^{-1}
\label{eq:grid_explanation}
\end{equation}
simultaneously for all $t\in [0,T]$ and $w_l\in \C\setminus\R$ with $\Re w_l\in\mathbf{B}_\kappa^{z_l}$ and $N^{-1+\epsilon}\le |\Im w_l|\le 1$, for any fixed $\xi,\kappa,\epsilon>0$. Then \eqref{eq:evectors1} follows from \eqref{eq:grid_explanation} as explained in the proof of Proposition~\ref{prop:overlap} presented in Supplementary Section~\ref{app:overlap}. By Proposition \ref{prop:2G_av} and \eqref{eq:M_12_beta_bound}, \eqref{eq:grid_explanation} holds for any $t\in [0,T]$ and $w_1,w_2$ satisfying conditions above with very high probability. The fact that \eqref{eq:grid_explanation} holds simultaneously for all values of these parameters with probability at least $1-N^{-D}$, for any fixed $D>0$, follows by a standard grid argument which proceeds by choosing sufficiently dense meshes for $t$, $w_1$, $w_2$, applying Proposition~\ref{prop:2G_av} to each of the mesh points, and then extending these result to the rest of the parameter values by H\"older continuity of the lhs. of \eqref{eq:grid_explanation}. This finishes the verification of~\eqref{eq:evectors1}.
\end{proof}

It remains to prove Theorem~\ref{eq:maintheoDBM1}.

\begin{proof}[Proof of Theorem~\ref{eq:maintheoDBM1}]
The proof of \eqref{eq:desbnew1} follows by combining \cite[Proposition~4.5,4.6]{bourgade2024fluctuations}, as done in the proof of \cite[Theorem~4.1]{bourgade2024fluctuations}. The only difference in this case is that we have a better control on the eigenvector overlaps in \eqref{eq:evectors1}, and thus obtaining an explicit dependence in $t$, $R$, and $|z_1-z_2|^2$ in the estimates in Proposition 4.5. We now explain the main steps of this proof, but omit several details that can be found in \cite{bourgade2024fluctuations}. To keep the presentation simpler we assume that the limiting eigenvalue density at $x=0$ is the same for $\lambda_i^{z_l}$ and $\mu_i^{(l)}$; if this is not the case, it can be achieved by a simple time re-scaling, we omit the details for brevity.

To study the closeness of $\lambda_i^{z_l}(t)$ and $\mu_i^{(l)}(t)$ in \eqref{eq:coupleDBM1} we introduce the interpolating process (see e.g.~\cite{fixed_E}) 
\begin{equation}
\label{eq:interdbm1}
\rd x_i^{z_l}(\alpha, t)=\frac{\rd b_i^{z_l}(\alpha,t)}{\sqrt{2N}}+\frac{1}{2N}\sum_{j\ne i}\frac{1}{x_i^{z_l}(\alpha, t)-x_j^{z_l}(\alpha, t)}\rd t,
\quad l=1,2, \quad \alpha\in [0,1],
\end{equation}
with initial data $x_i^{z_l}(\alpha, 0)=\alpha\lambda_i^{z_l}+(1-\alpha)\mu_i^{(l)}$, and $b_i^{z_l}(\alpha,t):=\alpha b_i^{z_l}(t)+(1-\alpha)\beta_i^{(l)}(t)$.

Since $z_1,z_2$ are fixed  and once the coupling in \eqref{eq:coupl1} is performed,
the index $l=1,2$ is also fixed, from now on we may drop the $z_l$ dependence of $x_i^{z_l}(\alpha, t)$ from the notation and use $x_i(\alpha, t)$ throughout this section. To further simplify the notation we consider $l=1$, so that in the computations below the only $(\sqrt{C(t)})_{ij}$ will appear are for $i,j\le R$. Additionally, since the whole analysis is performed for a fixed $\alpha\in [0,1]$, we may often omit the $\alpha$-dependence from the notation. 

The well-posedness of \eqref{eq:interdbm1} and the fact that its tangential flow $u_i(t)=u_i(\alpha, t):=\partial_\alpha x_i(\alpha,t)$ are well defined follows from \cite[Appendix~A]{macroCLT_real} and \cite[Appendix~B]{bourgade2024fluctuations}. In the following we will show that $|u_i(t)|\le rhs. \, \eqref{eq:desbnew1}$, which immediately gives the desired bound in \eqref{eq:desbnew1} by integrating in $\alpha$. By differentiating \eqref{eq:interdbm1} in $\alpha$, one sees that $u_i(t)$ is the solution of\footnote{Here ${\bm u}$ should not be confused with the singular vectors introduced below \eqref{eq:matDBM1}. A similar comment applies to ${\bm v}$ below.} (here $|i|, |j|\le N$)
\begin{equation}
\label{eq:derevol1}
\dif {\bm u}(t)=\big(\mathcal{B}{\bm u}\big)(t)\dif t+\frac{\dif {\bm \xi}}{\sqrt{N}}, \qquad\quad (\mathcal{B}{\bm u})_i:=\sum_j\mathcal{B}_{ij}(u_j-u_i), \quad \xi_i:=\partial_\alpha  b_i(\alpha,t),
\end{equation}
where $\partial_\alpha  b_i(\alpha,t)=b_i^{z_l}(t)-\beta_i^{(l)}(t)$ and
\[
\mathcal{B}_{ij}=\mathcal{B}_{ij}(\alpha, t):=\frac{1}{N(x_i(\alpha, t)-x_j(\alpha, t))^2}.
\]

Note that for the stochastic term in \eqref{eq:derevol1} we have (cf. with \cite[Assumption~4.2]{bourgade2024fluctuations})
\begin{equation}
\label{eq:newstoch}
\dif \big[\xi_i,\xi_j\big]_s=\dif \big[ b_i-\beta_i, b_j-\beta_j\big]_s\sim \begin{cases}
\big[(\sqrt{C(s)}-I)^2\big]_{i,j}\dif s & \mathrm{for}  |i|, |j|\le R, \\
\delta_{ij}\dif s & \mathrm{otherwise},
\end{cases}
\end{equation}
with
\begin{equation}
\label{eq:betterbstochterm}
\big|\big[(\sqrt{C(s)}-I)^2\big]_{i,j}\big|\lesssim \frac{R}{(N|z_1-z_2|^2)^2}.
\end{equation}
This last bound follows from
\begin{equation}
\label{eq:betterb}
\begin{split}
\big[(\sqrt{C(s)}-I)^2\big]_{ii}&=\sum_l (\sqrt{\lambda_l}-1)^2|\psi_l(i)|^2\le \sum_l (\lambda_l-1)^2|\psi_l(i)|^2 \\
&= \big[(C(s)-I)^2\big]_{ii}\le \sum_j \big|(C(s)-I)_{ij}\big|^2\lesssim \frac{R}{(N|z_1-z_2|^2)^2},
\end{split}
\end{equation}
and a simple Schwarz inequality for off-diagonal terms.

To show that $|u_i(t)|$ is bounded by the rhs. of \eqref{eq:desbnew1} we will follow two main steps. In particular, in the first step, we will study an analog of \eqref{eq:derevol1} without the stochastic term, as it was done in \cite[Proposition~4.6]{bourgade2024fluctuations}. Then, in the second step we will show that the stochastic term in \eqref{eq:derevol1} can in fact be removed if we assume that the $\lambda$ and $\mu$ flows have the same initial condition, i.e. for ${\bm u}_i(0)=0$; this will be achieved by two further intermediate steps by first estimating the long-range contribution and then studying the short-range one.

Let ${\bm v}(t)$ be the solution of the analog of \eqref{eq:derevol1} without the stochastic term, i.e.
\begin{equation}
\dif {\bm v}(t)=\big(\mathcal{B}{\bm v}\big)(t)\dif t,
\end{equation}
with two different initial datas ${\bm v}_i(0)=\mathrm{sgn}(i)(i^{1/2}/N)$ or ${\bm v}_i(0)=\mathrm{sgn}(i)(i/N)^2$ (see the paragraph after the proof of \cite[Lemma 4.21]{bourgade2024fluctuations} for more details). Then, by \cite[Proposition~4.6]{bourgade2024fluctuations}, it follows that\footnote{\label{ftn:DBM} We point out that in the current setting we could improve the rhs. of \eqref{eq:estv} by replacing $1/\sqrt{Nt}$ with $1/(Nt)$, as our initial conditions satisfy rigidity bounds (while in \cite{bourgade2024fluctuations} only a weaker condition was assumed). This confirms what was expected in \cite[Remark 4.19]{bourgade2024fluctuations}.}
\begin{equation}
\label{eq:estv}
\big|{\bm v}_i(t)\big|\lesssim N^\xi\frac{|i|}{N}\left(\frac{1}{\sqrt{Nt}}+\frac{|i|}{N}+t\right), \qquad\quad i\in [N],
\end{equation}
with very high probability. The combination of these two steps will give the desired result.

We are thus left with giving the details of the proof of the first step. Let ${\bm u}(t)$ be the solution of \eqref{eq:derevol1}, with ${\bm u}_i(0)=0$, and define the short- and long-range operators:
\begin{equation}
(\mathcal{S}{\bm u})_i:=\sum_{(i,j)\in \mathcal{I}_\mathcal{L}^c, \atop j\ne i}\mathcal{B}_{ij}(u_j-u_i), \qquad\quad (\mathcal{L}{\bm u})_i:=\sum_{(i,j)\in \mathcal{I}_\mathcal{L}}\mathcal{B}_{ij}(u_j-u_i).
\end{equation}
Here $0<\ell\le R/100$ is an $N$-dependent scale which we will choose later in the proof, and
\[
\mathcal{I}_\mathcal{L}:=\{(i,j):|i-j|>\ell, \min(|i|,|j|)\le cN\},
\]
for a small fixed $c>0$. We define the vector ${\bm w}\in \R^N$ as the solution of the following evolution with the short-range operator
\begin{equation}
\label{eq:shortstoch}
\dif {\bm w}(t)=\big(\mathcal{S}{\bm w}\big)(t)\dif t+\frac{\dif {\bm \xi}}{\sqrt{N}}, \qquad\quad {\bm w}(0)=0.
\end{equation}
We now first show that the contribution of the long-range part is negligible and then give a bound on ${\bm w}(t)$ using a weighted $\ell^2$ estimate. All the estimate below hold with very high probability, even if not stated explicitly. First, by \cite[Lemma~4.12]{bourgade2024fluctuations} we have
\begin{equation}
\lVert {\bm w}(t)\rVert_\infty\lesssim N^\xi\sqrt{\frac{t}{N}},
\end{equation}
with very high probability. Then, using Duhamel's formula we write (see \cite[Eq.(4.24)--(4.26)]{bourgade2024fluctuations})
\begin{equation}
\label{eq:apprduhamel}
{\bm u}(t)={\bm w}(t)+\int_0^t \mathcal{U}_\mathcal{S}(s,t)\mathcal{L}{\bm w}(s)\, \dif s={\bm w}(t)+\mathcal{O}\left(\frac{(\log N)^2\sqrt{Nt^3}}{\ell}\right),
\end{equation}
where $\mathcal{U}_\mathcal{S}$ is the propagator of $\partial_t {\bm f}_t=(\mathcal{S}{\bm f})_t$, and the error term is meant entry-wise. This shows that the long-range part can be neglected and so that it is enough to study ${\bm w}(t)$ to obtain the desired bound on ${\bm u}(t)$.

Next, to show that the effect of the stochastic term in \eqref{eq:shortstoch} is also negligible, we use a weighted $\ell^2$-bound. More precisely, one studies the evolution of (here $\zeta:=N/\ell$ and $\chi$ is a non-negative, symmetric, and smooth cut-off with $\int\chi=1$)
\begin{equation}
X_s:=\sum_i \phi_i(s){\bm w}(s), \qquad\quad \phi_i(s):=e^{-\zeta\psi(x_i(s))}, \qquad\quad \psi(x):=\zeta\int \min\{|x-y|, c\}\chi(\zeta y)\,\dif y
\end{equation}
as in \cite[Lemma~4.13]{bourgade2024fluctuations}. The main difference in this proof compared to \cite[Lemma~4.13]{bourgade2024fluctuations} is that the bound in \cite[Assumption~4.2]{bourgade2024fluctuations} can be improved to \eqref{eq:newstoch}--\eqref{eq:betterbstochterm} for $K=R$, giving a more explicit bound (which was already obtained in \cite{bourgade2024fluctuations}, but not exploited). Using the notation of \cite[Lemma~4.13]{bourgade2024fluctuations}, choosing $\ell=R/100$, $N^{\omega_K}=R$, and replacing $N^{-\tilde{\omega}}$ with $R/(N|z_1-~\!z_2|^2)^2$ from \eqref{eq:betterbstochterm}, we obtain
\begin{equation}
\begin{split}
X_t&\lesssim \frac{t\log N}{NR}+\frac{\sqrt{t}}{N^{3/2}R}+\frac{\sqrt{t R^2}}{N^{3/2} (N|z_1-z_2|^2)}+\frac{t}{N R^{1/4}}+\frac{t R^2}{N(N|z_1-z_2|^2)^2} \\
&\lesssim\frac{t}{N}\left(\frac{1}{R^{1/4}}+\frac{R^2}{(N|z_1-z_2|^2)^2}\right)+\frac{\sqrt{t}}{N^{3/2}}\left(\frac{1}{R}+\frac{R}{N|z_1-z_2|^2}\right).
\end{split}
\end{equation}
Next, using $R<N|z_1-z_2|^2$ and that $Nt\gg 1$, we obtain
\begin{equation}
\frac{t}{N}\left(\frac{1}{R^{1/4}}+\frac{R^2}{(N|z_1-z_2|^2)^2}\right)+\frac{\sqrt{t}}{N^{3/2}}\left(\frac{1}{R}+\frac{R}{N|z_1-z_2|^2}\right)\lesssim \frac{t}{N}\left(\frac{1}{R^{1/4}}+\frac{R}{N|z_1-z_2|^2}\right).
\end{equation}
We have thus concluded that for $|i|\le R$ we have
\begin{equation}
\label{eq:estw}
|w_i(t)|\lesssim X_t^{1/2}\lesssim \sqrt{\frac{t}{N}}\left(\frac{1}{R^{1/4}}+\frac{R}{N|z_1-z_2|^2}\right)^{1/2}.
\end{equation}

Finally, using that for $|i|\le R$ we have (schematically)
\begin{equation}
\big|\lambda_i^{z_l}(t)-\mu_i^{(l)}(t)\big|\lesssim \int_0^1\big(|v_i(t,\alpha)|+|w_i(t,\alpha)|\big)\,\dif \alpha,
\end{equation}
and combining \eqref{eq:estv}, \eqref{eq:apprduhamel}, and \eqref{eq:estw}, we obtain the bound \eqref{eq:desbnew1} and thus conclude the proof.
\end{proof}

\begin{remark}[Possible improvement on \eqref{eq:desbnew1}]\label{rem:DBM}
The main ingredient in the proof of Theorem~\ref{eq:maintheoDBM1} is an $\ell^2$-estimate to show that the stochastic term in \eqref{eq:derevol1} can be neglected, at the price of an explicit error depending on $t, R, |z_1-z_2|$. Inspired by \cite{cipolloni2022normal, marcinek2022high}, we believe, that this error can be improved complementing the $\ell^2$-estimate with a Nash-argument which gives the $\ell^2\to \ell^\infty$ contractivity of the kernel $\mathcal{B}$. Following \cite[Sections 6-7]{marcinek2022high}, it seems that the best possible error one could possibly obtain with this procedure is something of the form (for $T_1<T_2<T$)
\begin{equation}
\label{eq:firstb}
|\lambda_i(T_2)-\mu_i(T_2)|\prec \frac{1}{N}\left(\frac{NT_1}{R}+T_2-T_1\right)+\frac{1}{\sqrt{N(T_2-T_1)}}\times \mathrm{rhs. \,\, \eqref{eq:desbnew1}\, with \, t=T_1}.
\end{equation}
In addition, by Footnote~\ref{ftn:DBM}, the term $1/\sqrt{Nt}$ in \eqref{eq:desbnew1} can be replaced with $1/(Nt)$. This bound would clearly give a better error in Proposition~\ref{prop:EGG}. Carefully inspecting the proof of Theorem~\ref{theo:main1} one can see that with this better error $q_0$. Preliminary calculations indicate that with this better error it is possible to choose $q_0=1/3$. This is a significant improvement compared to $q_0=1/20$, however, it is still very far from the optimal choice $q_0=1$. For this reason we do not pursue the details of this proof and present all the details just to obtain $q_0=1/20$.

\end{remark}

\section{Proof of the multi-resolvent local laws from Section~\ref{sec:local_laws}}\label{sec:multiG_proof}

In this section we prove local laws from Propositions~\ref{prop:2G_av}, \ref{prop:2G_av_suboptimal}, and \ref{prop:multiG_oneG}. Our main focus is on Proposition~\ref{prop:2G_av} in the special case $b=0$, since this is our main novelty. Meanwhile, the general case $b\in [0,1]$ follows by a minor adjustment of this special case, and the proof of Proposition \ref{prop:2G_av_suboptimal} follows by a minor modification of the proof of \cite[Theorem~3.4]{nonHermdecay} using our new propagator bound from Lemma~\ref{lem:propag_bound2} (see also Lemma~\ref{lem:propag_bound} later), allowing to extend the result of \cite[Theorem~3.4]{nonHermdecay} from the imaginary axis to the entire bulk regime.

This section is structured in the following way. First, by the end of Section~\ref{sec:Zig_old} we prove Proposition~\ref{prop:2G_av} for $b=0$. Then, in Sections~\ref{sec:2G_av_subopt} and \ref{sec:general_b} we prove Proposition~\ref{prop:2G_av_suboptimal} and Proposition~\ref{prop:2G_av} for general $b\in [0,1]$.

Recall that $X_0$ is a complex i.i.d. matrix whose single-entry distribution does not need to be Gauss-divisible, whereas $X=\sqrt{1-\mathfrak{s}^2}X_0+\mathfrak{s}\widetilde{X}$ contains a Gaussian component. Consider the matrix-valued Ornstein-Uhlenbeck process starting from $X_0$:
\begin{equation}
\dif X_t = -\frac{1}{2}X_t + \frac{1}{\sqrt{N}}\dif B_t,
\label{eq:OU}
\end{equation}
where the entries of $B_t$ are independent standard complex-valued Brownian motions. Then
\begin{equation}
X_t\stackrel{d}{=}\ee^{-t/2}X_0 + \sqrt{1-\ee^{-t}}\widetilde{X},\quad \forall t\ge 0.
\label{eq:OU_distr}
\end{equation}
In particular, $X_t$ and $X$ follow the same distribution for $t=|\log(1-\mathfrak{s}^2)|$. Denote 
\begin{equation}
\label{eq:defres}
G_t^{z}(w):=\begin{pmatrix}
-w&X_t-z\\X_t^*-\overline{z}&-w
\end{pmatrix}^{-1},\quad \forall \, z\in\C,\, w\in\C\setminus\R. 
\end{equation}
With this notations the statement of Proposition~\ref{prop:2G_av} is equivalent to the 2-resolvent averaged local law for $G_t^{z_1}(w_1)$, $G^{z_2}_t(w_2)$ for a fixed $t>0$. From now on we adopt this formulation.

For convenience, we slightly modify the notations introduced in Proposition \ref{prop:2G_av} and use $X$ to denote a general $N\times N$ complex i.i.d. matrix satisfying Assumption \ref{ass:chi}(i) without requiring Gauss-divisibility. The original meaning of $X$ will appear only in the conclusion of the proof of Proposition \ref{prop:2G_av} around \eqref{eq:smaller_G_comp}--\eqref{eq:domain_inclusion}, so no ambiguity will arise. We also denote the resolvent associated to $X$ by $G^z(w)$ and set
\begin{equation}
G_j:=G^{z_j}(w_j),\quad \eta_j:=|\Im w_j|,\,\,\, j=1,2\quad \text{and}\quad \eta_*:=\eta_1\wedge\eta_2\wedge 1.
\end{equation}
for $z_j\in\C$ and $w_j\in\C\setminus\R$.

The proof of Proposition~\ref{prop:2G_av} proceeds in two steps. First we establish the global law, i.e. we prove \eqref{eq:2G_av} for $X_0$ and for spectral parameters $w_1, w_2$ with imaginary parts at least of order one. In the next step, we extend this result dynamically down to the real line in the bulk regime. This is achieved by complementing the process $\{X_t\}_{t\ge 0}$ by the deterministic evolution of $z_j$'s and $w_j$'s  along the \emph{characteristic flow} for $j=1,2$, and then tracking the resolvents under the combined time dependence arising both from $X_t$ and $z_j,w_j$. Our argument is a special case of the \emph{zigzag} strategy (introduced in \cite{cipolloni2024out}, see also \cite{Cipolloni_meso, cipolloni2023eigenstate}). Besides the global law and the \emph{zig} (characteristic flow) step described above, the full zigzag method includes also the \emph{zag} step which removes the Gaussian component from $X_t$ by the means of the Green function comparison. The zig and the zag steps are repeated alternatingly, gradually decreasing the imaginary part of the spectral parameter. This strategy allows to prove multi-resolvent local laws for matrices without a Gaussian component, and can be also applied in our set-up to prove \eqref{eq:2G_av} with $b=0$ for a general i.i.d. matrix~$X$. However, since Proposition~\ref{prop:2G_av} is formulated for Gauss-divisible ensembles, we omit the zag step and perform only one zig step.

The characteristic flow method does not allow to prove \eqref{eq:2G_av} by operating solely with the quantities of the form $\langle G_1B_1G_2B_2\rangle$ for some observables $B_1, B_2$, as it creates longer products of resolvents. To close the bounds on the products of two resolvents we also need to follow the evolution of products of three and four observables, albeit with less precision. Specifically, we consider resolvent chains of the following form:
\begin{equation}
\left\langle G_1B_1G_2B_2G_1^{(*)}B_3\right\rangle \quad \text{and}\quad \left\langle G_1B_1G_2B_2G_1^{(*)}B_3G_2^{(*)}B_4\right\rangle.
\label{eq:chain_types}
\end{equation}
Here $G_j^{(*)}$ indicates both choices $G_j$ and $G_j^*$. Throughout this section we will use the same convention for the spectral parameters, denoting by $w_j^{(*)}$ both choices $w_j$ and $\overline{w}_j$. Moreover, these choices for $G_j$ and $w_j$ each time are performed consistently due to the relation $(G^{z_j}(w_j))^*=G^{z_j}(\overline{w}_j)$. Note that the resolvent chains in \eqref{eq:chain_types} may contain simultaneously $G_1$ and $G_1^*$ (and similarly for $G_2$), since the first $G_1$ always comes without the conjugation, while the second one may be taken equal to $G_1^*$.

\begin{proposition}[Global law]\label{prop:global_law} Fix (small) $N$-independent constants $\delta, \varepsilon>0$. Let $X$ be a complex $N\times N$ i.i.d. matrix satisfying Assumption \ref{ass:chi}(i) and let $G^z(w)$ be defined as in \eqref{eq:def_resolvent}. For $z_1,z_2\in\C$ and $w_1,w_2\in\C\setminus\R$ denote $G_j:=G^{z_j}(w_j)$, $j=1,2$. We have
\begin{align}
\left\langle \left(G_1B_1G_2-M_{12}^{B_1}(w_1,w_2)\right)B_2\right\rangle&\prec \frac{1}{N},\label{eq:global_2G}\\
\left\langle \left(G_1B_1G_2B_2G_1^{(*)}-M_{121}^{B_1,B_2}\big(w_1,w_2,w_1^{(*)}\big)\right)B_3\right\rangle&\prec \frac{1}{N},\label{eq:global_3G}\\
\left\langle \left(G_1B_1G_2B_2G_1^{(*)}B_3G_2^{(*)}-M_{1212}^{B_1,B_2,B_3}\big(w_1,w_2,w_1^{(*)},w_2^{(*)}\big)\right)B_4\right\rangle&\prec \frac{1}{N},\label{eq:global_4G}
\end{align}
uniformly in $z_j\in(1-\delta)\bf{D}$, $w_j\in\C\setminus\R$ with $|\Im w_j|\ge \varepsilon$ for $j=1,2$, and $B_i\in\mathrm{span}\{E_{\pm},F^{(*)}\}$ for $i\in [4]$.
\end{proposition}

\begin{proof}[Proof of Proposition \ref{prop:global_law}] The two-resolvent global law \eqref{eq:global_2G} is a well-established result, e.g. this is a special case of \cite[Theorem 5.2]{macroCLT_complex}. The proof of \eqref{eq:global_3G}--\eqref{eq:global_4G} follows the same strategy as the global multi-resolvent laws for Wigner matrices \cite[Appendix B]{Multi_res_llaws}, with the simplification that we do not track the dependence of the bounds \eqref{eq:global_2G}--\eqref{eq:global_4G} on $\eta_1$ and $\eta_2$. In fact this dependence gives additional smallness in the regime $\eta_1,\eta_2\gg 1$, but is not needed here.

We present a sketch of the proof of \eqref{eq:global_3G}, referring to \cite[Appendix B]{Multi_res_llaws} for further details. Although \cite{Multi_res_llaws} concerns the case of Wigner matrices, where the self-energy operator $\mathcal{S}$ is given by the normalized trace rather than the sum of two traces as in \eqref{eq:def_S}, this does not affect the proof. The argument is based on the Gaussian renormalization technique (denoted by underline) analogous to the one used in the proof of Proposition \ref{prop:Cov}. Consider for definiteness the case where there is no star in \eqref{eq:global_3G}. Writing each observable simply by $B$ for brevity, we obtain from \eqref{eq:def_under_1} that
\begin{equation}
\begin{split}
\langle (G_1BG_2BG_1-M_{121})B\rangle &= \langle (G_2BG_1-M_{21})AB\rangle + \sigma \langle (G_1-M_1)E_\sigma\rangle \langle G_1BG_2BG_1AE_\sigma\rangle\\
& +\sigma \langle M_{12}E_\sigma\rangle \langle (G_2BG_1-M_{21})AE_\sigma\rangle+\sigma \langle M_{21}AE_\sigma\rangle \langle (G_1BG_2-M_{12})E_\sigma\rangle\\
&+ \sigma\langle (G_1-M_1)AE_\sigma\rangle\langle G_1BG_2BG_1E_\sigma\rangle - \langle A\underline{WG_1BG_2BG_1}\rangle,
\end{split}
\label{eq:global_exp}
\end{equation}
where $A:=(\mathcal{B}_{11}^{-1}[B^*])^*M_1$ (for more details see a similar calculation in \eqref{eq:G_exp}--\eqref{eq:Delta_exp}). In \eqref{eq:global_exp} we omitted superscripts of deterministic approximations as all of them are given by the appropriate number of $B$'s. Observe that $\|G_j\|\le \varepsilon^{-1}$ and that $\|\mathcal{B}_{11}^{-1}\|\lesssim \varepsilon^{-1}$ by \eqref{eq:eta_stab_bound} in the global regime (here we need this bound only for the one-body stability operator $\mathcal{B}_{11}$, but it holds for $\mathcal{B}_{12}$ as well), which implies that $\|A\|\lesssim 1$. Using the global law for one- and two-resolvent chains, and estimating longer chains by operator norm, we get from \eqref{eq:global_exp} that
\begin{equation}
\langle (G_1BG_2BG_1-M_{121})B\rangle = - \langle A\underline{WG_1BG_2BG_1}\rangle + \mathcal{O}_\prec (N^{-1}).
\label{eq:global_exp_simplified}
\end{equation}
Finally, using a minimalistic cumulant expansion as in \cite[Eq.(B.4)-(B.8)]{Multi_res_llaws}, we show that the first term in the rhs. of  \eqref{eq:global_exp_simplified} is stochastically dominated by $N^{-1}$ and conclude the proof of \eqref{eq:global_3G}. The proof of \eqref{eq:global_4G} follows the same strategy and relies on the global laws for shorter chains. 
\end{proof}

To propagate the global estimates \eqref{eq:global_2G}--\eqref{eq:global_4G} towards the real line, we consider the deterministic evolution of $z$ and $w$ along the \emph{characteristic flow}:
\begin{equation}
\frac{\dif}{\dif t} z_{t}=-\frac{1}{2} z_{t},\qquad\frac{\dif}{\dif t} w_{t} = -\frac{1}{2}w_{t} - \langle M^{z_t}(w_t)\rangle,
\label{eq:char_flow}
\end{equation} 
where $M^z(w)$ is defined in \eqref{eq:M}. While the precise initial conditions of this flow will be specified later, it is important to keep in mind that $(z_j,w_j)$  appearing in \eqref{eq:2G_av} is not the initial condition, but the \emph{target} of the characteristic flow, i.e. the value at the \emph{final time} $T$, which we choose to be a small fixed constant independent of $N$. An important feature of \eqref{eq:char_flow} is that $|\Im w_t|$ decreases with time, which allows us to propagate the local law estimates from the global regime $|\Im w|\gtrsim 1$ towards the real axis.

In fact, \eqref{eq:char_flow} is a special case of the characteristic flow for the deformed Wigner matrix model studied in \cite{eigenv_decorr}. Let $W^\#=(w^\#_{ab})_{a,b=1}^{2N}$ be a $2N\times 2N$ matrix from the \emph{Gaussian unitary ensemble} (GUE), i.e. $W^\#$ is Hermitian, its entries are independent up to the symmetry constraint, and for $a,b\in [2N]$ with $a\neq b$ (respectively, $a=b$) $w^\#_{ab}$ is a centered complex (respectively, real) Gaussian with variance $(2N)^{-1/2}$. One may relax the Gaussianity assumption on the entries of $W^\#$ in the way similar to Assumption \ref{ass:chi}(i), but since only the second order moment structure of $W^\#$ plays a role here, this generalization is not needed. For $z\in \C$, let $Z\in\C^{(2N)\times (2N)}$ be defined as in \eqref{eq:def_hermitization}. Then the matrix Dyson equation for the deformed Wigner matrix $H:=W^\# - Z$ is given by
\begin{equation}
-(M^\#(w))^{-1} = w+Z+\langle M^\#(w)\rangle,\quad \Im M^\#(w)\Im w>0,\quad w\in\C\setminus\R,
\label{eq:MDE_deformed}
\end{equation} 
see e.g. \cite[Eq.(1.4)]{eigenv_decorr}, where $M^\#(w)\in \C^{(2N)\times (2N)}$. We compare \eqref{eq:MDE_deformed} to the MDE for the Hermitization \eqref{eq:MDE}. Since all diagonal entries of $M^z(w)$ coincide by \eqref{eq:M}, we obtain from \eqref{eq:def_S} that $\mathcal{S}[M^z(w)]=\langle M^z(w)\rangle$. Therefore, \eqref{eq:MDE_deformed} and \eqref{eq:MDE} are identical, and by uniqueness of solutions to both equations we conclude that $M^\#(w)=M^z(w)$. In particular, the characteristic flows for these two models coincide. This allows us to borrow the characteristic flow analysis from \cite{eigenv_decorr}. However, we remark that the comparison between the deformed Wigner matrices with the deformation $-Z$ and the Hermitization model does not extend to the two-body analysis: the two-body stability operators for these models differ, since for deformed Wigner matrices the self-energy operator is given by the normalized trace, whereas in our Hermitization setting $\mathcal{S}$ is a sum of two traces (see \eqref{eq:def_S}).

To prove Proposition \ref{prop:2G_av} it is not sufficient to keep track of one trajectory for each spectral parameter, but instead we need to consider a family of trajectories. This leads us to the definition of the \emph{bulk-restricted spectral domains}, which is a special case of \cite[Definition 4.2]{eigenv_decorr} introduced in the context of deformed Wigner matrices.

\begin{definition}[Bulk-restricted spectral domains]\label{def:spec_dom} Fix (small) $\epsilon, \kappa>0$, the final time $T>0$ and $z_{T}\in \C$. Recall the definition of the bulk $\mathbf{B}^{z_T}_\kappa$ from \eqref{eq:def_bulk}. We define the \emph{bulk-restricted spectral domain} at time~$T$~by\footnote{In fact, the domain defined in \eqref{eq:def_domain_T} depends on $T$ only through $z_T$. However, we retain the $T$-dependence in this notation (and in $z_T$ itself) to keep it consistent with the notation for the time-dependent domain introduced in \eqref{eq:def_domain_t}.}
\begin{equation}
\Omega_{\kappa,\epsilon,T}^{z_{T}}:=\left\lbrace w\in\C\setminus\R\,:\, |\Im w|\ge \max \left\lbrace N^{\epsilon}(N\rho^{z_{T}})^{-1}, \mathrm{d}(\Re w, \mathbf{B}^{z_{T}}_\kappa)\right\rbrace\right\rbrace.
\label{eq:def_domain_T}
\end{equation}
Next, for $w_{T}\in \Omega_{\kappa,\epsilon,T}^{z_{T}}$ define the backward evolution operator associated to \eqref{eq:char_flow} by 
\begin{equation}
\mathfrak{F}_{t,T}^{z_{T}}(w_{T}):=w_{t},\quad \forall\, t\in [0,T].
\end{equation}
Finally, the bulk-restricted spectral domain at time $t\in [0,T)$ is given by
\begin{equation}
\Omega_{\kappa,\epsilon,t}^{z_{T}}:=\mathfrak{F}_{t,T}^{z_{T}}\left(\Omega_{\kappa,\epsilon,T}^{z_{T}}\right).
\label{eq:def_domain_t}
\end{equation}
\end{definition}

Note that for any $t\in [0,T]$ the upper index of $\Omega_{\kappa,\epsilon,t}^{z_{T}}$ always denotes the target value of the Hermitization parameter, that is $z_{T}$, rather than its value at the intermediate time $t$.  In the following lemma we collect the properties of the characteristic flow which will be later used in the proof of Proposition \ref{prop:2G_av}. The proof of  this lemma is elementary and thus omitted.

\begin{lemma}[Elementary properties of the characteristic flow]\label{lem:char_flow}
Fix $T>0$, $z_T\in \C$ and $w_T\in \C\setminus\R$. Let $z_t$, $w_t$ be the solution to \eqref{eq:char_flow}, $t\in [0,T]$, and denote $\eta_t:=|\Im w_t|$. Then the following holds. 
\begin{enumerate}
\item The map $t\mapsto \eta_t$ is monotone decreasing.
\item For any $t\in [0,T]$ we have $M^{z_t}(w_t)=\ee^{t/2}M^{z_0}(w_0)$.
\item The solution to \eqref{eq:char_flow} is explicitly given by
\begin{equation}
z_t = \ee^{-t/2}z_0,\quad w_t = \ee^{-t/2}w_0 -2m^{z_0}(w_0)\sinh \frac{t}{2},\quad \forall\, t\in [0,T].
\label{eq:char_flow_sol}
\end{equation}
\end{enumerate}
Assume additionally that $w_T\in \Omega_{\kappa,\epsilon,T}^{z_T}$ for some (small) $\kappa, \epsilon>0$.
\begin{enumerate}
\item[(4)] For any $s,t\in [0,T]$, $s\le t$, we have $\eta_s\sim \eta_t+|s-t|$.
\item[(5)] For any $a>1$ and $t\in [0,T]$ we have\footnote{A similar property holds also for general $w_T\in \C\setminus\R$ with the only modification that the integrands in \eqref{eq:eta_int} should be multiplied by $\rho_s:=\rho^{z_s}(w_s)$. For $w_T$ in the bulk-restricted spectral domain this additional factor is of order one whenever $\eta_s\lesssim 1$.} 
\begin{equation}
\int_0^t \frac{1}{(\eta_s\wedge 1)^a}\dif s\lesssim \frac{1}{(\eta_t\wedge 1)^{a-1}},\quad \int_0^t \frac{1}{\eta_s\wedge 1}\dif s\lesssim \log \frac{\eta_0\wedge 1}{\eta_t\wedge 1}.
\label{eq:eta_int}
\end{equation}
\end{enumerate}
\end{lemma}

To present the core ingredient for the proof of Proposition \ref{prop:2G_av}, which consists in propagation of the bounds \eqref{eq:global_2G}--\eqref{eq:global_4G} down to the real axis, we introduce some additional notation. We set
\begin{equation}
G_{j,t}:=G^{z_{j,t}}_t(w_{j,t}),\quad j=1,2,
\label{eq:def_Gt}
\end{equation} 
where $z_{j,t}$ and $w_{j,t}$ follow \eqref{eq:char_flow}. Here we use the convention that whenever the arguments of $G^z_t(w)$ are omitted, they are meant to be time-dependent. Using the same guiding principle, we denote the deterministic approximation of $G^{z_{1,t}}_t(w_1)B_1G^{z_{2,t}}_t(w_2)$ by $M_{12,t}^{B_1}(w_1,w_2)$, and abbreviate
\begin{equation}
M_{12,t}^{B_1}:=M_{12,t}^{B_1}(w_{1,t},w_{2,t}).
\end{equation}
The time-dependent deterministic approximations to the resolvent chains in \eqref{eq:chain_types} are defined in the same way. Finally, we denote by $\mathcal{B}_{12,t}(w_1,w_2)$ the 2-body stabilty operator \eqref{eq:def_B12} associated to $(z_{j,t},w_j)$, $j=1,2$, define $\widehat{\beta}_{12,t}(w_1,w_2)$ as in \eqref{eq:def_beta_hat} but with $\mathcal{B}_{12}(w_1^{(*)},w_2^{(*)})$ replaced by $\mathcal{B}_{12,t}(w_1^{(*)},w_2^{(*)})$, and set
\begin{equation}
\mathcal{B}_{12,t}:=\mathcal{B}_{12,t}(w_{1,t},w_{2,t}),\quad\widehat{\beta}_{12,t}:=\widehat{\beta}_{12,t}(w_{1,t},w_{2,t}),\quad\eta_{j,t}:=|\Im w_{j,t}|,\,j=1,2,\quad \eta_{*,t}:=\eta_{1,t}\wedge\eta_{2,t}\wedge 1.\label{eq:control_param_t}
\end{equation}

As usual, local laws come along with the bounds on the corresponding deterministic approximations. In the following statement we present bounds on the deterministic counterparts to the time-dependent resolvent chains.

\begin{proposition}[Bounds on the deterministic counterparts]\label{prop:M_bounds} Fix (small) $\delta, \epsilon, \kappa>0$ and the final time $T>0$. Let $\{\Omega_{\kappa,\epsilon,t}^{z_T}\}_{t\in[0,T]}$ be constructed as in Definition \ref{def:spec_dom} for $z_T\in\C$. It holds that
\begin{align}
\left\vert\left\langle M_{12,t}^{B_1}(w_{1,t},w_{2,t})\right\rangle\right\vert &\lesssim \frac{1}{\widehat{\beta}_{12,t}},\label{eq:M2_bound}\\
\left\vert\left\langle M_{121,t}^{B_1,B_2}\big(w_{1,t},w_{2,t},w_{1,t}^{(*)}\big)B_3\right\rangle\right\vert &\lesssim \frac{1}{\eta_{*,t}\widehat{\beta}_{12,t}},\label{eq:M3_bound}\\
\left\vert\left\langle M_{1212,t}^{B_1,B_2,B_3}\big(w_{1,t},w_{2,t},w_{1,t}^{(*)},w_{2,t}^{(*)}\big)B_4\right\rangle\right\vert&\lesssim \frac{1}{\eta_{*,t}(\widehat{\beta}_{12,t})^2},\label{eq:M4_bound}
\end{align}
uniformly in $t\in[0,T]$, $z_{j,T}\in \ee^{-T/2}(1-\delta)\mathbf{D}$, $w_{j,T}\in\Omega_{\kappa,\epsilon,T}^{z_T}$, $j=1,2$, and $B_i\in\mathrm{span}\{E_{\pm},F^{(*)}\}$, $i\in[4]$.
\end{proposition}

The proof of Proposition~\ref{prop:M_bounds} is presented in Section~\ref{sec:M_bounds}. Since the lhs. of \eqref{eq:M2_bound}--\eqref{eq:M4_bound} depend on~$t$ only through $z_{j,t}$ and $w_{j,t}$, there is no time-dependent nature in Proposition~\ref{prop:M_bounds}. Moreover, these bounds remain valid all the way down to the real line in the bulk regime, so the cut-off at the scale $N^{-1+\epsilon}$ in $\Omega_{\kappa,\epsilon,T}^{z_T}$ is not needed. However, we presented Proposition~\ref{prop:M_bounds} in this form to simplify the comparison with the time-dependent local laws, which we now state.

\begin{proposition}[Propagation of the local law bounds]\label{prop:Zig} Fix (small) $\delta, \epsilon, \kappa>0$. Let the final time $T>0$ be sufficiently small and let $\{\Omega_{\kappa,\epsilon,t}^{z_T}\}_{t\in [0,T]}$ be constructed as in Definition \ref{def:spec_dom} for $z_{T}\in \C$. Assume that for $t=0$ it holds that
\begin{align}
\left\vert\left\langle \left(G_{1,t}B_1G_{2,t}-M_{12,t}^{B_1}(w_{1,t},w_{2,t})\right)B_2\right\rangle\right\vert &\prec  \frac{1}{N\eta_{*,t}\widehat{\beta}_{12,t}},\label{eq:2G_bound_t}\\
\left\vert\left\langle \left(G_{1,t}B_1G_{2,t}B_2G_{1,t}^{(*)}-M_{121,t}^{B_1,B_2}\big(w_{1,t},w_{2,t},w_{1,t}^{(*)}\big)\right) B_3\right\rangle\right\vert &\prec \frac{1}{N\eta_{*,t}^3}\wedge\frac{1}{\sqrt{N\eta_{*,t}}\eta_{*,t}\widehat{\beta}_{12,t}},\label{eq:3G_bound_t}\\
\left\vert\left\langle \left(G_{1,t}B_1G_{2,t}B_2G_{1,t}^{(*)}B_3G_{2,t}^{(*)}-M_{1212,t}^{B_1,B_2,B_3}\big(w_{1,t},w_{2,t},w_{1,t}^{(*)}, w_{2,t}^{(*)}\big)\right) B_4\right\rangle\right\vert &\prec \frac{1}{N\eta_{*,t}^4}\wedge\frac{1}{\eta_{*,t}(\widehat{\beta}_{12,t})^2},\label{eq:4G_bound_t}
\end{align}
uniformly in $z_{j,T}\in \ee^{-T/2}(1-\delta)\mathbf{D}$, $w_{j,T}\in\Omega^{z_{j,T}}_{\kappa,\epsilon,T}$, $j=1,2$, and $B_i\in\{E_+,E_-,F,F^*\}$, $i\in [4]$. Then \eqref{eq:2G_bound_t}--\eqref{eq:4G_bound_t} hold uniformly in $t\in [0,T]$, $z_{j,T}$, $w_{j,T}$ and $B_i$ specified above.
\end{proposition}

Observe that \eqref{eq:2G_bound_t} improves upon the corresponding deterministic estimate \eqref{eq:M2_bound} by the optimal factor $(N\eta_*)^{-1}$. In contrast, the gain in \eqref{eq:3G_bound_t} compared to \eqref{eq:M3_bound} is suboptimal, amounting only to the factor $(N\eta_*)^{-1/2}$, while \eqref{eq:4G_bound_t} exhibits no improvement at all in some regimes, e.g. when $\widehat{\beta}_{12,t}\sim 1$ and $\eta_{*,t}\ll N^{-1/3}$. Nevertheless, this suboptimal estimates for three- and four-resolvent chains are sufficient to propagate \eqref{eq:2G_bound_t} in time, as we will show in Section \ref{sec:Zig_old}. We also remark that the bounds depending only on $\eta_{*,t}$ in the rhs. of \eqref{eq:3G_bound_t} and \eqref{eq:4G_bound_t} are needed for technical reasons and can be propagated without the more intricate estimates involving $\widehat{\beta}_{12,t}$, while the latter ones cannot be propagated on their own.

Given Propositions \ref{prop:global_law} and \ref{prop:Zig}, we now prove Proposition \ref{prop:2G_av}.

\begin{proof}[Proof of Proposition \ref{prop:2G_av} for $b=0$] Let $T_0>0$ be such a constant that Proposition \ref{prop:Zig} holds for any $T\le T_0$. Without loss of generality we may assume that $|\log(1-\mathfrak{s}^2)|\le T_0$, otherwise we take sufficiently small $\mathfrak{t}\in(0,\mathfrak{s})$ and observe that 
\begin{equation}
X=\sqrt{1-\mathfrak{s}^2}X_0+\mathfrak{s}\widetilde{X}\stackrel{d}{=}\sqrt{1-\mathfrak{t}^2}X_0'+\mathfrak{t}\widetilde{X}
\label{eq:smaller_G_comp}
\end{equation}
for some i.i.d. matrix $X_0'$ satisfying Assumption \ref{ass:chi}(i). Similarly we may also assume that $\ee^{T_0/2}(1-\delta)<(1-\delta/2)$. Take $T:=|\log(1-\mathfrak{s}^2)|$ and let $\kappa, \epsilon>0$ be as in Proposition \ref{prop:2G_av}. From Lemma \ref{lem:char_flow}(4) we have that 
\begin{equation}
\Omega^{z_T}_{\kappa,\epsilon,0}\subset \{w\in\C\,:\, |\Im w|\ge \varepsilon\},\quad \forall \, z_T\in (1-\delta)\mathbf{D},
\label{eq:domain_inclusion}
\end{equation}
for some $\varepsilon>0$. Therefore, \eqref{eq:2G_bound_t}--\eqref{eq:4G_bound_t} hold for $t=0$ by the global law from Proposition \ref{prop:global_law}. Finally, by Proposition \ref{prop:Zig} these estimates are valid for $t=T$, which in combination with \eqref{eq:OU_distr} completes the proof of Proposition \ref{prop:2G_av}.
\end{proof}

\subsection{Proof of Proposition \ref{prop:M_bounds}: bounds on the deterministic counterparts}\label{sec:M_bounds}

The bound \eqref{eq:M2_bound} on the two-resolvent deterministic approximation follows directly from \eqref{eq:M_12_beta_bound}, so we are left with \eqref{eq:M3_bound}--\eqref{eq:M4_bound}. First we establish these bounds in a weaker form, replacing $\widehat{\beta}_{12}$ by the smaller parameter $\eta_*$, and formulate them in the time-independent way. Then in the separate argument we improve these estimates by replacing back the appropriate number of $\eta_*$'s by $\widehat{\beta}_{12}$'s.

\begin{lemma}[Weak bounds on the deterministic counterparts]\label{lem:M_bounds_weak} Fix $\delta>0$. Uniformly in $z_j\in(1-\delta)\mathbf{D}$, $w_j\in\C\setminus\R$, $j=1,2$, and $B_i\in \mathrm{span}\{E_\pm,F^{(*)}\}$, $i\in [4]$ it holds that
\begin{align}
\left\vert\left\langle M_{121}^{B_1,B_2}\big(w_1,w_2,w_1^{(*)}\big)B_3\right\rangle\right\vert &\lesssim \frac{1}{\eta_*^2},\label{eq:M3_bound_weak}\\
\left\vert\left\langle M_{1212}^{B_1,B_2,B_3}\big(w_1,w_2,w_1^{(*)},w_2^{(*)}\big)B_4\right\rangle\right\vert&\lesssim \frac{1}{\eta_*^3},\label{eq:M4_bound_weak}
\end{align}
where $\eta_*:=|\Im w_1|\wedge|\Im w_2|\wedge 1$.
\end{lemma}

\begin{proof}[Proof of Lemma \ref{lem:M_bounds_weak}] 
The proof of \eqref{eq:M3_bound_weak}--\eqref{eq:M4_bound_weak} relies on the \emph{tensorization argument} (also known as the \emph{meta argument}), see e.g. \cite[Section~2.6]{Najim16} and \cite[Proof of Lemma~D.1]{non-herm_overlaps}. Although this argument is standard and well-known in the literature, we keep it here for the reader convenience, since it is frequently used in the proof of Proposition~\ref{prop:2G_av} and later in the proof of Proposition~\ref{prop:Cov}. Fix $N\in\N$, $z_j\in\C$ and $w_j\in\C\setminus\R$, $j=1,2$.  For $d\in\N$ let $X^{(d)}=(x^{(d)}_{ab})_{a,b=1}^{Nd}$ be an $Nd\times Nd$ i.i.d. matrix such that $x^{(d)}_{ab}\stackrel{d}{=} (Nd)^{-1/2}\chi$, where $\chi$ satisfies Assumption \ref{ass:chi}(i). Recall the definition of $Z_j\in \C^{(2N)\times (2N)}$ from \eqref{eq:def_hermitization} and denote
\begin{equation}
G_j^{(d)}:=\left(W^{(d)}-Z_j^{(d)}-w_j\right)^{-1},\quad \text{where}\quad Z_j^{(d)}:=Z_j\otimes I_d\in\C^{(2Nd)\times (2Nd)},\,\,  j=1,2.,
\end{equation} 
where $W^{(d)}$ is defined as in \eqref{eq:def_hermitization}, and $I_d$ is an operator acting identically on $\C^d$. As it immediately follows from this construction, the deterministic approximations to the products of $G_j^{(d)}$'s and observables from the set $\mathrm{span}\{E_{\pm},F^{(*)}\}$, which are viewed as $(2Nd)\times (2Nd)$ matrices, have $2\times 2$ block-constant structure, and only the size of blocks grows with $d$, while the scalars associated to these blocks are fixed. Thus from the global law \eqref{eq:global_3G} we have
\begin{equation}
\big\langle G_1^{(d)}B_1^{(d)}G_2^{(d)}B_2^{(d)}G_1^{(d)}B_3^{(d)}\big\rangle = \langle M_{121}^{B_1,B_2}B_3\rangle + \mathcal{O}\left(c(\eta_*)(Nd)^{-1}\right),\quad B_i^{(d)}:=B_i\otimes I_d,\, i\in[3],
\label{eq:meta_3G}
\end{equation}
where $c(\eta_*)$ is an implicit constant which depends on $\eta_*$, but not on $d$. For simplicity of presentation we further slightly abuse the notation and denote $B_i^{(d)}$ again by $B_i$. 

The Cauchy-Schwarz inequality along with the Ward identity imply
\begin{equation}
\left\vert\big\langle G_1^{(d)}B_1G_2^{(d)}B_2G_1^{(d)}B_3\big\rangle\right\vert \le \eta_*^{-1}\big\langle G_1^{(d)}B_1B_1^*(G_1^{(d)})^*B_3^*B_3\big\rangle^{1/2}\big\langle \Im G_1^{(d)} B_2^*\Im G_2^{(d)}B_2\big\rangle^{1/2}.
\label{eq:3G_Schwarz}
\end{equation}
Applying \eqref{eq:meta_3G} to the lhs. of \eqref{eq:3G_Schwarz}, using the same argument for both chains in the rhs. of \eqref{eq:3G_Schwarz} and taking the limit as $d$ goes to infinity, we upper bound $|\langle M_{121}^{B_1,B_2}B_3\rangle|$ by the product of two deterministic approximations to two-resolvent chains. By \eqref{eq:M2_bound} each of them has an upper bound of order $\eta_*^{-1}$. This finishes the proof of \eqref{eq:M3_bound_weak} with the choice, where there is no complex conjugation in the second $w_2$, while in the second case the proof is identical. The proof of \eqref{eq:M4_bound_weak} follows the same strategy, but instead of \eqref{eq:3G_Schwarz} makes use of
\begin{equation}
\left\vert \langle G_1B_1G_2B_2G_1B_3G_2B_4\rangle\right\vert \le \eta_*^{-1}\langle G_2B_1^*\Im G_1B_1G_2B_2B_2^*\rangle^{1/2}\langle G_2B_3^*\Im G_1B_3G_2B_4B_4^*\rangle^{1/2}, 
\label{eq:4G_Schwarz}
\end{equation} 
where we omitted index $(d)$ in $G_j$'s, and relies on \eqref{eq:M3_bound_weak} to bound the deterministic counterparts to the quantities in the rhs. of \eqref{eq:4G_Schwarz}.
\end{proof}

To prove \eqref{eq:M3_bound}--\eqref{eq:M4_bound}, we use the bounds \eqref{eq:M3_bound_weak}--\eqref{eq:M4_bound_weak} as an input and improve them using the characteristic flow method. For this purpose we differentiate the time-dependent deterministic approximations to the three- and four-resolvent chains along the characteristic flow. This idea of proving bounds on deterministic approximations via the characteristic flow was first used in \cite[Section~11.2]{RBM_Zigzag} in the context of random band matrices.

We remind the reader of the convention introduced above Lemma~\ref{lem:av_C_dif}: whenever index~$\sigma$ appears, it is meant to be summed over~$\sigma\in\{\pm\}$.

\begin{lemma}\label{lem:M_dif} Fix the final time $T>0$. For any $s\in [0,T]$ and $z_{j,T}\in\C$, $w_{j,T}\in\C\setminus\R$, $B_j\in\{E_+,E_-,F,F^*\}$ for $j\in[4]$, it holds that
\begin{align}
\partial_s \left\langle M_{123,s}^{B_1,B_2}B_3\right\rangle &= \frac{3}{2} \left\langle M_{123,s}^{B_1,B_2}B_3\right\rangle + \sum_{i=1}^3\sigma \left\langle M_{i,i+1,s}^{B_i}E_\sigma \right\rangle \left\langle M_{i+1,i+2,i+3,s}^{B_{i+1},B_{i+2}}E_\sigma\right\rangle,\label{eq:M3_dif}\\
\partial_s \left\langle M_{1234,s}^{B_1,B_2,B_3}B_4\right\rangle &= \frac{3}{2} \left\langle M_{1234,s}^{B_1,B_2,B_3}B_4\right\rangle + \sum_{i=1}^4\sigma \left\langle M_{i,i+1,s}^{B_i}E_\sigma \right\rangle \left\langle M_{i+1,i+2,i+3,i+4,s}^{B_{i+1},B_{i+2},B_{i+3}}E_\sigma\right\rangle\label{eq:M4_dif}\\
&+\sigma \left\langle M_{123,s}^{B_1,B_2}E_\sigma\right\rangle \left\langle M_{341,s}^{B_3,B_4}E_\sigma\right\rangle + \sigma \left\langle M_{234,s}^{B_2,B_3}E_\sigma\right\rangle \left\langle M_{412,s}^{B_4,B_1}E_\sigma\right\rangle.\nonumber
\end{align}
In \eqref{eq:M3_dif} we consider indices modulo 3 and in \eqref{eq:M4_dif} modulo 4.
\end{lemma}

In Lemma \ref{lem:M_dif} we considered the evolution of four pairs $(z_j,w_j)$ along the characteristic flow to make the formulation more transparent. In all applications below we will take $z_1=z_3$ and $z_2=z_4$. The proof of Lemma \ref{lem:M_dif} follows by a standard application of the meta-argument and the corresponding differential equations for the products of resolvents, where the fluctuation terms analogous to the lhs. of \eqref{eq:global_2G}--\eqref{eq:global_4G} are negligible as the tensorization parameter $d$ goes to infinity, see e.g. the proof of \cite[Lemma 4.8]{cipolloni2023eigenstate}.

As a final preparation for the proof of Proposition~\ref{prop:M_bounds}, we formulate several properties of $\widehat{\beta}_{12}$ analogous to the admissibility properties of the control parameter from \cite[Definition~4.4]{eigenv_decorr}. Our description of $\widehat{\beta}_{12}$ is slightly more extensive than what is required for the proof of Proposition~\ref{prop:M_bounds}, however the additional properties will be later used in the proof of Proposition~\ref{prop:Zig}, as it is explained in more detail below Lemma~\ref{lem:beta_hat}. We now remind the reader the definition of an auxiliary family of quantities introduced in \eqref{eq:def_beta_pm_main}. These quantities are used to analyze $\widehat{\beta}_{12}$, and will also naturally emerge on their own in the proof of Proposition~\ref{prop:Zig}:
\begin{equation}
\begin{split}
\beta_{12,\pm} =\beta_{12,\pm}(w_1,w_2) :&=1-\Re \left[z_1\overline{z}_2\right] u_1u_2\pm \sqrt{m_1^2m_2^2 - \left(\Im \left[z_1\overline{z}_2\right]\right)^2u_1^2u_2^2},\\
\beta_{12,*}=\beta_{12,*}(w_1,w_2):&=\min\{|\beta_{12,+}|,|\beta_{12,-}|\},
\end{split}
\label{eq:def_beta_pm}
\end{equation}
for any $z_j\in \C$ and $w_j\in \C\setminus\R$, where $m_j:=m^{z_j}(w_j)$ and $u_j:=u^{z_j}(w_j)$ are given by \eqref{eq:M}, \eqref{eq:m}. Note that
\begin{equation*}
\beta_{12,\pm}(w_1,w_2)=\beta_{21,\pm}(w_2,w_1)\quad\text{and}\quad \beta_{12,*}(w_1,w_2)=\beta_{21,*}(w_2,w_1).
\end{equation*}
We further define the time-dependent parameters $\beta_{12,\pm,t}$ and $\beta_{12,*,t}$ by replacing $m_j$ and $u_j$ in \eqref{eq:def_beta_pm} by $m_{j,t}:=m^{z_{j,t}}(w_{j,t})$ and $u_{j,t}:=u^{z_{j,t}}(w_{j,t})$, respectively.

\begin{lemma}[Properties of $\beta_{12,*}$ and $\widehat{\beta}_{12}$]\label{lem:beta_hat} Fix a (small) $\delta>0$.
\begin{enumerate}
\item \emph{[Relation between $\beta_{12,*}$ and $\widehat{\beta}_{12}$]} Uniformly in $z_j\in(1-\delta)\mathbf{D}$ and $w_j\in\C\setminus\R$, $j=1,2$, it holds that
\begin{equation}
\widehat{\beta}_{12}(w_1,w_2)\sim\min\left\lbrace\beta_{12,*}\big(w_1^{(*)}, w_2^{(*)}\big)\wedge 1\right\rbrace. 
\label{eq:beta_hat_to_beta}
\end{equation}
\item \emph{[Monotonicity in time]} Fix additionally (small) $\epsilon, \kappa>0$ and the final time $T>0$. Uniformly in $z_{j,T}\in \ee^{-T/2}(1-\delta)\mathbf{D}$ and $w_{j,T}\in \Omega^{z_{j,T}}_{\kappa,\epsilon,T}$ we have
\begin{equation}
\beta_{12,*,s}\sim \beta_{12,*,t}+|t-s|,\quad \widehat{\beta}_{12,s}\sim\widehat{\beta}_{12,t}+|t-s|,\quad\forall\, 0\le s\le t\le T.
\label{eq:gamma_monot}
\end{equation}
\item \emph{[Lipschitz continuity in space]} Fix additionally a (small) $\varepsilon>0$ and for any $z\in\C$ denote
\begin{equation}
\mathcal{D}^z_{\kappa,\varepsilon}:=\{w\in\C\setminus\R\,:\, \Re w\in \mathbf{B}^z_\kappa\}\cup\{w\in\C\setminus\R\,:\, |\Im w|\ge \varepsilon\}.
\end{equation}
Uniformly in $z_j\in (1-\delta)\mathbf{D}$, $w_j\in \mathcal{D}^{z_j}_{\kappa,\varepsilon}$ for $j=1,2$, and $w_2'\in\mathcal{D}^{z_2}_{\kappa,\varepsilon}$ it holds that
\begin{equation}
\widehat{\beta}_{12}(w_1,w_2) \lesssim \widehat{\beta}_{12}(w_1,w_2') + |w_2-w_2'|.
\label{eq:beta_hat_perturb}
\end{equation}
\item \emph{[Vague monotonicity in imaginary part]} Uniformly in $z_j\in (1-\delta)\mathbf{D}$, $w_j\in \mathcal{D}^{z_j}_{\kappa,\varepsilon}$, $j=1,2$, such that $\Im w_2>0$, and $x\ge 0$ it holds that
\begin{equation}
\widehat{\beta}_{12}(w_1,w_2)\lesssim \widehat{\beta}_{12}(w_1,w_2+\ii x).
\label{eq:beta_hat_vert_line}
\end{equation}
\end{enumerate}
\end{lemma}

The proof of Lemma \ref{lem:beta_hat} is presented in Section~\ref{app:beta_hat}. It is structured in such a way that first we establish \eqref{eq:beta_hat_to_beta} and then use this relation between $\widehat{\beta}_{12}$ and $\beta_{12,*}$ to transfer the properties of $\beta_{12,*}$, which follow from the explicit formula \eqref{eq:def_beta_pm}, to the less explicit quantity $\widehat{\beta}_{12}$. In particular, not only \eqref{eq:gamma_monot} holds both for $\widehat{\beta}_{12}$ and $\beta_{12,*}$, as currently stated, but also \eqref{eq:beta_hat_perturb} and \eqref{eq:beta_hat_vert_line} hold in this generality, although the $\beta_{12,*}$ counterparts of these two statements are not needed for the proof of Proposition~\ref{prop:2G_av}. Meanwhile, \eqref{eq:beta_hat_to_beta}, the second part of \eqref{eq:gamma_monot}, and \eqref{eq:beta_hat_perturb} are used in the proof of Proposition~\ref{prop:M_bounds}, whereas the remaining statements in Lemma~\ref{lem:beta_hat} are needed for the proof of Proposition~\ref{prop:Zig}.

Now we are ready to prove Proposition~\ref{prop:M_bounds}. For definiteness we consider only the case when there are no stars in \eqref{eq:M3_bound}--\eqref{eq:M4_bound}, while the argument is identical in the rest of the cases.

\medskip

\noindent\underline{Proof of \eqref{eq:M3_bound}.} First we prove \eqref{eq:M3_bound} in the case when $B_3\in\{E_+,E_-\}$. From \eqref{eq:char_flow_sol} and Definition~\ref{def:spec_dom} we have that there exist $\kappa',\varepsilon>0$ such that
\begin{equation}
\Omega^{z_T}_{\kappa,\epsilon,t}\subset \{w\in \C\setminus\R\,:\, \Re w\in \mathbf{B}^{z_t}_{\kappa'}\}\cup\{w\in \C\setminus\R\,:\, |\Im w|\ge \varepsilon\},\quad \forall z_T\in (1-\delta)\mathbf{D},\,\forall\, t\in [0,T].
\label{eq:new_domain}
\end{equation}
We use the time-independent formulation and consider any $z_j\in (1-\delta)\mathbf{D}$ and $w_j$ in the ambient domain in the rhs. of \eqref{eq:new_domain}. For $B_3=E_-$ we use that $G^{z_1}(w_1)E_-=-E_-G^{z_1}(-w_1)$ due to \eqref{eq:eigenvectors_symmetry}, which in combination with the resolvent identity gives that
\begin{equation}
\begin{split}
&\left\langle G^{z_1}(w_1)B_1G^{z_2}(w_2)B_2G^{z_1}(w_1)E_-\right\rangle = \frac{1}{2w_1}\left\langle \left(G^{z_1}(-w_1)-G^{z_1}(w_1)\right)B_1G^{z_2}(w_2)B_2E_-\right\rangle \\
&\quad = -\frac{1}{2w_1}\left(\left\langle G^{z_1}(w_1)E_-B_1G^{z_2}(w_2)B_2\right\rangle + \left\langle G^{z_1}(w_1)B_1G^{z_2}(w_2) B_2E_-\right\rangle\right).
\end{split}
\end{equation}
Together with the meta-argument and the bound on the deterministic approximation to the two-resolvent chain \eqref{eq:M2_bound} which is already established, this implies \eqref{eq:M3_bound} for $B_3=E_-$. For $B_3=E_+$ the argument is similar, but instead of the resolvent identity we employ the integral representation of $G_1^2$:
\begin{equation}
\left\langle (G^{z_1}(w_1))^2 B_1G^{z_2}(w_2)B_2\right\rangle = \frac{1}{2\pi\ii} \oint_{\mathscr{C}} \frac{1}{(\zeta-w_1)^2}\left\langle G^{z_1}(\zeta)B_1G^{z_2}(w_2)B_2\right\rangle\dif \zeta, 
\end{equation}
where $\mathscr{C}$ is a circle with the center at $w_1$ and radius $\eta_1/2$, and use that 
\begin{equation}
\widehat{\beta}_{12}(\zeta,w_2)\sim \widehat{\beta}_{12}(w_1,w_2),\quad \forall \zeta\in\mathscr{C},
\end{equation}
by \eqref{eq:eta_stab_bound} and \eqref{eq:beta_hat_perturb}.

Now we establish \eqref{eq:M3_bound} in the general case $B_3\in \mathrm{span}\{E_\pm, F^{(*)}\}$ in its time-dependent form. Multiplying \eqref{eq:M3_dif} by $\ee^{-3s/2}$ and integrating it in $s$ from 0 to $t$ we get
\begin{equation}
\begin{split}
\ee^{-3t/2}\left\langle M_{121,t}^{B_1,B_2}B_3\right\rangle &= \left\langle M_{121,0}^{B_1,B_2}B_3\right\rangle + \int_0^t \ee^{-3s/2} \sigma\left\langle M_{11,s}^{B_3}E_\sigma\right\rangle \left\langle M_{121,s}^{B_1,B_2}E_\sigma\right\rangle\dif s\\
&+\int_0^t \ee^{-3s/2} \sigma\left(\left\langle M_{12,s}^{B_1}E_\sigma\right\rangle \left\langle M_{121,s}^{E_\sigma,B_2}B_3\right\rangle+\left\langle M_{21,s}^{B_2}E_\sigma\right\rangle \left\langle M_{121,s}^{B_1,E_\sigma}B_3\right\rangle\right)\dif s,
\end{split}
\label{eq:M3_int}
\end{equation}
for any $t\in [0,T]$. To estimate the integral in the first line of \eqref{eq:M3_int}, we bound the $M_{11}$ term by $\eta_{*,s}^{-1}$ by means of \eqref{eq:M_12_beta_bound} and \eqref{eq:eta_stab_bound}, and the $M_{121}$ term by $(\eta_{*,s}\widehat{\beta}_{12,s})^{-1}$, since it corresponds to the already treated choice $B_3=E_\pm$. Performing the integration with the help of \eqref{eq:eta_int} and using that by \eqref{eq:gamma_monot} $\widehat{\beta}_{12,t}\lesssim \widehat{\beta}_{12,s}$ for any $s\in [0,t]$, we get
\begin{equation}
\int_0^t \left\vert\left\langle M_{11,s}^{B_3}E_\sigma\right\rangle \left\langle M_{121,s}^{B_1,B_2}E_\sigma\right\rangle\right\vert\dif s\lesssim \int_0^t \frac{1}{\eta_{*,s}}\frac{1}{\eta_{*,s}\widehat{\beta}_{12,s}}\dif s \lesssim \frac{1}{\eta_{*,t}\widehat{\beta}_{12,t}}.
\label{eq:M3_int_aux}
\end{equation}
For the integral in the second line of \eqref{eq:M3_int} the resulting bound is the same. In this case, however, we estimate the $M_{12}$ term by $\widehat{\beta}_{12,s}^{-1}$ using \eqref{eq:M2_bound} and the $M_{121}$ term by $\eta_{*,s}^{-2}$ using \eqref{eq:M3_bound_weak}. Finally, observing that the first term in the rhs. of \eqref{eq:M3_int} is bounded by order one by \eqref{eq:M3_bound_weak} and that $\eta_{*,0}\sim 1$, we complete the proof of \eqref{eq:M3_bound}.

\medskip

\noindent\underline{Proof of \eqref{eq:M4_bound}.} We prove \eqref{eq:M4_bound} in two steps. First we follow an argument similar to \eqref{eq:M3_int}--\eqref{eq:M3_int_aux} using the weak bound \eqref{eq:M4_bound_weak} on the four-resolvent deterministic approximation as an input, and improve \eqref{eq:M4_bound_weak} by a single $\eta_{*,t}/\widehat{\beta}_{12,t}$ factor. Next we use this intermediate estimate as a new input and partially rerun the previous part of the argument, gaining an additional factor $\eta_{*,t}/\widehat{\beta}_{12,t}$ and thereby proving \eqref{eq:M4_bound}.

As in \eqref{eq:M3_int}, we multiply \eqref{eq:M4_dif} by $\ee^{-2s}$ and integrate the resulting equation in $s$ from 0 to $t$. We treat the integral of the sum of four terms in the first line of \eqref{eq:M4_dif} analogously to the integral in the second line of \eqref{eq:M3_int} using \eqref{eq:M2_bound} and \eqref{eq:M4_bound_weak}. This gives an upper bound of order $(\eta_{*,t}^2\widehat{\beta}_{12,t})^{-1}$ on this integral. The integral of the first term in the second line of \eqref{eq:M4_dif} is estimated by
\begin{equation}
\int_0^t \left\vert \left\langle M_{121,s}^{B_1,B_2}E_\sigma\right\rangle \left\langle M_{121,s}^{B_3,B_4}E_\sigma\right\rangle\right\vert\dif s\lesssim \int_0^t \left(\frac{1}{\eta_{*,s}\widehat{\beta}_{12,s}}\right)^2\dif s\lesssim\frac{1}{\eta_{*,t}(\widehat{\beta}_{12,t})^2},
\end{equation}
where we used \eqref{eq:M3_bound}. For the second term in the second line of \eqref{eq:M4_dif} the bound is identical. Thus, we get 
\begin{equation}
\left\vert\left\langle M_{1212,t}^{B_1,B_2,B_3}\big(w_{1,t},w_{2,t},w_{1,t}^{(*)},w_{2,t}^{(*)}\big)B_4\right\rangle\right\vert\lesssim \frac{1}{\eta_{*,t}^2\widehat{\beta}_{12,t}}.\label{eq:M4_bound_subopt}
\end{equation}
To gain an additional $\eta_{*,t}/\widehat{\beta}_{12,t}$ improvement in the rhs. of \eqref{eq:M4_bound_subopt}, we return to the analysis on the integral of \eqref{eq:M4_dif}, and bound the sum of four terms in the first line of \eqref{eq:M4_dif} by \eqref{eq:M4_bound_subopt} instead of \eqref{eq:M4_bound_weak}. Keeping the rest of the analysis unchanged, we conclude the desired $(\eta_{*,t})^{-1}(\widehat{\beta}_{12,t})^{-2}$ bound on the lhs. of \eqref{eq:M4_bound_subopt}. This finishes the proof of Proposition \ref{prop:M_bounds}. \qed

\subsection{Proof of Proposition \ref{prop:Zig}}\label{sec:Zig_old} In this section we prove Proposition~\ref{prop:Zig} in several steps. We begin with the case when all observables in Proposition \ref{prop:Zig} are equal to $E_\pm$ and show that the bounds \eqref{eq:2G_bound_t}--\eqref{eq:4G_bound_t}, specialized to this setting, propagate in time. To this end, we differentiate the lhs. of \eqref{eq:2G_bound_t}--\eqref{eq:4G_bound_t} in time and analyze the structure of the resulting \emph{zig equations} in Section~\ref{sec:Zig_structure}. In Section~\ref{sec:Gronwall} we derive the desired bounds \eqref{eq:2G_bound_t}--\eqref{eq:4G_bound_t} from these equations using a Gr{\"o}nwall-type argument, deferring the necessary estimates on the error terms to Section~\ref{sec:err_bounds}. Finally, in Section~\ref{sec:general_B} we show how to adapt the argument presented in Sections~\ref{sec:Zig_structure}--\ref{sec:err_bounds} to the case when some of the observables are equal to~$F^{(*)}$.

\subsubsection{Structure of the zig equations}\label{sec:Zig_structure}

In this section we derive \eqref{eq:Y_syst_intro} and analyze the structure of the resulting system of SDEs. To make the presentation more transparent, we consider the general case of a $k$-resolvent chain, where $k\ge 2$ is an $N$-independent integer, and follow the evolution of $k$ resolvents $G_{1,t},\ldots, G_{k,t}$ defined as in \eqref{eq:def_Gt}. Later in the proof of  Proposition \ref{prop:Zig} we will need only the cases $k\in\{2,3,4\}$, but in Section \ref{sec:Zig_structure} we follow the most general set-up.

For any observables $B_i\in \mathrm{span}\{E_\pm,F^{(*)}\}$ and $t\in [0,T]$ we get from the It\^{o} calculus that
\begin{equation}
\begin{split}
\dif \left\langle G_{1,t}B_1G_{2,t}B_2\cdots G_{k,t}B_k\right\rangle &= \frac{k}{2} \left\langle G_{1,t}B_1G_{2,t}B_2\cdots G_{k,t}B_k\right\rangle\dif t+\dif\mathfrak{E}_{[k],t}(B_1,\ldots,B_k)\\
& + \sum_{1\le i<j\le k} \sigma \left\langle G_{i,t}B_i\cdots G_{j,t}E_\sigma\right\rangle\left\langle G_{j,t}B_j\cdots G_{i,t} E_\sigma\right\rangle\dif t\\
&+ \sum_{j=1}^k \sigma \left\langle (G_{j,t}-M_{j,t})E_\sigma\right\rangle\left\langle G_{1,t}B_1\cdots B_{j-1}G_{j,t}E_\sigma G_{j,t} B_j\cdots G_kB_k\right\rangle\dif t,
\end{split}
\label{eq:kG_dif}
\end{equation}
where the \emph{martingale term} $\dif \mathfrak{E}_{[k],t}$ is given by
\begin{equation}
\dif\mathfrak{E}_{[k],t}(B_1,\ldots,B_k):= \sum_{a,b=1}^{2N}\partial_{ab}\left\langle G_{1,t}B_1\cdots G_{k,t}B_k\right\rangle \frac{\dif B_{ab,t}}{\sqrt{N}}.
\label{eq:def_E}
\end{equation}
In \eqref{eq:kG_syst} we adopted the same notational convention as in Lemma \ref{lem:M_dif} that the free index $\sigma$ is meant to be summed over $\{\pm\}$. In the second line of \eqref{eq:kG_dif} the resolvents are indexed by $[i,j]$ in the first trace and by $[j,i]:=\{j,j+1,\ldots,k,1,\ldots,i\}$ in the second. Applying the meta-argument to \eqref{eq:kG_dif} in the same manner as in the proof of Lemma \ref{lem:M_bounds_weak}, we get the following generalization of Lemma \ref{lem:M_dif} for the deterministic approximation to the chain of arbitrary length:
\begin{equation}
\frac{\dif}{\dif t} \left\langle M_{[k],t}^{B_1,\ldots,B_{k-1}}B_k\right\rangle = \frac{k}{2} \left\langle M_{[k],t}^{B_1,\ldots,B_{k-1}}B_k\right\rangle + \sum_{1\le i<j\le k} \sigma \left\langle M_{[i,j],t}^{B_i,\ldots,B_{j-1}}E_\sigma\right\rangle \left\langle M_{[j,i],t}^{B_j,\ldots,B_{i-1}}E_\sigma\right\rangle.
\label{eq:kM_dif}
\end{equation}

To follow the time evolution of a fluctuation of a resolvent chain around its deterministic approximation, we introduce the following notation. For any $n\in \N$, set of indices $\mathcal{J}:=\{j_1,\ldots,j_n\}\in [k]^n$, observables $B_i\in \mathrm{span}\{E_\pm,F^{(*)}\}$, $i\in [n]$, and time $t\in[0,T]$, denote
\begin{equation}
Y_{\mathcal{J},t}(B_1,\ldots,B_n):=\left\langle \left(G_{j_1,t}B_1G_{j_2,t}\cdots G_{j_n,t} - M_{\mathcal{J},t}^{B_1,\ldots,B_{n-1}}\right)B_n\right\rangle.
\label{eq:def_Y_k}
\end{equation}
Subtracting \eqref{eq:kM_dif} from \eqref{eq:kG_dif} we get
\begin{equation}
\begin{split}
\dif Y_{[k],t}(B_1,\ldots,B_k)& = \frac{k}{2}Y_{[k],t}(B_1,\ldots,B_k) \dif t\\
&+ \sum_{j=1}^k \sigma \left\langle M_{[j,j+1],t}^{B_j}E_\sigma\right\rangle Y_{[k],t}(B_1,\ldots, B_{j-1},E_\sigma, B_{j+1},\ldots, B_k)\dif t\\
&+\mathfrak{F}_{[k],t}(B_1,\ldots,B_k)\dif t +\dif \mathfrak{E}_{[k],t}(B_1,\ldots,B_k),
\end{split}\label{eq:kY_dif}
\end{equation}
where in the second line of \eqref{eq:kY_dif} we denoted $[k,k+1]:=\{k,1\}$. Here the \emph{forcing term} $\mathfrak{F}_{[k],t}$ is given by
\begin{equation}
\begin{split}
&\mathfrak{F}_{[k],t}(B_1,\ldots,B_k):=\sum_{j=1}^k \sigma Y_{j,t}(E_\sigma)\left\langle G_{1,t}B_1\cdots B_{j-1}G_{j,t}E_\sigma G_{j,t} B_j\cdots G_kB_k\right\rangle\\
&+\sum_{1\le i<j\le k} Y_{[i,j],t}(B_i,\ldots,B_{j-1},E_\sigma)Y_{[j,i],t} (B_j,\ldots,B_{i-1},E_\sigma)\\
&+\sum_{\substack{1\le i<j\le k\\j\neq i+1}}\!\!\!\sigma\left( \!\left\langle\! M_{[i,j],t}^{B_i,\ldots,B_{j-1}}E_\sigma\!\right\rangle\! Y_{[j,i],t}(B_j,\ldots,B_{i-1},E_\sigma)\!+\! Y_{[i,j],t}(B_i,\ldots,B_{j-1},E_\sigma)\!\left\langle\! M_{[j,i],t}^{B_j,\ldots,B_{i-1}}E_\sigma\!\right\rangle\! \right)\!.
\end{split}
\label{eq:def_Fk}
\end{equation}
In the last line of \eqref{eq:def_Fk} the constraint $j\neq i+1$ is interpreted modulo $k$, i.e the pair $(i,j)=(k,1)$ is also removed from the summation. In particular, for $k=2$ no terms appear in the third line of \eqref{eq:def_Fk}. We remark that the terms in the first (respectively, second and third) line of \eqref{eq:def_Fk} originate from the third (respectively, second) line of \eqref{eq:kG_dif}. The terms with $j=i+1$ are not included into the third line of \eqref{eq:def_Fk}, instead, they are separated in the second line of \eqref{eq:kY_dif}, as they correspond to the fluctuations of the $k$-resolvent chains, while the terms in the third line of \eqref{eq:def_Fk} arise from shorter chains. Moreover, fluctuations of chains of length at least $k$ (namely, of length $k$ and $k+1$) appear in the forcing term only when multiplied by a fluctuation of another chain, which makes these terms small (for more details see the proof of Proposition~\ref{prop:err_bounds} later). For this reason we treat $\mathfrak{F}_{[k],t}$ as an error term, and refer to the terms in the first two lines in the rhs. of \eqref{eq:kY_dif} as the \emph{linear terms}.

Now we specialize \eqref{eq:kY_dif} to the case when $B_i=E_\pm$ for all $i\in [k]$. We identify the space of diagonal $2\times 2$ matrices with $\C^2$ by mapping $E_+$ and $E_-$ to the standard basis vectors in $\C^2$, which we denote by $\bm{e}_+$ and~$\bm{e}_-$, respectively. We further identify the space of $k$-tuples $(B_1,\ldots, B_k)$ of diagonal $2\times 2$ matrices with the tensor product $(\C^2)^{\otimes k}$, by mapping $(E_{\sigma_1},\ldots,E_{\sigma_k})$ to
\begin{equation}
\bm{e}_{\bm{\sigma}}=\bm{e}_{\sigma_1,\ldots,\sigma_k}:=\bm{e}_{\sigma_1}\otimes \cdots \otimes \bm{e}_{\sigma_k},\quad \text{where}\quad \bm{\sigma}:=(\sigma_1,\ldots,\sigma_k)\in\{\pm\}^k.
\end{equation}
Note that $\{\bm{e}_{\bm{\sigma}}\}_{\bm{\sigma}\in\{\pm\}^k}$ is an orthonormal basis in $(\C^2)^{\otimes k}$. Define
\begin{equation}
\mathcal{Y}_{[k],t}:=\sum_{\sigma_1,\ldots, \sigma_k\in\{\pm\}} Y_{[k],t}(E_{\sigma_1},\ldots, E_{\sigma_k})\bm{e}_{\sigma_1,\ldots,\sigma_k}\in (\C^2)^{\otimes k},\quad t\in [0,T],
\label{eq:def_Y_vector}
\end{equation}
the \emph{martingale term} $\dif \mathcal{E}_{[k],t}$ and the \emph{forcing term} $\mathcal{F}_{[k],t}$ are defined similarly to \eqref{eq:def_Y_vector} from \eqref{eq:def_E} and \eqref{eq:def_Fk}, respectively. For $\sigma\in\{\pm\}$, $i,j\in[k]$ and $t\in [0,T]$ denote further
\begin{equation}
\begin{split}
a_{ij}^{\sigma}=a_{ij}^{\sigma}(z_i,z_j;w_i,w_j)&:=\sigma \langle M_{ij}^{E_\sigma}E_\sigma\rangle,\quad\qquad\qquad\quad\,\,\,\,\,\,\, a_{ij,t}^\sigma:=a^\sigma_{ij} (z_{i,t},z_{j,t};w_{i,t},w_{j,t}),\\
d_{ij}=d_{ij}(z_i,z_j;w_i,w_j)&:=\langle M_{ij}^{E_+}E_-\rangle = -\langle M_{ij}^{E_-}E_+\rangle, \quad d_{ij,t}:=d_{ij}(z_{i,t},z_{j,t};w_{i,t},w_{j,t}).
\label{eq:propag_aux_notation}
\end{split}
\end{equation}
Here $\sigma\in\{\pm\}$ is fixed and we do not perform summation over $\sigma$. The identity in the last line of \eqref{eq:propag_aux_notation} is established for $w_i, w_j$ on the imaginary axes in \cite[Eq.(5.26)]{univ_extr}. Since the functions $\langle M_{ij}^{E_+}E_-\rangle$ and $-\langle M_{ij}^{E_-}E_+\rangle$ are analytic in $w_i,w_j$ in the upper and lower complex half-planes, this identity holds for all $w_i,w_j\in\C\setminus\R$. Finally, we define two families of linear operators acting on $(\C^2)^{\otimes k}$: the projections onto the basis vectors $\{P_{\bm{\sigma}}\}_{\bm{\sigma}\in \{\pm\}^k}$ and the \emph{spin-flip operators} $\{S_j\}_{j=1}^k$, where $S_j$ flips the sign of $j$-th spin, i.e.
\begin{equation}
P_{\bm{\sigma}}\bm{e}_{\bm{\omega}}:=\delta_{\bm{\sigma},\bm{\omega}}\bm{e}_{\bm{\sigma}}\quad\text{and}\quad S_j \bm{e}_{\omega_1,\ldots,\omega_j,\ldots,\omega_k} := \bm{e}_{\omega_1,\ldots,-\omega_i,\ldots,\omega_k},\quad \forall \bm{\sigma},\bm{\omega}\in \{\pm\}^k,\, j\in [k]. 
\end{equation}

Considering in \eqref{eq:kY_dif} all $2^k$ possible choices of observables $B_i=E_\pm$, $i\in [k]$, we get
\begin{equation}
\dif \mathcal{Y}_{[k],t}= \left(\frac{k}{2}I + \mathcal{A}_{[k],t}\right)\mathcal{Y}_{[k],t}\dif t + \mathcal{F}_{[k],t}\dif t +\dif\mathcal{E}_{[k],t},
\label{eq:kG_syst}
\end{equation}
where the \emph{generator} $\mathcal{A}_{[k],t}$ is given by
\begin{equation}
\mathcal{A}_{[k],t}: = \sum_{\sigma_1,\ldots,\sigma_k\in\{\pm\}} \bigg(\sum_{j=1}^k a_{j(j+1),t}^{\sigma_j}\bigg)P_{\sigma_1,\ldots,\sigma_k} - \sum_{j=1}^k d_{j(j+1),t}S_j, 
\label{eq:def_general_A}
\end{equation}
with the standard convention that the subscripts of $a$'s and $d$'s are viewed modulo $k$. Similarly to the discussion below the scalar equation \eqref{eq:kG_dif}, we view the forcing term $\mathcal{F}_{[k],t}\dif t$ and the martingale term $\dif\mathcal{E}_{[k],t}$ in \eqref{eq:kG_syst} as error terms and refer to the remaining first term in the rhs. of \eqref{eq:kG_syst} as to the \emph{linear term}. The terminology "error terms" is introduced purely for convenience: these terms are not negligible, but for the analysis of \eqref{eq:kG_syst} only upper bounds on them are needed, while the linear term will be treated with higher precision. We turn \eqref{eq:kG_syst} into an upper bound on $\|\mathcal{Y}_{[k],t}\|$ using the Gr{\" o}nwall inequality. To do so, one needs a control on the operator norm of the \emph{propagator} $\mathcal{P}_{[k],t,s}:(\C^2)^{\otimes k}\to (\C^2)^{\otimes k}$ for all $0\le s\le t\le T$, where for any $\bm{x}\in (\C^2)^{\otimes k}$, $\mathcal{P}_{[k],t,s}\bm{x}$ is defined by
\begin{equation}
\mathcal{P}_{[k],t,s} \bm{x}:=\bm{x}_t,\quad \text{with}\quad \frac{\dif}{\dif r} \bm{x}_r=\mathcal{A}_{[k],r}\bm{x}_r,\,\,\forall r\in [s,t],\,\, \bm{x}_s:=\bm{x}.
\end{equation}
Trivially, we have
\begin{equation}
\left\lVert \mathcal{P}_{[k],t,s}\right\rVert\le \exp\left\lbrace \int_s^t f_{[k],r}\dif r\right\rbrace,\quad \text{where}\quad f_{[k],r}:=\left[ \max \mathrm{spec}(\Re \mathcal{A}_{[k],r})\right]_+,\,\,\,\forall\, 0\le s\le r\le t\le T.
\label{eq:propagator}
\end{equation}
Here $[\cdot]_+$ stands for the positive part of a real number.

In the following statement which generalizes Lemma~\ref{lem:propag_bound2} we show that the off-diagonal entries of $\Re \mathcal{A}_{[k],t}$, namely $-\Re d_{ij,t}$, are negligible, and provide an exact calculation of the integral of the diagonal entries up to an error term of order~1. This allows us to get an optimal upper bound on integral of $f_{[k],r}$ introduced in \eqref{eq:propagator}, which we also present. The result is expressed in terms of quantities $\beta_{ij,*}$ for $i,j\in [k]$,  defined analogously to $\beta_{12,*}$ in \eqref{eq:def_beta_pm}, except that $z_1,z_2$ and $w_1,w_2$ are replaced by $z_i,z_j$ and $w_i,w_j$, respectively.

\begin{lemma}[Bound on the propagator]\label{lem:propag_bound} Fix (small) $\delta, \epsilon,\kappa>0$ and the final time $T>0$. Fix additionally an $N$-independent integer $k\ge 2$. The following statements hold uniformly in $z_{j,T}\in \ee^{-T/2}(1-\delta)\mathbf{D}$, $w_{j,T}\in \Omega_{\kappa,\epsilon,T}^{z_{j,T}}$, $j\in[k]$, and $s,t\in [0,T]$ with $s\le t$. First, for any $i,j\in [k]$ we have
\begin{equation}
\int_s^t \max\{\Re a_{ij,r}^+,\Re a_{ij,r}^-,0\}\dif r=\int_s^t \max\{\Re a_{ij,r}^+,\Re a_{ij,r}^-\}\dif r + \mathcal{O}(1) = \log  \frac{\beta_{ij,*,s}}{\beta_{ij,*,t}} + \mathcal{O}(1),\label{eq:a_bound}
\end{equation}
\begin{equation}
\int_s^t |\Re d_{ij,r}|\dif r=\mathcal{O}(1). \label{eq:d_bound}
\end{equation}
It further holds that
\begin{equation}
\exp\left\lbrace \int_s^t f_{[k],r}\dif r\right\rbrace \lesssim \prod_{j=1}^k \frac{\beta_{j(j+1),*,s}}{\beta_{j(j+1),*,t}},
\label{eq:propag_bound}
\end{equation}
where $\beta_{k(k+1),*,r}:=\beta_{k1,*,r}$ for any $r\in [0,T]$, and $f_{[k],r}$ is defined in \eqref{eq:propagator}.
\end{lemma}

The proofs of \eqref{eq:a_bound} and \eqref{eq:d_bound} were already presented in Appendix~\ref{app:propagator}, while \eqref{eq:propag_bound} immediately follows from \eqref{eq:a_bound}--\eqref{eq:d_bound}, similarly to the argument around \eqref{eq:f2_bound}. We omit further details.

\subsubsection{Application of the Gr{\"o}nwall inequality}\label{sec:Gronwall}

In this section we restrict the set-up of Section \ref{sec:Zig_structure} and consider only $k\in\{2,3,4\}$ in \eqref{eq:kG_syst}. Moreover, for $k=3$ we consider only the case when $z_{3,t}=z_{1,t}$ and $w_{3,t}=w_{1,t}^{(*)}$ for all $t\in [0,T]$, and for $k=4$ we additionally assume that $z_{4,t}=z_{2,t}$ and $w_{4,t}=w_{2,t}^{(*)}$. These quantities form a self-consistent family under the zig flow in the sense that no other constellations of $z_i$'s and $w_i$'s are needed to analyze~\eqref{eq:kG_syst}.

To control the error terms in the rhs. of \eqref{eq:kG_syst} for $k\in\{2,3,4\}$, we follow a stopping time argument. Denote
\begin{equation}
\alpha_{2,t}:=\frac{1}{N\eta_{*,t}\widehat{\beta}_{12,t}},\quad \alpha_{3,t}:=\frac{1}{N\eta^3_{*,t}}\wedge \frac{1}{\sqrt{N\eta_{*,t}}\eta_{*,t}\widehat{\beta}_{12,t}},\quad \alpha_{4,t}:=\frac{1}{N\eta^4_{*,t}}\wedge \frac{1}{\eta_{*,t}(\widehat{\beta}_{12,t})^2}
\label{eq:def_alpha}
\end{equation} 
for $t\in [0,T]$. These control parameters exactly match the rhs. of \eqref{eq:2G_bound_t}--\eqref{eq:4G_bound_t}, and depend on the trajectory of the characteristic flow, but as usual this dependence is omitted from notations. Fix small tolerance exponents $\xi_2,\xi_3,\xi_4$ satisfying
\begin{equation}
\xi_2,\xi_3,\xi_4\in (0,\epsilon/10),\quad \xi_2<\xi_3<\xi_4<2\xi_2.
\label{eq:tolerance_exp}
\end{equation}
Define the stopping times
\begin{equation}
\tau_k:= \inf\left\lbrace t\in [0,T]\,:\, \max_{s\in [0,t]}\max_{|z_{j,T}|\le \ee^{-T/2}(1-\delta)}\max_{w_{j,T}\in\Omega_{\kappa,\epsilon,T}^{z_{j,T}}} \alpha_{k,s}^{-1}\|\mathcal{Y}^{(k)}_s\|\ge N^{2\xi_k}\right\rbrace,\quad k=2,3,4,
\label{eq:def_tau_k}
\end{equation}
where we denoted
\begin{equation}
\mathcal{Y}^{(2)}_s:=\mathcal{Y}_{12,s},\quad \mathcal{Y}^{(3)}_s:=\mathcal{Y}_{121,s}\quad \text{and}\quad \mathcal{Y}^{(4)}_s:=\mathcal{Y}_{1212,s}.
\label{eq:Y_short}
\end{equation}
Here we slightly abuse the notation and abbreviate the quantities with stars in the lhs. of \eqref{eq:3G_bound_t}, \eqref{eq:4G_bound_t} using the same notation $\mathcal{Y}_{121}$ ($\mathcal{Y}_{1212}$, respectively). Finally, we set
\begin{equation}
\tau:=\min\left\lbrace \tau_2,\tau_3,\tau_4\right\rbrace.
\end{equation}
To prove Proposition \ref{prop:Zig} we need to show that $\tau=T$ almost surely.

We have the following bounds on the error terms in the rhs. of \eqref{eq:kG_syst}. The proof of this proposition is postponed to Section~\ref{sec:err_bounds}.

\begin{proposition}[Bounds on the error terms]\label{prop:err_bounds} Assume the set-up and conditions of Proposition \ref{prop:Zig}. Let $\mathcal{F}^{(k)}_s$ and $\dif\mathcal{E}^{(k)}_s$ be defined using the same convention as in \eqref{eq:Y_short} for $k=2,3,4$, and let $\mathcal{C}^{(k)}_s\dif s\in\C^{2^k\times 2^k}$ be the covariation process of $\dif\mathcal{E}^{(k)}_s$. Then for any $t\in [0,T]$ we have
\begin{equation}
\left(\int_0^{t\wedge\tau}\|\mathcal{F}_s^{(k)}\|\dif s\right)^2 +\int_0^{t\wedge\tau} \|\mathcal{C}^{(k)}_s\|\dif s \lesssim  N^{2\xi_k}\alpha_{k,t\wedge\tau}^2,\quad \text{for}\,\, k=2,3,4.
\label{eq:err_bounds}
\end{equation}
\end{proposition}
Having Proposition \ref{prop:err_bounds} in hand, we now prove Proposition \ref{prop:Zig}.

\begin{proof}[Proof of Proposition \ref{prop:Zig}] Throughout the proof $k$ is chosen from the range $\{2,3,4\}$ and the time parameters $r,s,t$ are in $[0,T]$. Similarly to \eqref{eq:Y_short}, we denote
\begin{equation*}
\mathcal{A}^{(2)}_s:=\mathcal{A}_{12},\quad \mathcal{A}^{(3)}_{s}:=\mathcal{A}_{121,s},\quad \mathcal{A}^{(4)}_s:=\mathcal{A}_{1212,s} \quad\text{and}\quad f_{s}^{(k)}:=\left[\max\mathrm{spec}\left(\Re \mathcal{A}^{(k)}_s\right)\right]_+.
\end{equation*}
From \eqref{eq:kG_syst}, Proposition \ref{prop:err_bounds} and the stochastic matrix-valued Gr{\"o}nwall inequality from \cite[Lemma 5.6]{univ_extr} we have
\begin{equation}
\begin{split}
\sup_{0\le s\le t\wedge \tau} \|\mathcal{Y}^{(k)}_s\|^2 &\lesssim \|\mathcal{Y}_0^{(k)}\|^2 + N^{2\xi_k+3\zeta}\alpha_{k,t\wedge \tau}^2\\
& + \int_0^{t\wedge \tau} \left(\lVert \mathcal{Y}_0^{(k)}\rVert^2 +N^{2\xi_k+3\zeta}\alpha_{k,s}^2\right)f_{s}^{(k)} \exp\left\lbrace 2(1+N^{-\zeta})\int_s^{t\wedge\tau} f_{r}^{(k)}\dif r\right\rbrace \dif s, 
\end{split}
\label{eq:stoch_Gronwall}
\end{equation}
for any (small) fixed $\zeta>0$.


We further focus on the case $k=3$, while for $k=2,4$ the argument is similar, and so omitted. It follows from \eqref{eq:propag_bound} that
\begin{equation}
\exp\left\lbrace\int_s^{t\wedge \tau}f_{r}^{(3)}\dif r\right\rbrace \lesssim \left(\frac{\beta_{12,*,s}}{\beta_{12,*,t\wedge\tau}}\right)^2 \frac{\beta_{11,*,s}}{\beta_{11,*,t\wedge\tau}}.
\label{eq:f_int}
\end{equation}
We plug \eqref{eq:f_int} into \eqref{eq:stoch_Gronwall}, estimate $\|\mathcal{Y}^{(3)}_0\|$ from the initial condition \eqref{eq:3G_bound_t} at time $t=0$, and bound $f_{s}^{(3)}\lesssim \eta_{*,s}^{-1}$ by \eqref{eq:M_12_beta_bound}, \eqref{eq:eta_stab_bound}, arriving to
\begin{equation}
\sup_{0\le s\le t\wedge \tau} \|\mathcal{Y}^{(3)}_s\|^2\lesssim N^{2\xi_3+3\zeta} \left(\alpha^2_{3,t\wedge\tau} + \int_0^{t\wedge\tau} \alpha_{3,s}^2 \frac{1}{\eta_{*,s}}\left(\frac{\beta_{12,*,s}}{\beta_{12,*,t\wedge\tau}}\right)^4 \left(\frac{\beta_{11,*,s}}{\beta_{11,*,t\wedge\tau}}\right)^2\dif s\right).
\label{eq:Gronwall_3G}
\end{equation}

We now show that
\begin{equation}
\int_0^{t\wedge\tau} \alpha_{3,s}^2 \frac{1}{\eta_{*,s}}\left(\frac{\beta_{12,*,s}}{\beta_{12,*,t\wedge\tau}}\right)^4 \left(\frac{\beta_{11,*,s}}{\beta_{11,*,t\wedge\tau}}\right)^2\dif s \lesssim \alpha_{3,t\wedge\tau}^2\log N.
\label{eq:int_propagation} 
\end{equation}
Recall from \eqref{eq:def_alpha} that $\alpha_{3,s}$ is a minimum of two control parameters: $(N\eta_{*,s}^3)^{-1}$ and $(\sqrt{N\eta_{*,s}}\eta_{*,s}\widehat{\beta}_{12,s})^{-1}$. To show that the lhs. of \eqref{eq:int_propagation} is smaller than the square of the first of them evaluated at the time $t\wedge \tau$, we observe that
\begin{equation}
\frac{\beta_{12,*,s}}{\beta_{12,*,t\wedge\tau}}\sim \frac{\beta_{12,*,t\wedge\tau}+|s-t\wedge\tau|}{\beta_{12,*,t\wedge\tau}}\lesssim  \frac{\eta_{*,t\wedge\tau}+|s-t\wedge\tau|}{\eta_{*,t\wedge\tau}}\sim \frac{\eta_{*,s}}{\eta_{*,t\wedge\tau}}
\label{eq:beta_ratio}
\end{equation}
for any $s\in [0,t\wedge\tau]$, where we used the first part of \eqref{eq:gamma_monot} in the first step, \eqref{eq:beta_hat_to_beta} and \eqref{eq:eta_stab_bound} in the second, and Lemma \ref{lem:char_flow}(4) in the third. Estimating $\beta_{11,*,s}/\beta_{11,*,t\wedge\tau}$ in the same way, we obtain that the lhs. of \eqref{eq:int_propagation} has an upper bound of order
\begin{equation}
\int_0^{t\wedge\tau} \left(\frac{1}{N\eta_{*,s}^3}\right)^2 \frac{1}{\eta_{*,s}}\left(\frac{\eta_{*,s}}{\eta_{*,t\wedge\tau}}\right)^6 \dif s = \left(\frac{1}{N\eta_{*,t\wedge\tau}^3}\right)^2 \int_0^{t\wedge\tau}\frac{\dif s}{\eta_{*,s}}\lesssim \left(\frac{1}{N\eta_{*,t\wedge\tau}^3}\right)^2\log N.
\label{eq:eta_propagation}
\end{equation}
In the last bound we used the second part of \eqref{eq:eta_int} and that $\eta_{*,t\wedge\tau}>N^{-1}$.

Now we prove that the lhs. of \eqref{eq:int_propagation} is smaller than the square of second component of $\alpha_{3,t\wedge \tau}$ (up to a $\log N$ factor). We again estimate the ratio of $\beta_{11,*}$'s as in \eqref{eq:beta_ratio} and further split the interval of integration into two regimes: $[0,\widetilde{t}]$ and $[\widetilde{t},T]$, where the random variable $\widetilde{t}$ is defined by
\begin{equation*}
\widetilde{t}:=\left[t\wedge\tau - \beta_{12,*,t\wedge\tau}\right]_+.
\end{equation*}
For any $s\in [0,\widetilde{t}]$ we have
\begin{equation}
\eta_{*,s}\lesssim \beta_{12,*,s}\sim \beta_{12,*,t\wedge\tau}+|t\wedge\tau-s|\sim |t\wedge\tau-s| \lesssim \eta_{*,t\wedge\tau} + |t\wedge\tau-s|\sim \eta_{*,s},
\label{eq:small_time}
\end{equation}
i.e. $\beta_{12,*,s}\sim \eta_{*,s}$, where we used the first part of \eqref{eq:gamma_monot} and Lemma \ref{lem:char_flow}(4). Therefore, the integral in the lhs. of \eqref{eq:int_propagation} restricted to $[0,\widetilde{t}]$ has an upper bound of order
\begin{equation}
\int_0^{\widetilde{t}} \left(\frac{1}{N\eta_{*,s}^3}\right)^2 \frac{1}{\eta_{*,s}}\left(\frac{\eta_{*,s}}{\widehat{\beta}_{12,t\wedge\tau}}\right)^4\left(\frac{\eta_{*,s}}{\eta_{*,t\wedge\tau}}\right)^2\dif s\lesssim \left(\frac{1}{N\eta_{*,t\wedge\tau}(\widehat{\beta}_{12,t\wedge\tau})^2}\right)^2\log N.
\label{eq:beta_propagation1}
\end{equation}
Here we used the second part of \eqref{eq:eta_int} and the bound $\beta_{12,*,t\wedge\tau}\ge \widehat{\beta}_{12,t\wedge\tau}$. In the regime $s\in [\widetilde{t},T]$, from the second step in \eqref{eq:small_time} we get that $\beta_{12,*,s}\sim \beta_{12,*,t\wedge\tau}$, so the integral in the lhs. of \eqref{eq:int_propagation} restricted to $[\widetilde{t}, T]$ is bounded by
\begin{equation}
\int_0^{t\wedge\tau} \left( \frac{1}{\sqrt{N\eta_{*,s}}\eta_{*,s}\widehat{\beta}_{12,s}}\right)^2 \frac{1}{\eta_{*,s}}\left(\frac{\eta_{*,s}}{\eta_{*,t\wedge\tau}}\right)^2\dif s\lesssim \left( \frac{1}{\sqrt{N\eta_{*,t\wedge\tau}}\eta_{*,t\wedge\tau}\widehat{\beta}_{12,t\wedge\tau}}\right)^2,
\label{eq:beta_propagation2}
\end{equation}
where we estimated $\widehat{\beta}_{12,s}\gtrsim \widehat{\beta}_{12,t\wedge\tau}$ from \eqref{eq:gamma_monot} and used \eqref{eq:eta_int} for $a=2$. Combining \eqref{eq:eta_propagation}, \eqref{eq:beta_propagation1} and \eqref{eq:beta_propagation2}, we finish the proof of \eqref{eq:int_propagation}.

Now we have all ingredients to complete the proof of Proposition \ref{prop:Zig}. Using \eqref{eq:Gronwall_3G} along with \eqref{eq:int_propagation} we~get
\begin{equation}
\sup_{0\le s\le t\wedge \tau} \|\mathcal{Y}^{(3)}_s\|^2\lesssim N^{2\xi_3+3\zeta}\alpha^2_{3,t\wedge\tau}\log N \lesssim N^{3\xi_3}\alpha^2_{3,t\wedge\tau},
\label{eq:tau3_saturation}
\end{equation}
provided that $\zeta>0$ is chosen to be sufficiently small. Recall the definition of the stopping time $\tau_3$ from \eqref{eq:def_tau_k}. From \eqref{eq:tau3_saturation} we have that either $\tau_3>\tau$ or $\tau_3=T$. Arguing similarly for $\tau_2$ and $\tau_4$ and using that $\tau=\min\{\tau_1,\tau_2,\tau_3\}$, we conclude that $\tau=T$. This finishes the proof of Proposition \ref{prop:Zig}.
\end{proof}

\subsubsection{Proof of Proposition \ref{prop:err_bounds}: bounds on the error terms}\label{sec:err_bounds}

For simplicity, we establish the bounds only on the error terms arising from the differentiation of resolvent chains, where the complex conjugates are not involved. The remaining cases are completely analogous and thus are omitted. 

Recall from \eqref{eq:def_alpha} that each of the control parameters $\alpha_{k,t}$, $k\in\{2,3,4\}$, is a minimum of two components: $(N\eta_{*,t}^{k})^{-1}$ and another term containing $\widehat{\beta}_{12,t}$. The case $k=2$ is not an exception as it may seem from the definition of $\alpha_{2,t}$, since $(N\eta_{*,t}\widehat{\beta}_{12,t})^{-1}$ is always smaller than $(N\eta_{*,t}^2)^{-1}$, due to \eqref{eq:eta_stab_bound}. While the proof of the $(N\eta_{*,t\wedge\tau}^{k})^{-2}$ bounds on the lhs. of \eqref{eq:err_bounds} is fairly standard, we first prove the upper bounds involving $\widehat{\beta}_{12,t\wedge\tau}$ and later in the end of Section \ref{sec:err_bounds} explain how to adjust the proof to obtain the bounds purely in terms of $\eta_{*,t\wedge\tau}$.

We start with the estimates on the second term in the lhs. of \eqref{eq:err_bounds}, which arises from the martingale term $\dif\mathfrak{E}^{(k)}$ and show that
\begin{equation}
\int_0^{t\wedge\tau} \|\mathcal{C}^{(k)}_s\| \dif s \lesssim 
\begin{cases}
N^{2\xi_2}(N\eta_{*,t\wedge\tau}\widehat{\beta}_{12,t\wedge\tau})^{-2},&k=2,\\
N^{2\xi_3}(\sqrt{N\eta_{*,t\wedge\tau}}\eta_{*,t\wedge\tau}\widehat{\beta}_{12,t\wedge\tau})^{-2},&k=3,\\
N^{2\xi_4}(\eta_{*,t\wedge\tau}(\widehat{\beta}_{12,t\wedge\tau})^2)^{-2},&k=4.
\end{cases}
\label{eq:QV_beta_bound}
\end{equation}
Recall from the formulation of Proposition \ref{prop:err_bounds} that $\mathcal{C}^{(k)}_s\dif s$ is the covariation process of $\dif\mathfrak{E}^{(k)}_s$. Thus, $\mathcal{C}^{(k)}_s$ is a matrix of size $2^k\times 2^k$ with entries indexed by $\bm{\sigma},\bm{\omega}\in\{\pm\}^k$.
Since the size of $\mathcal{C}^{(k)}_s$ does not depend on $N$, we have
\begin{equation}
\|\mathcal{C}^{(k)}_s\|\lesssim \sum_{\bm{\sigma}, \bm{\omega}} \Big\vert\big(\mathcal{C}^{(k)}_s\big)_{\bm{\sigma}, \bm{\omega}}\Big\vert,
\label{eq:QV_matrix_bound}
\end{equation}
where the implicit constant in the inequality depends on $k$, but not on $N$. Therefore, it is sufficient to prove \eqref{eq:QV_beta_bound} separately for each of the terms in the rhs. of \eqref{eq:QV_matrix_bound}. To simplify the presentation, we present the bounds only for diagonal entries, i.e. for $\bm{\sigma}=\bm{\omega}$, while for the off-diagonal entries the proof is identical. We also denote $B_i:=E_{\sigma_i}$ for $i\in [k]$.

Now we consider separately each of the cases $k\in\{2,3,4\}$.

\medskip

\noindent\underline{$k=2$.} We have from \cite[Eq.(5.28),(5.29)]{univ_extr} the following simple estimate on the entries of $\mathcal{C}^{(k)}_s$ based on the explicit calculation and the Cauchy-Schwarz inequality:
\begin{equation}
\Big\vert\big(\mathcal{C}^{(2)}_s\big)_{\bm{\sigma}, \bm{\sigma}}\Big\vert \lesssim \frac{1}{N^2\eta_{1,s}^2}\left\langle \Im G_{1,s}B_1G_{2,s} B_2 \Im G_{1,s}B_2^*G_{2,s}^* B_1^*\right\rangle + \frac{1}{N^2\eta_{2,s}^2} \left\langle \Im G_{2,s}B_2G_{1,s}B_1\Im G_{2,s} B_1^* G_{1,s}^*B_2^*\right\rangle.
\label{eq:QV_2G}
\end{equation} 
Both terms in the rhs. of \eqref{eq:QV_2G} are the four-resolvent chains exactly of the type which is controlled by the stopping time $\tau_4$ defined in \eqref{eq:def_tau_k}. Decomposing each of these chains into the deterministic approximation and the fluctuation around it, we get from \eqref{eq:M4_bound} and \eqref{eq:def_tau_k} that the rhs. of \eqref{eq:QV_2G} has an upper bound of order
\begin{equation}
\frac{1}{N^2\eta_{*,s}^2}\left( \frac{1}{\eta_{*,s}(\widehat{\beta}_{12,s})^2} + \frac{N^{\xi_4}}{N\eta_{*,s}^4}\wedge\frac{N^{\xi_4}}{\eta_{*,s}(\widehat{\beta}_{12,s})^2}\right)\lesssim \frac{N^{\xi_4}}{N^2\eta_{*,s}^3 (\widehat{\beta}_{12,s})^2},\quad \forall\, s\in [0,t\wedge\tau].
\label{eq:QV_2G_bound}
\end{equation}
Estimating $\widehat{\beta}_{12,s}\gtrsim \widehat{\beta}_{12,t\wedge\tau}$ from \eqref{eq:gamma_monot}, integrating the rhs. of \eqref{eq:QV_2G_bound} in $s\in [0,t\wedge\tau]$ by the means of \eqref{eq:eta_int} applied to $a=3$, and recalling from \eqref{eq:tolerance_exp} that $\xi_4<2\xi_2$, we finish the proof of \eqref{eq:QV_beta_bound} for $k=2$.

\medskip

\noindent\underline{$k=3$.} Similarly to \eqref{eq:QV_2G}, we upper bound $|(\mathcal{C}^{(3)}_s)_{\bm{\sigma},\bm{\sigma}}|$ by the sum of three terms, two of which are given~by
\begin{equation}
\begin{split}
&\frac{1}{N^2\eta_{1,s}^2} \left\langle \Im G_{1,s} B_1 G_{2,s} B_2 G_{1,s} B_3 \Im G_{1,s} B_3^* G_{1,s}^* B_2^* G_{2,s}^* B_1^*\right\rangle,\\
&\frac{1}{N^2\eta_{2,s}^2} \left\langle \Im G_{2,s} B_2 G_{1,s} B_3 G_{1,s} B_1\Im G_{2,s} B_1^* G_{1,s}^*B_3^* G_{1,s}^* B_2^*\right\rangle,
\end{split}
\label{eq:QV_3G}
\end{equation}
and the third term is analogous to the one in the first line of \eqref{eq:QV_3G}. Each of these terms contains six resolvents and thus it is not directly controlled by the stopping time $\tau$, in contrast to the terms in the rhs. of \eqref{eq:QV_2G}. This forces us to invoke the so-called \emph{reduction inequalities} (as introduced in \cite{Multi_res_llaws}), which bound traces of longer resolvent chains in terms of traces of shorter ones.
 To estimate the term in the second line of \eqref{eq:QV_3G}, we observe that for all $R,S\in \C^{(2N)\times (2N)}$ such that each of the matrices $S,T$ is either positive or negative semi-definite, it holds that
\begin{equation}
|\langle RS\rangle|\le 2N|\langle R\rangle\langle S\rangle|.
\label{eq:reduction_simple}
\end{equation}
We apply this elementary reduction inequality to
\begin{equation}
R:=G_{1,s}^*B_2^*\Im G_{2,s}B_2G_{1,s}\quad \text{and}\quad S:= B_3G_{1,s}B_1 \Im G_{2,s} B_1^* G_{1,s}^*B_3^*,
\end{equation}
and note that the trace in the second line of \eqref{eq:QV_3G} equals to $\langle RS\rangle$. Since $\langle R\rangle$ and $\langle S\rangle$ are controlled by the stopping time $\tau_3$, we further estimate each of these traces in the way similar to \eqref{eq:QV_2G_bound}. We omit further details concerning integration of the resulting bound in time, since this procedure is similar to the one discussed below \eqref{eq:QV_2G_bound}.

Now we estimate the term in the first line of \eqref{eq:QV_3G} and instead of \eqref{eq:reduction_simple} use the following reduction bound from \cite[Eq.(5.27)]{eigenv_decorr}:  
\begin{equation}
\left\vert\left\langle G_{2,s}RG_{2,s}^*S\right\rangle\right\vert \lesssim N\left\vert \left\langle |G_{2,s}|R\right\rangle \left\langle |G_{2,s}|S\right\rangle\right\vert,
\label{eq:reduction}
\end{equation}
uniformly in the same set of matrices $R,S$ as in \eqref{eq:reduction_simple}. We take
\begin{equation}
R:= B_2G_{1,s}B_3\Im G_{1,s}B_3^*G_{1,s}^*B_2^*\quad \text{and}\quad S:= B_1^*\Im G_{1,s} B_1.
\label{eq:def_RS_QV3}
\end{equation}
However, $\langle |G_{2,s}|R\rangle$ and $\langle |G_{2,s}|S\rangle$ are not resolvent chains in the usual sense, since they involve the absolute value of a resolvent. To address this issue, we use the following integral representation from \cite[Eq.(5.4)]{Multi_res_llaws}
\begin{equation}
|G^z(E+\ii\eta)| = \frac{2}{\pi} \int_0^\infty \frac{\Im G^z(E+\ii\sqrt{\eta^2+x^2})}{\sqrt{\eta^2+x^2}}\dif x,\quad \forall\,z\in\C,\, E\in\R,\, \eta>0. 
\label{eq:abs_int_rep}
\end{equation}
We complement \eqref{eq:abs_int_rep} by the following geometric property of the bulk-restricted domains, essential for the application of \eqref{eq:abs_int_rep}. This property is a special case of \cite[Lemma 5.4]{eigenv_decorr} in the sense discussed above Definition \ref{def:spec_dom}.

\begin{lemma}[Ray property of the bulk-restricted spectral domains]\label{lem:ray_property} Fix (small) $\kappa, \epsilon ,\delta>0$ and the final time $T>0$. There exists $t_*\in [0,T]$ such that $T-t_*\sim 1$ and we have the following. For any $t\in [t_*,T]$, $z_T\in \ee^{-T/2}(1-\delta)\bf{D}$, $w\in \Omega^{z_t}_{\kappa,\epsilon,t}$, and $x\ge 0$ such that $|\Im w|+x\le N^{100}$, it holds that
\begin{equation*}
w+\mathrm{sgn}(\Im w)\ii x\in \Omega^{z_t}_{\kappa,\epsilon,t}.
\end{equation*}
\end{lemma}

That is, for $\Im w>0$ ($\Im w<0$) the vertical ray which starts at $w$ and goes up, leaves the spectral domain only after reaching points with imaginary part greater than $N^{100}$ (smaller than $N^{-100}$). Since $T-t_*\sim 1$, we may assume that $t_*=0$. Otherwise it suffices to redefine $T$, taking the final time equal to $T-t_*$.

Applying \eqref{eq:abs_int_rep} to $|G_{2,s}|$, we get
\begin{equation}
\left\langle |G_{2,s}|S\right\rangle = \frac{2}{\pi}\int_0^\infty\frac{1}{\sqrt{\eta_{2,t}^2+x^2}}\left\langle \Im G^{z_{2,t}}_t\left(\Re w_{2,t}+\ii\sqrt{\eta_{2,t}^2+x^2}\right)B_1^*\Im G_{1,t}B_1\right\rangle\dif x,
\label{eq:abs_int_rep_2G}
\end{equation}
where $\eta_{2,t}:=|\Im w_{2,t}|$. In the regime $\sqrt{\eta_{2,t}^2+x^2}\ge N^{100}$ we trivially estimate the trace in the rhs. of \eqref{eq:abs_int_rep_2G} by the product of operator norms of the involved resolvents, which integrates to the bound of order $N^{-100}\eta_{1,t}^{-1}$. In the complementary regime $\sqrt{\eta_{2,t}^2+x^2}< N^{100}$, the spectral parameter $w_2':=\Re w_{2,t}+\ii\sqrt{\eta_{2,t}^2+x^2}$ lies in $\Omega_{\kappa,\epsilon,t}^{z_{2,t}}$ by Lemma \ref{lem:ray_property}, so the trace is controlled by the stopping time $\tau_2$ and has an upper bound of order $(\widehat{\beta}_{12,t}\left(w_{1,t},w_2'\right))^{-1}$, where we also employed \eqref{eq:M2_bound}. Using additionally the vague monotonicity of $\widehat{\beta}_{12}$ in imaginary part from Lemma \ref{lem:beta_hat}(4), we conclude that
\begin{equation}
\left\vert\left\langle |G_{2,s}|B_1^*\Im G_{1,s} B_1\right\rangle\right\vert \lesssim \frac{\log N}{\widehat{\beta}_{12,s}}
\label{eq:2G_abs_bound}
\end{equation} 
with very high probability, for more details see e.g. \cite[Eq.(5.29)--(5.31)]{Multi_res_llaws}.

Recall the definition of $R$ from \eqref{eq:def_RS_QV3}. We further bound
\begin{equation}
\left\vert\left\langle |G_{2,s}|R\right\rangle\right\vert = \left\vert\left\langle \Im G_{1,s} B_3^*G_{1,s}^*B_2^* |G_{2,s}|B_2 G_{1,s}B_3\right\rangle\right\vert\le \lVert \Im G_{2,s}\rVert \left\langle B_3^*G_{1,s}^*B_2^* |G_{2,s}|B_2 G_{1,s}B_3\right\rangle.
\label{eq:QV_3G_GR}
\end{equation}
Finally, we estimate the trace in the rhs. of \eqref{eq:QV_3G_GR} from above by $(\eta_{*,s}\widehat{\beta}_{12,s})^{-1}$ similarly to \eqref{eq:abs_int_rep_2G}--\eqref{eq:2G_abs_bound}, bound trivially $\lVert \Im G_{2,s}\rVert \le \eta_{*,s}^{-1}$, and combine the resulting bound with \eqref{eq:reduction} and \eqref{eq:2G_abs_bound}. We thus get
\begin{equation*}
\frac{1}{N^2\eta_{1,s}^2} \left\vert\left\langle \Im G_{1,s} B_1 G_{2,s} B_2 G_{1,s} B_3 \Im G_{1,s} B_3^* G_{1,s}^* B_2^* G_{2,s}^* B_1^*\right\rangle\right\vert \lesssim \frac{(\log N)^2}{N\eta_{*,s}^4(\widehat{\beta}_{12,s})^2}.
\end{equation*}
Integrating this bound over $s\in [0,t\wedge\tau]$, we finish the proof of \eqref{eq:QV_beta_bound} for $k=3$.

\medskip

\noindent\underline{$k=4$.} Similarly to \eqref{eq:QV_2G}, we upper bound $|(\mathcal{C}^{(4)}_s)_{\bm{\sigma},\bm{\sigma}}|$ by the sum of four terms of the following type:
\begin{equation}
\frac{1}{N^2\eta_{1,s}^2} \left\langle \Im G_{1,s} B_1G_{2,s}B_2G_{1,s}B_3G_{2,s}B_4\Im G_{1,s} B_4^* G_{2,s}^* B_3^* G_{1,s}^* B_2^* G_{2,s}^* B_1^*\right\rangle.
\label{eq:QV_4G}
\end{equation}
We apply \eqref{eq:reduction} with $G_{1,s}$ instead of $G_{2,s}$ and choose
\begin{equation}
R:=B_3G_{2,s}B_4\Im G_{1,s} B_4^*G_{2,s}^*B_3^*,\quad S:=B_2^*G_{2,s}^*B_1^*\Im G_{1,s} B_1G_{2,s} B_2.
\end{equation}
This yields that \eqref{eq:QV_4G} is bounded by
\begin{equation}
\frac{1}{N\eta_{1,s}^2} \left\vert \left\langle |G_{1,s}| B_3G_{2,s}B_4\Im G_{1,s} B_4^*G_{2,s}^*B_3^*\right\rangle \left\langle |G_{1,s}| B_2^*G_{2,s}^*B_1^*\Im G_{1,s} B_1G_{2,s} B_2\right\rangle\right\vert.
\label{eq:QV_4G_bound1}
\end{equation}
Since the traces in the rhs. of \eqref{eq:QV_4G_bound1} are analogous, we further focus on the first of them. We use the integral representation \eqref{eq:abs_int_rep} for $|G_{1,s}|$ and denote $\widetilde{G}_{1,s}:=G^{z_{1,s}}_s \left(\Re w_{1,s} +\sqrt{\eta_{1,s}^2 +x^2}\right)$, omitting $x$ from this notation. We obtain the four-resolvent chain with Hermitization parameters alternating between $z_{1,s}$ and $z_{2,s}$, however the spectral parameters of $\widetilde{G}_{1,s}$ and $G_{1,s}$ differ more than just by a complex conjugation. To remove this discrepancy, we observe that for any $X,Y\in \C^{(2N)\times (2N)}$ it holds that
\begin{equation}
\left\vert\left\langle G_{2,s}XG_{2,s}^*Y\right\rangle\right\vert \le \left\langle |G_{2,s}|X|G_{2,s}|X^*\right\rangle^{1/2} \left\langle |G_{2,s}|Y|G_{2,s}|Y^*\right\rangle^{1/2}.
\label{eq:reduction_Schwarz}
\end{equation}
We stress that the semi-definiteness of $X$ and $Y$ is not required in \eqref{eq:reduction_Schwarz}. This bound immediately follows from the fact that $G_{2,s}$ is diagonalizable and from the Cauchy-Schwarz inequality. Applying \eqref{eq:reduction_Schwarz} for $X=B_4\Im G_{1,s}B_4^*$ and $Y:=B_3^*\Im \widetilde{G}_{1,s}B_3$, we get
\begin{equation}
\begin{split}
&\left\vert\left\langle \Im \widetilde{G}_{1,s} B_3G_{2,s}B_4\Im G_{1,s} B_4^*G_{2,s}^*B_3^*\right\rangle\right\vert\\
&\quad \le \left\langle |G_{2,s}| B_4\Im G_{1,s}B_4^*|G_{2,s}| B_4^*\Im G_{1,s}B_4\right\rangle^{1/2}\left\langle |G_{2,s}| B_3\Im \widetilde{G}_{1,s}B_3^*|G_{2,s}| B_3^*\Im \widetilde{G}_{1,s}B_3\right\rangle^{1/2}.
\end{split}
\label{eq:QV_4G_bound2}
\end{equation}
We further focus again on the first factor in the rhs. of \eqref{eq:QV_4G_bound2} and use the integral representation \eqref{eq:abs_int_rep} for both matrices $|G_{2,s}|$ appearing in the product. Denoting the resolvents arising from these representations by $\widetilde{G}_{2,s}$ and $\widetilde{G}_{2,s}'$, we reduce the first trace in the rhs. of \eqref{eq:QV_4G_bound2} to
\begin{equation}
\begin{split}
&\left\vert\left\langle \Im\widetilde{G}_{2,s} B_4\Im G_{1,s}B_4^*\Im\widetilde{G}_{2,s}' B_4^*\Im G_{1,s}B_4\right\rangle\right\vert\\
&\quad \le \left\langle \Im G_{1,s}B_4^*\Im\widetilde{G}_{2,s}' B_4^*\Im G_{1,s} B_4 \Im\widetilde{G}_{2,s}'B_4\right\rangle^{1/2} \left\langle \Im G_{1,s}B_4^*\Im\widetilde{G}_{2,s} B_4^*\Im G_{1,s} B_4 \Im\widetilde{G}_{2,s}B_4\right\rangle^{1/2}, 
\end{split}
\label{eq:QV_4G_bound3}
\end{equation}
where we used the Cauchy-Schwarz inequality to go from the first to the second line. From the definition of the stopping time $\tau_4$ \eqref{eq:def_tau_k}, the fact that $s\le \tau_4$, and \eqref{eq:beta_hat_vert_line}, we get that each of the traces in the rhs. of \eqref{eq:QV_4G_bound3} has an upper bound of order $N^{2\xi_4}\eta_{*,s}^{-1}(\widehat{\beta}_{12,s})^{-2}$. Finally, combining this bound with \eqref{eq:QV_4G_bound1}--\eqref{eq:QV_4G_bound3} and performing the integration similarly to \eqref{eq:abs_int_rep_2G}--\eqref{eq:2G_abs_bound} in the representations \eqref{eq:abs_int_rep} which we used along the way, we obtain the following bound on the integral of \eqref{eq:QV_4G} in $s\in [0,t\wedge\tau]$:
\begin{equation}
\int_0^{t\wedge\tau}\frac{1}{N\eta_{*,s}^2}\left(\frac{N^{2\xi_4}}{\eta_{*,s}\widehat{\beta}_{12,s}}\right)^2\dif s\lesssim \frac{N^{4\xi_4}}{N\eta_{*,t\wedge\tau}}\left(\frac{1}{\eta_{*,t\wedge\tau}(\widehat{\beta}_{12,t\wedge\tau})^2}\right)^2.
\label{eq:QV_4G_bound4}
\end{equation}
Here we used the second part of \eqref{eq:gamma_monot} and \eqref{eq:eta_int} for $a=4$. Since $N\eta_{*,t\wedge\tau}\ge N^{\epsilon}$ and $\xi_4< \epsilon/10$ by \eqref{eq:tolerance_exp}, \eqref{eq:QV_4G_bound4} finishes the proof of \eqref{eq:QV_beta_bound} for $k=4$.

Next, we prove that
\begin{equation}
\int_0^{t\wedge\tau} \|\mathcal{F}^{(k)}_s\| \dif s \lesssim 
\begin{cases}
N^{\xi_2}(N\eta_{*,t\wedge\tau}\widehat{\beta}_{12,t\wedge\tau})^{-1},&k=2,\\
N^{\xi_3}(\sqrt{N\eta_{*,t\wedge\tau}}\eta_{*,t\wedge\tau}\widehat{\beta}_{12,t\wedge\tau})^{-1},&k=3,\\
N^{\xi_4}(\eta_{*,t\wedge\tau}(\widehat{\beta}_{12,t\wedge\tau})^2)^{-1},&k=4.
\end{cases}
\label{eq:F_beta_bound}
\end{equation}
Since the proof of \eqref{eq:F_beta_bound} does not require any additional ideas compared to the proof of \eqref{eq:QV_beta_bound}, we only comment on the involved reduction inequalities. The proof of \eqref{eq:F_beta_bound} for $k=2$ does not require any reductions of resolvent chains to shorter ones, and only utilizes the bounds provided by the stopping time $\tau$, single-resolvent local law \eqref{eq:av_ll_simple}, and Proposition \ref{prop:M_bounds}. For $k=3$ reductions are needed only for the terms in the first line of \eqref{eq:def_Fk}. There are two different types of such terms, which we present in the lhs. of the following two lines and immediately demonstrate the corresponding reductions performed by the means of the Cauchy-Schwarz inequality:
\begin{equation}
\begin{split}
&\left\vert\left\langle G_{1,s}E_\sigma G_{1,s}B_1G_{2,s}B_2G_{1,s}B_3\right\rangle\right\vert \le \frac{1}{\eta_{*,s}^{3/2}}\left\vert\left\langle \Im G_{1,s}B_3^*B_3\right\rangle\right\vert^{1/2}\left\vert\left\langle \Im G_{1,s}B_1G_{2,s}B_2\Im G_{1,s}B_2^*G_{2,s}^*B_1^*\right\rangle\right\vert^{1/2}\\
&\left\vert\left\langle G_{1,s}B_1G_{2,s}E_\sigma G_{2,s}B_2G_{1,s}B_3\right\rangle\right\vert \le \frac{1}{\eta_{*,s}^{3/2}}\left\vert\left\langle \Im G_{1,s}B_1\Im G_{2,s}B_1^*\right\rangle\right\vert^{1/2}\left\vert\left\langle G_{1,s}^*B_2^*\Im G_{2,s}B_2G_{1,s}B_3B_3^*\right\rangle\right\vert^{1/2}.
\end{split}
\label{eq:F_3G}
\end{equation}
Here $\sigma\in \{\pm\}$ is fixed, i.e. we do not perform the summation over this index. When the lhs. of both lines of \eqref{eq:F_3G} appear in the first line of \eqref{eq:def_Fk}, they are multiplied by $\left\langle (G_{j,s}-M_{j,s})E_\sigma\right\rangle$ for $j=1$ and $j=2$ respectively, which gives an additional small factor $(N\eta_{*,s})^{-1}$. Estimating the rhs. of both lines of \eqref{eq:F_3G} using that $s\le \tau$ and employing the deterministic bounds from Proposition \ref{prop:M_bounds}, we prove \eqref{eq:F_beta_bound} for $k=3$.

We are now left with \eqref{eq:F_beta_bound} for $k=4$, where the reduction is needed again only for the terms in the first line of \eqref{eq:def_Fk}. All of these terms are of the same form, so we consider only one representative:
\begin{equation}
\left\langle (G_{1,s}-M_{1,s})E_\sigma\right\rangle \left\langle G_{1,s}E_\sigma G_{1,s}B_1G_{2,s}B_2G_{1,s}B_3G_{2,s}B_4\right\rangle.
\label{eq:F_4G}
\end{equation}
Similarly to \eqref{eq:reduction_Schwarz} we estimate the second trace in \eqref{eq:F_4G} from above by
\begin{equation}
\frac{1}{\eta_{*,s}} \left\vert\left\langle \Im G_{1,s}B_1G_{2,s}B_2|G_{1,s}|B_2^*G_{2,s}^*B_1^*\right\rangle\right\vert^{1/2}\left\vert\left\langle \Im G_{1,s} B_4^*G_{2,s}^*B_3^*|G_{1,s}|B_3G_2B_4\right\rangle\right\vert^{1/2}.
\label{eq:F_4G_bound}
\end{equation}
Arguing analogously to \eqref{eq:QV_4G_bound2}--\eqref{eq:QV_4G_bound3}, we get that each of the traces in \eqref{eq:F_4G_bound} has an upper bound of order $N^{2\xi_4}\eta_{*,s}^{-1}(\widehat{\beta}_{12,s})^{-2}$, up to some irrelevant $\log N$ factor. Together with the single-resolvent local law \eqref{eq:av_ll_simple}, this yields an upper bound of order
\begin{equation*}
\frac{N^{2\xi_4}}{N\eta_{*,s}^3(\widehat{\beta}_{12,s})^2}
\end{equation*}
on the absolute value of the term in \eqref{eq:F_4G}. Integrating this bound over $s\in [0,t\wedge\tau]$, we finish the proof of \eqref{eq:F_beta_bound} for $k=4$.

Finally, we explain how to adapt the proof of \eqref{eq:QV_beta_bound} and \eqref{eq:F_beta_bound} to prove the bounds of order 
\begin{equation*}
N^{2\xi_k}(N\eta_{*,t\wedge \tau}^k)^{-2}\quad\text{and}\quad N^{\xi_k}(N\eta_{*,t\wedge \tau}^k)^{-1}
\end{equation*}
on the lhs. of \eqref{eq:QV_beta_bound} and \eqref{eq:F_beta_bound}, respectively. First, throughout the proof we replace the use of Proposition \ref{prop:M_bounds} for estimating deterministic approximations, by Lemma \ref{lem:M_bounds_weak}, which provides a weaker control, but purely in terms of $\eta_*$. Similarly, for all times $s\le t\wedge\tau$ we use only the bounds involving solely $\eta_*$ which are provided by \eqref{eq:def_tau_k}, e.g. for four-resolvent chains we apply
\begin{equation}
\left\vert\left\langle \left(G_{1,s}B_1G_{2,s}B_2 G_{1,s}B_3G_{2,s}-M_{1212,s}^{B_1,B_2,B_3}\right)B_4\right\rangle\right\vert \le \frac{N^{2\xi_4}}{N\eta_{*,s}^4}\quad \text{for}\quad 0\le s\le t\wedge\tau.
\end{equation}
A further modification is required for the estimates where the reduction inequalities \eqref{eq:reduction_simple} and \eqref{eq:reduction} were used, since these bounds introduce an additional factor $N$, which is not affordable now. These reductions were used only in three instances: for both terms in \eqref{eq:QV_3G} and for the term in \eqref{eq:QV_4G}. Now we show how to control these quantities in terms of $\eta_*$ without incurring factor $N$. We demonstrate this argument only for the first term in \eqref{eq:QV_3G}, while for the second term in \eqref{eq:QV_3G} and for \eqref{eq:QV_4G} the estimates are analogous. We have
\begin{equation}
\begin{split}
&\left\vert\left\langle \Im G_{1,s} B_1 G_{2,s} B_2 G_{1,s} B_3 \Im G_{1,s} B_3^* G_{1,s}^* B_2^* G_{2,s}^* B_1^*\right\rangle\right\vert\\
&\quad \le \lVert B_1^*\Im G_{1,s}B_1\rVert \left\vert \left\langle G_{2,s} B_2 G_{1,s} B_3 \Im G_{1,s} B_3^* G_{1,s}^* B_2^* G_{2,s}^* \right\rangle\right\vert\\
&\quad \le \lVert B_1^*\Im G_{1,s}B_1\rVert \lVert B_3\Im G_{1,s}B_3^*\rVert \left\vert \left\langle G_{2,s} B_2 G_{1,s}G_{1,s}^* B_2^* G_{2,s}^* \right\rangle\right\vert\le \eta_{*,s}^{-4}\left\vert \left\langle \Im G_{2,s} B_2 \Im G_{1,s} B_2^* \right\rangle\right\vert\lesssim \eta_{*,s}^{-5}.
\end{split}
\label{eq:QV_3G_eta}
\end{equation}
To go from the first to the second line of \eqref{eq:QV_3G_eta} we used that the resolvent chain in the second factor in the second line of \eqref{eq:QV_3G_eta} is either positive or negative semi-definite, depending on the sign of $\Im w_{1,s}$. Similarly, to go from the second to the third line we used that $G_{2,s} B_2 G_{1,s}G_{1,s}^* B_2^* G_{2,s}^*\ge 0$. Additionally, in the first inequality in the third line of \eqref{eq:QV_3G_eta} we estimated $\|\Im G_{1,s}\|\le \eta_{*,s}^{-1}$ and used the Ward identity twice. Finally, it the last estimate \eqref{eq:M2_bound} and the concentration bound provided by $\tau_2$ were used. This completes the proof of Proposition \ref{prop:err_bounds}. \qed

\subsubsection{Modifications for the case of general $B_i\in\mathrm{span}\{E_\pm, F^{(*)}\}$}\label{sec:general_B}

In Sections \ref{sec:Zig_structure}--\ref{sec:err_bounds} we proved that the bounds \eqref{eq:2G_bound_t}--\eqref{eq:4G_bound_t} propagate in time in the special case when $B_i\in\{E_\pm\}$, $i\in [4]$. Notably, in this argument we do not assume \eqref{eq:2G_bound_t}--\eqref{eq:4G_bound_t} at time $t=0$ for the remaining choices of observables. In the current section we treat the general case $B_i\in\mathrm{span}\{E_\pm,F^{(*)}\}$ and complete the proof of Proposition \ref{prop:Zig}. To do so, we use the already established concentration bounds \eqref{eq:2G_bound_t}--\eqref{eq:4G_bound_t} for all $t\in [0,T]$ in the case when $B_i=E_\pm$, $i\in [4]$, and slightly adjust the argument presented in Sections \ref{sec:Zig_structure}--\ref{sec:err_bounds} to cover this more general case.

Throughout this section we consider only four choices for each observable, specifically $B_i\in \{E_\pm, F^{(*)}\}$, $i\in [4]$. This does not decrease the generality of the set-up, since any other observable in Proposition \ref{prop:Zig} is a linear combination of these four, and the lhs. of \eqref{eq:2G_bound_t}--\eqref{eq:4G_bound_t} are linear in each of the observables. In fact, when some of the observables in the resolvent chain equal to $F^{(*)}$, the fluctuation of this chain around its deterministic approximation is expected to be smaller, see e.g. \cite[Theorem~4.4]{non-herm_overlaps} and \cite[Theorem~3.5]{univ_extr}. However, this improvement is not needed for our purposes, and we do not pursue it further. This makes the proof of Proposition \ref{prop:Zig} for $B_i\in \{E_\pm, F^{(*)}\}$ fairly analogous to the proof in the special case $B_i\in \{E_\pm\}$. These two arguments can be even combined into one, but we keep them separate for clarity of the presentation.

Throughout the proof we freely use the notations introduced in Section \ref{sec:Zig_structure} without mentioning this further. We first consider a general set-up introduced above \eqref{eq:kG_dif}, i.e. we follow the time evolution of $k\ge 2$ resolvents $G_{j,t}=G^{z_{j,t}}_t(w_{j,t})$, $j\in [k]$. The starting point of our analysis is \eqref{eq:kY_dif}. We again treat the terms in the third line of \eqref{eq:kY_dif} as error terms and focus on the structure of the so-called linear terms in the second line of \eqref{eq:kY_dif}. For any choice of $B_j\in\{E_\pm, F^{(*)}\}$, $j\in [k]$, the number of observables equal to $F^{(*)}$ in each of these terms does not exceed the number of $F^{(*)}$ observables among $\{B_j\}_{j=1}^k$. Moreover, the terms in the second line of \eqref{eq:kY_dif} corresponding to the summation index $j$ with $B_j\in\{F^{(*)}\}$, contain strictly less $F^{(*)}$ observables. This suggests the following structural decomposition of the terms in the second line of \eqref{eq:kY_dif} into two groups. 
The sum of terms in the first group is given by
\begin{equation}
\mathrm{L}_{[k],t}:=\sum_{j=1}^k \sigma\mathbf{1}\left\lbrace B_j\in\{F^{(*)}\}\right\rbrace \left\langle M_{[j,j+1],t}^{B_j}E_\sigma\right\rangle Y_{[k],t}(B_1,\ldots, B_{j-1},E_\sigma, B_{j+1},\ldots, B_k)\dif t,
\label{eq:def_L}
\end{equation}
while the remaining terms contain the same number of $F^{(*)}$ observables as the lhs. of \eqref{eq:kY_dif}. We treat \eqref{eq:def_L} as an error term, while the rest of the terms in the second line of \eqref{eq:kY_dif} contribute to the propagator of the system of equations analogous to \eqref{eq:kG_syst}, which we now construct.

Unlike in Section \ref{sec:Zig_structure}, instead of one system for each $k\ge 2$, we construct several systems of stochastic differential equations analogous to \eqref{eq:kG_syst} to accommodate observables equal to $F^{(*)}$. These systems are labeled  by the set of indices $\mathcal{J}\subset [k]$, which indicates the indices of observables equal to $F^{(*)}$ in $Y_{[k],t}(B_1,\ldots,B_k)$, and by the exact choice of $\{B_j\}_{j\in \mathcal{J}} \in \{F^{(*)}\}^{|\mathcal{J}|}$. Denote
\begin{equation*}
f:=|\mathcal{J}|,\quad \mathcal{J}:=\{j_1,\ldots,j_f\},\quad [k]\setminus\mathcal{J}:=\{i_1,\ldots,i_{k-f}\}.
\end{equation*}
We fix $B_{j}\in \{F^{(*)}\}$ for all $j\in \mathcal{J}$ and for the remaining indices $i\in [k]\setminus\mathcal{J}$ consider all $2^{k-f}$ possible choices $B_i\in\{E_\pm\}$. Similarly to \eqref{eq:def_Y_vector} define
\begin{equation}
\mathcal{Y}_{[k],\mathcal{J},t}=\mathcal{Y}_{[k],\mathcal{J},t}(B_{j_1},\ldots,B_{j_f}):=\sum_{\sigma_{i_1},\ldots,\sigma_{i_{k-f}}\in\{\pm\}} Y_{[k],t}(\cdots) \bm{e}_{\sigma_{i_1},\ldots,\sigma_{i_{k-f}}}\in (\C^2)^{\otimes (k-f)},
\label{eq:def_YJ_vector}
\end{equation}
where the arguments of $Y_{[k],t}$ in the rhs. of \eqref{eq:def_YJ_vector} are given by $\{B_j'\}_{j=1}^k$ with $B_j'=B_j$ for $j\in \mathcal{J}$ and $B_j=E_{\sigma_j}$ for $j\notin \mathcal{J}$. Define further the $(\C^2)^{\otimes (k-f)}$ vectors $\mathcal{L}_{[k],\mathcal{J},t}$, $\mathcal{F}_{[k],\mathcal{J},t}$ and $\dif\mathcal{E}_{[k],\mathcal{J},t}$ analogously to \eqref{eq:def_YJ_vector}, from \eqref{eq:def_L}, \eqref{eq:def_Fk} and \eqref{eq:def_E}, respectively. With these notations we get from \eqref{eq:kG_dif} the following analogue of \eqref{eq:kG_syst}:
\begin{equation}
\dif \mathcal{Y}_{[k],\mathcal{J},t} = \left(\frac{k}{2} + \mathcal{A}_{[k],\mathcal{J},t}\right)\mathcal{Y}_{[k],\mathcal{J},t} \dif t + \mathcal{L}_{[k],\mathcal{J},t}\dif t + \mathcal{F}_{[k],\mathcal{J},t} \dif t + \dif\mathcal{E}_{[k],\mathcal{J},t},
\label{eq:kJ_syst}
\end{equation}
where the time-dependent linear operator $\mathcal{A}_{[k],\mathcal{J},t}:(\C^2)^{\otimes (k-f)}\to (\C^2)^{\otimes (k-f)}$ is given by
\begin{equation}
\mathcal{A}_{[k],\mathcal{J},t}: = \sum_{\sigma_{i_1},\ldots,\sigma_{i_{k-f}}\in\{\pm\}} \bigg(\sum_{l=1}^{k-f} a_{{i_l}(i_{l}+1),t}^{\sigma_{i_l}}\bigg)P_{\sigma_{i_1},\ldots,\sigma_{i_{f-k}}} - \sum_{l=1}^{k-f} d_{i_l(i_l+1),t}S_{l}.
\label{eq:def_general_AJ}
\end{equation}

Importantly, no additional analysis of the operator $\mathcal{A}_{[k],\mathcal{J},t}$ is required compared to the one performed in Lemma \ref{lem:propag_bound} in the special case $f=0$. Analogously to \eqref{eq:propagator} we denote
\begin{equation*}
f_{[k],\mathcal{J},r}:=\left[ \max \mathrm{spec}(\Re \mathcal{A}_{[k],\mathcal{J},r})\right]_+,\quad \forall\, r\in[0,T]. 
\end{equation*}
Then we have
\begin{equation}
\exp\left\lbrace\int_s^t f_{[k],\mathcal{J},r} \dif r\right\rbrace \lesssim \prod_{l=1}^{k-f} \frac{\beta_{i_l(i_l+1),*,s}}{\beta_{i_l(i_l+1),*,t}}\lesssim \prod_{i=1}^{k} \frac{\beta_{i(i+1),*,s}}{\beta_{i(i+1),*,t}}.
\label{eq:propag_F}
\end{equation}
Here the first bound is derived from \eqref{eq:a_bound} and \eqref{eq:d_bound} similarly to \eqref{eq:propag_bound} (see also \eqref{eq:f2_bound}), while the second one follows from the first part of \eqref{eq:gamma_monot}. In principle, \eqref{eq:propag_F} is an overestimate, and instead of working with the rightmost expression in \eqref{eq:propag_F} one could perform the refined analysis using the middle expression in \eqref{eq:propag_F}, which would lead to the improvement mentioned in the second paragraph of this section.

By definition, $\mathcal{L}_{[k],\mathcal{J},t}$ in the rhs. of \eqref{eq:kJ_syst} is constructed from the fluctuations of multi-resolvent chains with strictly less than $f$ observables. It is easy to see from \eqref{eq:def_Fk} that the same holds for the forcing term $\mathcal{F}_{[k],\mathcal{J},t}$. However, the quadratic variation of the martingale term $\dif \mathcal{E}_{[k],\mathcal{J},t}$ typically contains resolvent chains with more than $f$ observables equal to $F^{(*)}$. This prevents us from analysing \eqref{eq:kJ_syst} inductively, obtaining at first the upper bound on $Y_{[k],t}(B_1,\ldots,B_k)$ with only one $F^{(*)}$ observable, and then concluding the bound in the case of $f$ observables $F^{(*)}$, relying on the bounds for all smaller numbers of the off-diagonal observables. The resolution of this issue is to control \eqref{eq:kJ_syst} at the same time for all $\mathcal{J}\subset [k]$ with $|\mathcal{J}|\ge 1$ and for all choices $\{B_j\}_{j\in\mathcal{J}}\in \{F^{(*)}\}^{|\mathcal{J}|}$.

Now we restrict the setting to the one introduced in the beginning of Section \ref{sec:Zig_old}. In particular, from now on $k\in \{2,3,4\}$, $z_{3,t}=z_{1,t}$ and $z_{4,t}=z_{2,t}$. Similarly to \eqref{eq:Y_short} we denote
\begin{equation}
\mathcal{Y}^{(2)}_{\mathcal{J},t}:=\mathcal{Y}_{12,\mathcal{J},t},\quad \mathcal{Y}^{(3)}_{\mathcal{J},t}:=\mathcal{Y}_{121,\mathcal{J},t},\quad \mathcal{Y}^{(4)}_{\mathcal{J},t}:=\mathcal{Y}_{1212,\mathcal{J},t}, \quad \forall \mathcal{I}\subset [k],
\end{equation}
and use the same notational convention for $\mathcal{A}_{[k], \mathcal{J},t}$, $f_{[k],\mathcal{J},t}$, $\mathcal{L}_{[k],\mathcal{J},t}$, $\mathcal{F}_{[k],\mathcal{J},t}$ and $\mathcal{E}_{[k],\mathcal{J},t}$. We fix (small) tolerance exponents $\xi_{k,m}$ for $k\in \{2,3,4\}$, $m\in [k]$, such that
\begin{equation}
\xi_{k,m}\in (0,\epsilon/10),\quad \xi_{k,1}< \xi_{k,2}<\ldots< \xi_{k,k} \quad\text{and}\quad \xi_{4,4}/2< \xi_{2,1}< \xi_{2,2} < \xi_{4,1},
\label{eq:xi_new}
\end{equation}
for all $k\in \{2,3,4\}$ and $m\in [k]$, and for the same range of $(k,m)$ define the stopping times $\tau_{k,m}$ by
\begin{equation}
\tau_{k,m}\!:= \inf\left\lbrace \!t\!\in\! [0,T]\,:\, \max_{|\mathcal{J}|=m}\,\max_{B_j\in \{F^{(*)}\},j\in\mathcal{J}}\,\max_{s\in [0,t]}\max_{|z_{j,T}|\le \ee^{-T/2}(1-\delta)}\max_{w_{j,T}\in\Omega_{\kappa,\epsilon,T}^{z_{j,T}}} \alpha_{k,s}^{-1}\|\mathcal{Y}^{(k)}_{\mathcal{J},s}\|\ge N^{2\xi_{k,m}}\!\right\rbrace\!,
\label{eq:def_tau_k_new}
\end{equation}
where $\mathcal{J}=\{j_1,\ldots,j_m\}\subset [k]$ and $\mathcal{Y}^{(k)}_{\mathcal{J},s} = \mathcal{Y}^{(k)}_{\mathcal{J},s}(B_{j_1},\ldots,B_{j_m})$. The control parameters $\alpha_{k,s}$ for $k\in \{2,3,4\}$ are defined in \eqref{eq:def_alpha} and do not depend on $\mathcal{J}$. Finally, we set
\begin{equation}
\widetilde{\tau}:=\min_{k\in \{2,3,4\}}\min_{m\in [k]} \tau_{k,m}.
\end{equation}
We claim that
\begin{equation}
\left(\int_0^{t\wedge\widetilde{\tau}}\left(\|\mathcal{L}_{\mathcal{J},s}^{(k)}\| + \|\mathcal{F}_{\mathcal{J},s}^{(k)}\|\right)\dif s\right)^2 +\int_0^{t\wedge\widetilde{\tau}} \|\mathcal{C}^{(k)}_{\mathcal{J},s}\|\dif s \lesssim  N^{4\xi_{k,|\mathcal{J}|}-\xi}\alpha_{k,t\wedge\widetilde{\tau}}^2,
\label{eq:err_bounds_new}
\end{equation}
uniformly in $t\in [0,T]$, $k\in \{2,3,4\}$, $\mathcal{J}\subset [k]$ with $|\mathcal{J}|\ge 1$ for some exponent $\xi>0$ which depends only on the exponents introduced in \eqref{eq:xi_new}. In \eqref{eq:err_bounds_new} we denoted by $\mathcal{C}^{(k)}_{\mathcal{J},s}\dif s\in\C^{2^{k-m}\times 2^{k-m}}$, $m:=|\mathcal{J}|$, the covariation process of $\dif\mathcal{E}^{(k)}_{\mathcal{J},s}$. We also omitted for brevity the arguments of $\mathcal{L}_{\mathcal{J},s}^{(k)}$, $\mathcal{F}_{\mathcal{J},s}^{(k)}$ and $\mathcal{C}^{(k)}_{\mathcal{J},s}$, which are given by  $B_{j_1},\ldots, B_{j_m}\in\{F^{(*)}\}$, where $\mathcal{J}=\{j_1,\ldots,j_m\}$. The bound \eqref{eq:err_bounds_new} is uniform in these arguments.

The proof of \eqref{eq:err_bounds_new} is completely analogous to the proof of Proposition \ref{prop:err_bounds}, since there we never use that $B_j\in \{E_\pm\}$ for $j\in [k]$. The only difference lies in the additional term $\mathcal{L}_{\mathcal{J},s}^{(k)}$, which is directly controlled by the stopping time $\widetilde{\tau}$ for $s\in [0,\widetilde{\tau}]$. Having \eqref{eq:err_bounds_new} in hand, we finish the proof of Proposition~\ref{prop:Zig} similarly to the argument presented in Section \ref{sec:Zig_old} relying on the stochastic Gr{\"o}nwall inequality and on~\eqref{eq:propag_F}.

\subsection{Proof of Proposition~\ref{prop:multiG_oneG}}
\label{sec:multiG_oneG}

To keep the presentation short we only prove this result for matrices with an order one Gaussian component. Then, this additional Gaussian component can be easily removed by a standard GFT argument. Since this proof is similar to the proof of Proposition~\ref{prop:Zig}, and in fact much simpler, we will omit several details and only present the main steps of the proof (we also refer to \cite{sparse_universality_Osman} for an analogous proof in the special case when all the $w_i$'s are on the imaginary axis). For $\eta_*\gtrsim 1$, the local laws in \eqref{eq:multiG_av}--\eqref{eq:multiG_iso} follow analogously to the proof of Proposition~\ref{prop:global_law}. We thus now focus on the proof that these local laws can be propagated down to $N\eta_*\ge N^\epsilon$. For this purpose we consider the flow \eqref{eq:OU}, define $G_t^z(w)$ as in \eqref{eq:defres}, and use the short-hand notation $G_{i,t}:=G_t^{z_i}(w_i)$. Denote (here we suppress the dependence on the $B_i$'s from the notation)
\begin{equation}
	{{G}}_{[{i}, {j}],t}:=\begin{cases}
		G_{i,t}B_i\dots B_{j-1}G_{j,t} & \mathrm{if}\quad  i<j \\
		G_{i,t} & \mathrm{if}\quad i=j \\
		G_{i,t}B_{i,t}\dots G_{k,t}B_k G_{1,t}B_1\dots B_{j-1} G_{j,t} &\mathrm{if}\quad i>j\,,
	\end{cases}
\end{equation}
and by $M_{[i,j],t}$ its deterministic approximation (see e.g. \cite[Lemma D.1]{non-herm_overlaps} for a recursive relation  for this deterministic term). For this deterministic approximation we have (this follows by a simple \emph{meta argument} as in the proof of Lemma~\ref{lem:M_bounds_weak})
\begin{equation}
\label{eq:Mboundeas}
\lVert M_{[i,j],t} \rVert\lesssim \frac{1}{\eta_{*,t}^{j-i}}.
\end{equation}
Additionally, we define $G_{[{i}, {j}],t}^{(l)}$ exactly as $G_{[{i},{j}],t}$ but with the $l$--th factor $G_{l,t}$ being
replaced with $G_{l,t}^2$.  Then, by It\^{o}'s formula and \eqref{eq:kM_dif} we have
\begin{equation}
	\begin{split}
		\label{eq:flowka}
		\dif \langle (G_{[1,k],t}-M_{{[1,k]},t})B_k\rangle&=\frac{1}{\sqrt{N}}\sum_{a,b=1}^N \partial_{ab}  \langle G_{[1,k],t}B_k\rangle\dif B_{ab,t}+\frac{k}{2}\langle (G_{[1,k],t}-M_{{[1,k]},t})B_k\rangle\dif t \\
		&\quad+\sum_{i,j=1\atop i< j}^k\sigma \langle G_{[i,j],t}E_\sigma\rangle\langle G_{[j,i],t} E_\sigma\rangle\dif t-\sum_{i,j=1\atop i< j}^k\sigma \langle M_{[i,j],t}E_\sigma\rangle\langle M_{[j,i],t} E_\sigma\rangle\dif t \\
		&\quad+\sum_{i=1}^k \langle G_{i,t}-m_{i,t}\rangle \langle G^{(i)}_{[1,k],t}B_k\rangle\dif t.
	\end{split}
\end{equation}
For the special choice $B_k=N{\bm x}^* {\bm y}$, \eqref{eq:flowka} gives the evolution for isotropic quantities. To keep the presentation short, from now on we only focus on proving the averaged local law \eqref{eq:multiG_av}, the proof for the isotropic local law \eqref{eq:multiG_iso} being completely analogous.

Define the stopping time
\begin{equation}
\tau:=\inf\left\{t\in [0,1]: \big|\langle (G_{[i,j],t}-M_{{[i,j]},t})B\rangle\big|=\frac{N^{2\xi_{j-i+1}}}{N\eta_{*,t}^{j-i+1}}\,\, \forall i<j\le k\right\},
\end{equation}
for small $\xi_l\le \epsilon/10$ with $\xi_{l-1}\le \xi_{l/2}$. From now on, we always assume that $t\le \tau$. Additionally, throughout this proof we will use \eqref{eq:eta_int}, even if not stated explicitly. We now discuss the estimate in the rhs. of \eqref{eq:flowka} one by one. The second term in the rhs. of \eqref{eq:flowka} can be easily neglected, as it amounts to a rescaling $e^{-kt/2}$, by a simple change of variables.

For the term in the third line of \eqref{eq:flowka} we estimate (assume $i\ne 1$ for simplicity of notation)
\[
\big|\langle G_{i,s}-m_{i,s}\rangle\big|\le \frac{N^\xi}{N\eta_{i,s}}, \qquad\quad \big|\langle G^{(i)}_{[1,k],s}B_k\rangle\big|\le \langle G_1B_1B_1^*G_1^*\rangle^{1/2} \langle G^{(i)}_{[2,k],s}B_kB_k^*(G^{(i)}_{[2,k],s})^*\rangle^{1/2}\lesssim \frac{1}{\eta_{*,s}^k},
\]
where we also used the deterministic norm bound $\lVert G_{i,t}\rVert\le \eta_{i,t}^{-1}$ (and that the $B_i$'s have bounded norms). We thus obtain
\begin{equation}
\left|\int_0^t\sum_{i=1}^k \langle G_{i,s}-m_{i,s}\rangle \langle G^{(i)}_{[1,k],s}B_k\rangle\, \dif s\right|\lesssim N^\xi \int_0^t \frac{1}{N\eta_{i,s}}\frac{1}{\eta_{*,s}^k}\, \dif s\lesssim \frac{N^\xi}{N\eta_{*,t}^k}.
\end{equation}

For the terms in the second line of \eqref{eq:flowka}, using \eqref{eq:Mboundeas}, we estimate (for simplicity of notation we consider $i=1$)
\begin{equation}
\begin{split}
\left|\sigma\int_0^t\left\langle \big(G_{[1,j],t}-M_{[1,j],t}\big)E_\sigma\right\rangle\langle M_{[j,1],t} E_\sigma\rangle\,\dif s\right|&\lesssim \int_0^t \frac{N^{\xi_j}}{N\eta_{*,s}^j} \frac{1}{\eta_{*,s}^{k+1-j}}\,\dif s\lesssim \frac{N^{\xi_j}}{N\eta_{*,t}^k}, \\
\left|\sigma\int_0^t\langle M_{[1,j],t}E_\sigma\rangle\left\langle \big(G_{[j,1],t}-M_{[j,1],t}\big) E_\sigma\right\rangle\,\dif s\right|&\lesssim \frac{N^{\xi_{k+1-j}}}{N\eta_{*,t}^k},\\
\left|\sigma\int_0^t\left\langle \big(G_{[1,j],t}-M_{[1,j],t}\big)E_\sigma\right\rangle \left\langle \big(G_{[j,1],t}-M_{[j,1],t}\big) E_\sigma\right\rangle\,\dif s\right|&\lesssim \int_0^t \frac{N^{\xi_j}}{N\eta_{*,s}^j} \frac{N^{\xi_{k+1-j}}}{N\eta_{*,s}^{k+2-j}}\,\dif s\\
&\lesssim \frac{N^{\xi_j+\xi_{k+1-j}}}{N\eta_{*,t}} \frac{1}{N\eta_{*,t}^k}.
\end{split}
\end{equation}

Finally, we consider the first term in the rhs. of \eqref{eq:flowka}. Similarly to the computations above, it is easy to see that its quadratic variation process is bounded by $(N^2\eta_{*,s})^{-2k-1}$, with very high probability. We thus obtain
\begin{equation}
\sup_{0\le s\le t}\left|\int_0^s \frac{1}{\sqrt{N}}\sum_{a,b=1}^N \partial_{ab}  \langle G_{[1,k],r}B_k\rangle\dif B_{ab,r}\right|\lesssim \log N\left(\int_0^t \frac{1}{N^2\eta_{*,s}^{2k+1}}\,\dif \right)^{1/2}\lesssim  \frac{\log N}{N\eta_{*,s}^k}.
\end{equation}

Combining all these estimates  we obtain
\begin{equation}
\big|\langle (G_{[i,j],t}-M_{{[i,j]},t})B\rangle\big|\lesssim\frac{N^{\xi_{j-i+1}}}{N\eta_{*,t}^{j-i+1}},
\end{equation}
with very high probability for any $i<j\le k$, thus concluding the proof of Proposition~\ref{prop:multiG_oneG}.\hfill$\qed$

\subsection{Proof of Proposition \ref{prop:2G_av} for general $b\in [0,1]$}\label{sec:general_b} 
The statement of Proposition \ref{prop:2G_av} for $b=1$ immediately follows from the averaged local law in Proposition~\ref{prop:multiG_oneG} with $k=2$, so we further focus on the case $b\in (0,1)$. Since we proceed by a minor modification of the proof of Proposition \ref{prop:2G_av} for $b=0$, we first remind the reader the overall structure of that proof, consisting of two main steps. The first step is the global law from Proposition~\ref{prop:global_law} which establishes the desired bounds on the fluctuations of resolvent chains for spectral parameters at distance of order one from the real axis. In the second step (see Proposition~\ref{prop:Zig}) we propagate these estimates down to the real axis via the characteristic flow. The proof of Proposition \ref{prop:2G_av} for $b\in (0,1)$ has the same structure up to the following adjustments, which we now describe.

We fix the final time $T:=N^{-b}$, which is consistent with the choice $T\sim 1$ made for $b=0$, and keep the definition of the bulk truncated spectral domains $\{\Omega_{\kappa,\epsilon,t}^{z_T}\}_{t\in [0,T]}$ from Definition~\ref{def:spec_dom}. By Lemma~\ref{lem:char_flow}(4), for any $|z_T|\le 1-\delta$ and $w\in \Omega_{\kappa,\epsilon,0}^{z_T}$ it holds that $|\Im w|\gtrsim N^{-b}$. In particular, the global law from Proposition~\ref{prop:global_law} cannot be applied in this regime, so we use Proposition~\ref{prop:multiG_oneG} instead for the lhs. of \eqref{eq:global_2G}--\eqref{eq:global_4G}. This estimate is then used as the initial condition, which we propagate towards the real axis similarly to Proposition~\ref{prop:Zig}. To adapt this proposition to the current setting, we introduce the new control parameter:
\begin{equation}
\widehat{\beta}^{[b]}_{12}(w_1,w_2):=\widehat{\beta}_{12}(w_1,w_2)\wedge N^{-b} + |\Im w_1|\wedge|\Im w_2|\wedge 1,
\label{eq:beta_hat_b}
\end{equation}
and define its time-dependent version $\widehat{\beta}^{[b]}_{12,t}$ for $t\in [0,T]$ similarly to \eqref{eq:control_param_t}. Then Proposition~\ref{prop:Zig} holds with $\widehat{\beta}_{12,t}$ replaced by $\widehat{\beta}_{12,t}^{[b]}$. To see this we observe that the proof of Proposition~\ref{prop:Zig} does not directly rely on the definition of $\widehat{\beta}_{12}$ given in \eqref{eq:def_beta_hat}, but it only relies on the properties of this control parameter stated in \eqref{eq:M_12_beta_bound}, \eqref{eq:beta_hat_trivial}, and Lemma~\ref{lem:beta_hat}. Trivially, $\widehat{\beta}_{12}^{[b]}$ satisfies these properties as well.

It remains to observe that $\widehat{\beta}^{[b]}_{12,0}\sim \eta_{*,0}$, so \eqref{eq:2G_bound_t}--\eqref{eq:4G_bound_t} folds for $t=0$ with $\widehat{\beta}_{12,t}$ replaced by $\widehat{\beta}_{12,t}^{[b]}$, by Proposition~\ref{prop:multiG_oneG}. Propagating these bounds up to the final time $T$, we finish the proof of Proposition~\ref{prop:2G_av} for $b\in (0,1)$.\hfill$\qed$

\subsection{Proof of Proposition \ref{prop:2G_av_suboptimal}}\label{sec:2G_av_subopt} It is known from \cite[Theorem~3.4]{nonHermdecay} that \eqref{eq:2G_av_suboptimal} holds for spectral parameters on the imaginary axis, i.e. for $\Re w_1=\Re w_2=0$. Thus, we only need to extend this result to the bulk regime $\Re w_i\in\mathbf{B}_\kappa^{z_i}$, $i=1,2$, which we do by a minor modification of the proof of \cite[Theorem~3.4]{nonHermdecay}. The zag step (removal of the Gaussian component) remains unchanged compared in the one in \cite[Section~5]{nonHermdecay}. In the zig step (adding a Gaussian component), we propagate the local law estimate down to the real axis by applying Gronwall inequality to the analogue of \eqref{eq:kG_syst}, instead of integrating these equations directly as it was done in \cite{nonHermdecay} (see e.g. Eq.~(4.28) therein). This allows us to avoid extension of the deterministic estimates \cite[Eq.(4.27)]{nonHermdecay} to the bulk regime, and use Lemma~\ref{lem:propag_bound2} instead.


We note that \cite{nonHermdecay} fully covers the regime $|z_i|\le 1+\delta$, $i=1,2$, for a small $\delta>0$, while Proposition~\ref{prop:2G_av_suboptimal} concerns only $|z_i|\le 1-\delta$. This allows to simplify the proof of \cite[Theorem~3.4]{nonHermdecay} in our set-up by neglecting all $\rho^z$ factors appearing in the estimates, since they are of order~1.

\section{The real case}
\label{sec:real}

In this section we prove Theorem~\ref{theo:main2} by following a similar proof to Theorem~\ref{theo:main1}. The overall  proof structure in the real case remains the same as in the complex case, we thus primarily focus on explaining the main  technical differences and new difficulties. We frequently use the convention that a free index $\sigma$ appearing in formulas is meant to be summed over $\sigma\in\{\pm\}$ if the opposite is not stated. We also use $\xi>0$ to denote an $N$-independent exponent which can be taken arbitrarily small and whose exact value may change from line to line. 

Denote the self-energy operator introduced in the complex case in \eqref{eq:def_S} by $\mathcal{S}_2$, where the subscript reflects the symmetry type $\beta=2$. In the real case we define
\begin{equation}
\mathcal{S}_1[R] := \mathcal{S}_2[R] +\frac{\sigma}{N}E_\sigma R^\mathfrak{t}E_\sigma,\quad \forall R\in\C^{(2N)\times (2N)}.
\label{eq:real_def_S}
\end{equation}
Here and everywhere further in this section we denote the transpose of a matrix $R$ by $R^\mathfrak{t}$. Another quantity closely related to \eqref{eq:real_def_S} which changes compared to the complex case is the bilinear form
\begin{equation}
\E \langle WR\rangle \langle WS\rangle = \frac{\sigma}{4N^2}\langle RE_\sigma SE_\sigma \rangle + \frac{\sigma}{4N^2}\langle RE_\sigma S^\mathfrak{t}E_\sigma \rangle,
\label{eq:real_EWW_id}
\end{equation}
cf. with \eqref{eq:EWW_id} in the complex case, where only the first term in the rhs. of \eqref{eq:real_EWW_id} is present. Though $\mathcal{S}_1$ differs from $\mathcal{S}_2$, we still use the MDE \eqref{eq:MDE} with $\mathcal{S} = \mathcal{S}_2$ and the notions of deterministic approximations to resolvent chains introduced in Section~\ref{sec:local_laws}, without adapting them to the real case. As our proof shows, this slight conceptual mismatch is still affordable as the new term in \eqref{eq:real_def_S} only amounts to a negligible contribution. The advantage, in particular, is that the two-body stability operator $\mathcal{B}_{12}$ defined in \eqref{eq:def_B12} is unchanged, so we do not need to redo the stability analysis from Proposition~\ref{prop:stab_bound}, as well as the rest of the analysis of deterministic approximations along the proof of Theorem~\ref{theo:main1}.

Due to the real symmetry of the model, we often encounter transposes of resolvents in the proof of Theorem~\ref{theo:main2}. These can be written again as resolvents, but with conjugated Hermitization parameter:
\begin{equation}
\left(G^z(w)\right)^\mathfrak{t} = G^{\overline{z}}(w).
\label{eq:G_symmetry}
\end{equation}
Thus, it is convenient to denote the deterministic approximation to $G_1B_1G_2^\mathfrak{t}$ with $G_j=G^{z_j}(w_j)$ for $j=1,2$, by $M_{1\overline{2}}^{B_1}(w_1,w_2)$. In other words, a bar over a subscript $j$ means that we are considering the deterministic approximation of the product of two resolvents when $G^{z_j}(w_j)$ is replaced by $G^{\overline{z}_j}(w_j)$. We will use this convention also for $\widehat{\beta}_{12}$, which is defined in \eqref{eq:def_beta_hat}, and for the deterministic approximations to longer resolvent chains.

This section is structured as follows. First, we state the real case analogues of the technical results listed in Sections~\ref{sec:large_eta}--\ref{sec:interm_eta}. Next, in Section~\ref{sec:theo2} we prove Theorem~\ref{theo:main2} relying on these results. Finally, in Sections~\ref{sec:real_local_laws}--\ref{sec:real_crit_scale} we adapt the proofs of the technical ingredients presented in Sections~\ref{sec:proofprop}--\ref{sec:Cov} to the real case.

We start with stating the analogue of Proposition \ref{prop:Cov} for real i.i.d. matrices, postponing the proof to Section~\ref{sec:real_chaos_exp}.

\begin{proposition}\label{prop:real_Cov} Let $X$ be a real $N\times N$ i.i.d. matrix satisfying Assumption~\ref{ass:chi}. Denote $\kappa_4:=\E|\chi|^4-3$. Fix (small) $\delta,\delta_r, \epsilon, \kappa, \xi>0$. Uniformly in $z_l\in (1-\delta)\mathbf{D}$ with $|\Im z_l|\ge N^{-1/2+\delta_r}$, and $w_l\in\C\setminus\R$ with $E_l:=\Re w_l\in \mathbf{B}_\kappa^{z_l}$ and $\eta_l\in [N^{-1+\epsilon}, 1]$, $l=1,2$, it holds that
\begin{equation}
\begin{split}
\Cov \left(\langle G^{z_1}(w_1)\rangle,\langle G^{z_2}(w_2)\rangle\right) &= \frac{1}{N^2}\cdot \frac{\widehat{V}_{12}+\kappa_4U_1U_2}{2}\\
& + \mathcal{O}\left(\left(\frac{1}{N\widehat{\gamma}}+\frac{1}{N|\Im z_1|^2}+\frac{1}{N|\Im z_2|^2} +N^{-1/4}\right)\frac{N^\xi}{N^2\eta_1\eta_2}\right),
\end{split}
\label{eq:real_main_cov}
\end{equation}
where $\widehat{V}_{12}=\widehat{V}_{12}(z_1,z_2,w_1,w_2)$ and $\widehat{\gamma}=\widehat{\gamma}(z_1,z_2,w_1,w_2)$ are defined as
\begin{equation}
\begin{split}
\widehat{\gamma}(z_1,z_2,w_1,w_2):=&\min_{z_1'\in \{z_1,\overline{z}_1\}}\min_{z_2'\in\{z_2,\overline{z}_2\}} \gamma(z_1',z_2',w_1,w_2),\\
\widehat{V}_{12}(z_1,z_2,w_1,w_2):=&V_{12}(z_1,z_2,w_1,w_2)+V_{12}(z_1,\overline{z}_2,w_1,w_2),
\end{split}
\label{eq:real_def_UV}
\end{equation} 
while $V_{12}$, $U_l$ for $l=1,2$, and $\gamma$ are defined in \eqref{eq:def_UV}.
\end{proposition}

The leading order term in the rhs. of \eqref{eq:real_main_cov} was initially computed in \cite[Proposition~3.3]{macroCLT_real}, and our contribution consists in the improvement of the bound on the error term. Compared to \eqref{eq:main_cov}, the bound on the error term in the rhs. of \eqref{eq:real_main_cov} contains the additional terms $(N|\Im z_j|^2)^{-1}$ which blow up as $z_j$ approaches the real line. In fact, these terms are present only for technical reasons and they can be removed by extending the hierarchy of covariances introduced in the proof of Proposition~\ref{prop:Cov}. Specifically, one would need to include into consideration the quantities \eqref{eq:def_delta}, where some of the resolvents $G_1$ are replaced by $G_1^\mt$. However, in order to  keep the adjustments  in the real case  minimal, we use the same hierarchy of covariances, which leads to the deterioration of the estimate on the error term. 

To analyze the sub-microscopic regime ($\eta\ll N^{-1}$) in the Girko's formula \eqref{eq:Girko}, we establish the following real version of Proposition~\ref{prop:tail_bound}.

\begin{proposition}[Left tail of the least positive eigenvalue distribution; real case]\label{prop:real_tail_bound} Let $X$ be a real $N\times N$ i.i.d. matrix satisfying Assumption~\ref{ass:chi}. Fix (small) $\delta, \delta_r, \xi>0$ and set $\nu_1:=1/10$ as in Proposition~\ref{prop:tail_bound}. Uniformly in the parameters
\begin{equation*}
z\in (1-\delta)\mathbf{D}, \,|\Im z|\ge N^{-1/2+\delta_r}\quad \text{and}\quad N^\xi\left(N^{-\nu_1} + \left(N^{1/2}|\Im z|\right)^{-1/2}\right)\le x\le 1,
\end{equation*}
it holds that
\begin{equation}
\mathbf{P}[\lambda_1^z\le N^{-1}x]\lesssim \log N \cdot x^2.
\label{eq:real_tail_bound}
\end{equation}
\end{proposition}

In other words, Proposition~\ref{prop:real_tail_bound} asserts that
\begin{equation}
\mathbf{P}[\lambda_1^z\le N^{-1}x]\lesssim \log N \cdot x^2 + N^\xi\left(N^{-2\nu_1} + \left(N^{1/2}|\Im z|\right)^{-1}\right),
\label{eq:real_tail_bound_alternative}
\end{equation}
for any fixed $\xi>0$, uniformly in $0<x\le 1$. This bound immediately follows from \eqref{eq:real_tail_bound} and the fact that the lhs. of \eqref{eq:real_tail_bound} monotonically increases in $x$. Note that compared to the complex analogue~\eqref{eq:tail_bound}, the estimate \eqref{eq:real_tail_bound_alternative} contains the additional term $\left(N^{1/2}|\Im z|\right)^{-1}$, which makes \eqref{eq:real_tail_bound_alternative} ineffective for $|\Im z|$ just slightly above $N^{-1/2}$.

The real version of Proposition~\ref{prop:EG} required for the analysis of the microscopic regime ($\eta\sim N^{-1}$) in the Girko's formula \eqref{eq:Girko} is as follows.

\begin{proposition}\label{prop:real_EG} Let $X$ be a real $N\times N$ i.i.d. matrix satisfying Assumption~\ref{ass:chi}, and let $\widetilde{X}$ be an $N\times N$ GinOE matrix. Let $\widetilde{H}^z$ and $\widetilde{G}^z$ be defined as in Proposition \ref{prop:EG}. Fix (small) $\delta,\delta_r,\epsilon, \omega_*>0$. Then, uniformly in $z\in (1-\delta)\mathbf{D}$ with $|\Im z|\ge N^{-1/2+\delta_r}$ and $\eta\in [N^{-3/2+\epsilon},1]$, it holds that 
\begin{equation}
\E\langle G^z(\ii\eta)\rangle = \E \langle \widetilde{G}^z(\ii\eta)\rangle + \mathcal{O}(N^\xi\Phi_1^r(\eta,z))
\label{eq:real_EG_comp}
\end{equation}
for any fixed $\xi>0$, where $\Phi_1^r=\Phi_1^r(\eta,z)$ is defined by
\begin{equation}
\begin{split}
\Phi_1^r:=&\min_{N^{-1+\omega_*}\le t_2\le t_1\le N^{-\omega_*}}\bigg(N\mathcal{E}_0(t_2)\left(1+N\eta + \frac{N\mathcal{E}_0(t_2)}{N\eta}\right)\\
 &+Nt_2\bigg(1+\left(\frac{1}{N\eta}\right)^3\bigg)\left(\frac{1}{N|\Im z|^2} +\frac{1}{Nt_1}\right) +\sqrt{N}t_1\left(1+\frac{1}{N\eta}\right)^4\bigg),
\end{split}
\label{eq:real_def_Phi}
\end{equation}
with $\mathcal{E}_0$ defined in \eqref{eq:def_Phi}.
\end{proposition}

The main difference of Proposition~\ref{prop:real_EG} from Proposition~\ref{prop:EG} lies in the first term in the second line of \eqref{eq:real_def_Phi}, while the remaining terms in the definition of the control parameter \eqref{eq:real_def_Phi} already appeared in \eqref{eq:def_Phi}. This additional term arises from the fact that we do not derive the analogue of the DBM relaxation from \cite[Proposition~4.6]{bourgade2024fluctuations} in the real case but instead directly use \cite[Proposition~4.6]{bourgade2024fluctuations}. This proof strategy necessitates adding a small complex Gaussian component to the real matrix $X$. Removing this component by GFT leads to the aforementioned deterioration of the bound on the error term. For more details see Section~\ref{sec:real_crit_scale}. Alternatively, one could establish Proposition~\ref{prop:real_EG} by developing a real analogue of \cite[Proposition~4.6]{bourgade2024fluctuations}. However, the proof of this result would not follow by a simple modification of the argument in \cite{bourgade2024fluctuations} and would require a substantial separate analysis, which we avoid here for the sake of brevity.

The second ingredient required for the analysis of the microscopic regime is the real analogue of Proposition~\ref{prop:EGG}, which we now state.

\begin{proposition}\label{prop:real_EGG} Let $X$ be a real $N\times N$ i.i.d. matrix satisfying Assumption~\ref{ass:chi}. Fix (small) $\delta,\delta_r,\epsilon, \omega_*~\!>~\!0$. Then uniformly in $z_l\in (1-\delta)\mathbf{D}$ with $|\Im z_l|\ge N^{-1/2+\delta_r}$ and $\eta_l\in[ N^{-3/2+\epsilon},1]$, $l=1,2$, it holds that
\begin{equation}
\left\vert\mathrm{Cov}(\langle G^{z_1}(\ii\eta_1)\rangle,\langle G^{z_2}(\ii\eta_2)\rangle)\right\vert \lesssim N^\xi \Phi_2^r(\eta_1,\eta_2,z_1,z_2), 
\label{eq:real_EGG}
\end{equation}
for any fixed $\xi>0$, where $\Phi_2^r=\Phi_2^r(\eta_1,\eta_2,z_1,z_2)$ is defined by
\begin{equation}
\begin{split}
\Phi_2^r:=&\min_{N^{-1+\omega_*}\le t_2\le t_1\le N^{-\omega_*}}\min_{0\le R\le t_2|z_1-z_2|^2}\bigg( N\left(\mathcal{E}_1(t_2,R)+\mathcal{E}_0(t_2)\right)\left(N\eta_1+\frac{1}{N\eta_1}\right)\left(N\eta_2+\frac{1}{N\eta_2}\right)\\
&+ Nt_2 \left(\frac{1}{N|\Im z_1|^2}+\frac{1}{N|\Im z_2|^2}+\frac{1}{N|z_1-\overline{z}_2|^2}+\frac{1}{Nt_1}\right)\left(1+\frac{1}{N\eta_1}\right)^2\left(1+\frac{1}{N\eta_2}\right)^2\\
&+\sqrt{N}t_1 \left(1+\frac{1}{N\eta_*}\right)^3 \frac{1}{N^2\eta_1\eta_2}\bigg),
\end{split}
\label{eq:real_def_Phi2}
\end{equation}
with $\eta_*:=\eta_1\wedge\eta_2$ and $\mathcal{E}_0$, $\mathcal{E}_1$ defined in \eqref{eq:def_Phi} and \eqref{eq:def_Phi2}, respectively.
\end{proposition}

The source of the term in the second line of \eqref{eq:real_def_Phi2} is the same as in the discussion below Proposition~\ref{prop:real_EG}: we do not prove the real version of Theorem~\ref{eq:maintheoDBM1_main}, but instead use  the complex result directly, which yields an additional loss in the GFT owing to the slight mismatch in the second moments. Eventually, this leads to the additional term in the rhs. of \eqref{eq:main2} compared to the \eqref{eq:main1} in the complex case. We also remark that due to this term, the control parameter $\Phi_2^r$ becomes ineffective in the regime when $|z_1-\overline{z}_2|\lesssim N^{-1/2}$, in the sense that the local law \eqref{eq:1G_ll_av} implies a better bound in this regime compared to \eqref{eq:EGG}. Moreover, in the critical case $z_1=\overline{z}_2$, $\Phi_2^r$ blows up. This is a manifestation of the symmetry of the spectrum of a real matrix $X$ with respect to the real axis.

Finally, we formulate the analogue of the local law from Proposition~\ref{prop:2G_av} in the real setting. While Proposition~\ref{prop:2G_av} addresses only the two-resolvent chains, along its proof, which is presented in Section~\ref{sec:multiG_proof}, we also obtained suboptimal local laws for three- and four-resolvent chains as a byproduct. Though these local laws were not needed for the proof of Theorem~\ref{theo:main1}, now they are required as ingredients for the GFT in Propositions~\ref{prop:real_EG} and~\ref{prop:real_EGG}. Additionally, for the same purpose we need to cover the case when~$X$ is a real i.i.d. matrix with a small complex Gaussian component, however with much less precision: instead of local laws only the size bounds on the two- and four-resolvent chains are required. This is the second statement of the following proposition.

\begin{proposition}\label{prop:real_2G_av} (i) Assume the set-up and conditions of Proposition \ref{prop:2G_av} modulo the replacement of the complex $X_0$ and $\widetilde{X}$ matrices by real matrices. Recall the definition of $\widehat{\beta}_{12}^{[b]}=\widehat{\beta}^{[b]}_{12}(w_1,w_2)$ from \eqref{eq:beta_hat_b}. Denote $G_l:=G^{z_l}(w_l)$ for $l=1,2$. Then, for any fixed $\delta,\epsilon,\kappa>0$, we have
\begin{align}
\left\vert \left\langle \left(G_1B_1G_2-M_{12}^{B_1}(w_1,w_2)\right)B_2\right\rangle\right\vert &\prec \frac{1}{N\eta_*\widehat{\beta}_{12}^{[b]}},\label{eq:real_2G_av}\\
\left\vert\left\langle \left(G_1B_1G_2B_2G_1^{(*)}-M_{121}^{B_1,B_2}\big(w_1,w_2,w_1^{(*)}\big)\right) B_3\right\rangle\right\vert &\prec \frac{1}{N\eta_*^3}\wedge\frac{1}{\sqrt{N\eta_*}\eta_*\widehat{\beta}_{12}^{[b]}},\label{eq:real_3G_av}\\
\left\vert\left\langle \left(G_1B_1G_2B_2G_1^{(*)}B_3G_2^{(*)}-M_{1212}^{B_1,B_2,B_3}\big(w_1,w_2,w_1^{(*)}, w_2^{(*)}\big)\right) B_4\right\rangle\right\vert &\prec \frac{1}{N\eta_*^4}\wedge\frac{1}{\eta_*\big(\widehat{\beta}_{12}^{[b]}\big)^2},\label{eq:real_4G_av}
\end{align}  
uniformly in $B_l\in\mathrm{span}\{E_\pm,F^{(*)}\}$, $l\in[4]$, $|z_i|\le 1-\delta$, $\Re w_i\in\mathbf{B}_\kappa^{z_i}$, and $N^{-1+\epsilon}\le \eta_i\le N^{100}$ for $i\in[2]$.

(ii) Fix further $\omega_*>0$ and consider $\widehat{X}:=\sqrt{1-\widehat{\mathfrak{s}}^2}X + \widehat{\mathfrak{s}}\widetilde{X}_c$, where $\widetilde{X}_c$ is an $N\times N$ complex Ginibre matrix, independent of $X$, $\widehat{\mathfrak{s}}\in [0,1]$, and $X$ satisfies the conditions of Proposition~\ref{prop:real_2G_av}(i). Denoting again by $G_l=G^{z_l}(w_l)$ the resolvent associated to the Hermitization of $\widehat{X}$, we have 
\begin{align}
\left\vert\langle G_1B_1G_2^\mt B_2\rangle\right\vert &\prec \frac{1}{\widehat{\beta}^{[b]}_{1\overline{2}}},\label{eq:2G_bound_cr}\\
\left\vert \langle G_1B_1G_2^\mt B_2G_1^{(*)}B_3\rangle\right\vert&\prec \frac{1}{\eta_*\widehat{\beta}^{[b]}_{1\overline{2}}},\label{eq:3G_bound_cr}\\
\left\vert \langle G_1B_1G_2^\mt B_2G_1^{(*)}B_3(G_2^\mt)^{(*)}\rangle\right\vert&\prec \frac{1}{\eta_*\big(\widehat{\beta}^{[b]}_{1\overline{2}}\big)^2},\label{eq:4G_bound_cr}
\end{align}
uniformly in $\widehat{\mathfrak{s}}\in [0,N^{-\omega_*}]$ and in the rest of the parameters as stated below \eqref{eq:real_4G_av}.
\end{proposition}

The proof of Proposition~\ref{prop:real_2G_av} is presented in Section~\ref{sec:real_local_laws}. Apart from Proposition~\ref{prop:real_2G_av}, the proof of Theorem~\ref{theo:main2} requires also the real versions of the local laws from Proposition~\ref{prop:2G_av_suboptimal} and \ref{prop:multiG_oneG}. These results hold without any changes in the formulations, and their proofs do not require any additional adjustments apart from the ones discussed in Section~\ref{sec:real_local_laws}. We omit further details.  

\subsection{Proof of Theorem \ref{theo:main2}}\label{sec:theo2} In this section we derive the analogues of Propositions~\ref{prop:E_omega} and~\ref{prop:Var_main} in the real case from Propositions~\ref{prop:real_Cov}--\ref{prop:real_EGG} and conclude the proof of Theorem~\ref{theo:main2}. We omit most of the calculations as they are analogous to the ones presented in Section \ref{sec:Girko_calculations}.

\subsubsection{Analysis of the expectation} First, fix $\delta_r>0$ and assume the set-up and conditions of Proposition~\ref{prop:E_omega} (modulo the replacement of the complex case by the real one). Denote $\nu_2:=1/11$ and recall from Proposition~\ref{prop:E_omega} that $\nu_0=1/14$. We claim that
\begin{equation}
\E\mathcal{L}_N\big(\omega_{a,N}^{(z_0)}\big) = \frac{N}{\pi}\left(1+\mathcal{O}(N^{-c})\right) + N^{2a+\xi}\mathcal{O}\left(\left(N|\Im z_0|^2\right)^{-1/3} + N^{-\nu_2}\right),
\label{eq:E_omega_real}
\end{equation}
for any $a\in [1/2,1/2+\nu_0)$, uniformly in $|z|\le 1-\delta$ with $|\Im z|\ge N^{-1/2+\delta_r}$. Here $c>0$ is some small constant independent of $N$. The main difference of \eqref{eq:E_omega_real} from its complex counterpart \eqref{eq:E_omega} is that the error term in \eqref{eq:E_omega_real} may be larger than the leading one for sufficiently small $|\Im z_0|$, while in \eqref{eq:E_omega} the error term is always smaller. 

To prove \eqref{eq:E_omega_real}, we follow the proof of Proposition \ref{prop:E_omega}. Denote 
\begin{equation*}
p:=N|\Im z_0|^2\in [N^{2\delta_r},N].
\end{equation*}
Since the diameter of support of $\omega_{a,N}^{(z_0)}$ is of order $N^{-a}\ll |\Im z_0|$, it holds that $N|\Im z|^2\sim p$ for every $z\in\mathrm{supp}(\omega_{a,N}^{(z_0)})$. We take $\eta_c:=N^{-1+\delta_1}$ for a small $\delta_1>0$, and $\eta_0:=N^{-1}p^{-1/6}$ to optimize the bounds. The estimates in the mesoscopic regime $\eta\ge \eta_c$ do not change, while instead of \eqref{eq:E_cut_tails} we have
\begin{equation}
|\E J_T| + |\E I_{0}^{\eta_L}|+|\E I_{\eta_L}^{\eta_0}|\lesssim ((N\eta_0)^2+N^{-2\nu_1} + p^{-1/2})N^{2a+\xi}.
\label{eq:real_E_cut_tails}
\end{equation}
Here the additional $p^{-1/2}$ term arises from \eqref{eq:real_tail_bound_alternative} employed in the analysis of $|\E I_{\eta_L}^{\eta_0}|$. In the intermediate regime $\eta\in [\eta_0,\eta_c]$ we apply Proposition~\ref{prop:real_EG} and optimize \eqref{eq:real_def_Phi} by choosing 
\begin{equation*}
t_1:=N^{-3/4}(N\eta)^{3/2}p^{1/3}\quad\text{and}\quad t_2:= N^{-1}(N\eta)^2 p^{2/3}.
\end{equation*}
For $p\lesssim N^{3/11}$ this gives
\begin{equation*}
\Phi_1^r(\eta)\lesssim  (N\eta)^{-1}p^{-1/3} + (N\eta)^{-3}p^{-2/3} + N^{-1/4}(N\eta)^{-5/2}p^{1/3}.
\end{equation*}
Arguing similarly to the proof of Proposition \ref{prop:E_omega} and taking into account the error term from \eqref{eq:real_E_cut_tails}, we conclude that
\begin{equation*}
\E\mathcal{L}_N\big(\omega_{a,N}^{(z_0)}\big) = \frac{N}{\pi}\left(1+\mathcal{O}(N^{-c})\right) +N^{2a+\xi}\mathcal{O}\left(p^{-1/3} + N^{-2\nu_1}+N^{-1/4} p^{7/12} \right),
\end{equation*}
which immediately implies \eqref{eq:E_omega_real} for $p\le N^{3/11}$. It remains to notice that the bounds in Propositions~\ref{prop:real_tail_bound} and~\ref{prop:real_EG} improve as $p$ increases, so for $p>N^{3/11}$ \eqref{eq:E_omega_real} holds with $p$ replaced by $N^{3/11}$, which gives the error term of order $N^{-1/11}$. This finishes the proof of \eqref{eq:E_omega_real}.

Recall the definition of $f_{a,N}^\pm$ from \eqref{eq:def_f_pm}. We show that
\begin{equation}
\left\vert\E \mathcal{L}_N\left(f_{a,N}^+-f_{a,N}^-\right)\right\vert\lesssim N^{1/2-\alpha+\xi} \left( N^{1/2-a} + N^{a-1/2-\nu_2} + N^{a-1/2-2(1/2-\alpha)/3}\right).
\label{eq:real_E_contribution}
\end{equation}
To establish \eqref{eq:real_E_contribution}, we argue as in \eqref{eq:E_dif_f} and use the definition of $\Omega_N^\pm$ from \eqref{eq:Omega_pm}. We get
\begin{equation}
\E \mathcal{L}_N\left(f_{a,N}^+-f_{a,N}^-\right) = \int_{\Omega_N^+\setminus\Omega_N^-} \E\mathcal{L}_N\big(\omega_{a,N}^{(z)}\big) \dif z.
\label{eq:real_E_rep}
\end{equation}
For $|\Im z|\ge N^{-1/2+\delta_r}$ we estimate the integrand in \eqref{eq:real_E_rep} by \eqref{eq:E_omega_real}, while for $|\Im z|<N^{-1/2+\delta_r}$ we use the rough upper bound of order $N^{2a+\xi}$, which immediately follows from the Girko's formula \eqref{eq:Girko} and the single-resolvent local law \eqref{eq:1G_ll_av}. However, in the second regime the smallness comes from the volume factor. Indeed, by Assumption~\ref{ass:Omega_real} we have
\begin{equation*}
\left\vert \left(\Omega_N^+\setminus\Omega_N^-\right) \cap \left\lbrace z\,:\, |\Im z|<N^{-1/2+\delta_r}\right\rbrace\right\vert \lesssim N^{-1/2-a+\delta_r}.
\end{equation*}
Combining these inputs and choosing sufficiently small $\delta_r>0$, we obtain \eqref{eq:real_E_contribution} by an elementary calculation.

\subsubsection{Analysis of the variance} Assume the set-up and conditions of Proposition \ref{prop:Var_main} (modulo the replacement of the complex case by the real one). Denote $\nu_3:=1/53$ and recall that $q_0=1/20$. We claim that
\begin{equation}
\Var\left[ \Lin_N(f_{a,N})\right]\lesssim N^{2(a-\alpha) - 2q_0(1/2-\alpha)+\xi} +  N^{2(a-\alpha)-\nu_3+\xi}.
\label{eq:Var_main_real}
\end{equation}
The proof of \eqref{eq:Var_main_real} closely follows the one of Proposition~\ref{prop:Var_main}. First, we upper bound $\Var[I_{\eta_c}^T]$ where $\eta_c$ is chosen as before. Fix a (small) $\delta_r>0$. We represent $\Var[I_{\eta_c}^T]$ in the form \eqref{eq:Var_I_ll} and distinguish between the two regimes of $(z_1,z_2)$-integration:  $|\Im z_1|\wedge |\Im z_2|\le N^{-1/2+\delta_r}$ and the complementary regime. In the first regime we use the trivial bound
\begin{equation}
\left\vert\mathrm{Cov}\left(\langle G^{z_1}(\ii\eta_1)\rangle,\langle G^{z_2}(\ii\eta_2)\rangle\right)\right\vert \lesssim \frac{N^\xi}{N^2\eta_1\eta_2},
\label{eq:Cov_trivial}
\end{equation}
which immediately follows from \eqref{eq:1G_ll_av}. The small volume of the $z_1,z_2$ integration regime will compensate
for the crudeness of the bound~\eqref{eq:Cov_trivial}. In the regime $|\Im z_1|\wedge |\Im z_2|\ge N^{-1/2+\delta_r}$ we employ Proposition~\ref{prop:real_Cov} and argue as in \eqref{eq:Cov_err_int}--\eqref{eq:tails_UV}. The only difference is that now the regime $|\Im z_1|\wedge |\Im z_2|\le N^{-1/2+\delta_r}$ is missing in \eqref{eq:expl_UV_calc}, so it should be added back and taken into account in \eqref{eq:tails_UV}. This is done by simply estimating $|V_{12}|\lesssim (\eta_1\eta_2)^{-1}$ and $|U_j|\lesssim 1$, $j=1,2$, as it follows from \eqref{eq:def_UV}. Using these inputs, we get that $\Var[I_{\eta_c}^T]$ admits the bound \eqref{eq:Var_ll_bound} in the real case, i.e. the upper bound on the mesoscopic regime does not deteriorate compared to the complex case.

Next, we derive an upper bound on 
\begin{equation}
\Var[I_{\eta_L}^{\eta_c}]= \left(\frac{N}{2\pi\ii}\right)^2 \int_\C\int_\C \Delta f_N(z_1)\overline{\Delta f_N(z_2)} \int_{\eta_L}^{\eta_c}\int_{\eta_L}^{\eta_c} \Cov\left(\langle G^{z_1}(\ii\eta_1)\rangle, \langle G^{z_2}(\ii\eta_2)\rangle\right)\dif\eta_1\dif\eta_2\dif^2 z_1\dif^2 z_2.
\label{eq:Var_Girko_real}
\end{equation}
Denote for short
\begin{equation}
p=p(z_1,z_2):=\min\left\lbrace N|\Im z_1|^2, N|\Im z_2|^2, N|z_1-z_2|^2, N|z_1-\overline{z}_2|^2\right\rbrace.
\end{equation} 
In the regime $p\le N^{2\delta_r}$ we again use the bound \eqref{eq:Cov_trivial} together with the small volume effect. Hence, we can further focus on the complementary regime. For each pair of parameters $z_1,z_2\in \mathrm{supp}(f_N)$ with $p>N^{2\delta_r}$ we take $\eta_0=\eta_0(z_1,z_2)\in (\eta_L,N^{-1})$. This intermediate scale will be chosen at the end to optimize the bound on the rhs. of \eqref{eq:Var_Girko_real}. Note that unlike in the complex case, $\eta_0$ depends on $z_1, z_2$, so  the splitting of the integration regime $[\eta_L,\eta_c]$ into $[\eta_L,\eta_0]$ and $[\eta_0,\eta_c]$ must be done inside the $z_1, z_2$-integration in~\eqref{eq:Var_Girko_real}. 

Fix $z_1,z_2$ and consider the cross-regime when $\eta_1\in [\eta_L,\eta_0]$ and $\eta_2\in [\eta_0,\eta_c]$. We have
\begin{equation}
\begin{split}
\left\vert\int_{\eta_L}^{\eta_0}\int_{\eta_0}^{\eta_c} \mathrm{Cov}\left(\langle G^{z_1}(\ii\eta_1)\rangle, \langle G^{z_2}(\ii\eta_2)\rangle\right)\dif\eta_1\dif\eta_2\right\vert \lesssim &N^{-1+\xi} \E \left\vert\int_{\eta_L}^{\eta_0}\langle G^{z_1}(\ii\eta_1)\rangle\dif \eta_1\right\vert \\
\lesssim &N^{-2+\xi}\left((N\eta_0)^2+N^{-2\nu_1} + \left(N|\Im z_1|^2\right)^{-1/2}\right).
\end{split}
\label{eq:Var_cross_term}
\end{equation}
Here in the first line we used \eqref{eq:av_ll_simple} and to go from the first to the second line we combined Proposition~\ref{prop:real_tail_bound} with the argument from the proof of Lemma~\ref{lem:Girko_reduction}. Similarly we get that the same bound holds in the case when $\eta_1,\eta_2\in [\eta_L,\eta_0]$. In the regime $\eta_1,\eta_2\in [\eta_0,\eta_c]$ we apply Proposition~\ref{prop:real_EGG} and optimize the error term $\Phi_2^r$ in \eqref{eq:real_def_Phi2} by taking
\begin{equation}
R:=p^{4/5},\quad t_1:=N^{-3/4}p^{1/20}(N\eta_*)^{1/2},\quad t_2:=N^{-1}p^{1/10}.
\label{eq:real_Phi2_optimization}
\end{equation}
For sufficiently small $\eta_*$ this choice is invalid, since $t_1$ becomes smaller than $t_2$. However, in this situation taking \eqref{eq:real_Phi2_optimization} in the rhs. of \eqref{eq:real_def_Phi2} would give an upper bound on 
$\Phi_2^r$ which is much larger than $(N^2\eta_1\eta_2)^{-1}$. Since this latter bound holds by the local law \eqref{eq:1G_ll_av}, one can make a choice \eqref{eq:real_Phi2_optimization} for any $\eta_*\le \eta_c$. 
Integrating the resulting bound for $\Phi_2^r$ over $\eta_1,\eta_2\in [\eta_0,\eta_c]$ and using \eqref{eq:Var_cross_term}, we get
\begin{equation}
\begin{split}
&N^2\left\vert\int_{\eta_L}^{\eta_c}\int_{\eta_L}^{\eta_c} \mathrm{Cov}\left(\langle G^{z_1}(\ii\eta_1)\rangle,\langle G^{z_2}(\ii\eta_2)\rangle\right)\dif \eta_1\dif\eta_2\right\vert\\
&\quad\lesssim N^{8\delta_1+\xi}\left(p^{-1/20} + p^{-9/10} (N\eta_0)^{-2} + N^{-1/4} p^{1/20}(N\eta_0)^{-5/2} +(N\eta_0)^2\right).
\end{split}
\label{eq:real_Var_ll_bound}
\end{equation}
To optimize \eqref{eq:real_Var_ll_bound} over $\eta_0$, we take $\eta_0:=N^{-1}p^{-9/40}$.
We use this bound for $p\le N^{20/53}$, while in the complementary regime employ the same monotonicity idea as above \eqref{eq:real_E_contribution}. It remains to perform the integration over $z_1,z_2$ by the means of Lemma~\ref{lem:int}. This finishes the proof of~\eqref{eq:Var_main_real}.

\subsubsection{Conclusion of the proof of Theorem \ref{theo:main2}} Combining \eqref{eq:real_E_contribution} with \eqref{eq:Var_main_real} as in Section \ref{sec:proof_thm1}, we get
\begin{equation}
\begin{split}
\Var\left[\mathcal{L}_N(\phi_N)\right] \lesssim &N^{2(a-\alpha)+\xi} \left( N^{-2q_0(1/2-\alpha)} + N^{-\nu_3}\right)\\
 + &N^{1-2\alpha+\xi}\left(N^{1-2a} +N^{2a-1-2\nu_2} + N^{2a-1-4(1/2-\alpha)/3}\right).
\end{split}
\label{eq:real_Var_prefin}
\end{equation}
We optimize the rhs. of \eqref{eq:real_Var_prefin} in $a\in [1/2,1/2+\nu_0)$, which leads to the choice
\begin{equation}
a:=\frac{1}{2}+\frac{1}{4}\min\left\lbrace 2q_0\left(\frac{1}{2}-\alpha\right), \nu_3\right\rbrace.
\label{eq:real_a}
\end{equation}
Substituting \eqref{eq:real_a} into \eqref{eq:real_Var_prefin}, we finish the proof of Theorem~\ref{theo:main2}.\hfill $\qed$

\subsection{Proof of Proposition \ref{prop:real_2G_av}}\label{sec:real_local_laws}

First we prove that \eqref{eq:real_2G_av}--\eqref{eq:real_4G_av} hold for the Gauss-divisible real matrix $X$ satisfying the assumptions of Proposition~\ref{prop:real_2G_av}. We closely follow the proof of Proposition \ref{prop:2G_av} presented in Section~\ref{sec:multiG_proof} and only outline the differences arising in the real set-up compared to the complex one. In particular, we focus on the case $b=0$ which corresponds to $\mathfrak{s}\sim 1$, while the general case $b\in [0,1]$ follows by a minor adjustment, similarly to Section~\ref{sec:general_b}. Moreover, we consider only observables equal to $E_\pm$, for the treatment of off-diagonal observables see Section~\ref{sec:general_B}.

We replace the complex-valued Brownian motion in the Ornstein-Uhlenbeck process \eqref{eq:OU} with the real-valued one and keep the characteristic flow \eqref{eq:char_flow} unchanged. This does not affect the estimates in Proposition~\ref{prop:Zig}, and there is only a minor modification in the zig equations in Section~\ref{sec:Zig_structure}. Specifically, the generator $\mathcal{A}_{[k],t}$ in \eqref{eq:kG_syst} as well as the analysis of the martingale term $\dif\mathfrak{E}_{[k],t}$ remain unchanged, while the following new terms appear in the forcing term \eqref{eq:def_Fk} due to the real symmetry of the model:
\begin{equation}
\frac{\sigma}{N}\sum_{i\le j}\left\langle G_{1,t}B_1\cdots G_{i,t}E_\sigma \left(G_{i,t}B_i\cdots G_{j,t}\right)^\mathfrak{t}E_\sigma G_{j,t}B_{j}\cdots G_{k,t}B_k\right\rangle.
\label{eq:real_forcing}
\end{equation}
The emergence of these terms does not change the bounds in Proposition \ref{prop:err_bounds}. Moreover, we estimate all terms in \eqref{eq:real_forcing} in the same way by upper bounding the absolute value of the term indexed by $1\le i\le j\le k$ in \eqref{eq:real_forcing} (together with the prefactor $N^{-1}$) as follows:
\begin{equation}
N^{-1}\left\langle \left\vert G_{i,t}B_i\cdots G_{j,t}\right\vert^2\right\rangle^{1/2}\left\langle \left\vert G_{j,t}B_j\cdots G_{k,t}B_kG_{1,t}B_1\cdots G_{j,t}\right\vert^2\right\rangle^{1/2},
\end{equation}
by the means of Cauchy-Schwarz inequality. This reduces the quantities in \eqref{eq:real_forcing} to the resolvent chains containing only $G_{1,t}^{(*)}$ and $G_{2,t}^{(*)}$, which are then further reduced to shorter chains by the reduction inequalities \eqref{eq:reduction_simple}, \eqref{eq:reduction}, and estimated by the stopping time \eqref{eq:def_tau_k}.

Now we prove Proposition~\ref{prop:real_2G_av}(ii). Let $X$ be a real matrix satisfying assumptions of Proposition~\ref{prop:real_2G_av}(i). Define $X_t$ by the flow \eqref{eq:OU} starting from $X$, where $B_t$ is the complex-valued Brownian motion. Additionally, we evolve $z_j,w_j$ for $j=1,2$, along the characteristic flow \eqref{eq:char_flow}. Since \eqref{eq:2G_bound_cr}--\eqref{eq:4G_bound_cr} hold for $t=0$ by Propositions~\ref{prop:real_2G_av}(i) and~\ref{prop:M_bounds}, it suffices to show that these bounds propagate along the flow up to the time $t\le N^{-\omega_*}$ in the same sense as in Proposition~\ref{prop:Zig}. We need two inputs to prove this propagation. First, one needs to show that $X_t$ satisfies the local low \eqref{eq:1G_ll_av}. This follows from \cite[Lemma~B.7]{cipolloni2024maximum}, which states this result for spectral parameters on the imaginary axis, though its proof immediately extends to the bulk regime. The second input is the two-resolvent local law
\begin{equation}
\left\vert\left\langle \left(G_{1,t}B_1G_{1,t}^{(*)}-M_{11,t}^{B_1}(w_{1,t},w_{2,t}^{(*)})\right)B_2\right\rangle \right\vert\prec \frac{1}{N\eta_{1,t}^2},
\label{eq:mixed_2G_simple}
\end{equation} 
for $B_1,B_2\in\mathrm{span}\{E_\pm,F^{(*)}\}$ and $t\le N^{-\omega_*}$. Since the resolvent chain in \eqref{eq:mixed_2G_simple} does not involve $G_1^\mt$, the proof of \eqref{eq:mixed_2G_simple} is identical to the zig step in the proof of \eqref{eq:multiG_av} for $k=2$.

Let us differentiate the time-dependent analogue of \eqref{eq:2G_bound_cr} along the flow. We get, similarly to \eqref{eq:kG_dif}, that
\begin{equation}
\begin{split}
\dif \langle G_{1,t}B_1G_{2,t}^\mathfrak{t}B_2\rangle = &\langle G_{1,t}B_1G_{2,t}^\mathfrak{t}B_2\rangle\dif t +\frac{\sigma}{N} \langle G_{1,t}(G_{1,t}B_1G_{2,t}^{\mathfrak{t}})^\mathfrak{t}G_{2,t}^\mathfrak{t}B_2\rangle\dif t\\
 +&\langle G_{1,t}-M_{1,t}\rangle \langle G_{1,t}^2B_1G_{2,t}^\mathfrak{t}B_2\rangle\dif t + \langle G_{2,t}-M_{2,t}\rangle \langle G_{1,t}B_1(G_{2,t}^\mathfrak{t})^2B_2\rangle\dif t +\dif\mathcal{E}_{[2],t},
\end{split}
\label{eq:2G_dif_cr}
\end{equation}
where $\dif\mathcal{E}_{[2],t}$ is a martingale term. Unlike in Section~\ref{sec:Zig_structure}, we do not subtract the corresponding differential equation for the deterministic approximation from \eqref{eq:2G_dif_cr}. Observe that the only linear term in \eqref{eq:2G_dif_cr}, in the terminology introduced below \eqref{eq:def_Fk}, is the first term in the rhs. of \eqref{eq:2G_dif_cr}. This means that for $k=2$ there is no analogue of the generator $\mathcal{A}_{[2],t}$ from \eqref{eq:kG_syst} in the current setting, so no Duhamel formula or propagators are needed. To estimate the remaining terms in the rhs. of \eqref{eq:2G_dif_cr}, one needs only the size bounds of the form \eqref{eq:2G_bound_cr}--\eqref{eq:2G_bound_t}, which are available by a stopping time argument, and the single-resolvent local law mentioned above. For $k=3$ and $k=4$ we obtain similar equations with drift terms, which can be estimated relying only the size bounds for two, three and four-resolvent chains, \eqref{eq:mixed_2G_simple} and the single-resolvent local law. The remaining details of the proof are standard and thus are omitted.

\subsection{Proof of Proposition \ref{prop:real_Cov}}\label{sec:real_chaos_exp}

We proceed by a minor modification of the proof of Proposition~\ref{prop:Cov} and focus on the case when $X$ has a Gaussian component of order 1. Similarly to the proof of \eqref{eq:init_expansion} presented in Section~\ref{sec:Cov_reduction}, we derive the following initial expansion in the real case:
\begin{equation}
\begin{split}
\mathrm{Cov}\left(\langle G_1\rangle,\langle G_2\rangle\right) =& \mathrm{Cov}\left(\langle G_1-M_1\rangle \langle (G_1-M_1)A\rangle, \langle G_2\rangle\right) + \frac{\sigma}{2N}\mathrm{Cov}\left(\langle G^\mathfrak{t}_1E_\sigma G_1AE_\sigma\rangle,\langle G_2\rangle\right)\\
+&\frac{1}{4N^2}\left\langle G_1AE_\sigma \left(G_2^2 + (G_2^\mathfrak{t})^2\right)E_\sigma\right\rangle - \E\left[\underline{\langle WG_1A\rangle \left(\langle G_2\rangle-\E\langle G_2\rangle\right)}\right],
\end{split}
\label{eq:real_init_expansion}
\end{equation}
where $A$ is the same as in \eqref{eq:init_expansion}. Applying \eqref{eq:G_symmetry} to $G_2$ and recalling \eqref{eq:V_comp3}, we get that the deterministic approximation to the first term in the second line of \eqref{eq:real_init_expansion} equals to $\widehat{V}_{12}/(2N^2)$. Arguing further similarly to the analysis of \eqref{eq:init_expansion} in Section~\ref{sec:Cov_reduction} and using \eqref{eq:real_2G_av} instead of \eqref{eq:2G_av}, we get that the sum of the terms in the second line of \eqref{eq:real_init_expansion} equals to 
\begin{equation}
\frac{1}{N^2}\cdot \frac{\widehat{V}_{12}+\kappa_4U_1U_2}{2} + \mathcal{O}\left(\left(\frac{1}{N\widehat{\gamma}}+N^{-1/2}\right)\frac{N^\xi}{N^2\eta_1\eta_2}\right),
\end{equation} 
We do not perform iterative expansions for the last term in the first line of \eqref{eq:real_init_expansion}, but bound it from above directly:
\begin{equation}
\begin{split}
N^{-1}\left\vert\mathrm{Cov}\left(\langle G^\mathfrak{t}_1E_\sigma G_1AE_\sigma\rangle,\langle G_2\rangle\right)\right\vert \le &N^{-1}\E \left[ \big\vert \big\langle \big( G^\mathfrak{t}_1E_\sigma G_1 -M_{\overline{1}1}^{E_\sigma}\big)AE_\sigma\big\rangle\big\vert\cdot \big\vert \langle G_2-M_2\rangle\big\vert\right]\\
\lesssim &\frac{1}{N|\Im z_1|^2}\cdot\frac{N^\xi}{N^2\eta_1\eta_2},
\end{split}
\end{equation}
for $\sigma\in\{\pm\}$. Here we used \eqref{eq:real_2G_av} for $G_1$ and $G_1^\mathfrak{t}$, and \eqref{eq:1G_ll_av} for $G_2$. Thus, to conclude the proof of Proposition \ref{prop:real_Cov} it remains to show that
\begin{equation}
\left\vert\mathrm{Cov}\left(\langle G_1-M_1\rangle \langle (G_1-M_1)A\rangle, \langle G_2\rangle\right)\right\vert \lesssim \left(\frac{1}{N\widehat{\gamma}}+\frac{1}{N|\Im z_1|^2}+\frac{1}{N|\Im z_2|^2} +N^{-1/2}\right)\frac{N^\xi}{N^2\eta_1\eta_2},
\label{eq:real_psi22}
\end{equation}
which is the analogue of Proposition~\ref{prop:psi_22}.

The proof of \eqref{eq:real_psi22} follows the lines of the proof of Proposition~\ref{prop:psi_22} and it relies on the same hierarchy of covariances. The only difference is that the covariances containing both $G_1$ and $G_1^\mathfrak{t}$ emerge along the chaos expansion, as we have seen in \eqref{eq:real_init_expansion}. Similarly to the treatment of the last term in the first line of \eqref{eq:real_init_expansion}, we do not expand these quantities further, but treat them as error terms and estimate them using the local law \eqref{eq:real_2G_av}. This increases the error term $\Upsilon$ in the rhs. of \eqref{eq:master_S}, so in the real case it is given by
\begin{equation}
\Upsilon^r=\Upsilon^r(z_1,z_2;\eta_1,\eta_2):=\left(\sqrt{N}+\frac{1}{\widehat{\gamma}} + \frac{1}{|\Im z_1|^2} + \frac{1}{|\Im z_2|^2}\right)\eta_* N^\xi,
\end{equation}
compared with \eqref{eq:def_Upsilon} in the complex case. However, this increase is affordable, since the estimate on the error term in \eqref{eq:real_main_cov} is weaker than the one in \eqref{eq:main_cov}. This finishes the proof of Proposition~\ref{prop:real_Cov}.\hfill$\qed$

\subsection{Proof of Propositions \ref{prop:real_tail_bound}, \ref{prop:real_EG} and \ref{prop:real_EGG}}\label{sec:real_crit_scale}

The proofs of Propositions~\ref{prop:real_tail_bound}, \ref{prop:real_EG}, and \ref{prop:real_EGG} closely follow the proofs of Propositions~\ref{prop:tail_bound}, \ref{prop:EG}, and \ref{prop:EGG}, respectively, presented in Section~\ref{sec:proofprop}. We start with outlining the main technical differences. Recall that the key inputs in Section~\ref{sec:proofprop} are the relaxation of DBM from \cite[Proposition~4.6]{bourgade2024fluctuations} and the quantitative decorrelation of a pair of DBMs from Theorem~\ref{eq:maintheoDBM1}. These results are stated in the complex case, which in the framework of Theorem~\ref{eq:maintheoDBM1} means that the driving Brownian motion in \eqref{eq:matDBM1} is complex-valued and the initial condition $X_0$ is a complex i.i.d. matrix. Although the analogue Theorem~\ref{eq:maintheoDBM1} can be established in the real case, corresponding to the real $B_t$ and $X_0$ in \eqref{eq:matDBM1}, this would lead to additional technical difficulties, as the analysis of DBM in the real case is more delicate. Instead, we use Theorem \ref{eq:maintheoDBM1} with a real-valued initial condition $X_0$ and a complex driving Brownian motion (see Section~\ref{sec:DBM} for the complex case). We also use \cite[Proposition~4.6]{bourgade2024fluctuations} in this mixed setting instead of proving this result in the real case.

Recall from \eqref{eq:distr_t}--\eqref{eq:strat_concl} that the proofs of Propositions~\ref{prop:tail_bound}, \ref{prop:EG}, and~\ref{prop:EGG} proceed in two steps: first, the results are established for a matrix with a small Gaussian component, and second, this component is removed by a GFT. Since we now use the same inputs \cite[Proposition~4.6]{bourgade2024fluctuations} and Theorem~\ref{eq:maintheoDBM1} as in Section~\ref{sec:proofprop}, the first step of this strategy does not change. However, this does not come for free, as the second step (GFT) becomes more involved, which we now explain. To prove Propositions~\ref{prop:real_tail_bound}, \ref{prop:real_EG}, and~\ref{prop:real_EGG} we must remove a small \emph{complex} Ginibre component added to a \emph{real} i.i.d. matrix. This is achieved using the flow \eqref{eq:OU_GFT} with real~$X_0$ and complex driving Brownian motion. Unlike in Section~\ref{sec:proofprop}, this flow does not preserve the second-order correlation structure of~$X^t$, since the first two moments of~$X_0$ and~$B_t$ are not matched. This leads to a non-negligible contribution from the second moments in the GFT, which for instance in the set-up of Proposition~\ref{prop:EG} means that the second-order terms ($\ell=2$) in the cumulant expansion \eqref{eq:cum_exp_first} do not cancel fully with the second term in the rhs. of \eqref{eq:dG}. As a consequence, the error term coming from the GFT step deteriorates, which should be properly taken into account.

We control the contribution from the second-order cumulants in the GFT in the proofs of Propositions~\ref{prop:real_tail_bound}, \ref{prop:real_EG}, and~\ref{prop:real_EGG}, by the means of the following elementary consequence of Proposition~\ref{prop:real_2G_av}(ii) proved in Supplementary Section~\ref{app:G_bounds_extension}. This statement gives an upper bound on the two specific types of resolvent chains. As we will see later in this section, no other new terms arise from the GFT.

\begin{lemma}\label{lem:G_bounds_extension} Fix a (small) $\epsilon>0$ and assume the set-up and conditions of Proposition~\ref{prop:real_2G_av}(ii). Then for $\sigma\in\{\pm\}$ we have
\begin{align}
\left\vert\langle G_1^2 E_\sigma (G_2^\mt)^2E_\sigma \rangle\right\vert &\prec \frac{N^{2\epsilon}}{N^2\eta_1\eta_2}\cdot \frac{1}{\widehat{\beta}^{[b]}_{1\overline{2}}},\label{eq:real_size_bound1}\\
\left\vert\langle (G_1^2 E_\sigma G_2^\mt E_\sigma \rangle\right\vert &\prec \frac{N^{2\epsilon}}{N^2\eta_1\eta_2}\cdot\frac{1}{\eta_1\widehat{\beta}^{[b]}_{1\overline{2}}},\label{eq:real_size_bound2}
\end{align}
uniformly in $|z_i|\le 1-\delta$, $\Re w_i\in\mathbf{B}_\kappa^{z_i}$ and $0<\eta_i<N^{-1+\epsilon}$ for $i\in[2]$. 
\end{lemma}

We remind the reader that in the set-up of Proposition \ref{prop:real_2G_av}(ii) and Lemma \ref{lem:G_bounds_extension}, the iid matrix $\widetilde{X}$ contains a real Ginibre component of order $N^{-b}$ and a complex Ginibre component of order $\widehat{\mathfrak{s}}<N^{-\omega_*}$. As the size of the real Ginibre component in $X$ decreases, the rhs. of the bounds \eqref{eq:real_size_bound1} and \eqref{eq:real_size_bound2} deteriorate. As explained below Proposition \ref{prop:2G_av}, this deterioration is purely technical: one could set $b:=0$ after performing the GFT for the local laws in Proposition~\ref{prop:real_2G_av}(i). We do not pursue this approach here, and the price is that crude estimates in  Lemma \ref{lem:G_bounds_extension} enable us to remove only a small complex Gaussian component added to a real matrix which already has sufficiently large real Gaussian component. 

With the discussion above in mind, we outline the proof strategy for Propositions~\ref{prop:real_tail_bound}, \ref{prop:real_EG} and~\ref{prop:real_EGG}. Let $X$ be the original real i.i.d. matrix from Propositions \ref{prop:real_tail_bound}, \ref{prop:real_EG}, and \ref{prop:real_EGG}. We take two ($N$-dependent) times $N^{-1+\omega_*}<t_2\le t_1<N^{-\omega_*}$ for a small fixed $\omega_*>0$. All estimates will be optimized over $t_1,t_2$ at the end.

\medskip

\noindent\textbf{Step 1.1.} We embed $X$ into the real-valued Ornstein-Uhlenbeck process
\begin{equation}
\dif X^{t}=-\frac{1}{2}X^{t}+\frac{\dif B_t^{(r)}}{\sqrt{N}}, \quad X^{0}:=X,\quad t\in [0,t_1],
\label{eq:OU_real}
\end{equation}  
where $B_t^{(r)}$ is an $N\times N$ matrix composed of $N^2$ independent real-valued Brownian motions. By \eqref{eq:OU_distr}, $X^{t_1}$ has a real Ginibre component of order $t_1$. This step is only needed to introduce a real Ginibre component in $X$ so that the local laws from Proposition~\ref{prop:real_2G_av} become better, which will be used in the following steps. In particular, we do not follow any estimates along the flow \eqref{eq:OU_real}.

\medskip

\noindent\textbf{Step 1.2.} Add a complex Ginibre component to $X^{[t_1]}$ running the following flow for time $t_2$:
\begin{equation}
\dif X^{[t_1]}_t=\frac{\dif B_t^{(c)}}{\sqrt{N}},\quad X^{[t_1]}_0:=X^{t_1},\quad t\in [0,t_2],
\label{eq:BM_for_real}
\end{equation}
where $B_t^{(r)}$ is an $N\times N$ matrix composed of $N^2$ independent complex-valued Brownian motions. In \eqref{eq:BM_for_real} we use the square brackets around $t_1$ to indicate that this time parameter comes from the initial condition and is fixed along the flow.  This step is the same as the first step of the strategy outlined in Section~\ref{sec:proofprop}, and relies on the properties of the complex DBM from \cite[Proposition~4.6]{bourgade2024fluctuations} and Theorem~\ref{eq:maintheoDBM1}. It establishes the desired result for $X^{[t_1]}_{t_2}$.

\medskip

\noindent\textbf{Step 2.1.} The complex Ginibre component is removed from $X^{[t_1]}_{t_2}$ using the flow 
\begin{equation}
\dif X^{[t_1],t}=-\frac{1}{2}X^{[t_1],t}+\frac{\dif B_t^{(c)}}{\sqrt{N}},\quad X^{[t_1],0}:=X^{t_1},\quad t\in [0,s_2],\quad s_2:=\log(1+t_2),
\label{eq:OU_compl_for_real}
\end{equation}
and \eqref{eq:strat_concl} to match this step with the previous one. We note that one does not need to assume anything about the joint distribution of driving Brownian motions in \eqref{eq:BM_for_real} and \eqref{eq:OU_compl_for_real}, since we only need to match $X^{[t_1]}_{t_2}$ with $X^{[t_1],s_2}$ in distribution by \eqref{eq:strat_concl}. This step contains a technical novelty compared to the complex case, since we need to estimate the contribution from the second order cumulants in the GFT, as discussed above.

\medskip

\noindent\textbf{Step 2.2.} The real Ginibre component is removed from $X^{t_1}$ using the flow \eqref{eq:OU_real}. Since \eqref{eq:OU_real} preserves the second-order correlation structure, this step is analogous to the second step of the strategy outlined in Section~\ref{sec:proofprop} and does not require any additional technical inputs. 

\medskip

We stress once again that the proof strategy presented above is shaped by the two technical resolutions, independent from each other. First, we avoid real DBM and transfer difficulties to the GFT, which is the content of Steps~1.2 and~2.1. Second, we do not perform the GFT in the local law in Proposition~\ref{prop:real_2G_av}, but remove the real Gaussian component directly from the main quantities in Propositions~\ref{prop:real_tail_bound}, \ref{prop:real_EG}, and~\ref{prop:real_EGG}. This necessitates performing Steps~1.1 and~2.2. A different distribution of technical difficulties would lead to a different strategy.

The only step in our strategy which requires an additional argument compared to the complex case, is Step~2.1. In contrast, Step~1.1 does not require following any estimates along the flow, as noted above; Step~1.2 coincides with the first step of the strategy introduced in the complex case in Section~\ref{sec:proofprop}; and Step~2.2 is analogous to the second step of that strategy. In the remaining part of this section we prove Propositions~\ref{prop:real_tail_bound} and~\ref{prop:real_EGG} focusing on the novel Step~2.1. The adjustments required for the proof of Proposition~\ref{prop:real_EG} compared to the proof of Proposition~\ref{prop:EG} are similar to the ones required for the proof of Proposition~\ref{prop:real_EGG}. Thus, we omit the proof of Proposition~\ref{prop:real_EG}.

\subsubsection{Proof of Proposition \ref{prop:real_tail_bound}}  
Denote the least positive eigenvalue of the Hermitization of $X^{t}-z$ by $\lambda_1^{t,z}$. Since \eqref{eq:OU_real} preserves the first two moments of matrix entries, arguing similarly to the GFT in the proof of Proposition~\ref{prop:tail_bound} (see in particular \eqref{eq:F_GFT_propag}--\eqref{eq:F_iteration}) we get that 
\begin{equation}
\mathbf{P}[\lambda^{0,z}_{1}\le N^{-1}x] \lesssim \mathbf{P}[\lambda^{t_1,z}_1\le N^{-1}2x]+N^{-D},
\label{eq:real_tail_step1}
\end{equation}
for any fixed (small) $\xi>0$, (large) $D>0$, $N^{-1/6+\xi}\le x\le 1$ and $t_1\le N^{-1/2-\xi}x^3$, where the implicit constant in \eqref{eq:real_tail_step1} depends on $\xi$ and~$D$. We take $t_1:=N^{-1/2+\xi}x^3$ and complete the first and the last steps of the strategy explained around \eqref{eq:OU_real}--\eqref{eq:OU_compl_for_real}. To simplify the notation, we further assume that $X$ has a real Gaussian component of order~$t_1$ and drop the superscript $[t_1]$ in the notations introduced in \eqref{eq:BM_for_real} and \eqref{eq:OU_compl_for_real}. We also use the notations introduced in the proof of Proposition \ref{prop:tail_bound} without explicitly mentioning this further.

It is easy to see that \eqref{eq:tail_t} holds also in the real case, so we are left with the GFT removing a small complex Ginibre component of order $t_2$ added to a real i.i.d. matrix. We perform a calculation similar to \eqref{eq:dF} and get in the rhs. of \eqref{eq:dF} the following additional term coming from the second order cumulants:
\begin{equation}
\E\sum_{\alpha,\beta} \left(\E[w_{\alpha,s}w_{\beta,s}] - \E [\widetilde{w}_\alpha\widetilde{w}_\beta]\right)\partial_\alpha\partial_\beta F\left(\mathrm{Tr}\bm{1}_{E+l}*\theta_\eta(H^{s,z})\right).
\label{eq:tail_GFT_2nd_order}
\end{equation}
Here $\alpha,\beta\in ([N]\times [N+1,2N])\cup ([N+1,2N]\times [N])$, $\widetilde{W}=(w_{ab})_{a,b\in[2N]}$ is a Hermitization of an independent GinUE matrix and $w_{\alpha,s}$, $w_{\beta,s}$ are the entries of the Hermitization of $X^s$. Performing the differentiation in \eqref{eq:tail_GFT_2nd_order}, we arrive to two types of terms. The first type is obtained by differentiating $F$ only once, and once the derivative of the internal function $\mathrm{Tr}\bm{1}_{E+l}*\theta_\eta(H^{s,z})$, while in the terms of the second type $F$ is differentiated twice. Summing up the terms of the first type over $\alpha,\beta$, we get
\begin{equation}
\E\left[F'\left(\mathrm{Tr}\bm{1}_{E+l}*\theta_\eta(H^{s,z})\right) \Im \int_{-(E+l)}^{E+l} \sigma \langle G^2(y+\ii\eta)E_\sigma G^T(y+\ii\eta)E_\sigma\rangle\dif y\right],
\label{eq:tail_GFT_summed1}
\end{equation}
up to an $s$-dependent factor of order one. We also get several terms which differ from \eqref{eq:tail_GFT_summed1} only by conjugating some of the resolvents, but their treatment is identical to the analysis of \eqref{eq:tail_GFT_summed1} below and thus is omitted. Define the exponent $b$ by $t_1=N^{-b}$ and observe that
\begin{equation}
\widehat{\beta}_{1\overline{1}}^{[b]}(y+\ii\eta,y+\ii\eta) \ge  \widehat{\beta}_{1\overline{1}}(y+\ii\eta,y+\ii\eta)\wedge N^{-b}\gtrsim |\Im z|^2\wedge t_1,
\label{eq:real_beta}
\end{equation}
where we used \eqref{eq:beta_hat_b} and Proposition~\ref{prop:stab_bound}. We apply \eqref{eq:real_size_bound2} from Lemma~\ref{lem:G_bounds_extension} to the three-resolvent chain in \eqref{eq:tail_GFT_summed1} and obtain the following upper bound on the absolute value of~\eqref{eq:tail_GFT_summed1}:
\begin{equation}
\E\left[F\left(\mathrm{Tr}\bm{1}_{E+3l}*\theta_\eta(H^{s,z})\right)\right]\frac{E+l}{N^2\eta^3}\left(\frac{1}{|\Im z|^2}+\frac{1}{t_1}\right)\lesssim N^{6\epsilon} x^{-2}\!\!\left(\frac{1}{|\Im z|^2}\!+\!\frac{1}{t_1}\right)\E\left[F\left(\mathrm{Tr}\bm{1}_{E+3l}*\theta_\eta(H^{s,z})\right)\right],
\label{eq:tail_GFT_bound1}
\end{equation}
up to an $N^\xi$ factor. Here we bounded
\begin{equation}
\E\left[ \left\vert F'\left(\mathrm{Tr}\bm{1}_{E+l}*\theta_\eta(H^{s,z})\right)\right\vert\right]\lesssim\E\left[F\left(\mathrm{Tr}\bm{1}_{E+3l}*\theta_\eta(H^{s,z})\right)\right]+\mathcal{O}(N^{-D}),
\label{eq:F_prime}
\end{equation}
arguing similarly to \cite[Eq.(2.43)--(2.46)]{small_dev}, and used the relation between $x,E,l$ and $\eta$ from~\eqref{eq:F_GFT_param}.

Summing up the terms in \eqref{eq:tail_GFT_2nd_order} where $F$ is differentiated twice, we get
\begin{equation}
\E\left[F''\left(\mathrm{Tr}\bm{1}_{E+l}*\theta_\eta(H^{s,z})\right)\int_{-(E+l)}^{E+l}\int_{-(E+l)}^{E+l} \sigma \left\langle \Im G^2(y_1+\ii\eta)E_\sigma \left(\Im G^2(y_2+\ii\eta)\right)^TE_\sigma\right\rangle\dif y_1\dif y_2\right],
\label{eq:tail_GFT_summed2}
\end{equation}
up to an $s$-dependent factor of order one. Using the analogue of \eqref{eq:F_prime} for $F''$, \eqref{eq:F_GFT_param}, \eqref{eq:real_beta}, and \eqref{eq:real_size_bound2} to bound the four-resolvent chain in \eqref{eq:tail_GFT_summed2}, we obtain that the absolute value of \eqref{eq:tail_GFT_summed2} has an upper bound of order
\begin{equation}
\E\left[F\left(\mathrm{Tr}\bm{1}_{E+3l}\!*\!\theta_\eta(H^{s,z})\right)\right]\! \frac{(E+l)^2}{N^2\eta^4}\!\left(\frac{1}{|\Im z|^2}\! + \!\frac{1}{t_1}\right)\!\lesssim\! N^{8\epsilon} x^{-2}\!\!\left(\frac{1}{|\Im z|^2}\!+\!\frac{1}{t_1}\right)\E\left[F\left(\mathrm{Tr}\bm{1}_{E+3l}*\theta_\eta(H^{s,z})\right)\right]\!.
\label{eq:tail_GFT_bound2}
\end{equation}
up to an $N^\xi$ factor. Combining \eqref{eq:tail_GFT_bound1} and \eqref{eq:tail_GFT_bound2}, and estimating the contribution from the third and higher-order cumulants as in the proof of Proposition~\ref{prop:tail_bound}, we get the following analogue of~\eqref{eq:dF}:
\begin{equation}
\begin{split}
&\frac{\dif}{\dif s}\E\left[F\left(\mathrm{Tr}\bm{1}_{E+l}*\theta_\eta(H^{s,z})\right)\right]\lesssim\\
&\qquad N^{8\epsilon+\xi}t_2\left(x^{-2}\left(\frac{1}{|\Im z|^2}+\frac{1}{t_1}\right)+N^{1/2}x^{-3}\right)\E\left[F\left(\mathrm{Tr}\bm{1}_{E+3l}*\theta_\eta(H^{s,z})\right)\right] +N^{-D},
\end{split}
\label{eq:real_F_dif}
\end{equation}
for any fixed $D>0$.

Iterating \eqref{eq:real_F_dif} as explained in \eqref{eq:F_iteration}, we obtain
\begin{equation}
\mathbf{P}[\lambda^{z}_1\le  N^{-1}x] \le (\log N) x^2 + N^\xi (N\mathcal{E}_0(t_2))^2
\label{eq:real_tail_t}
\end{equation}
for any fixed $\xi>0$ and $t_2$ satisfying conditions
\begin{equation}
N^{-1+\xi}\le t_2\le N^{-8\epsilon -\xi}\left(x^{-2}\left(\frac{1}{|\Im z|^2}+\frac{1}{t_1}\right)+N^{1/2}x^{-3}\right)^{-1},\quad\text{where}\quad t_1:=N^{-1/2-\xi}x^3.
\label{eq:t2_range}
\end{equation}
Finally, we optimize \eqref{eq:real_tail_t} over $t_2$ by choosing the maximal possible $t_2$ in \eqref{eq:t2_range}, take $\epsilon>0$ to be sufficiently small and complete the proof of Proposition~\ref{prop:real_tail_bound}.\hfill $\qed$

\subsubsection{Proof of Proposition \ref{prop:real_EGG}} The term in the first line of \eqref{eq:real_def_Phi2} arises from Step~1.2 introduced around \eqref{eq:BM_for_real}, while the last term in the second line is picked up from Step~2.1. Since these estimates are completely analogous to the ones in the proof of Proposition~\ref{prop:EGG}, we further focus on the analysis of the second order cumulant terms in the GFT in Step~2.1. As in the proof of Proposition~\ref{prop:real_tail_bound}, we may assume that $X$ contains a real Ginibre component of size $t_1$, so we further drop the superscript $[t_1]$ in \eqref{eq:OU_compl_for_real}. 

Denote $G_j:=G^{z_j,s}(\ii\eta_j)$, for $j=1,2$, and differentiate $\mathrm{Cov}(\langle G_1\rangle,\langle G_2\rangle)$ along the flow \eqref{eq:OU_compl_for_real}. Performing the cumulant expansion similarly to \eqref{eq:cum_exp_first}, we get the following contribution from the second order cumulant terms:
\begin{equation}
\frac{\sigma}{2N}\E\left[ \left\langle G_1^2E_\sigma G_1^\mt E_\sigma\right\rangle\left(\langle G_2\rangle -\E\langle G_2\rangle\right)+\left(\langle G_1\rangle -\E\langle G_1\rangle\right)\left\langle G_2^2E_\sigma G_2^\mt E_\sigma \right\rangle + \frac{1}{N}\left\langle G_1^2 E_\sigma (G_2^2)^\mt E_\sigma\right\rangle\right],
\label{eq:real_Cov_2nd_order}
\end{equation}
up to an $s$-dependent factor of order one. By \cite[Lemma~B.7]{cipolloni2024maximum} and a standard monotonicity argument (see e.g. \cite[Lemma~5.3]{cusp_univ}), we have 
\begin{equation*}
\left\vert\langle G_j\rangle -\E\langle G_j\rangle\right\vert\prec 1+\frac{1}{N\eta_j},\quad j=1,2.
\end{equation*}
Using Lemma~\ref{lem:G_bounds_extension} to estimate the rest of the terms in \eqref{eq:real_Cov_2nd_order} and employing \eqref{eq:real_beta}, we get that the absolute value of \eqref{eq:real_Cov_2nd_order} has an upper bound of order
\begin{equation*}
N^{1+\xi}\left(\frac{1}{N|\Im z_1|^2}+\frac{1}{N|\Im z_2|^2}+ \frac{1}{N|z_1-\overline{z}_2|^2} +\frac{1}{Nt_1}\right)\left(1+\frac{1}{N\eta_1}\right)^2\left(1+\frac{1}{N\eta_2}\right)^2,
\end{equation*}
where the term containing $|z_1-\overline{z}_2|$ comes from the upper bound on the last term in the rhs. of \eqref{eq:real_Cov_2nd_order}. Therefore, the contribution from the second order cumulants in the GFT in Step~2.1 is bounded by the first term in the second line of \eqref{eq:real_def_Phi2}, while the contribution from the higher order cumulants is dominated by the last term in the second line of \eqref{eq:real_def_Phi2}. This finishes the proof of Proposition~\ref{prop:real_EGG}.\hfill$\qed$

\section{Proofs of the additional technical results}

\subsection{Relaxation of Assumption \ref{ass:Omega}}\label{app:Omega} In this section we show that Theorem~\ref{theo:main1} holds for an open simply connected domain $\Omega_N$ given by
\begin{equation}
\Omega_N:=z_N + N^{-\alpha}\widetilde{\Omega}\subset (1-\delta)\mathbf{D}
\label{eq:Omega_relaxed}
\end{equation}
for some fixed $\delta>0$ and $\alpha\in [0,1/2)$, where $z_N\in \C$ and $\widetilde{\Omega}\subset \C$ is a domain with a piecewise $C^2$ boundary. The model \eqref{eq:Omega_relaxed} extends the class of domains introduced in Assumption~\ref{ass:Omega} by allowing $\partial\Omega_N$ to have angles. In particular, $\partial\Omega_N$ can be a polygon, which is not compatible with the smoothness condition from Assumption~\ref{ass:Omega}. Strictly speaking, \eqref{eq:Omega_relaxed} does not include all domains satisfying Assumption~\ref{ass:Omega}, since $\widetilde{\Omega}$ does not depend on~$N$ in \eqref{eq:Omega_relaxed}. This assumption can be further relaxed, but we prefer to fix $\widetilde{\Omega}$ to avoid further technical assumptions.

We start with discussion of modifications needed for the proof of Proposition~\ref{prop:Var_main} in the set-up~\eqref{eq:Omega_relaxed}. There are two main aspects in the proof of these results where the smoothness of $\partial\Omega_N$ is used. The first is the bound on the area of the tubular neighborhood $\Delta\Omega_N$ of $\partial\Omega_N$ in \eqref{eq:Delta_Omega_area}, while the second is the bound on the regularized singular integrals over $(\Delta\Omega_N)^2$ in Lemma~\ref{lem:int}. In fact, both of these results hold under much weaker assumption on $\partial\Omega_N$ compared to smoothness. One only needs to assume that for any $x>0$ there exists a constant $C(x)>0$ such that
\begin{equation}
\left\vert \mathfrak{N}_h(\partial\Omega_N)\cap B_r(z)\right\vert \le C(x)hr\quad \text{for}\quad h:=xN^{-a},\, \forall r\in [h, \mathrm{diam}(\Omega_N)],\, \forall z\in \partial\Omega_N,\, \forall N\in\N,
\label{eq:reg_boundary}
\end{equation} 
where $\mathrm{diam}$, $\mathfrak{N}_h$ and $B_r$ are defined in \eqref{eq:diam}, \eqref{eq:def_Nh} and below \eqref{eq:def_Nh}, respectively. Since \eqref{eq:reg_boundary} holds for any~$C^2$ curve, it also holds for a finite collection of such curves, in particular for $\partial \Omega_N$ in the case when it is a piecewise~$C^2$ curve. Therefore, Proposition~\ref{prop:Var_main} holds for the model~\eqref{eq:Omega_relaxed}.

Now we discuss the modifications needed for the proof of Theorem~\ref{theo:main1}, were the smoothness of $\partial\Omega_N$ was additionally used in the construction of domains $\Omega_N^\pm$ in \eqref{eq:Omega_pm}. We keep this construction and note that for $\Omega_N$ satisfying \eqref{eq:Omega_relaxed}, the domains $\Omega_N^\pm$ are not necessarily of the form \eqref{eq:Omega_relaxed}. However, it is easy to see that $\Omega_N^\pm$ satisfy \eqref{eq:reg_boundary}, so Proposition~\ref{prop:Var_main} can be applied to these domains in the same way as discussed in Section~\ref{sec:proof_thm1}. This finishes the proof of Theorem~\ref{theo:main1} for $\Omega_N$ satisfying~\eqref{eq:Omega_relaxed}.

\subsection{Proof of \eqref{eq:tails_UV}}\label{app:cut_tails} Recall that the non-zero contribution in the lhs. of~\eqref{eq:tails_UV} comes only from $z_l\in \Delta\Omega_N$, $l=1,2$, and that $\Delta\Omega_N\subset (1-\delta/2)\mathbf{D}$. From now on, we consider only $|z_l|\le 1-\delta/2$. Our proof relies on \cite[Eq. (4.24)]{macroCLT_complex}, which states that
\begin{equation}
|V_{12}|\lesssim \frac{\left[(\eta_1+\rho_1)(\eta_2+\rho_2)\right]^{-2}}{|z_1-z_2|^4 + (\eta_1+\eta_2)^2 \min\{\rho_1,\rho_2\}^4},\quad |U_l|\lesssim \frac{1}{\rho_l^2 +\eta_l^3},
\label{eq:UV_bound}
\end{equation}
where $\rho_l=\rho^{z_l}(\ii\eta_l)$ for $l=1,2$. We start with estimating the integral of $U_1U_2$ term in~\eqref{eq:tails_UV}. Since $|z_l|\le 1-\delta/2$, it holds that $\rho_l\sim 1$ for $\eta_l\in [0,\eta_c]$, so~\eqref{eq:UV_bound} implies
\begin{equation*}
\left(\int_0^{\eta_c}+\int_T^\infty\right)U_i\dif\eta_i\lesssim \eta_c +T^{-2}\lesssim \eta_c.
\end{equation*}
Therefore, the integral of $U_1U_2$ in \eqref{eq:tails_UV} has an upper bound of order
\begin{equation}
\lVert \Delta f_N\rVert_1^2\eta_c^2\lesssim N^{2(a-\alpha-1)-2\delta_c},
\label{eq:app_U_bound}
\end{equation}
where we additionally used \eqref{eq:Delta_f_L1}.

Next, we estimate the integral of $V_{12}$ in \eqref{eq:tails_UV} in the regime $\eta_1,\eta_2\in [0,\eta_c]$, where the first bound in \eqref{eq:UV_bound} simplifies to
\begin{equation}
|V_{12}|\lesssim \frac{1}{|z_1-z_2|^4 + (\eta_1+\eta_2)^2},
\label{eq:V_bound_small_eta}
\end{equation}
since $\rho_1\sim\rho_2\sim 1$. Integrating \eqref{eq:V_bound_small_eta} over $z_1,z_2\in (\Delta \Omega_N)^2$ by the means of Lemma~\ref{lem:int} applied to $q=4$ and $r_N:=(\eta_1+\eta_2)^2$, we get
\begin{equation}
\begin{split}
\int_{\Delta \Omega_N}\int_{\Delta \Omega_N} |V_{12}|\dif^2 z_1\dif^2 z_2 &\lesssim N^{-\alpha-2a}(\eta_1+\eta_2)^{-3/2} |\log (\eta_1+\eta_2)|\lesssim N^{-\alpha-2a}(\eta_1+\eta_2)^{-3/2-\xi}\\
&\lesssim N^{-\alpha-2a} \eta^{-3/4-\xi/2}_1\eta^{-3/4-\xi/2}_2
\end{split}
\label{eq:int_V_bouns_small_eta}
\end{equation}
for any fixed $\xi>0$. Here we estimated $|\log x|\lesssim x^{-\xi}$ for $x\in (0,2\eta_c]$. We integrate \eqref{eq:int_V_bouns_small_eta} over $\eta_1,\eta_2\in [0,\eta_c]$, recall that $\|\Delta f_N\|_\infty\lesssim N^{2a}$, and conclude that 
\begin{equation}
\int_{\C}\!\int_{\C}\! |\Delta f_N(z_1)||\Delta f_N(z_2)|\!\left(\int_0^{\eta_c}\!\!\!\int_0^{\eta_c}\!|V_{12}|\dif\eta_1\dif\eta_2\right) \dif^2 z_1\dif^2 z_2\lesssim N^{2a-\alpha}\eta_c^{1/2-\xi}\lesssim N^{2a-\alpha -1/2-\delta_c/2+\xi}.
\label{eq:app_V_bound}
\end{equation}

In order to deal with the remaining regime $\max\{\eta_1,\eta_2\}\ge T$, recall from \eqref{eq:def_UV} that
\begin{equation*}
V_{12}=\frac{1}{2}\partial_{\eta_1}\partial_{\eta_2} \log A,\quad \text{where}\quad A=1-u_1u_2\left( 1-|z_1-z_2|^2 + (1-u_1) |z_1|^2 + (1-u_2)|z_2|^2\right).
\end{equation*}
Since $A$ goes to $1$ when $\max\{\eta_1,\eta_2\}$ goes to infinity, we have
\begin{equation*}
\begin{split}
\int_0^{\eta_c}\int_T^\infty V_{12}\dif\eta_1\dif\eta_2  &= - \frac{1}{2}\int_0^{\eta_c} \partial_{\eta_1} \log A(z_1,z_2,\eta_1,T)\dif\eta_1\\
&\lesssim \left\vert \log A(z_1,z_2,0,T)\right\vert + \left\vert \log A(z_1,z_2,\eta_c,T)\right\vert\lesssim T^{-1},
\end{split} 
\end{equation*}
where we used that $|u_l|\lesssim 1/(1+\eta_l)$, $l=1,2$. The bound in the integration regime $\eta_1,\eta_2\in [T,\infty)$ follows similarly. Together with \eqref{eq:app_U_bound} and \eqref{eq:app_V_bound} this finishes the proof of \eqref{eq:tails_UV}. \qed

\subsection{Proof of \eqref{eq:multiG_size}}\label{app:multiG_size} Since $\overline{\bm{B}}^{[1]}$ consists of $n-1$ blocks containing $S-k_1$ resolvents in total, the averaged multi-resolvent local law \eqref{eq:multiG_av} implies that 
\begin{equation*}
\left\vert\Delta\left(\overline{\bm{B}}^{[1]}\right)\right\vert \prec N^{-(n-1)} \eta_1^{-S+k_1}.
\end{equation*}
Thus, to prove \eqref{eq:multiG_size} it suffices to show that
\begin{equation}
\left\vert\langle G_1 B_1G_1\cdots G_1 B_{k-1}G_1B_kG_2^2 B_{k+1}\rangle\right\vert \prec \frac{1+\bm{1}_{k>1}(N\eta_*)^{1/2}}{\eta_*^{k-1}\eta_2\gamma},
\label{eq:multiG_size_aux}
\end{equation}
where we denoted $k:=k_1$, $B_j:=B_j^{(1)}$ for $j\in[k-1]$, $B_{k}:=AE_\sigma$ and $B_{k+1}:=E_\sigma$. Consider first the case $k=1$, where~\eqref{eq:multiG_size_aux} is of the form
\begin{equation}
\left\vert \langle G_1B_1G_2^2 B_2\rangle\right\vert \prec \frac{1}{\eta_2\gamma}.
\end{equation}
This bound immediately follows from \eqref{eq:stab_LT} and the 2-resolvent averaged local law from Proposition~\ref{prop:2G_av} after representing $G_2^2$ in terms of a contour integral as it was done in~\eqref{eq:G^2_rep}.

For $k=2$ we apply Cauchy-Schwarz inequality and get
\begin{equation}
\begin{split}
\left\vert \langle G_1B_1G_1B_2G_2^2B_3\rangle\right\vert &\le \left\langle  |G_1B_2G_2|^2\right\rangle^{1/2}\left\langle |G_2B_3G_1B_1|^2\right\rangle^{1/2} \le \lVert B_1\rVert  \left\langle  |G_1B_2G_2|^2\right\rangle^{1/2}\left\langle |G_2B_3G_1|^2\right\rangle^{1/2}\\
&= \frac{\lVert B_1\rVert}{\eta_1\eta_2}\langle \Im G_1B_2\Im G_2 B_2^*\rangle^{1/2}\langle \Im G_1B_3^*\Im G_2 B_3\rangle^{1/2}\prec \frac{1}{\eta_1\eta_2\gamma},
\end{split}
\label{eq:app_chain_k=2}
\end{equation}
i.e. \eqref{eq:multiG_size_aux} holds for $k=2$. In \eqref{eq:app_chain_k=2} we used Ward identity to go from the first to the second line and employed Proposition~\ref{prop:2G_av} in the last bound.

Now we focus on the case $k\ge 3$ and present a unified argument for all these values of $k$. Estimating the lhs. of \eqref{eq:multiG_size_aux} via a Cauchy-Schwarz inequality followed by Ward identity we get that it admits the following upper bound
\begin{equation}
\eta_1^{-1}\langle G_2^2 B_{k+1}\Im G_1 B_{k+1}^* (G_2^*)^2 B_k^*\Im G_1 B_k \rangle^{1/2}\langle B_1G_1\cdots G_1B_{k-1}B_{k-1}^* G_1^*\cdots G_1^* B_1^*\rangle^{1/2}.
\label{eq:reduction_Schwarz} 
\end{equation}
In the second trace in \eqref{eq:reduction_Schwarz} only one type of resolvents appears, and their number equals to $2k-4$. Thus, \eqref{eq:multiG_av} along with \eqref{eq:multiM_bound} give an upper bound of order $\eta_1^{-k+5/2}$ on the aquare root of this trace. From \cite[Eq.(5.27)]{eigenv_decorr} we have that the first trace in \eqref{eq:reduction_Schwarz} has an upper bound of order
\begin{equation}
N\langle |G_2|^2 B_{k+1}\Im G_1 B_{k+1}^*\rangle \langle |G_2|^2 B_k^*\Im G_1 B_k\rangle.
\end{equation}
Let $B$ be equal either to $B_{k+1}$ or to $B_k^*$. Since $B\Im G_1 B^*\ge 0$, we have
\begin{equation}
\langle |G_2|^2 B\Im G_1 B^*\rangle \le \frac{1}{\eta_2}\langle |G_2|B\Im G_1B^*\rangle\lesssim \frac{\log N}{\eta_2\gamma},
\end{equation}
with very high probability, where in the last bound we argued as in \eqref{eq:2G_abs_bound}. Collecting all the bounds presented above we finish the proof of \eqref{eq:multiG_size_aux}.

\subsection{Proof of Lemma \ref{lem:beta_hat}: properties of $\hat{\beta}_{12}$} \label{app:beta_hat} 

\underline{Proof of \eqref{eq:beta_hat_to_beta}}. Throughout the proof we consider the action of $\mathcal{B}_{12}$ and $\mathcal{B}_{12}^{-1}$ only on the $2\times 2$ block-constant matrices, i.e. on the elements of $\mathrm{span}\{E_\pm,F^{(*)}\}$, without explicitly mentioning this further. First, we show that
\begin{equation}
\|\left(\mathcal{B}_{12}^{-1}(w_1,w_2)\right)[R]\|\lesssim \frac{1}{\beta_{12,*}(w_1,w_2)\wedge 1},\quad \forall R\in\mathrm{span}\{E_\pm,F^{(*)}\}.
\label{eq:beta_bound} 
\end{equation}
Once \eqref{eq:beta_bound} is obtained, the proof of \eqref{eq:beta_hat_to_beta} goes as follows. We have from \eqref{eq:beta_bound} that
\begin{equation}
\min_{\|R\|=1} \left\lVert\left(\mathcal{B}_{12}(w_1,w_2)\right)[R]\right\rVert\gtrsim \beta_{12,*}(w_1,w_2)\wedge 1.
\label{eq:B12_lower_bound}
\end{equation}
On the other hand, $\beta_{12,*}(w_1,w_2)\wedge 1$ is the absolute value of the closest to zero eigenvalue of $\mathcal{B}_{12}(w_1,w_2)$ (see the discussion above Lemma \ref{lem:beta_hat}). Therefore, the inequality in \eqref{eq:B12_lower_bound} holds also in the reversed direction, i.e. both sides of \eqref{eq:B12_lower_bound} are of the same order. Taking the minimum of \eqref{eq:B12_lower_bound} over all choices of $w_j$ and $\overline{w}_j$, and recalling the definition of $\widehat{\beta}_{12}(w_1,w_2)$ from \eqref{eq:def_beta_hat}, we conclude the proof of \eqref{eq:beta_hat_to_beta}. Thus, it remains to prove \eqref{eq:beta_bound}.

From now on and up to the end of the proof of \eqref{eq:beta_bound} we drop the $(w_1,w_2)$-dependence from the notations $\mathcal{B}_{12}(w_1,w_2)$ and $\beta_{12,\pm}(w_1,w_2)$. We prove \eqref{eq:beta_bound} by explicitly inverting the two-body stability operator. For any $R\in\mathrm{span}\{E_\pm,F^{(*)}\}$ denote $X:=\mathcal{B}_{12}^{-1}[R]$. Using the explicit form of $\mathcal{S}$ from \eqref{eq:def_S} we get
\begin{equation}
\mathcal{P}
\begin{pmatrix}
\langle X\rangle\\ \langle XE_-\rangle
\end{pmatrix}
=
\begin{pmatrix}
\langle R\rangle\\ \langle RE_-\rangle
\end{pmatrix},\quad \mathcal{P}:= \begin{pmatrix}
1-\langle M_1M_2\rangle&\langle M_1E_-M_2\rangle\\ 
-\langle M_2E_-M_1\rangle & 1+\langle M_1E_-M_2E_-\rangle
\end{pmatrix},
\label{eq:B_inverse_aux}
\end{equation}
where $M_j:=M^{z_j}(w_j)$, $j=1,2$. Denote the eigenvectors of $\mathcal{B}_{12}$ associated to $\beta_{12,\pm}$ by $X_\pm$. It immediately follows from the eigenvector equations for $X_\pm$ that $\mathcal{S}[X_+]$ and $\mathcal{S}[X_-]$ are not colinear. Therefore, $\beta_{12,\pm}$ are the eigenvalues of $\mathcal{P}$, while the corresponding eigenvectors are given by $(\langle X_\pm \rangle,\langle X_\pm E_-\rangle)\in\C^2$. Next, observe that
\begin{equation}
\max\left\lbrace |\beta_{12,+}|,|\beta_{12,-}|\right\rbrace\sim 1.
\label{eq:beta_bulk}
\end{equation}
Indeed, from \cite[Eq. (3.5)]{macroCLT_real} we have that $|u_{j,r}|<1$, so
\begin{equation}
\Re\left[ 1-\Re [z_1\overline{z}_2]u_1u_2\right] \ge 1-|z_1z_2|\ge \delta,
\label{eq:Re_beta}
\end{equation}
where in the last step we estimated $|z_j|\le 1-\delta$ for $j=1,2$. Therefore, either $\Re \beta_{12,+}\ge \delta$ or $\Re \beta_{12,-}\ge \delta$, depending on the sign of the real part of the square root in the rhs. of \eqref{eq:def_beta_pm}, i.e. \eqref{eq:beta_bulk} holds. We conclude from \eqref{eq:beta_bulk} that
\begin{equation}
\left\vert\det \mathcal{P}\right\vert = \left\vert \beta_{12,+}\beta_{12,-}\right\vert \sim \beta_{12,*}.
\label{eq:detP}
\end{equation}
Together with the fact that the entries of $\mathcal{P}$ have an upper bound of order one in absolute value, this yields $\|\mathcal{P}^{-1}\|\lesssim \beta_{12,*}^{-1}$. Combining this bound with \eqref{eq:B_inverse_aux} we get $\|\mathcal{S}[X]\|\lesssim \lVert \mathcal{P}^{-1}\rVert\lVert \mathcal{S}[R]\rVert\lesssim \beta_{12,*}^{-1}\|R\|$, which yields
\begin{equation}
\lVert\mathcal{B}_{12}^{-1}[R]\rVert=\lVert X\rVert = \lVert R+M_1\mathcal{S}[X]M_2\rVert \lesssim \lVert R\rVert + \lVert \mathcal{S}[X]\rVert\lesssim (1 +\beta_{12,*}^{-1})\|R\|.
\end{equation}
This completes the proof of \eqref{eq:beta_bound}.

\medskip

\noindent\underline{Proof of \eqref{eq:gamma_monot}.} Note that the second part of \eqref{eq:gamma_monot} immediately follows from \eqref{eq:beta_hat_to_beta} and the first part of \eqref{eq:gamma_monot}, so we focus on \eqref{eq:gamma_monot} for $\beta_{12,*}$. From Lemma \ref{lem:char_flow} we have that $z_{j,r}, m_{j,r}, u_{j,r}$ scale with time as
\begin{equation}
z_{j,r}=\ee^{-r/2}z_{j,0},\quad m_{j,r}=\ee^{r/2} m_{j,0},\quad u_{j,r}=\ee^r u_{j,0},\qquad \forall \,r\in [0,T], \, j=1,2.
\label{eq:time_scaling}
\end{equation}
Together with \eqref{eq:def_beta_pm} this gives
\begin{equation}
\beta_{12,\sigma,r} = 1-\ee^{r}(1-\beta_{12,\sigma,0}),\quad \sigma\in\{\pm\},\,\, r\in [0,T].
\label{eq:beta_time_dep}
\end{equation}
From \eqref{eq:beta_time_dep} and the trivial bound $|\beta_{12,\sigma,r}|\lesssim 1$ we have that
\begin{equation}
|\beta_{12,\sigma,s}-\beta_{12,\sigma,t}|\lesssim |s-t|,\quad \forall\, \sigma\in\{\pm\},\, s,t\in [0,T].
\label{eq:beta_dif}
\end{equation}
Combining the stability bound $|\beta_{12,\sigma,s}|\gtrsim \eta_{*,s}$ which follows from \eqref{eq:eta_stab_bound} and \eqref{eq:beta_hat_to_beta}, with the estimate $\eta_{*,s}\gtrsim T-s$ from Lemma \ref{lem:char_flow}(4), we get $|\beta_{12,\sigma,s}|\gtrsim |T-s|$. Together with \eqref{eq:beta_dif} and \eqref{eq:beta_time_dep} this immediately implies that 
\begin{equation}
|\beta_{12,\sigma,s}|\sim |\beta_{12,\sigma,t}| + |t-s|,\quad \forall\, \sigma\in\{\pm\},\, s,t\in [0,T],\, s\le t.
\label{eq:beta_monot}
\end{equation}
Finally, taking the minimum over $\sigma\in\{\pm\}$ in \eqref{eq:beta_monot}, we complete the proof of \eqref{eq:gamma_monot}.

\medskip

\noindent\underline{Proof of \eqref{eq:beta_hat_perturb}.} Since $|\partial_w m^z(w)|\lesssim 1$ for any $z\in (1-\delta)\mathbf{D}$ and $w\in \mathcal{D}^z_{\kappa,\varepsilon}$, from the explicit equations for $\beta_{12, \pm}$ given in \eqref{eq:def_beta_pm} we get that
\begin{equation}
|\beta_{12,\sigma}(w_1,w_2)|\lesssim |\beta_{12,\sigma}(w_1,w_2')|+|w_2-w_2'|,\quad \sigma\in\{\pm\},
\label{eq:beta_sigma_perturb}
\end{equation}
under the assumptions of Lemma \ref{lem:beta_hat}(3). By taking minimum of \eqref{eq:beta_sigma_perturb} over $\sigma\in\{\pm\}$ we obtain that \eqref{eq:beta_sigma_perturb} holds for $\beta_{12,*}$ as well, which together with \eqref{eq:beta_hat_to_beta} yields \eqref{eq:beta_hat_perturb}.

\medskip

\noindent\underline{Proof of \eqref{eq:beta_hat_vert_line}.} It is easy to see that \eqref{eq:beta_hat_vert_line} immediately follows from \eqref{eq:beta_hat_perturb} and \eqref{eq:eta_stab_bound} applied to the rhs. of \eqref{eq:beta_hat_vert_line}.\hfill$\qed$

\subsection{Proof of Lemma \ref{lem:G_bounds_extension}}\label{app:G_bounds_extension}

We derive Lemma~\ref{lem:G_bounds_extension} from the following result for deterministic matrices.

\begin{lemma}\label{lem:size_bound_extension} For $n\in\N$ consider $D_j=D_j^*\in \C^{n\times n}$ for $j=1,2$, and for $w\in\C\setminus\R$ denote $G_j(w):=(D_j-w)^{-1}$. Assume that for some $B=B^*\in\C^{n\times n}$, $\eta_0>0$ and $\gamma>0$ it holds that
\begin{align}
\left\vert\langle G_1(\pm\ii\eta_0)BG_2(\pm\ii\eta_0)B\rangle\right\vert &\le \frac{1}{\gamma},\label{eq:size_input2}\\
\left\vert \langle G_1(\pm\ii\eta_0)BG_2(\pm\ii\eta_0)BG_1(\pm\ii\eta_0)BG_2(\pm\ii\eta_0)B\rangle\right\vert &\le \frac{1}{\eta_0\gamma^2},\label{eq:size_input4} 
\end{align}
where in \eqref{eq:size_input2}--\eqref{eq:size_input4} all choices of $\pm$ are considered independently from each other. Then for any $\eta_1,\eta_2\in (0,\eta_0]$ we have
\begin{align}
\left\vert\langle (G_1(\ii\eta_1))^2 B (G_2(\ii\eta_2))^2B\rangle\right\vert &\le \frac{\eta_0^2}{(\eta_1\eta_2)^2}\frac{1}{\gamma},\label{eq:size_bound1}\\
\left\vert\langle (G_1(\ii\eta_1))^2 B G_2(\ii\eta_2)B\rangle\right\vert &\le\left(\frac{\eta_0^2}{\eta_1\eta_2}+\frac{\eta_0}{\eta_1} \langle \Im G_1(\ii\eta_0)\rangle^{1/2}\right)\frac{1}{\eta_1\gamma}.\label{eq:size_bound2}
\end{align}
\end{lemma}

\begin{proof}[Proof of Lemma~\ref{lem:G_bounds_extension}] Fix a small $\xi>0$. We apply Lemma~\ref{lem:size_bound_extension} with
\begin{equation*}
n:=2N, \,\,D_1:=W-Z_1-\Re w_1,\,\, D_2:=(W-Z_2-\Re w_2)^\mt,\,\, B:=E_\sigma,\,\,\eta_0:=N^{-1+\epsilon}\,\,\text{and}\,\,\gamma:=N^{-\xi}\widehat{\beta}_{1\overline{2}}^{[b]},
\end{equation*}
where $W$ is the Hermitization of $\widehat{X}$. Then \eqref{eq:2G_bound_cr} and \eqref{eq:4G_bound_cr} imply that the assumptions \eqref{eq:size_input2} and \eqref{eq:size_input4} hold with very high probability. Using that $|\langle \Im G_1(\ii\eta_0)\rangle|\lesssim 1$ by \cite[Lemma~B.7]{cipolloni2024maximum} and applying Lemma~\ref{lem:size_bound_extension}, we finish the proof of Lemma~\ref{lem:G_bounds_extension}.
\end{proof}

\begin{proof}[Proof of Lemma \ref{lem:size_bound_extension}] We start with the proof of \eqref{eq:size_bound1}. By the Cauchy-Schwarz inequality and the Ward identity we have
\begin{equation}
\left\vert\langle (G_1(\ii\eta_1))^2 B (G_2(\ii\eta_2))^2B\rangle\right\vert \le \frac{1}{\eta_1\eta_2}\langle \Im G_1(\ii\eta_1)B \Im G_2(\ii\eta_2)B\rangle.
\label{eq:determ_Schwarz}
\end{equation}
Denote the eigenvalues of $D_j$ by $\{\lambda_i^{(j)}\}_{i=1}^n$ and let $\{\bm{w}_i^{(j)}\}_{i=1}^n$ be the associated normalized eigenvectors. Since $B=B^*$, we get
\begin{equation}
\begin{split}
&\langle \Im G_1(\ii\eta_1)B \Im G_2(\ii\eta_2)B\rangle = \frac{1}{N}\sum_{i,j=-N}^N \frac{\eta_1}{\big(\lambda_i^{(1)}\big)^2+\eta_1^2}\frac{\eta_2}{\big(\lambda_j^{(2)}\big)^2+\eta_2^2}\left\vert\langle \bm{w}^{(1)}_i,B\bm{w}^{(2)}_j\rangle\right\vert^2\\
&\quad \le \frac{\eta_0^2}{N\eta_1\eta_2}\sum_{i,j=-N}^N \frac{\eta_0}{\big(\lambda_i^{(1)}\big)^2+\eta_0^2}\frac{\eta_0}{\big(\lambda_j^{(2)}\big)^2+\eta_0^2}\left\vert\langle \bm{w}^{(1)}_i,B\bm{w}^{(2)}_j\rangle\right\vert^2 = \frac{\eta_0^2}{\eta_1\eta_2} \langle \Im G_1(\ii\eta_0)B \Im G_2(\ii\eta_0)B\rangle.
\end{split}
\label{eq:determ_monot}
\end{equation}
Combining \eqref{eq:determ_Schwarz}, \eqref{eq:determ_monot} and~\eqref{eq:size_input2}, we complete the proof of~\eqref{eq:size_bound1}.

To prove \eqref{eq:size_bound2}, we represent the lhs. of \eqref{eq:size_bound2} as follows:
\begin{equation}
\langle (G_1(\ii\eta_1))^2 B G_2(\ii\eta_2)B\rangle = \langle (G_1(\ii\eta_1))^2 B G_2(\ii\eta_0)B\rangle - \int_{\eta_2}^{\eta_0}\partial_y \langle (G_1(\ii\eta_1))^2 B G_2(\ii y)B\rangle\dif y.
\label{eq:size_bound2_aux1}
\end{equation}
Next we upper bound the integral in the rhs. of \eqref{eq:size_bound2_aux1} using \eqref{eq:size_bound1}:
\begin{equation}
\left\vert\int_{\eta_2}^{\eta_0}\!\!\partial_y \langle (G_1(\ii\eta_1))^2 B G_2(\ii y)B\rangle\dif y\right\vert \le \int_{\eta_2}^{\eta_0}\!\! \left\vert \langle (G_1(\ii\eta_1))^2 B (G_2(\ii y))^2 B\rangle\right\vert\dif y\le \int_{\eta_2}^{\eta_0}\!\! \frac{\eta_0^2}{\eta_1^2y^2\gamma}\dif y\le \frac{\eta_0^2}{\eta_1^2\eta_2\gamma}.
\end{equation}
For the first term in the rhs. of \eqref{eq:size_bound2_aux1} we analogously to \eqref{eq:size_bound2_aux1} write
\begin{equation}
\langle (G_1(\ii\eta_1))^2 B G_2(\ii\eta_0)B\rangle = \langle (G_1(\ii\eta_0))^2 B G_2(\ii\eta_0)B\rangle-\int_{\eta_1}^{\eta_0}\partial_y \langle (G_1(\ii y))^2 B G_2(\ii\eta_0)B\rangle\dif y
\label{eq:size_bound2_aux2}
\end{equation}
We compute the derivative in the rhs. of \eqref{eq:size_bound2_aux2} and apply the Cauchy-Schwarz inequality and the Ward identity as in \eqref{eq:determ_Schwarz}:
\begin{equation}
\begin{split}
&\frac{1}{2}\left\vert \partial_y \langle (G_1(\ii y))^2 B G_2(\ii\eta_0)B\rangle\right\vert =\left\vert \langle (G_1(\ii y))^3 B G_2(\ii\eta_0)B\rangle\right\vert \\
&\quad\le  \frac{1}{\sqrt{y}}\langle\Im G_1(\ii y)\rangle^{1/2}\frac{1}{y}\langle \Im G_1(\ii y)BG_2(\ii \eta_0)B\Im G_1(\ii y)B G_2(-\ii\eta_0)B\rangle^{1/2}\\
&\quad \le \frac{\eta_0^{3/2}}{y^3} \langle \Im G_1(\ii\eta_0)\rangle^{1/2}\langle \Im G_1(\ii \eta_0)BG_2(\ii \eta_0)B\Im G_1(\ii \eta_0)B G_2(-\ii\eta_0)B\rangle^{1/2}\le \frac{\eta_0}{y^3\beta} \langle \Im G_1(\ii\eta_0)\rangle^{1/2}.
\end{split}
\label{eq:determ_Schwarz_monot}
\end{equation}
Here to go from the second to the third line we used the spectral decomposition of $G_1, G_2$ similarly to \eqref{eq:determ_monot}, and in the last bound applied \eqref{eq:size_input4}. Combining \eqref{eq:size_bound2_aux1}--\eqref{eq:determ_Schwarz_monot}, we finish the proof of \eqref{eq:size_bound2}.
\end{proof}

\end{document}